\documentclass[reqno]{amsart}

\usepackage{mathrsfs}
\usepackage{amsmath}
\usepackage{amssymb}
\usepackage{cite}
\usepackage{latexsym}
\usepackage{graphicx}

\usepackage{color}
\usepackage{comment}

\theoremstyle{plain}
\newtheorem{thm}{Theorem}[section]

\newtheorem{cor}[thm]{Corollary}

\newtheorem{prop}[thm]{Proposition}
\newtheorem{rmk}[thm]{Remark}

\def\D{\mathrm{D}}
\def\L{\mathscr{L}}

\def\U{\mathscr{U}}

\def\d{\mathrm{d}}

\def\Nset{\mathbb{N}}

\def\Rset{\mathbb{R}}
\def\Sset{\mathbb{S}}
\def\Zset{\mathbb{Z}}

\def\epsilon{\varepsilon}

\makeatletter
 \@addtoreset{equation}{section}
\makeatother
\def\theequation{\arabic{section}.\arabic{equation}}

\begin{document}


\title[Feedback control of the Kuramoto model]%
{Feedback control of twisted states in the Kuramoto model
 on nearest neighbor and complete simple graphs}

\author{Kazuyuki Yagasaki}

\address{Department of Applied Mathematics and Physics, Graduate School of Informatics,
Kyoto University, Yoshida-Honmachi, Sakyo-ku, Kyoto 606-8501, JAPAN}
\email{yagasaki@amp.i.kyoto-u.ac.jp}

\date{\today}

\begin{abstract}
We study feedback control of twisted states in the Kuramoto model (KM)
 of identical oscillators defined on deterministic nearest neighbor graphs
 containing complete simple ones when it may have phase-lag.
Bifurcations of such twisted solutions in the continuum limit (CL)
 for the uncontrolled KM defined on nearest neighbor graphs
 that may be deterministic dense, random dense or random sparse
 were discussed very recently by using the center manifold reduction,
 which is a standard technique in dynamical systems theory.
In this paper we analyze the stability and bifurcations 
 of twisted solutions in the CL for the KM subjected to feedback control.
In particular, {\color{black}
it is shown that} the twisted solutions exist and can be stabilized
 not only for nearest neighbor graphs but also for complete simple graphs.
Moreover, the CL is shown to suffer bifurcations
 at which the twisted solution becomes unstable
 and a stable one-parameter family of modulated or oscillating twisted solutions
 is born, depending on whether the phase-lag is zero or not.
We demonstrate the theoretical results by numerical simulations
 for the feedback controlled KM
 on deterministic nearest neighbor and complete simple graphs.
\end{abstract}

\maketitle


\section{Introduction}

\subsection{\color{black}Background}
Coupled oscillators in complex networks
 have recently attracted {\color{black}significant} attention
 and have been studied {\color{black}with rapidly increasing intensity
 \cite{RPJK16,S16,BAGLLSWZ16,EPW19,BGLM20,WL20,ZSZK21,AA22,MSD23,FM24,MRSKG24,EN25}.}
They provide many mathematical models in various fields
 such as physics, chemistry, biology, social sciences and engineering,
{\color{black}and exhibit several collective dynamics
 including synchronization, chimeras and chaos.}
Among them, the Kuramoto model (KM) \cite{K75,K84}
 is one of the most representative models and has been generalized in several directions.
{\color{black}
It was originally proposed by Kuramoto \cite{K75,K84} half a century ago,
 and has continued to be the subject of enormous research,
 especially to discuss synchronization phenomena in diverse fields, since then.
The range of its direct applications now spreads
 to power grids \cite{FNP08,DB12,SA15,LM17,SKWLMP19,LZ20,GZLWY21},
 neuroscience \cite{BHD10,FCPL12,TE14,ZLQWWL22,PDS24,YLYWQZH25},
 machine learning \cite{WHKY21,ABMH22,SCOK23} and so on.
See the above surveys and \cite{S00,PRK01,ABPRS05,BLMCH06,ADKMZ08,DB14,PR15}}
 for the reviews of vast literature on coupled oscillators in complex networks
 including the KM and its generalizations.

The control problem of nonlinear oscillator networks
 is important not only in theoretical interest but also in applications,
 and has drawn much attention \cite{CKM13,M15,DNM22}.
Feedback control of synchronized states different from twisted states
 in the KM on deterministic dense, random dense and random sparse graphs
 was studied numerically or theoretically
 in \cite{SA15,SA16,SNKJZBPPB20,DNM22,IY23,Y24d,KY24}.
In particular, {\color{black}
the asymptotic stability of even orbits that do not coincide with the desired orbit
 but approach it as the feedback gain tends to infinity} has been discussed theoretically
 only in \cite{Y24d,KY24}.

\subsection{\color{black}Feedback controlled Kuramoto model}
In this paper we consider feedback control of the KM
 consisting of identical oscillators on a deterministic dense graph
 $G_{n}=\langle V(G_n),E(G_n),W(G_n)\rangle$,
\begin{align}
\frac{\d}{\d t} u_k^n (t)
=& \omega
 +\frac{1}{n}\sum_{j=1}^{n}w_{kj}^n\sin\left( u_j^n(t)-u_k^n(t)+\sigma\right)\notag\\
& +b_1(\hat{u}_k^n(t)-u_k^n(t))+b_3(\hat{u}_k^n(t)-u_k^n(t))^3,\quad
k\in[n]:=\{1,2,\ldots,n\},
\label{eqn:dsys}
\end{align}
where $u_k^n:\Rset\to\Sset^1:=\Rset/2\pi\Zset$
 stands for the phase of oscillator at the node $k\in [n]$;
 $\hat{u}_k^n(t)$, $k\in[n]$, represent the target orbit;
 $\omega$ is the natural frequency;
 $\sigma\in(-\tfrac{1}{2}\pi,\tfrac{1}{2}\pi)$ is the phase-lag parameter;
and  $b_1,b_3>0$ are the linear and nonlinear feedback gains.
{\color{black}
We use the convenient notation $[n]$, which represents the set $\{1,2,\ldots,n\}$,
 throughout this paper.
In \cite{SA15,SA16,SNKJZBPPB20,IY23,Y24d,KY24},
 a different nonlinear feedback control input,
\[
\tilde{b}_1\sin\left({\color{black}\hat{u}(t)}-u_i^n(t)\right)+\tilde{b}_0,
\]
where {\color{black}the target orbit $\hat{u}(t)$ is independent of $n$ and $k\in[n]$,}
 and $\tilde{b}_1,\tilde{b}_0$ are constants with $\tilde{b}_1>0$,
 was treated.

On the other hand,}
 $V(G_{n})=[n]$ and $E(G_{n})$ represent the sets of nodes and edges, respectively,
 and $W(G_{n})$ is an $n\times n$ weight matrix given by
\begin{equation*}
(W(G_{n}))_{kj}=
\begin{cases}
w_{kj}^{n} & \mbox{if $(k,j)\in E(G_{n})$};\\
0 &\rm{otherwise}.
\end{cases}
\end{equation*}
So we express
\[ 
E(G_{n})=\{(k,j)\in[n]^2\mid (W(G_{n}))_{kj}\neq 0\}, 
\] 
where each edge is represented by an ordered pair of nodes $(k,j)$, 
 which is also denoted by $j\to k$, and a loop is allowed.
If $W(G_{n})$ is symmetric,
 then $G_{n}$ represents an undirected weighted graph 
 and each edge is also denoted by $k\sim j$ instead of $j\to k$. 
When $G_{n}$ is a simple graph, 
 $W(G_{n})$ is a matrix whose elements are $\{0,1\}$-valued. 
Moreover, the weight matrix $W(G_{n})$ is assumed to be given as follows.
Let $I=[0,1]$ and let $W^n\in L^2 (I^2)$, $n\in\Nset$, be nonnegative functions.
We have
\begin{equation*}
w_{kj}^{n} = \langle W^n\rangle_{kj}^{n}
:= n^2 \int_{I_k^n \times I_j^n}W^n(x,y) \d x\d y,
\end{equation*}
where
\[
I_k^n:=
\begin{cases}
  [(k-1)/n,k/n) & \mbox{for $k<n$};\\
  [(n-1)/n,1] & \mbox{for $k=n$}.
\end{cases}
\]
Such a function as $W^n(x,y)$ is usually called a \emph{graphon} \cite{L12}. 
We also assume that there exists a measurable function $W\in L^2(I^2)$ such that 
\begin{equation*}
\|W(x,y)-W^n(x,y)\|_{L^2(I^2)}=\int_{I^2}|W(x,y)-W^n(x,y)|^2\d x\d y\to 0
\label{eqn:W}
\end{equation*}
as $n\to\infty$.
{\color{black}We only assume here that $G_n$ is deterministic,
 but can similarly treat the case in which it is random dense or sparse.}
 
\subsection{\color{black}Continuum Limit (CL)}
In \cite{IY23}, a general coupled oscillator network
\begin{equation}
\frac{\d}{\d t} u_k^n (t)
=f(u_k^n(t),t)
 +\frac{1}{n}\sum_{l=1}^{m}\sum_{j=1}^{n}w_{kj}^{ln}D_l\left( u_j^n(t)-u_k^n(t)\right),\quad
k\in[n],
\label{eqn:dsys0}
\end{equation}
defined on $m$ $(\ge 2)$ graphs $G_{ln}$, $l\in[m]$,
 which may be not only deterministic dense but also random dense or sparse,
 was studied and shown to be well approximated
 by the corresponding continuum limit,
\begin{equation*}
\frac{\partial}{\partial t}u(t,x)
=f(u(t,x),t)+\sum_{l=1}^{m}\int_I W_l(x,y)D_l(u(t,y)-u(t,x))\d y,\quad x \in I,
\end{equation*}
where $u_k^n:\Rset\to\Rset$, $k\in[n]$,
 $f(u,t)$ is Lipschitz continuous in $u$ and continuous in $t$
 and $D_l(u)$, $l\in[m]$, are bounded and Lipschitz continuous.
The result is also applicable
 when the natural frequency at each node is different,
{\color{black}and it was improved in \cite{Y24a},
 so that relationships between such a coupled oscillator network and its CL
 on the stability of their solutions were further developed.
See Section~2 for more details on these results.
Moreover,  in \cite{Y24e},
 they were modified for the random natural frequency case
 by introducing a random permutation.}
Similar results for such networks defined on single graphs
 and having the same {\color{black}or equivalently zero,} natural frequency at each node
 were obtained earlier in \cite{KM17,M14a,M14b,M19}
 {\color{black}although neither the random natural frequencies nor stability of solutions
 were treated}.
Such a CL was introduced for the classical KM,
 which is defined on a single complete simple graph
 but may have natural frequencies depending on each oscillator,
 without a rigorous mathematical guarantee very early in \cite{E85},
{\color{black}
 and it was fully discussed very recently in \cite{Y24a}
 when the natural frequencies are uniformly spaced.
In \cite{Y24d} and \cite{KY24}, respectively,
 the results of \cite{IY23,Y24a,Y24e} were used or extended successfully
 to discuss feedback control of the KM on uniform graphs
 which may be complete, random dense or sparse
 when the natural frequencies are uniformly spaced and random.}
Similar CLs were utilized for the KM
 with nonlocal coupling and a single or zero natural frequency
 in \cite{GHM12,M14c,MW17,WSG06}.
However, the feedback controlled KM \eqref{eqn:dsys}
 does not have the form \eqref{eqn:dsys0}.

{\color{black}Here we} take as the graphons $W^n(x,y)$ and $W(x,y)$
\[
W^n(x,y)=\begin{cases}
1 & \mbox{if $(x,y)\in I_k^n\times I_j^n$ with $|k-j|\le n\kappa$ or $|k-j|\ge n(1-\kappa)$};\\
0 & \mbox{otherwise},
\end{cases}
\]
and
\[
W(x,y)=\begin{cases}
1 & \mbox{if $|x-y|\le\kappa$ or $|x-y|\ge1-\kappa$};\\
0 & \mbox{otherwise},
\end{cases}
\]
with $0<\kappa\le\tfrac{1}{2}$,
 which correspond to nearest {\color{black}(more specifically, $\lfloor n\kappa\rfloor$-nearest)}
 neighbor graphs{\color{black},
 where $\lfloor z\rfloor$ represents the maximum integer that is not greater than $z\in\Rset$.}
For $\kappa=\tfrac{1}{2}$,
 they become $W^n(x,y),W(x,y)\equiv 1$ and correspond to complete simple graphs.
Moreover, we choose as the target orbit
\begin{equation}
\hat{u}_k^n(t)=\frac{2\pi qk}{n}+\Omega_\D^n t,\quad
k\in[n],\quad 
q\in\Nset,
\label{eqn:ts}
\end{equation}
where $\Omega_\D^n$ is a constant given by
\begin{equation*}
\Omega_\D^n=\omega
 +\frac{1}{n}\sum_{|j|\le n\kappa}\sin\left(\frac{2\pi qj}{n}+\sigma\right).
\label{eqn:OmegaD*}
\end{equation*}
Note that Eq.~\eqref{eqn:ts}  provides a particular solution to \eqref{eqn:dsys}
 even if $\kappa=\tfrac{1}{2}$, i.e., the graph $G_n$ is complete simple,
 whether $b_1,b_3=0$ or not.
Letting $v_k^n(t)=u_k^n(t)-\hat{u}_k^n(t)$, $k\in[n]$,
 we rewrite \eqref{eqn:dsys} as 
\begin{align}
\frac{\d}{\d t} v_k^n (t)
=& \omega
 +\frac{1}{n}\sum_{j=1}^{n}w_{kj}^n\cos\frac{2\pi q(j-k)}{n}
  \sin\left(v_j^n(t)-v_k^n(t)+\sigma\right)\notag\\
&
-\frac{1}{n}\sum_{j=1}^{n}w_{kj}^n\sin\frac{2\pi q(j-k)}{n}
  \cos\left(v_j^n(t)-v_k^n(t)+\sigma\right)\notag\\
&
-b_1v_k^n(t)-b_3v_k^n(t)^3,\quad
k\in[n],
\label{eqn:dsys1}
\end{align}
which has the form \eqref{eqn:dsys0}.
Using the result{\color{black}s of \cite{IY23,Y24a}},
 we see that the coupled oscillator network \eqref{eqn:dsys1} is well approximated 
 by its CL
\begin{align}
\frac{\partial}{\partial t}v(t,x)
=&\omega+\int_I W(x,y)\cos2\pi q(y-x)\sin(v(t,y)-v(t,x)+\sigma)\d y\notag\\
&-\int_I W(x,y)\sin2\pi q(y-x)\cos(v(t,y)-v(t,x)+\sigma)\d y,\notag\\
& -b_1v(x,t)-b_3v(x,t)^3,\quad x \in I,
\label{eqn:csys1}
\end{align}
so that the KM \eqref{eqn:dsys} is well approximated by the CL,
\begin{align}
\frac{\partial}{\partial t}u(t,x)
=& \omega+\int_I W(x,y) \sin(u(t,y)-u(t,x)+\sigma)\d y\notag\\
& +b_1(\hat{u}(t,x)-u(t,x))+b_3(\hat{u}(t,x)-u(t,x))^3,\quad x \in I,
\label{eqn:csys}
\end{align}
where
\begin{equation}
\hat{u}(t,x)=2\pi qx+\Omega t,\quad
q\in\Nset,
\label{eqn:tsol}
\end{equation}
and
\begin{equation}
\Omega=\omega+\int_{x-\kappa}^{x+\kappa}\sin(2\pi q(y-x)+\sigma)\d y
=\omega+\frac{\sin 2\pi q\kappa\sin\sigma}{\pi q}.
\label{eqn:Omega}
\end{equation}
Note that Eq.~\eqref{eqn:tsol} also provides a particular solution to \eqref{eqn:csys},
 even if $\kappa=\tfrac{1}{2}$, i.e., $W(x,y)\equiv 1$, whether $b_1,b_3=0$ or not,
 although it is unstable for $\kappa=\tfrac{1}{2}$ and $b_1=0$
 (see Remark~\ref{rmk:3a}(iii) below),
 and that $\Omega=\lim_{n\to\infty}\Omega_\D^n$.
Moreover,
\begin{equation}
\Bigl\|\hat{u}(t,x)-\sum^{n}_{i=1}\hat{u}_k^n(t) \mathbf{1}_{I_k^n}\Bigr\|_{L^2(I)}
 =\int_I\Bigl|\hat{u}(t,x)-\sum^{n}_{i=1}\hat{u}_k^n(t) \mathbf{1}_{I_k^n}\Bigr|^2\d x\to 0
\label{eqn:conv}
\end{equation}
as  $n\to\infty$,
 where $\mathbf{1}_{I_j^n}$ represents the characteristic function of $I_k^n$, $i\in[n]$.

\subsection{\color{black}Outline of the results}
In the previous work \cite{Y24c},
 the uncontrolled CL \eqref{eqn:csys} with $b_1,b_3=0$
 was studied and bifurcations of the $q$-twisted solutions
\begin{equation}
u(t,x)=2\pi qx+\Omega t+\theta\in\Sset^1,\quad
\theta\in\Sset^1,\quad
q\in\Nset,
\label{eqn:tsol0}
\end{equation}
when $\kappa$ is taken as a control parameter were analyzed
 by using the center manifold reduction \cite{GH83,HI11,K04},
 which is a standard technique in dynamical systems theory.
In particular, it was shown {\color{black}in \cite{Y24c}} that at the bifurcation points,
 the twisted solutions change their stability from stable to unstable, 
 and stable or unstable modulated or oscillating twisted solutions are born,
 depending on whether $\sigma=0$ or not,
 when $\kappa$ is increased.
The theoretical results for $q\in[2]$ or $[4]$
 {\color{black}(recall that $[2]$ and $[4]$ represent the sets $\{1,2\}$ and $\{1,2,3,4\}$,
 respectively)}
 were also demonstrated in numerical simulations
 for the uncontrolled KM \eqref{eqn:dsys} with $b_1,b_3=0$
 defined on deterministic dense, random dense and random sparse nearest neighbor graphs
 although the observation of the modulated and oscillating twisted solutions
 born at the bifurcations
 was impossible or very subtle since they are unstable or considered to disappear
 near the bifurcations.
A similar bifurcation analysis was performed
 for synchronized solutions in the CL
 of the KM {\color{black}with two-mode interaction} defined on two graphs in \cite{Y24b}.

Here we take the linear feedback gain $b_1$ as a control parameter
 and analyze bifurcations of the $q$-twisted solution \eqref{eqn:tsol}
 in the controlled CL \eqref{eqn:csys} with $b_1,b_3\neq 0$ for $q\in[4]$,
 using the center manifold reduction as in \cite{Y24b,Y24c}.
We show that at the bifurcation points,
 the twisted solutions change their stability from stable to unstable,
 and stable modulated or oscillating twisted solutions
 depending on whether $\sigma=0$ or not,
 when $b_1$ is decreased and $b_3>0$ is sufficiently large.
This is similar to the uncontrolled case of $b_1,b_3 =0$,
 but a little different since modulated and oscillating twisted solutions
 born at the bifurcations are stable in a wide range of the control parameter $b_1$.
The CL \eqref{eqn:csys} can be controlled
 to exhibit the desired $q$-twisted solution \eqref{eqn:tsol},
 which is the same as \eqref{eqn:tsol0} with $\theta=0$, for $q\in[4]$,
 by choosing appropriate values of $b_1,b_3$,
 even if the one-parameter family \eqref{eqn:tsol0} of $q$-twisted solutions
 is unstable when $b_1,b_3=0$.
It follows from the results of \cite{IY23,Y24a} that 
 such bifurcations are also expected to occur in the KM \eqref{eqn:dsys}
 (see Remarks~\ref{rmk:4a}(iii) and \ref{rmk:4c}(iii) below).
In particular, if the $q$-twisted solution \eqref{eqn:tsol}
 is asymptotically stable in the CL \eqref{eqn:csys},
 then our feedback control is considered to be successful
 since the $q$-twisted state \eqref{eqn:ts}
 is expected to be asymptotically stable in the KM \eqref{eqn:dsys}.
We demonstrate our theoretical results
 by numerical simulations for the KM \eqref{eqn:dsys}
 on deterministic nearest neighbor and complete simple graphs.
Stable modulated and oscillating twisted states
 are observed clearly in wide ranges of $b_1$ in the numerical simulations.
 
The remainder of this paper is as follows:
In Section~2 we briefly review the previous fundamental results
 of \cite{IY23,Y24a} {\color{black}on relationships
 between coupled oscillators and their CLs
 in the context of \eqref{eqn:dsys} and \eqref{eqn:csys}.
In particular, we will see that
 an asymptotically stable and unstable solutions to the CL \eqref{eqn:csys}, respectively,
 behave as if they are asymptotically stable and unstable in the KM \eqref{eqn:dsys}.}
We analyze the associated linear eigenvalue problem
 for $q$-twisted solutions to the CL \eqref{eqn:csys} in Section~3,
 and discuss their bifurcations for $q\in[4]$ in Section~4.
Numerical simulation results of the KM \eqref{eqn:dsys}
 on the nearest neighbor and complete simple graphs
 are given in Sections~5 and 6, respectively.


\section{\color{black}
Relationships between the KM \eqref{eqn:dsys} and CL \eqref{eqn:csys}}

We first review the results of {\color{black}\cite{IY23,Y24a}
 on relationships between coupled oscillator networks and their CLs}
 in the context of the KM \eqref{eqn:dsys} and CL \eqref{eqn:csys}.
See Section~2 and Appendices~A and B of \cite{IY23}
 and Section~2 of \cite{Y24a} for more details,
 including the proofs of the theorems stated below.
These results for \eqref{eqn:dsys} and \eqref{eqn:csys}
 follow from application of the results of \cite{IY23,Y24a}
 to \eqref{eqn:dsys1} and \eqref{eqn:csys1}.

Let $g(x)\in L^2(I)$
 and let $\mathbf{u}:\Rset\to L^2(I)$ stand for an $L^2(I)$-valued function.
We have the following on the existence and  uniqueness of solutions
 to the initial value problem (IVP) of the CL \eqref{eqn:csys}
 (see Theorem~2.1 of \cite{IY23}).
 
\begin{thm}
\label{thm:2a}
There exists a unique solution $\mathbf{u}(t)\in C^1(\Rset,L^2(I))$
 to the IVP of \eqref{eqn:csys} with
\begin{equation*}
u(0,x)=g(x).
\end{equation*}
Moreover, the solution depends continuously on $g$.
\end{thm}

We next consider the IVP of the KM \eqref{eqn:dsys}
 and turn to the issue on convergence of solutions in \eqref{eqn:dsys}
 to those in the CL \eqref{eqn:csys}.
Since the right-hand side of \eqref{eqn:dsys} is Lipschitz continuous in $u_k^n$, $i\in[n]$,
 we see by a fundamental result of ordinary differential equations
 (e.g., Theorem~2.1 of Chapter~1 of \cite{CL55})
 that the IVP of \eqref{eqn:dsys} has a unique solution.
Given a solution $u_n(t)=(u_1^n(t),\ldots, u_n^n(t))$ to the IVP of \eqref{eqn:dsys},
 we define an $L^2(I)$-valued function $\mathbf{u}_n:\Rset\to L^2(I)$ as
\begin{equation*}
\mathbf{u}_n(t) = \sum^{n}_{j=1} u_j^n(t) \mathbf{1}_{I_j^n}.
\end{equation*}
Recall that $\mathbf{1}_{I_j^n}$ represents the characteristic function of $I_j^n$, $j\in[n]$.
Let $\|\cdot\|$ denote the norm in $L^2(I)$.
In our setting stated in Section~1,
 we slightly modify the arguments given in the proof of Theorem~2.3 of \cite{IY23}
 to obtain the following
 (see also Theorem~2.2 of \cite{Y24a}).

\begin{thm}
\label{thm:2b}
If $\mathbf{u}_{n}(t)$ is the solution to the IVP of \eqref{eqn:dsys}
 with the initial condition
\[
\lim_{n\to\infty}\|\mathbf{u}_n(0)-\mathbf{u}(0)\|=0,
\]
then for any $T > 0$ we have
\[
\lim_{n \rightarrow \infty}\max_{t\in[0,T]}\|\mathbf{u}_n(t)-\mathbf{u}(t)\|=0,
\]
where $\mathbf{u}(t)$ represents the solution
 to the IVP of the CL \eqref{eqn:csys}. 
\end{thm}

We also obtain the following result,
 slightly modifying the proof of Theorem~2.5 in \cite{IY23}
 (see also Theorem~2.3 of \cite{Y24a}).

\begin{thm}
\label{thm:2c}
Suppose that the KM \eqref{eqn:dsys} and CL \eqref{eqn:csys}
 have solutions $\bar{\mathbf{u}}_n(t)$ and $\bar{\mathbf{u}}(t)$, respectively, such that
\begin{equation}
\lim_{n\to\infty}\|\bar{\mathbf{u}}_n(t)-\bar{\mathbf{u}}(t)\|=0
\label{eqn:2c}
\end{equation}
for any $t\in[0,\infty)$.
Then the following hold$:$
\begin{enumerate}
\setlength{\leftskip}{-1.8em}
\item[\rm(i)]
{\color{black}
If for any $\epsilon>0$, there exist $\delta_1>0$
 such that for $n>0$ sufficiently large,
 any solution $\bar{u}_j^n(t)$, $j\in[n]$, to the KM \eqref{eqn:dsys} with
\[
|u_j^n(0)-\bar{u}_j^n(0)|<\delta_1,
\quad j\in[n],
\] 
satisfies
\[
|u_j^n(t)-\bar{u}_j^n(t)|<\epsilon,
\quad j\in[n],
\] 
then $\bar{\mathbf{u}}(t)$ is stable.
Moreover, If for any $\epsilon>0$, there exists $\delta_2>0$
 such that for $n>0$ sufficiently large,
 any solution $\bar{u}_j^n(t)$, $j\in[n]$, to the KM \eqref{eqn:dsys}  with
\[
|u_j^n(t)-\bar{u}_j^n(t)-\theta|<\delta_2,\quad j\in[n],
\] 
converges to $\bar{u}_j^n(t)$, $j\in[n]$,
 then $\bar{\mathbf{u}}(t)$ is asymptotically stable.}
\item[\rm(ii)]
{\color{black}
If $\bar{\mathbf{u}}(t)$ is stable, then for any $\epsilon,T>0$ there exists $\delta>0$
 such that for $n>0$ sufficiently large,
 if $\bar{\mathbf{u}}_n(t)$ is any solution to the KM \eqref{eqn:dsys} satisfying
\begin{equation*}
\|\mathbf{u}_n(0)-\bar{\mathbf{u}}_n(0)\|<\delta,
\end{equation*}
then
\begin{equation*}
\|\mathbf{u}_n(t)-\bar{\mathbf{u}}_n(t)\|<\epsilon.
\end{equation*}
Moreover, if $\bar{\mathbf{u}}(t)$ is asymptotically stable, then
\begin{equation*}
\lim_{t\to\infty}\lim_{n\to\infty}\|\mathbf{u}_n(t)-\bar{\mathbf{u}}_n(t)\|=0,
\end{equation*}
where $\mathbf{u}_n(t)$ is any solution to \eqref{eqn:dsys}
 such that $\mathbf{u}_n(0)$ is contained in the basin of attraction for $\bar{\mathbf{u}}(t)$.}
\end{enumerate}
\end{thm}

\begin{rmk}\
\label{rmk:2a}
\begin{enumerate}
\setlength{\leftskip}{-1.8em}
\item[\rm(i)]
The solution $\bar{\mathbf{u}}_n(t)$ may not be asymptotically stable
 in the KM \eqref{eqn:dsys}
 for $n>0$ sufficiently large even if so is $\bar{\mathbf{u}}(t)$  in the CL \eqref{eqn:csys}.
In the definition of stability and asymptotic stability of solutions to the CL \eqref{eqn:csys},
 we cannot distinguish two solutions that are different only in a set with the Lebesgue measure zero.
\item[\rm(ii)]
From the proof of Theorem~$2.5$ in {\rm\cite{IY23}}
 the stability stated in Theorem~$\ref{thm:2c}$
 contains not only the Lyapunov meaning but also the orbital one.
\end{enumerate}
\end{rmk}

We have the following as a corollary of Theorem~\ref{thm:2c},
 without assuming the existence
 of the solution $\bar{\mathbf{u}}_n(t)$ to the KM \eqref{eqn:dsys} satisfying \eqref{eqn:2c}
 (see the proof of Theorem~2.4(ii) and Corollary~2.6 of \cite{Y24a}).

\begin{cor}
\label{cor:2a}
Suppose that the CL \eqref{eqn:csys} has a stable solution $\bar{\mathbf{u}}(t)$.
Then for any $\epsilon,T>0$ there exists $\delta>0$ such that for $n>0$ sufficiently large,
 if $\mathbf{u}_n(t)$ is any solution to the KM \eqref{eqn:dsys} satisfying
\[
\|\mathbf{u}_n(0)-\bar{\mathbf{u}}(0)\|<\delta,
\]
then
\[
\max_{t\in[0,T]}\|\mathbf{u}_n(t)-\bar{\mathbf{u}}(t)\|<\epsilon.
\]
Moreover, if $\bar{\mathbf{u}}(t)$ is asymptotically stable, then
\begin{equation*}
\lim_{t\to\infty}\lim_{n\to\infty}\|\mathbf{u}_n(t)-\bar{\mathbf{u}}(t)\|=0,
\end{equation*}
where $\mathbf{u}_n(t)$ is any solution to \eqref{eqn:dsys}
 such that $\mathbf{u}_n(0)$ is contained in the basin of attraction for $\bar{\mathbf{u}}(t)$.
\end{cor}

\begin{rmk}\
\label{rmk:2b}
\begin{enumerate}
\setlength{\leftskip}{-1.6em}
\item[\rm(i)]
In Corollary~$2.6$ of {\rm\cite{Y24a}} only complete simple graphs were treated
 but Corollary~{\rm\ref{cor:2a}} can be proven similarly
 since its proof relies only on Theorem~$2.2$ of {\rm\cite{Y24a}},
 of which extension to \eqref{eqn:dsys} and \eqref{eqn:csys} is Theorem~{\rm\ref{thm:2b}}.
\item[\rm(ii)]
Corollary~$\ref{cor:2a}$ implies that $\bar{\mathbf{u}}(t)$ behaves
 as if it is an $($asymptotically$)$ stable solution in the KM \eqref{eqn:dsys}.
{\color{black}
\item[\rm(iii)]
The statements of Theorem~$\ref{thm:2c}$ and Corollary~$\ref{cor:2a}$
 hold for a one-parameter family of solutions
  in the KM~$\eqref{eqn:dsys}$ and CL~$\eqref{eqn:csys}$.
See Theorem~$2.3$ and Corollary~$2.5$ of {\rm\cite{Y24c}}. 
}
\end{enumerate}
\end{rmk}

Finally, modifying the arguments given
 in the proof{\color{black}s of Theorems~2.7 and 2.9} of \cite{Y24a} slightly,
 we obtain the following results.
 
\begin{thm}
\label{thm:2d}
Suppose that the hypothesis of Theorem~$\ref{thm:2c}$ holds.
Then the following hold$:$
\begin{enumerate}
\setlength{\leftskip}{-1.8em}
\item[\rm(i)]
If $\bar{\mathbf{u}}_n(t)$ is unstable for $n>0$ sufficiently large
 and {\color{black}no stable solution to the KM \eqref{eqn:dsys} converges to $\bar{\mathbf{u}}(t)$}
 as $n\to\infty$,
 then $\bar{\mathbf{u}}(t)$ is unstable.
\item[\rm(ii)]
If $\bar{\mathbf{u}}(t)$ is unstable,
 then so is $\bar{\mathbf{u}}_n(t)$ for $n>0$ sufficiently large.
\end{enumerate}
\end{thm}

{\color{black}
\begin{thm}
\label{thm:2e}
If $\bar{\mathbf{u}}(t)$ is unstable,
 then for any $\epsilon,\delta>0$
 there exists $T>0$ such that for $n>0$ sufficiently large
\[
\max_{t\in[0,T]}\|\mathbf{u}_n(t)-\bar{\mathbf{u}}(t)\|>\epsilon,
\]
where $\mathbf{u}_n(t)$ is a solution to the KM \eqref{eqn:dsys} satisfying
\[
\|\mathbf{u}_n(0)-\bar{\mathbf{u}}(0)\|<\delta.
\]
\end{thm}}

\begin{rmk}\
\label{rmk:2c}
\begin{enumerate}
\setlength{\leftskip}{-1.6em}
\item[\rm(i)]
Only under the hypothesis of {\color{black}Corollary~$\ref{cor:2a}$},
 $\mathbf{u}(t)$ is not necessarily unstable
  even if $\mathbf{u}_n(t)$ is unstable for $n>0$ sufficiently large.
Moreover, $\mathbf{u}(t)$ may be asymptotically stable
 even if $\mathbf{u}_n(t)$ is unstable for $n>0$ sufficiently large.
See {\rm\cite{Y24a}} for such an example.
\item[\rm(ii)]
In Theorem~$2.9$ of {\rm\cite{Y24a}} only complete simple graphs were treated
 but Theorem~{\rm\ref{thm:2e}} can be proven similarly,
 like Corollary~{\rm\ref{cor:2a}}, as stated in Remark~{\rm\ref{rmk:2b}(i)}.
\item[\rm(iii)]
Theorem~$\ref{thm:2e}$ implies that $\bar{\mathbf{u}}(t)$ behaves
 as if it is an unstable solution in the KM \eqref{eqn:dsys}.
{\color{black}
\item[\rm(iv)]
The statements of Theorems~$\ref{thm:2d}$ and $\ref{thm:2e}$
 hold for a one-parameter family of solutions
  in the KM~$\eqref{eqn:dsys}$ and CL~$\eqref{eqn:csys}$.
See Theorems~$2.6$ and $4.7$ of {\rm\cite{Y24c}}. 
}

\end{enumerate}
\end{rmk}

Thus, the relationship between the KM \eqref{eqn:dsys} and CL \eqref{eqn:csys} is subtle.
However, under the hypothesis of {\color{black}Corollary~\ref{cor:2a}},
 if $\bar{\mathbf{u}}(t)$ is asymptotically stable in the CL \eqref{eqn:csys},
 then a solution to the KM \eqref{eqn:dsys}
 starting in the basin of attraction of $\bar{\mathbf{u}}(t)$
 stays near $\bar{\mathbf{u}}(t)$ for $n,t>0$ sufficiently large.
This conclusion indicates that
 the ``asymptotic stability'' of $\mathbf{u}_n(t)$ is observed in numerical simulations
 since they can be performed only for large values of $n,t>0$.
We will observe this behavior in {\color{black}numerical simulations} in Sections~5 and 6.


\section{Linear Stability}

We {\color{black}now turn to the CL \eqref{eqn:csys}
and first} determine the linear stability of the solution \eqref{eqn:tsol}
 to the CL \eqref{eqn:csys}.
Following the approach of Section~3 in \cite{Y24c},
 we consider the eigenvalue problem
 for the linear operator $\L:L^2(I)\to L^2(I)$ given by
\begin{align*}
\L\phi(x)
=& \int_I W(x,y)\cos(2\pi q(y-x)+\sigma)(\phi(y)-\phi(x))\d y-b\phi(x)\notag\\
=& \int_{x-\kappa}^{x+\kappa}\cos(2\pi q(y-x)+\sigma)\phi(y)\d y
-\left(\frac{\cos\sigma\sin 2\pi q\kappa}{\pi q}+b_1\right)\phi(x)
\end{align*}
for the linearization of \eqref{eqn:csys} around $u(t,x)=2\pi q x$.

Obviously, $\phi(x)=1$ is an eigenfunction for the eigenvalue $\lambda=-b_1<0$.
Moreover,  if $\sigma=0$, then
\[
\phi(x)=\cos 2 \pi\ell x,\
\sin 2\pi\ell x
\]
are eigenfunctions for the eigenvalue
\[
\lambda=\chi_1(\ell, q)-b_1
\]
for each $\ell\in\Nset$, and if $\sigma\neq 0$, then
\[
\phi(x)=\cos 2 \pi\ell x\pm i\sin 2\pi\ell x
\]
are eigenfunctions for the eigenvalue
\[
\lambda=\chi_1(\ell, q)\cos\sigma-b_1\mp i\chi_2(\ell,q)\sin\sigma
\]
for each $\ell\in\Nset$, where the upper or lower signs are taken simultaneously,
\begin{align*}
\chi_1(\ell, q)=\begin{cases}
\displaystyle
\kappa+\frac{\sin 4\pi q\kappa}{4\pi q}-\frac{\sin 2\pi q\kappa}{\pi q} & \mbox{if $\ell= q$;}\\[2ex]
\displaystyle
\frac{\sin 2\pi(\ell-q)\kappa}{2\pi(\ell-q)}+\frac{\sin 2\pi(\ell+q)\kappa}{2\pi(\ell+q)}
 -\frac{\sin 2\pi q\kappa}{\pi q} & \mbox{otherwise}
\end{cases}
\end{align*}
and
\begin{align*}
\chi_2(\ell, q)=\begin{cases}
\displaystyle
\kappa-\frac{\sin 4\pi q\kappa}{4\pi q} & \mbox{if $\ell= q$;}\\[2ex]
\displaystyle
\frac{\sin 2\pi(\ell-q)\kappa}{2\pi(\ell-q)}-\frac{\sin 2\pi(\ell+q)\kappa}{2\pi(\ell+q)}
  & \mbox{otherwise}.
\end{cases}
\end{align*}
See also Section~3 of \cite{Y24c}.
These eigenvalues are the only ones of $\L$
 since the Fourier expansion of any function in $L^2(I)$ converges a.e.
 by Carleson's theorem \cite{C66}.
Thus, if
\begin{equation}
b_1>\chi_1(\ell, q)\cos\sigma
\label{eqn:ls}
\end{equation}
for any $\ell\in\Nset$, then the $q$-twisted solution \eqref{eqn:tsol} is linearly stable.
In addition, if
\begin{equation}
b_1=\chi_1(\ell, q)\cos\sigma
\label{eqn:ezero}
\end{equation}
then $\L$ has a zero eigenvalue of geometric multiplicity two when $\sigma=0$
 and a pair of purely imaginary eigenvalues
 when $\sigma\neq 0$ and $\chi_2(\ell,q)\neq 0$,
 so that a bifurcation may occur in the CL \eqref{eqn:csys}.

Let
\[
\varphi(\zeta)=\frac{2\sin\zeta}{\zeta}-\frac{\sin 2\zeta}{2\zeta}
=\frac{\sin\zeta}{\zeta}(2-\cos\zeta).
\]
We show that the equation $\varphi(\zeta)=1$ has a unique root at
\[
\zeta_0=2.1391\ldots
\]
in $(0,\pi)$ (see Section~3 of \cite{Y24c}).
The following properties on $\chi_1(\ell,q)$ hold,
 as proven in Proposition~3.2 of \cite{Y24c}.

\begin{prop}\
\label{prop:3a}
\begin{enumerate}
\setlength{\leftskip}{-1.6em}
\item[\rm(i)]
$\chi_1(q, q)=\tfrac{1}{2}$
 and $\chi_1(\ell,q)=0$ for $\ell\neq q$ at $\kappa=\tfrac{1}{2}$,
 while $\chi_1(\ell, q)\to0$ as $\kappa\to+0$ for any $\ell\in\Nset$.
\item[\rm(ii)]
$\chi_1( q, q)<0$ for $\kappa\in(0,\kappa_q)$
 and $\chi_1( q, q)>0$ for  $\kappa\in(\kappa_q,\tfrac{1}{2})$,
 where $\kappa_q=\zeta_0/2\pi q$.
\item[\rm(iii)]
For any $\ell\in\Nset$, $\chi_1(\ell, q)<0$ when $\kappa>0$ is sufficiently small.
\item[\rm(iv)]
If $\ell\ge 2q$ and $\kappa\le\kappa_q$,
 then $\chi_1(\ell, q)<\chi_1(q, q)$.
\end{enumerate}
\end{prop}

\begin{rmk}\
\label{rmk:3a}
\begin{enumerate}
\setlength{\leftskip}{-1.8em}
\item[\rm(i)]
From Proposition~{\rm\ref{prop:3a}(i)} we see that for each $q\in\Nset$,
 the $q$-twisted solution \eqref{eqn:tsol} is linearly stable
 if $b_1>\chi(\kappa;q,q)$, near $\kappa=\tfrac{1}{2}$,
 and especially at $\kappa=\tfrac{1}{2}$,
 i.e., when the graph $G_n$ is complete simple.
\item[\rm(ii)]
It follows from Proposition~{\rm\ref{prop:3a}(ii)-(iv)} that for $q\in\Nset$
 the $q$-twisted solution \eqref{eqn:tsol} is asymptotically stable for $b_1>0$
 when $\kappa>0$ is sufficiently small
 or when $\kappa<\kappa_q$ and $\chi_1(\kappa;\ell,q)<0$ for $\ell<2q$.
\item[\rm(iii)]
When $b_1=0$ and $\kappa=\tfrac{1}{2}$,
 the $q$-twisted solution \eqref{eqn:tsol} is unstable
 for any $q\in\Nset$, by Proposition~{\rm\ref{prop:3a}(i)}.
\end{enumerate}
\end{rmk}

\begin{figure}
\includegraphics[scale=0.54]{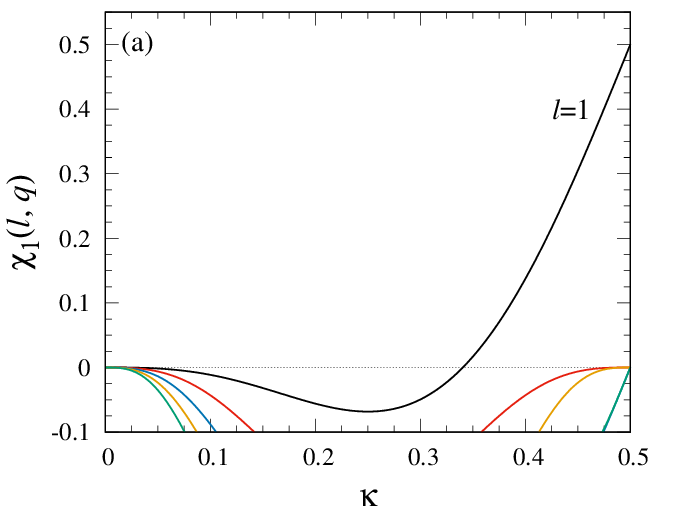}\
\includegraphics[scale=0.54]{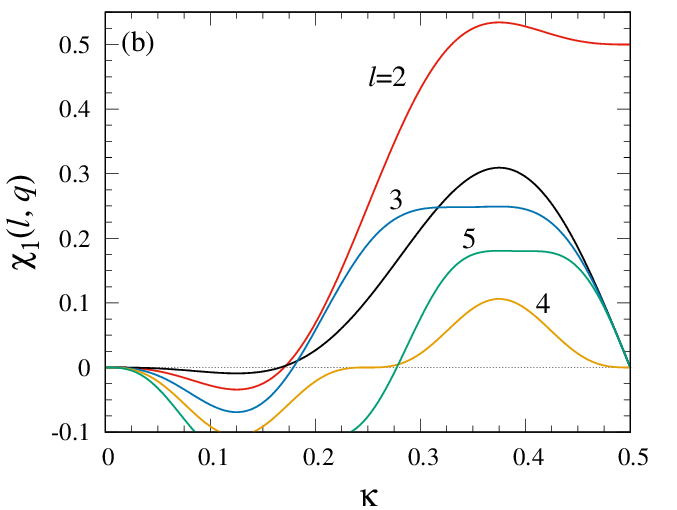}\\[1ex]
\includegraphics[scale=0.54]{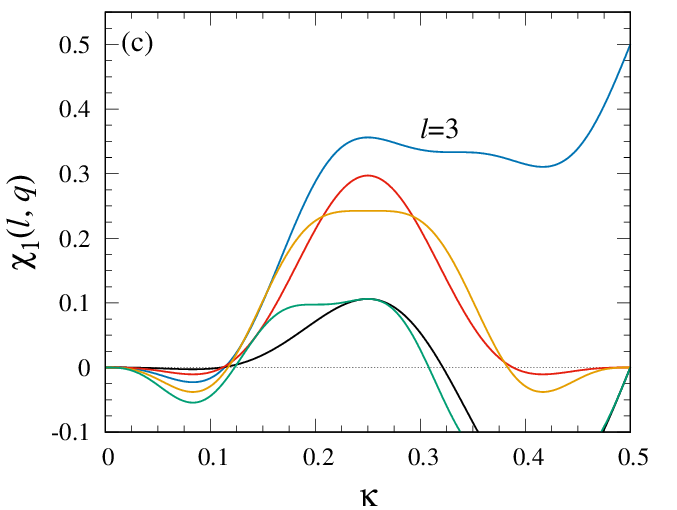}\
\includegraphics[scale=0.54]{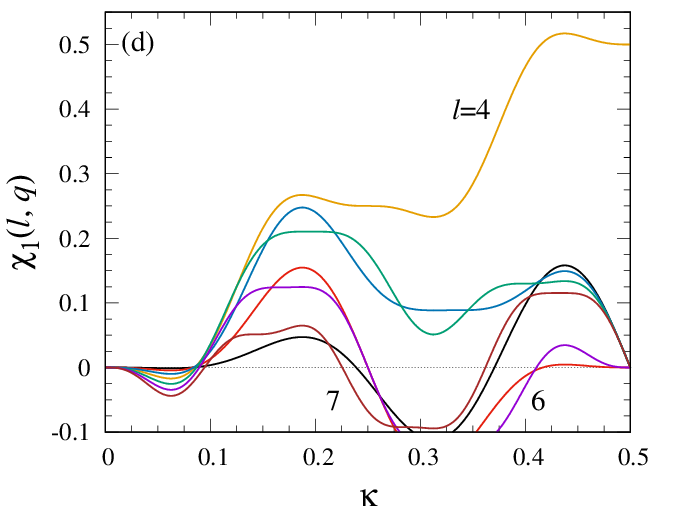}
\caption{Dependence of $\chi_1(l,q)$ on $\kappa$ for $l=1$-$7$:
(a) $q=1$; (b) $q=2$; (c) $q=3$; (d) $q=4$.
It is plotted as the line of which color is black for $l=1$, red for $l=2$, blue for $l=3$,
 orange for $l=4$, green for $l=5$, purple for $l=6$ and brown for $l=7$.
\label{fig:3a}}
\end{figure}

Figure~\ref{fig:3a}
 displays the dependence of $\chi_1(\ell,q)$ on $\kappa$
 for $q\in[4]$ and some values of $\ell$.
{\color{black}If we choose such a sufficient large value of $b_1$
 that Eq.~\eqref{eqn:ls} holds, i.e., $\chi_1(\ell,q)<b_1$, for all $\ell\in\Nset$,
 then $q$-twisted solution \eqref{eqn:tsol} is asymptotically stable.
In particular, we can stabilize the $q$-twisted solution \eqref{eqn:tsol}
 for complete simple graphs.}
We {\color{black}also} observe that $\chi_1( q, q)\to\tfrac{1}{2}$
 and $\chi_1(\ell,q)\to0$ for $\ell\neq q$ as $\kappa\to\tfrac{1}{2}-0$
 while $\chi_1(\ell,q)<0$ near $\kappa=0$,
 as stated in Proposition~\ref{prop:3a}(i) and (iii).
For $q\in[4]$,
 the $q$-twisted solution \eqref{eqn:tsol} is asymptotically stable
 if $b_1>0$ and $\kappa<\kappa_q$, by Remark~\ref{rmk:3a}(ii).
Moreover, $\chi_1(q,q)>\chi_1(\ell,q)$ in wide ranges of $\kappa$ for $q\in[4]$.
In particular, $\chi_1(1,1)>\chi_1(\ell,1)$ on $\kappa\in(0,\tfrac{1}{2})$ for any $\ell>1$,
 by Proposition~\ref{prop:3a}(iv).

Since $\chi_2(q,q)$ is monotonically increasing
 as a function of $\kappa$ on $(0,\tfrac{1}{2}]$,
 we easily prove the following.

\begin{prop}
\label{prop:3b}
$\chi_1(q,q)\to0$ as $\kappa\to+0$,
 $\chi_2(q, q)=\tfrac{1}{2}$ at $\kappa=\tfrac{1}{2}$,
 and $\chi_2(q, q)>0$ for $\kappa\in(0,\tfrac{1}{2}]$.
\end{prop}


\section{Bifurcations}

We now take the linear feedback gain $b_1$ as a control parameter
 and analyze bifurcations of the $q$-twisted solution \eqref{eqn:tsol}
 in the CL \eqref{eqn:csys} for $q\in[4]$.
Our approach is similar to that of \cite{Y24c}
 for the uncontrolled case of $b_1,b_3=0$
 but some modifications are required.

Let $b_{1q}$ denote the value of $b_1$
 satisfying \eqref{eqn:ezero} for $\ell=q$, i.e.,
\begin{equation}
b_{1q}=\chi_1(q,q)\cos\sigma.
\label{eqn:b1q}
\end{equation}
From the analysis of Section~3 we see that
 a bifurcation may occur at $b_1=b_{1q}$ for each $q\in\Nset$.
In the following
 we analyze the bifurcation which may occur at \eqref{eqn:b1q} for $q\in[4]$.
We assume for each $q\in[4]$ that $b_1\approx b_{1q}$
 and condition \eqref{eqn:ls} holds for $\ell\neq q$,
 and introduce a parameter $\mu=b_1-b_{1q}\approx 0$.
Moreover, we write solutions to the CL \eqref{eqn:csys}
 near the $q$-twisted solution \eqref{eqn:tsol} as
\begin{equation}
u(t,x)=2\pi qx+\Omega t+\xi_0(t)+\sum_{j=1}^\infty(\xi_j(t)\cos 2\pi jx+\eta_j(t)\sin 2\pi jx)
\label{eqn:solex}
\end{equation}
and regard $\mu$ as a state variable.

\subsection{Center manifold reduction}
We substitute \eqref{eqn:solex} into \eqref{eqn:csys}
 and integrate the resulting equation from $x=0$ to $1$ directly
 or after multiplying it with $\cos 2\pi jx$ or $\sin 2\pi jx$, $j\in\Nset$, to obtain
\begin{equation}
\begin{split}
\dot{\xi}_q=& -\mu\xi_q-\nu_q\eta_q
 -(\tfrac{3}{4}b_3+\beta_1\cos\sigma)(\xi_q^2+\eta_q^2)\xi_q
+\delta_1\cos\sigma(\xi_q\eta_{2q}-\xi_{2q}\eta_q)\\
&
+\sin\sigma(-\beta_2(\xi_q^2+\eta_q^2)\eta_q
 +\delta_2(\xi_q\xi_{2q}+\eta_q\eta_{2q}))+\cdots,\\
\dot{\eta}_q=&\nu_q\xi_q-\mu\eta_q
-(\tfrac{3}{4}b_3+\beta_1\cos\sigma)(\xi_q^2+\eta_q^2)\eta_q
 -\delta_1\cos\sigma(\xi_q\xi_{2q}+\eta_q\eta_{2q})\\
& +\sin\sigma(\beta_2(\xi_q^2+\eta_q^2)\xi_q
 +\delta_2(\xi_q\eta_{2q}-\xi_{2q}\eta_q))+\cdots,\\
\dot{\xi}_{2q}=&\mu_{2q}\xi_{2q}-\nu_{2q}\eta_{2q}
 -2\rho_1\xi_q\eta_q\cos\sigma
 +\rho_2\sin\sigma(\xi_q^2-\eta_q^2)+\cdots,\\
\dot{\eta}_{2q}=&\nu_{2q}\xi_{2q}+\mu_{2q}\eta_{2q}
 +\rho_1\cos\sigma(\xi_q^2-\eta_q^2)
 +2\rho_2\xi_q\eta_q\sin\sigma+\cdots,\\
\dot{\xi}_j=&\mu_j\xi_j-\nu_j\eta_j+\cdots,\quad
\dot{\eta}_j=\nu_j\xi_j+\mu_j\eta_j+\cdots,\quad
j\in\Nset\setminus\{q,2q\}\\
\dot{\xi}_0=&-b_{1q}\xi_0+\cdots,\quad
\dot{\mu}=0,
\end{split}
\label{eqn:ifex}
\end{equation}
for $q\in[4]$,  where  `$\cdots$' represents higher-order terms of
\[
O\left(\xi_q^4+\eta_q^4+\xi_0^2+\sum_{j=1,j\neq q}^\infty(\xi_j^2+\eta_j^2)+\mu^2\right)
\]
for the first and second equations and
\[
O\left(\xi_0^2+\sum_{j=1}^\infty(\xi_j^2+\eta_j^2)+\mu^2\right)
\]
for the other equations, and
\begin{align*}
&
\beta_1=\tfrac{3}{8}a_2(q,0)-\tfrac{1}{2}a_2(q,q)+\tfrac{1}{8}a_2(q,2q),\quad
\beta_2=\tfrac{1}{4}a_1(q,q)-\tfrac{1}{8}a_1(q,2q),\\
&
\delta_1=a_1(q,q)-\tfrac{1}{2}a_1(q,2q),\quad
\delta_2=\tfrac{1}{2}a_2(q,0)-\tfrac{1}{2}a_2(q,2q),\\
&
\rho_1=\tfrac{1}{2}a_1(q,q)-\tfrac{1}{4}a_1(q,2q),\quad
\rho_2=\tfrac{1}{4}a_2(q,0)-\tfrac{1}{2}a_2(q,q)+\tfrac{1}{4}a_2(q,2q),\\
&
\mu_j=-b_{1q}+\chi_1(j,q)\cos\sigma,\quad
j\in\Nset\setminus\{1\},\\
&
\nu_j=\chi_2(j,q)\sin\sigma,\quad
j\in\Nset,
\end{align*}
with
\begin{align*}
&
a_1(q,j)=\begin{cases}
\displaystyle
\frac{\sin(4\pi q\kappa)}{4\pi q}-\kappa
& \mbox{for $j=q$};\\[2ex]
\displaystyle
\frac{q\sin(2\pi j\kappa)\cos(2\pi q\kappa)-j\cos(2\pi j\kappa)\sin(2\pi q\kappa)}{\pi(q^2-j^2)}
& \mbox{for $j\neq q$},
\end{cases}\\
&
a_2(q,j)=\begin{cases}
\displaystyle
-\frac{\sin4\pi q\kappa}{4\pi q}-\kappa
& \mbox{for $j=q$};\\[2ex]
\displaystyle
\frac{j\sin(2\pi j\kappa)\cos(2\pi q\kappa)-q\cos(2\pi j\kappa)\sin(2\pi q\kappa)}{\pi(q^2-j^2)}
& \mbox{for $j\neq q$}.
\end{cases}
\end{align*}
See Appendix~A for the derivation of \eqref{eqn:ifex}.

Henceforth we assume that $\mu_j<0$ for any $j\neq q\in[4]$.
Actually, this assumption holds near $\kappa=\tfrac{1}{2}$ by Remark~\ref{rmk:3a}(i),
 and in a wide range of $\kappa$ containing $(0,\kappa_q]$
 as seen from Fig.~\ref{fig:3a} and Proposition~\ref{prop:3a}(iv).
The origin in the infinite-dimensional system \eqref{eqn:ifex} is an equilibrium
 having a three-dimensional center manifold $W^\mathrm{c}$, even if $\sigma\neq 0$.
Using the standard approach \cite{GH83,HI11,K04},
 we obtain the following.

\begin{prop}
\label{prop:4a}
The center manifold is expressed as 
\begin{align*}
W^{\mathrm{c}}=\{
\xi_{2q}=&\bar{\xi}_{2q}(\xi_q,\eta_q)+O(3),
\eta_{2q}=\bar{\eta}_{2q}(\xi_q,\eta_q)+O(3),\\
&
\xi_0=O(3),\xi_j=O(3),\eta_j=O(3),j\neq q,2q\}
\end{align*}
near the origin, where $O(k)$ represents higher-order terms
 of $O\left(\sqrt{\xi_q^{2k}+\eta_q^{2k}+\mu^4}\right)$, and
\[
\bar{\xi}_{2q}(\xi_q,\eta_q)
=c_1(\xi_q^2-\eta_q^2)+2c_2\xi_q\eta_q,\quad
\bar{\eta}_{2q}(\xi_q,\eta_q)
=-c_2(\xi_q^2-\eta_q^2)+2c_1\xi_q\eta_q
\]
with
\begin{align*}
&
c_1
=\frac{(2\nu_q-\nu_{2q})\rho_1\cos\sigma-\mu_{2q}\rho_2\sin\sigma}
{\mu_{2q}^2+(2\nu_q-\nu_{2q})^2},\\
&
c_2
=\frac{\mu_{2q}\rho_1\cos\sigma-(2\nu_q-\nu_{2q})\rho_2\sin\sigma}
{\mu_{2q}^2+(2\nu_q-\nu_{2q})^2}.
\end{align*}
\end{prop}

Based on Proposition~\ref{prop:4a},
 we apply the center manifold reduction \cite{HI11} to \eqref{eqn:ifex}.
 and obtain
\begin{equation}
\begin{split}
\dot{\xi}_q=& -\mu\xi_q-\nu_q\eta_q\\
& -(\tfrac{3}{4}b_3+\beta_1\cos\sigma)(\xi_q^2+\eta_q^2)\xi_q
 +\delta_1\cos\sigma(\bar{\eta}_{2q}(\xi_q,\eta_q)\xi_q-\bar{\xi}_{2q}(\xi_q,\eta_q)\eta_q)\\
&
+\sin\sigma(-\beta_2(\xi_q^2+\eta_q^2)\eta_q
 +\delta_2(\bar{\xi}_{2q}(\xi_q,\eta_q)\xi_q+\bar{\eta}_{2q}(\xi_q,\eta_q)\eta_q))+O(4),\\
\dot{\eta}_q=&\nu_q\xi_q-\mu\eta_q\\
& -(\tfrac{3}{4}b_3+\beta_1\cos\sigma)(\xi_q^2+\eta_q^2)\eta_q
 +\delta_1\cos\sigma(\bar{\xi}_{2q}(\xi_q,\eta_q)\xi_q+\bar{\eta}_{2q}(\xi_q,\eta_q)\eta_q)\\
& +\sin\sigma(\beta_2(\xi_q^2+\eta_q^2)\xi_q
 +\delta_2(\bar{\eta}_{2q}(\xi_q,\eta_q)\xi_q-\bar{\xi}_{2q}(\xi_q,\eta_q)\eta_q))+O(4),\\
\dot{\mu}=&0
\end{split}
\label{eqn:cmex}
\end{equation}
on $W^{\mathrm{c}}$.
See Appendix~B of \cite{Y24b}
 for the validity of application of the center manifold theory
 on infinite-dimensional dynamical systems \cite{HI11}.
The origin $(\xi_q,\eta_q,\mu)=(0,0,0)$ is always an equilibrium in \eqref{eqn:cmex}.
This is because the twisted solution \eqref{eqn:tsol}
 necessarily satisfies the CL \eqref{eqn:csys}. 

\subsection{Case of $\sigma=0$}

We set $\sigma=0$,
 so that $c_1=0$ and $c_2=-\rho_1/\mu_{2q}$.
We remark that
\[
\mu_{2q}=-b_{1q}+\chi_1(2q,q)=-\kappa+\frac{\sin 2\pi q\kappa}{2\pi q}<0
\]
for $\kappa>0$.
Letting $r=\sqrt{\xi_q^2+\eta_q^2}\ge 0$,
 we rewrite \eqref{eqn:cmex}.as
\begin{equation}
\dot{r}=-\mu r-\beta_0 r^3
 +O(\sqrt{r^8+\mu^4}),\quad
\dot{\mu}=0,
\label{eqn:r}
\end{equation}
where
\begin{equation}
\beta_0
=\tfrac{3}{4}b_3+\bar{\beta}_1,\quad
\bar{\beta}_1=\beta_1+\frac{\delta_1\rho_1}{\mu_{2q}}.
\label{eqn:beta1}
\end{equation}
In particular, $\bar{\beta}_1=0$ when $\kappa=\tfrac{1}{2}$.
Here by the translation symmetry {\color{black}(see Eq.~\eqref{eqn:csys1})},
 the first equation of \eqref{eqn:r} must depend only on $r$ and $\mu$,
  even if the higher-order terms are included.
We take such a sufficiently large value for $b_3$ as $\beta_0>0$.
We easily show the following for \eqref{eqn:r}:

\begin{figure}
\includegraphics[scale=0.55]{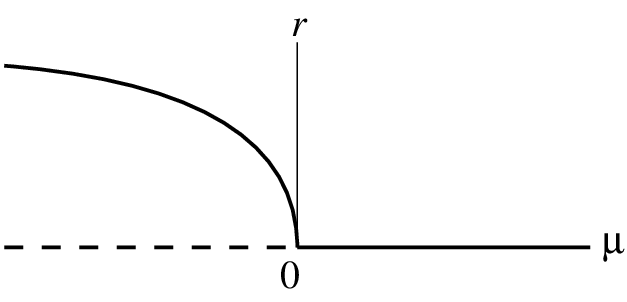}
\caption{Bifurcation diagrams for \eqref{eqn:r}.
\label{fig:4a}}
\end{figure}

\begin{enumerate}
\setlength{\leftskip}{-1.8em}
\item[(i)]
The equilibrium $r=0$ is stable for $\mu>0$ and unstable for $\mu<0$;
\item[(ii)]
There exists another stable equilibrium at
\[
r=\sqrt{-\frac{\mu}{\beta_0}}
\]
for $\mu<0$.
\end{enumerate}
See Fig.~\ref{fig:4a} for the bifurcation diagrams for \eqref{eqn:r}.
From this result,
 we obtain the following for the CL \eqref{eqn:csys}.

\begin{thm}
\label{thm:4a}
Let $q\in[4]$ and suppose that $\beta_0>0$, i.e., $b_3>-\tfrac{4}{3}\bar{\beta}_1$,
 and $\mu_j=-b_{1q}+\chi_1(j,q)<0$ for any $j\neq q$.
Then the following bifurcation of the twisted solution \eqref{eqn:tsol} occurs
 at $b_1=b_{1q}$ in the CL \eqref{eqn:csys} with $\sigma=0:$
\begin{enumerate}
\setlength{\leftskip}{-1.8em}
\item[(i)]
The twisted solution \eqref{eqn:tsol}
 is stable for $b_1>b_{1q}$ and unstable for $b_1<b_{1q};$
\item[(ii)]
There exists a stable one-parameter family of modulated twisted solutions
\begin{align}
\U^q
=\biggl\{u=2\pi qx+\sqrt{-\frac{b_1-b_{1q}}{\beta_0}}&\sin(2\pi qx+\psi)\notag\\
+&\Omega t+O(b_1-b_{1q})\,\bigg|\,\psi\in\Sset^1\biggr\}
\label{eqn:thm4a}
\end{align}
for $b_1<b_{1q}$ near  $b_1=b_{1q}$,
 where $\beta_0=O(1)$ is given in \eqref{eqn:beta1}.
\end{enumerate}
\end{thm}
 
\begin{rmk}\
\label{rmk:4a}
\begin{enumerate}
\setlength{\leftskip}{-1.6em}
\item[\rm(i)]
A bifurcation similar to one detected in Theorem~$\ref{thm:4a}$
 also occurs at $b_1=b_{1q}$ even if $b_{3}<-\tfrac{4}{3}\bar{\beta}_1$
 or $\mu_j>0$ for some $j\neq q$,
 although the one-parameter family $\U^q$ of modulated twisted solutions born there
 is unstable.
\item[\rm(ii)]
We suspect for any $q\in\Nset$ that Eq.~\eqref{eqn:ifex} is valid
 and the statements of Theorem~$\ref{thm:4a}$ also hold.
However, it is very hard to derive \eqref{eqn:ifex} for any $q\in\Nset$,
 so that we restrict ourselves to $q\in[4]$ in the above analysis.
\item[\rm(iii)]
{\color{black}
Noting the relation \eqref{eqn:conv} and using the theory of Section~$2$
 $($see Theorems~{\rm\ref{thm:2c}(ii), {\rm\ref{thm:2d}(ii)} and $\ref{thm:2e}$,
 Corollary~$\ref{cor:2a}$ and Remarks~{\rm\ref{rmk:2b}} and {\rm\ref{rmk:2c})},
 we see that the target orbit \eqref{eqn:ts}} $($resp. $\U^q)$ behaves
 as if it is an asymptotically stable solution or it is actually unstable
 $($resp. as if it is an asymptotically stable family of solutions$)$
 in the KM~\eqref{eqn:dsys} near $b_1=b_{1q}$ for $n>0$ sufficiently large.}
Thus, the KM \eqref{eqn:dsys} suffers a ``bifurcation''
 similar to one detected in Theorem~$\ref{thm:4a}$ for the CL \eqref{eqn:csys}.
\end{enumerate}
\end{rmk}

\begin{figure}
\includegraphics[scale=0.6]{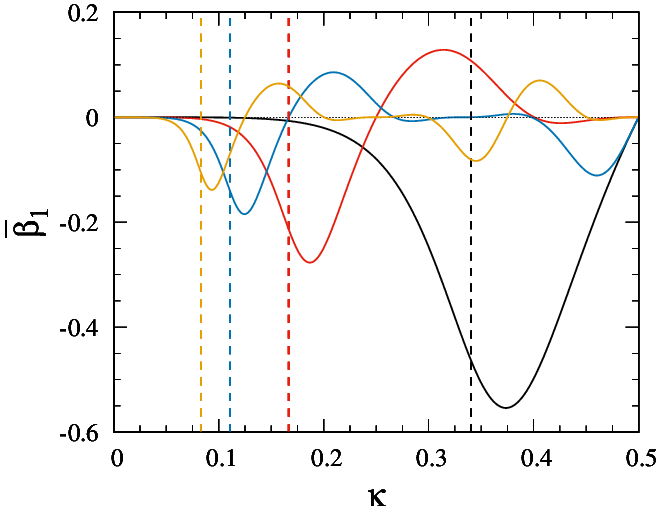}
\caption{Dependence of $\bar{\beta}_1$ on $\kappa$.
The black, red, blue and orange lines
 represent the cases of $q=1,2,3$ and $4$, respectively,
The dashed lines with the same colors
 represent $\kappa=\kappa_q$, on which $b_{1q}=0$, for $q\in[4]$.
\label{fig:4b}}
\end{figure}

Figure~\ref{fig:4b} displays the dependence of $\bar{\beta}_1$ given by \eqref{eqn:beta1}
 on $\kappa$ for $q\in[4]$.
In particular, $\bar{\beta}_1=0$ at $\kappa=\tfrac{1}{2}$.
{\color{black}
Theorem~\ref{thm:4a} requires $b_3>-\tfrac{4}{3}\bar{\beta}_1$
 and $b_{1q}=\chi_1(q,q)>\chi_1(j,q)$ for any $j\neq q$ as its hypotheses.
The nonlinear feedback gain $b_{3}$ has to be positive at least
 for the former to hold in the range of $\kappa$ where $\bar{\beta}_1$ is negative,
 especially near $\kappa=\tfrac{1}{2}$ and for $q=1$.
From Fig.~\ref{fig:3a} and Proposition~\ref{prop:3a}(iv)
 we see that the latter holds for $q\in[4]$ when $\kappa\le\kappa_q$.}

\subsection{Case of $\sigma\neq 0$}
 
We next consider the case of $\sigma\neq 0$.
Letting $\xi_q=r\cos\psi$ and $\eta_q=r\sin\psi$,
 we rewrite \eqref{eqn:cmex} as
\begin{equation}
\begin{split}
\dot{r}=-\mu r-\beta_\sigma r^3+O(\sqrt{r^8+\mu^4}),\quad
\dot{\psi}=\nu_q+O(\sqrt{r^2+\mu^2}),\quad
\dot{\mu}=0,
\end{split}
\label{eqn:polar}
\end{equation}
where
\begin{equation}
\beta_\sigma
=\tfrac{3}{4}b_3+\bar{\beta}_{1 \sigma}
\label{eqn:beta1s}
\end{equation}
with
\begin{align}
\bar{\beta}_{1\sigma}
=&\beta_1\cos\sigma+\frac{1}{2(\mu_{2q}^2+(2\nu_q-\nu_{2q})^2)}
 (\mu_{2q}(\delta_1\rho_1+\delta_2\rho_2)\notag\\
& +\mu_{2q}(\delta_1\rho_1-\delta_2\rho_2)\cos2\sigma
 +(2\nu_q-\nu_{2q})(\delta_1\rho_2-\delta_2\rho_1)\sin2\sigma).
\label{eqn:beta1s'}
\end{align}
Here by the translation symmetry,
 Eq.~ \eqref{eqn:polar} must depend only on $r$ and $\mu$,
 even if the higher-order terms are included, {\color{black}like \eqref{eqn:r}}.
We take such a sufficiently large value for $b_3$ as $\beta_\sigma>0$.
Noting that $\nu_q>0$ for $\kappa>0$ by Proposition~\ref{prop:3b},
 we easily show that a Hopf bifurcation \cite{GH83,HI11,K04} occurs
 in \eqref{eqn:polar} as follows {\color{black}(cf. Fig.~\ref{fig:4a})}:

\begin{enumerate}
\setlength{\leftskip}{-1.8em}
\item[(i)]
The equilibrium $r=0$ is stable for $\mu>0$ and unstable for $\mu<0$;
\item[(ii)]
There exists a stable periodic orbit given by
\begin{equation}
r=\sqrt{-\frac{\mu}{\beta_\sigma}}+O(\mu),\quad
\psi=\nu_q t+O(\sqrt{\mu}),
\label{eqn:ic2}
\end{equation}
for $\mu<0$.
\end{enumerate}
From this result,
 we obtain the following for the CL \eqref{eqn:csys}.

\begin{thm}
\label{thm:4b}
Let $q\in[4]$ and suppose that $\beta_\sigma>0$,
 i.e., $b_3>-\tfrac{4}{3}\bar{\beta}_{1\sigma}$,
 and $\mu_j=-b_{1q}+\chi_1(j,q)\cos\sigma<0$ for any $j\neq q$.
Then the following bifurcation of the twisted solution \eqref{eqn:tsol} occurs
 at $b_1=b_{1q}$ in the CL \eqref{eqn:csys} with $\sigma\neq 0:$
\begin{enumerate}
\setlength{\leftskip}{-1.8em}
\item[(i)]
The twisted solution \eqref{eqn:tsol} is stable for $b_1>b_{1q}$
 and unstable for $b_1<b_{1q};$
\item[(ii)]
There exists a stable one-parameter family
 of oscillating twisted solutions
\begin{align}
\tilde{\U}^q
=\biggl\{u=2\pi qx+\sqrt{-\frac{b_1-b_{1q}}{\beta_\sigma}}&
 \sin(2\pi qx+\tilde{\psi}(t)+\psi)\notag\\
+&\Omega t+O(b_1-b_{1q})
 \,\Big|\,\psi\in\Sset^1\biggr\}
\label{eqn:thm4b}
\end{align}
for $b_1<b_{1q}$ near  $b_1=b_{1q}$,
 where $\tilde{\psi}(t)\in\Sset^1$ is a periodic function
  whose period is approximately $2\pi/\nu_q$.
Here $\Omega$ and $\beta_\sigma=O(1)$ are given
 in \eqref{eqn:Omega} and \eqref{eqn:beta1s}, respectively.
\end{enumerate}
\end{thm}
 
\begin{rmk}\
\label{rmk:4c}
\begin{enumerate}
\setlength{\leftskip}{-1.6em}
\item[\rm(i)]
A bifurcation similar to one detected in Theorem~$\ref{thm:4b}$ also occurs
 at $b_1=b_{1 q}$ even if $b_{3}<-\tfrac{4}{3}\bar{\beta}_{1\sigma}$ or $\mu_j>0$
 although the one-parameter family $\tilde{\U}^q$
 of oscillating twisted solutions born there is unstable
 $($cf. Remark~{\rm\ref{rmk:4a}(i))}.
\item[\rm(ii)]
We suspect for any $q\in\Nset$ that
 the statements of Theorem~$\ref{thm:4b}$ also hold
 $($cf. Remark~{\rm\ref{rmk:4a}(ii))}.
\item[\rm(iii)]
As in Remark~{\rm\ref{rmk:4a}(iii)},
{\color{black}
 the target orbit \eqref{eqn:ts} $($resp. $\tilde{\U}^q)$ behaves
 as if it is an asymptotically stable solution or it is actually unstable
 $($resp. as if it is an} asymptotically stable families of solutions$)$ 
 in the KM~\eqref{eqn:dsys} near $b_1=b_{1q}$ for $n>0$ sufficiently large.
Thus, the KM \eqref{eqn:dsys} suffers a ``bifurcation''
 similar to one detected in Theorem~$\ref{thm:4b}$ for the CL \eqref{eqn:csys}.
\end{enumerate}
\end{rmk}

\begin{table}
\caption{Constants appearing in Eq.~\eqref{eqn:beta1s'} for $\kappa=0.4,0.5$.
The numbers are rounded up to the fifth decimal point.	
\label{tbl:4a}}

\begin{tabular}{c|c|c|c|c|c}
\hline
$\kappa$ & \multicolumn{4}{c|}{$0.4$} & $0.5$\\
\hline
$q$  & $1$ & $2$ & $3$ & $4$ & $[4]$\\
\hline
$\displaystyle\frac{b_{1q}}{\cos\sigma}$
 & $0.13722$ & $0.52798$ & $0.31468$ & $0.46570$ & $0.5$\\[2ex]
$\beta_1$ & $0.07400$ & $0.25258$ & $0.16495$ & $0.23150$ & $0.25$\\
$\delta_1$ & $-0.45414$ & $-0.46902$ & $-0.35398$ & $-0.38647$ & $-0.5$\\
$\delta_2$& $-0.02155$ & $0.04564$ & $-0.03042$ & $0.00539$ & $0$\\
$\rho_1$ & $-0.22707$ & $-0.23451$ & $-0.17699$ & $-0.19323$ & $-0.25$\\
$\rho_2$ & $-0.12616$ & $-0.20333$ & $-0.19778$ & $-0.21846$ & $-0.25$\\
$\displaystyle\frac{\mu_{2q}+b_{1q}}{\cos\sigma}$
& $-0.04309$ & $0.09127$ & $-0.06085$ & $0.01077$ & $0$\\[2ex]
$\displaystyle\frac{\nu_q}{\sin\sigma}$
& $0.47568$ & $0.423387$ & $0.38441$ & $0.38108$ & $0.5$\\[2ex]
$\displaystyle\frac{\nu_{2q}}{\sin\sigma}$
& $0.04309$ & $-0.09127$ & $0.06085$ & $-0.01077$ & $0$\\[1.5ex]
\hline
\end{tabular}
\end{table}

\begin{figure}
\includegraphics[scale=0.6]{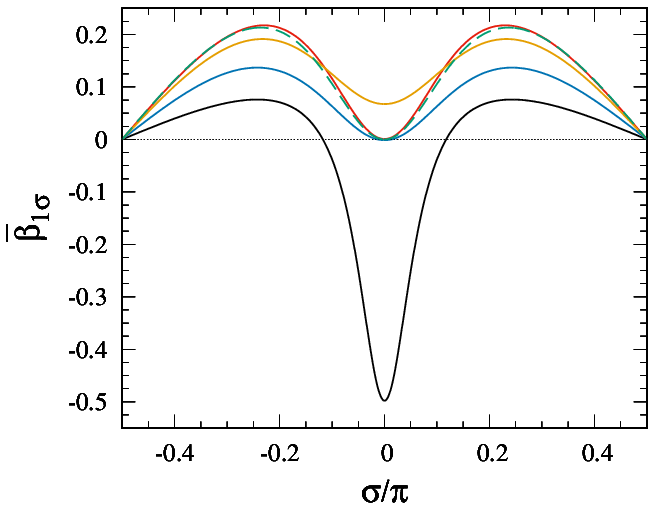}\quad
\caption{Dependence of $\bar{\beta}_{1\sigma}$ on $\sigma$ for $\kappa=0.4,0.5$:
The black, red, blue and orange lines represent the cases of $q=1,2,3$ and $4$, respectively,
 for $\kappa=0.4$,
 while the green dashed line represents the case of $q\in[4]$ for $\kappa=0.5$.
 \color{black}
 Note that $\bar{\beta}_{1\sigma}$ is independent of $q$ when $\kappa=0.5$.
\label{fig:4c}}
\end{figure}

The values of constants appearing in \eqref{eqn:beta1s'} 
 are provided in Table~\ref{tbl:4a} for $q\in[4]$ and $\kappa=0.4,0.5$.
In particular, when $\kappa=\tfrac{1}{2}$,
\[
\bar{\beta}_{1\sigma}=\frac{3\cos\sigma\sin^2\sigma}{2(4-3\cos^2\sigma)}.
\]
Figure~\ref{fig:4c} shows the dependence of $\beta_{1\sigma}$ on $\sigma$
 for $q\in[4]$ and $\kappa=0.4,0.5$.
{\color{black}
Note that 
Theorem~\ref{thm:4b} requires $b_3>-\tfrac{4}{3}\bar{\beta}_{1\sigma}$
 as its hypotheses.}
In particular, the one-parameter family $\tilde{\U}^q$ born at the bifurcation
 is stable for $\sigma\neq0,\pm\tfrac{1}{2}\pi$ if $q=2$-$4$ or $\kappa=0.5$,
 and in some range  of $\sigma$ if $q=1$ and $\kappa=0.4$,  even when $b_3=0$.
 

\section{Numerical Simulations: Nearest Neighbor Graphs}

In this and the next sections,
we give numerical simulation results for the KM \eqref{eqn:dsys} with
{\color{black}the phase-lag $\sigma=0$ or $\pi/3$
 defined on deterministic {\color{black}$\lfloor n\kappa\rfloor$}-nearest neighbor
 and complete simple graphs, respectively.
Here  we assume 
\begin{equation}
\omega=-\frac{\sin(2\pi q\kappa)\sin\sigma}{\pi q}
\label{eqn:omega}
\end{equation}
without loss of generality, considering an adequate rotational frame if necessary.}
Note that $\omega=0$ for $\sigma=0$ or $\kappa=\tfrac{1}{2}$
 and that the $q$-twisted solution \eqref{eqn:tsol} in the CL \eqref{eqn:csys}
 has $\Omega=0$ by \eqref{eqn:Omega}
 for any $\sigma\in(-\tfrac{1}{2}\pi,\tfrac{1}{2}\pi)$.

\begin{figure}[t]
\begin{minipage}[t]{0.495\textwidth}
\begin{center}
\includegraphics[scale=0.265]{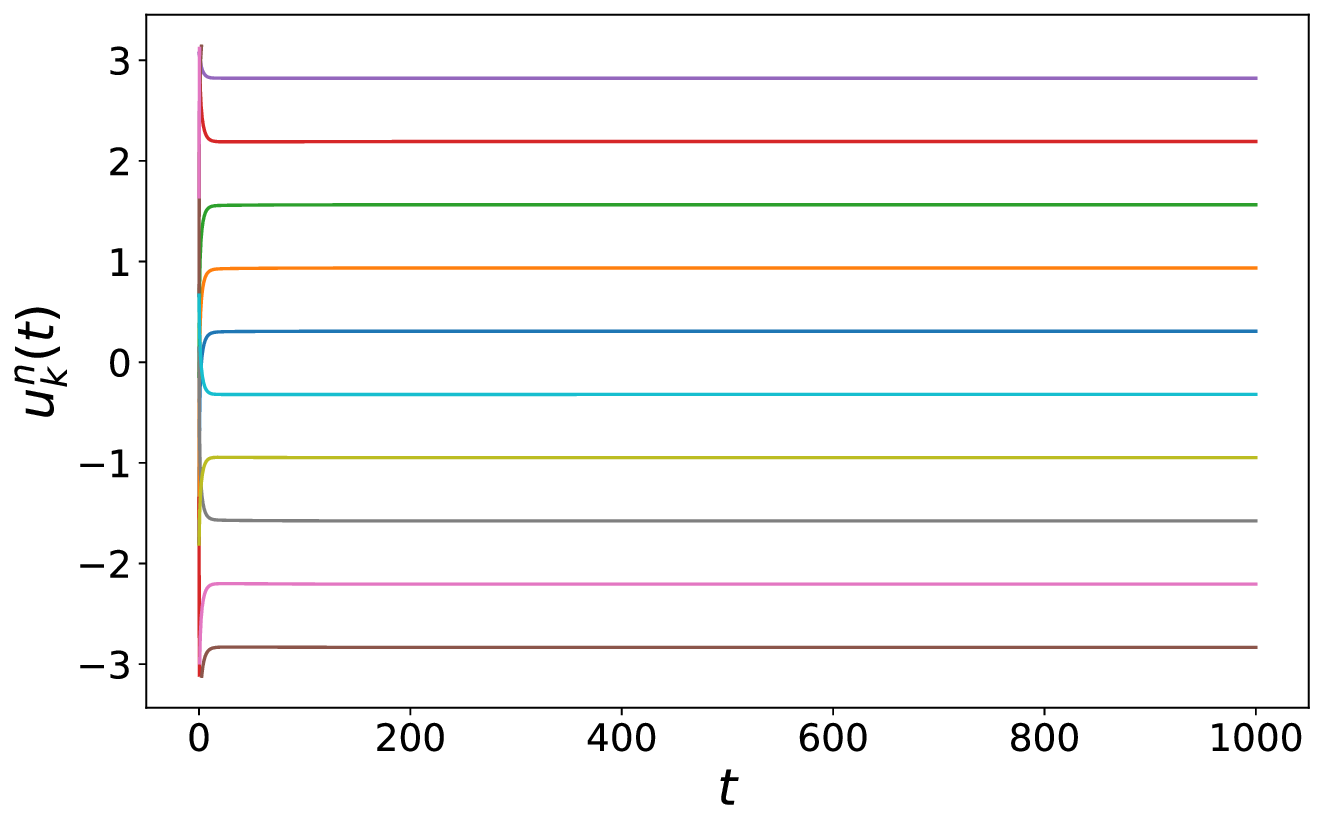}\\[-1ex]
{\footnotesize(a)}
\end{center}
\end{minipage}
\begin{minipage}[t]{0.495\textwidth}
\begin{center}
\includegraphics[scale=0.265]{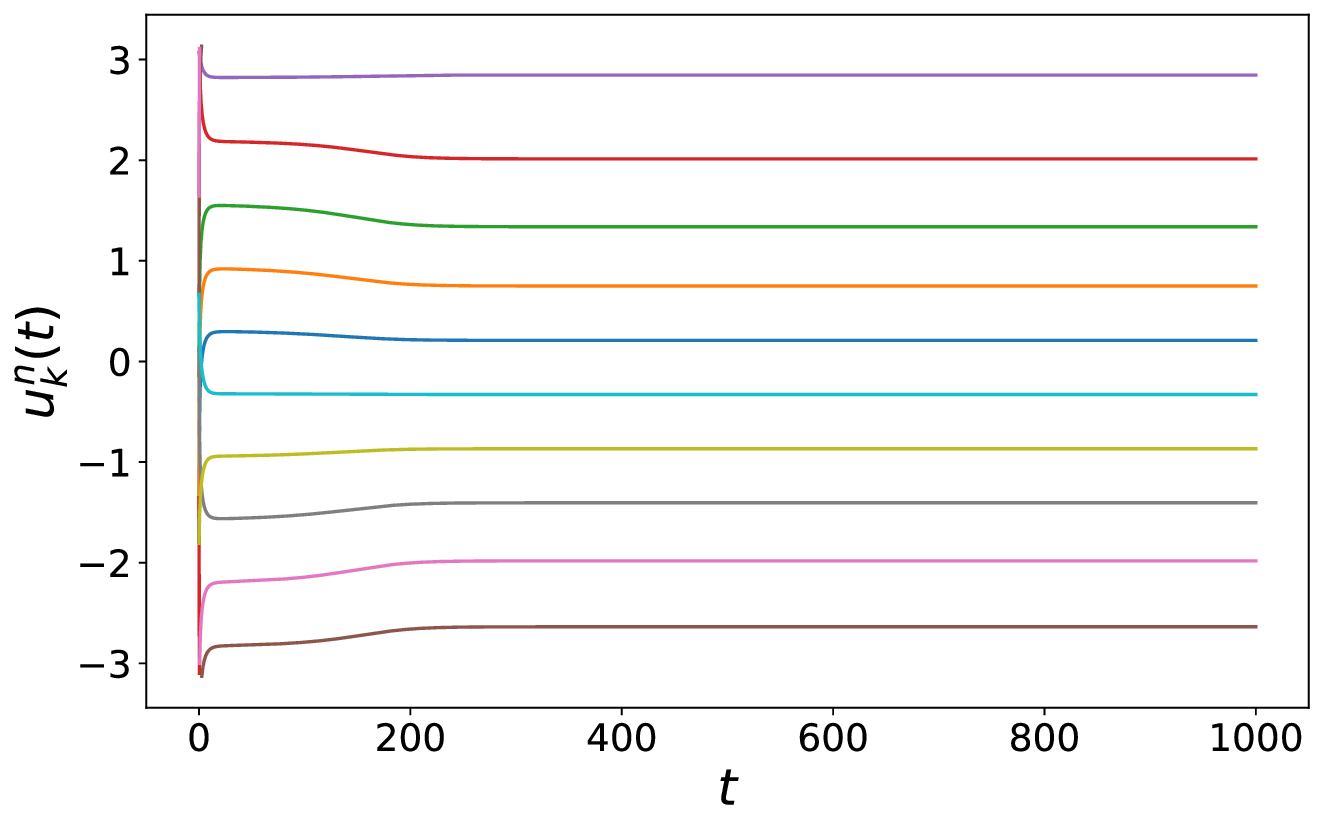}\\[-1ex]
{\footnotesize(b)}
\end{center}
\end{minipage}
\vspace*{0.5ex}

\begin{minipage}[t]{0.495\textwidth}
\begin{center}
\includegraphics[scale=0.265]{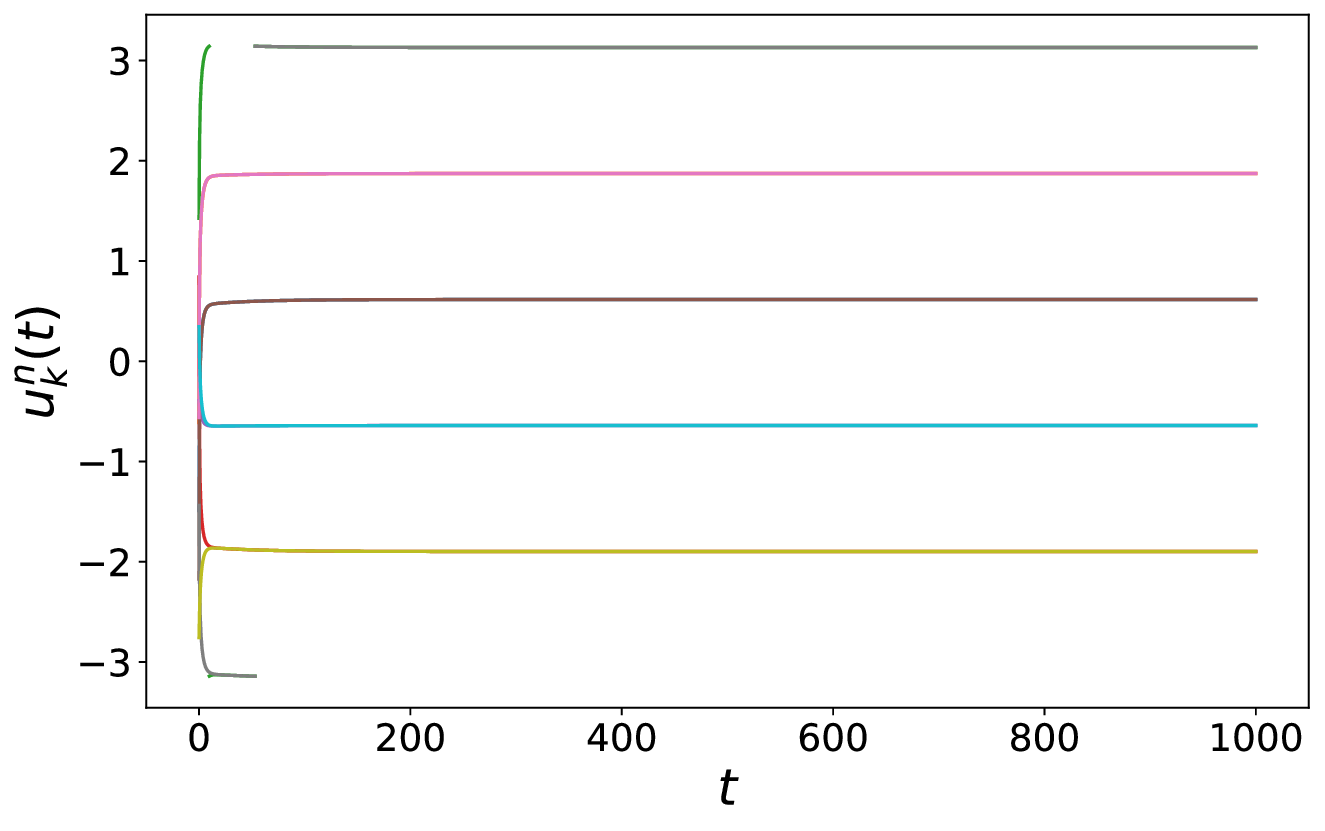}\\[-1ex]
{\footnotesize(c)}
\end{center}
\end{minipage}
\begin{minipage}[t]{0.495\textwidth}
\begin{center}
\includegraphics[scale=0.265]{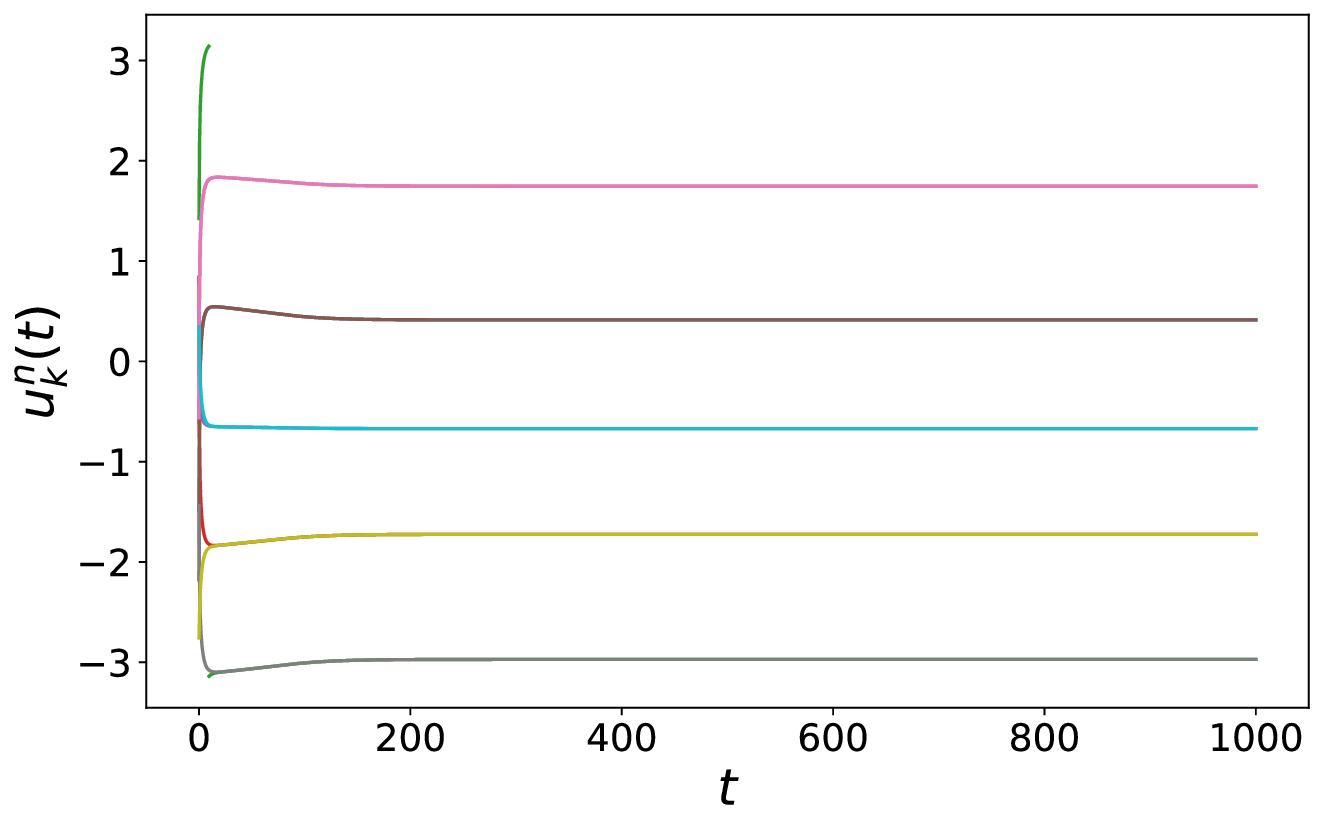}\\[-1ex]
{\footnotesize(d)}
\end{center}
\end{minipage}
\vspace*{0.5ex}

\begin{minipage}[t]{0.495\textwidth}
\begin{center}
\includegraphics[scale=0.265]{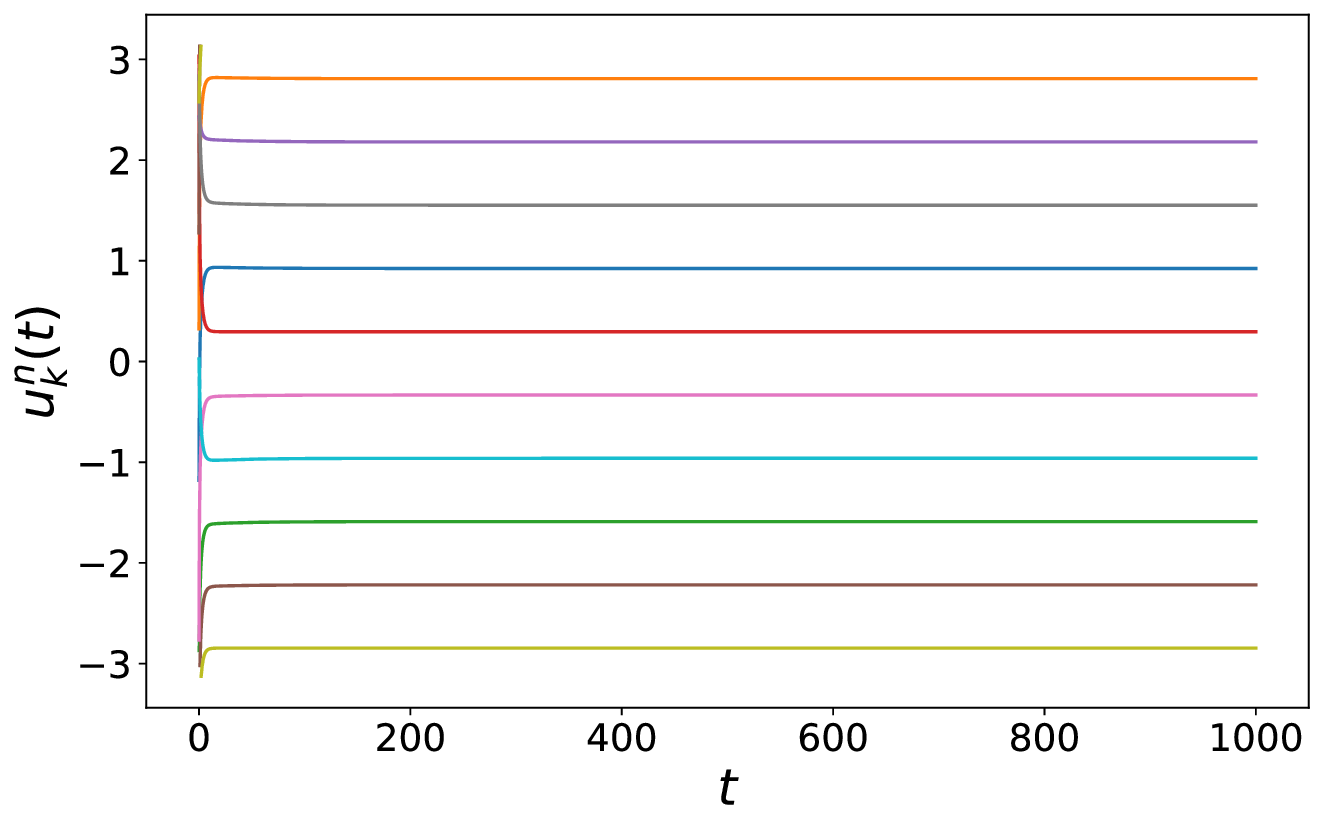}\\[-1ex]
{\footnotesize(e)}
\end{center}
\end{minipage}
\begin{minipage}[t]{0.495\textwidth}
\begin{center}
\includegraphics[scale=0.265]{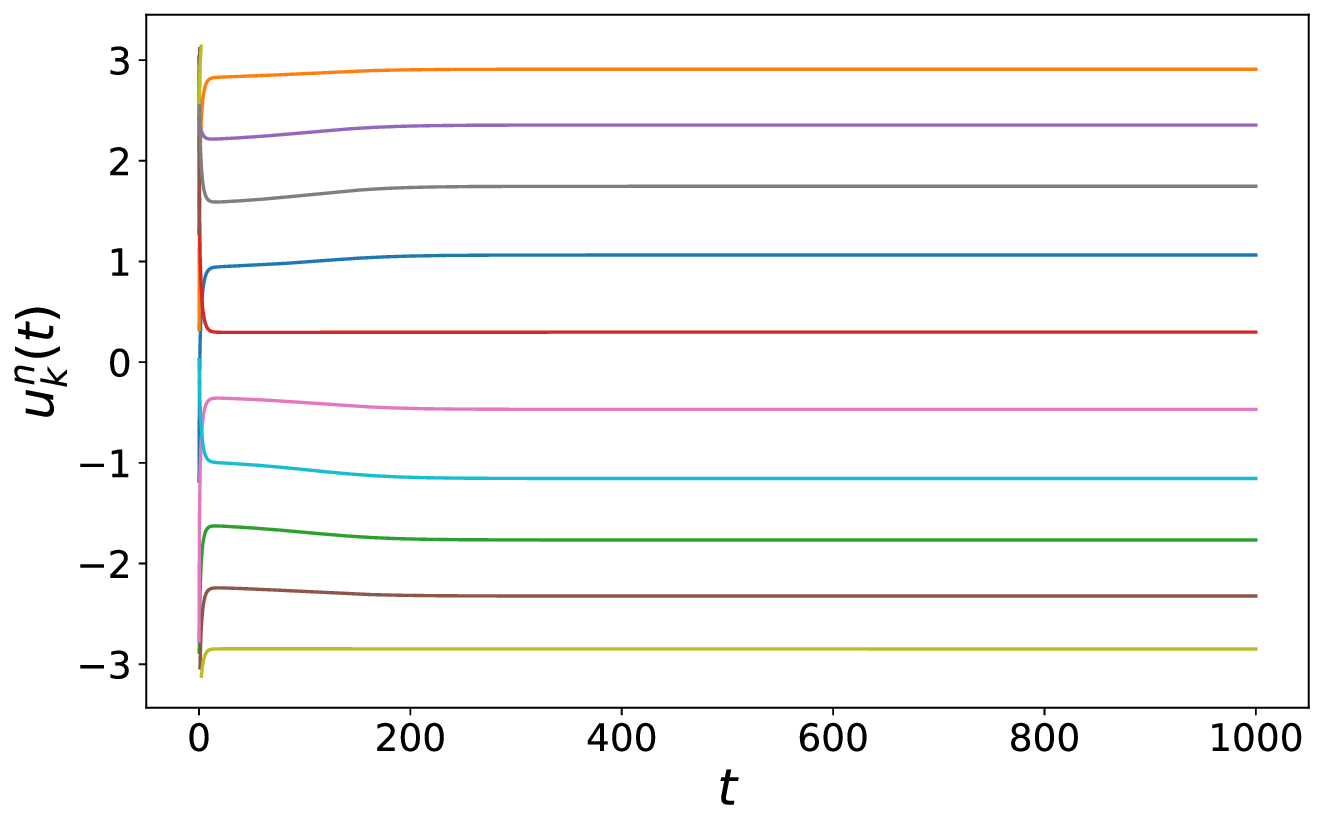}\\[-1ex]
{\footnotesize(f)}
\end{center}
\end{minipage}
\vspace*{0.5ex}

\begin{minipage}[t]{0.495\textwidth}
\begin{center}
\includegraphics[scale=0.265]{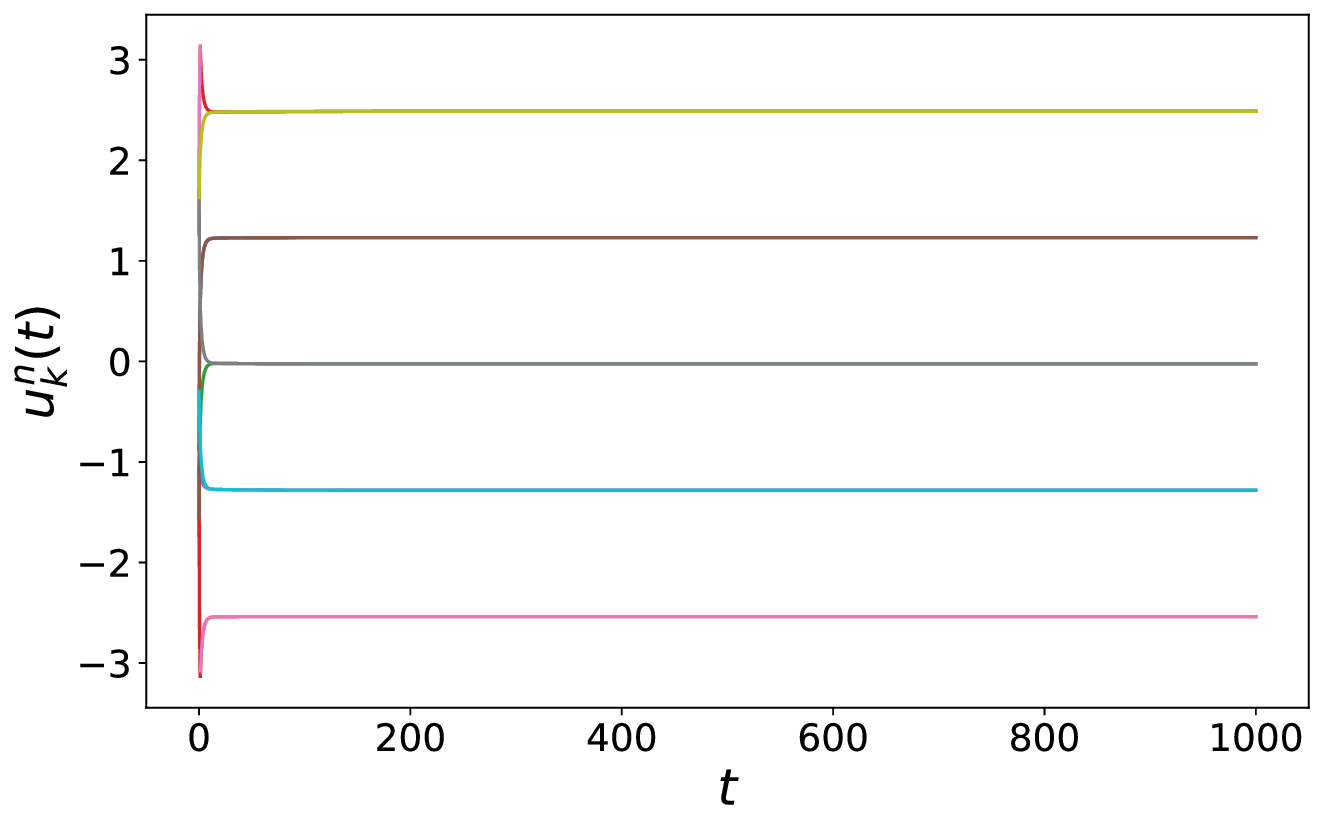}\\[-1ex]
{\footnotesize(g)}
\end{center}
\end{minipage}
\begin{minipage}[t]{0.495\textwidth}
\begin{center}
\includegraphics[scale=0.265]{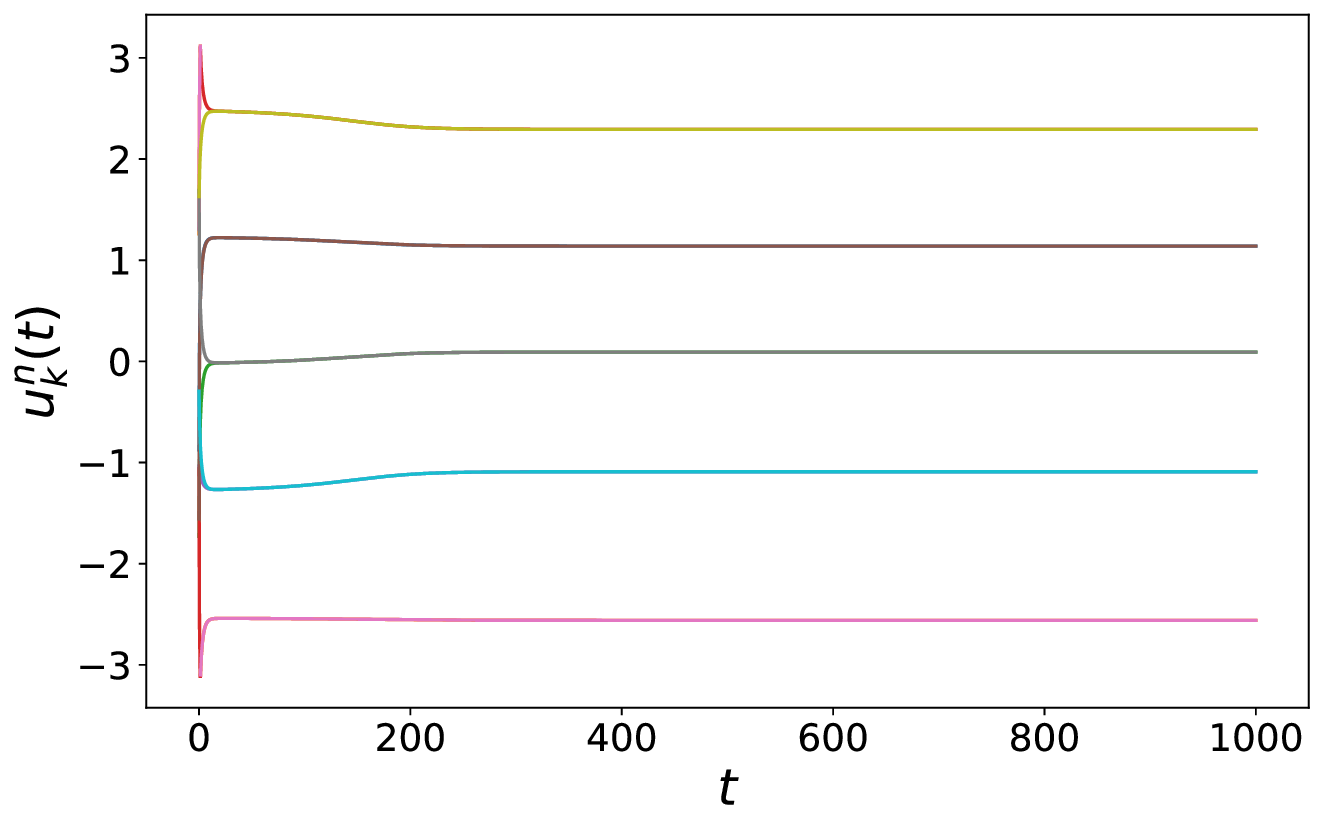}\\[-1ex]
{\footnotesize(h)}
\end{center}
\end{minipage}
\caption{Numerical simulation results for the KM \eqref{eqn:dsys}
 with $n=1000$, $\kappa=0.4$ and $\sigma=0$:
(a) $(q,b_1,b_3)=(1,0.16,1)$; (b) $(1,0.12,1)$;
(c) $(2,0.55,0.5)$; (d) $(2,0.51,0.5)$;
(e) $(3,0.34,0.5)$; (f) $(3,0.3,0.5)$;
(g) $(4,0.49,0.5)$; (h) $(4,0.45,0.5)$.
The values of $u_k^n(t)\mod 2\pi$, $k\in[n]$, are plotted as the ordinates.
The five pairs of two lines coincide almost completely in Figs.~(c), (d), (g) and (h).}
\label{fig:5a1}
\end{figure}

\begin{figure}[t]
\begin{minipage}[t]{0.495\textwidth}
\begin{center}
\includegraphics[scale=0.265]{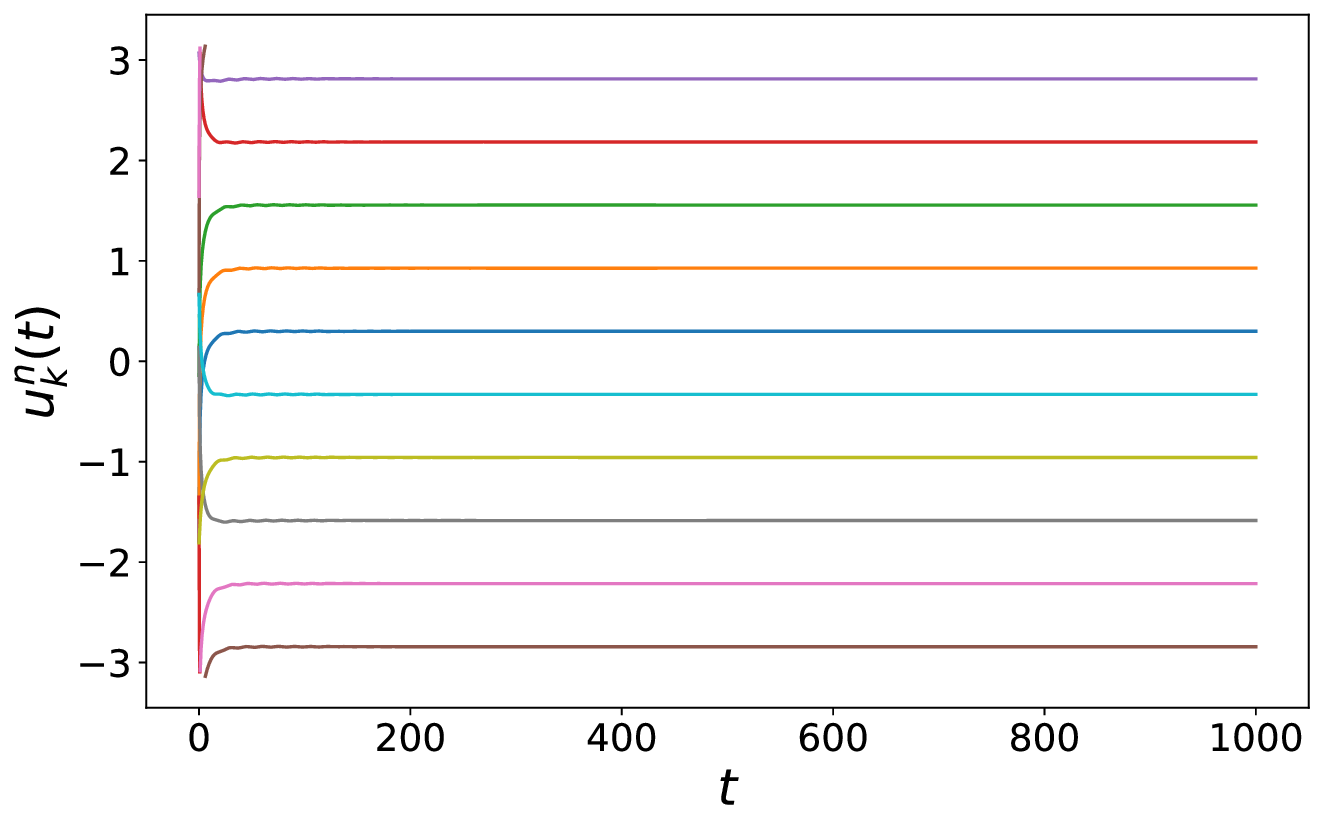}\\[-1ex]
{\footnotesize(a)}
\end{center}
\end{minipage}
\begin{minipage}[t]{0.495\textwidth}
\begin{center}
\includegraphics[scale=0.265]{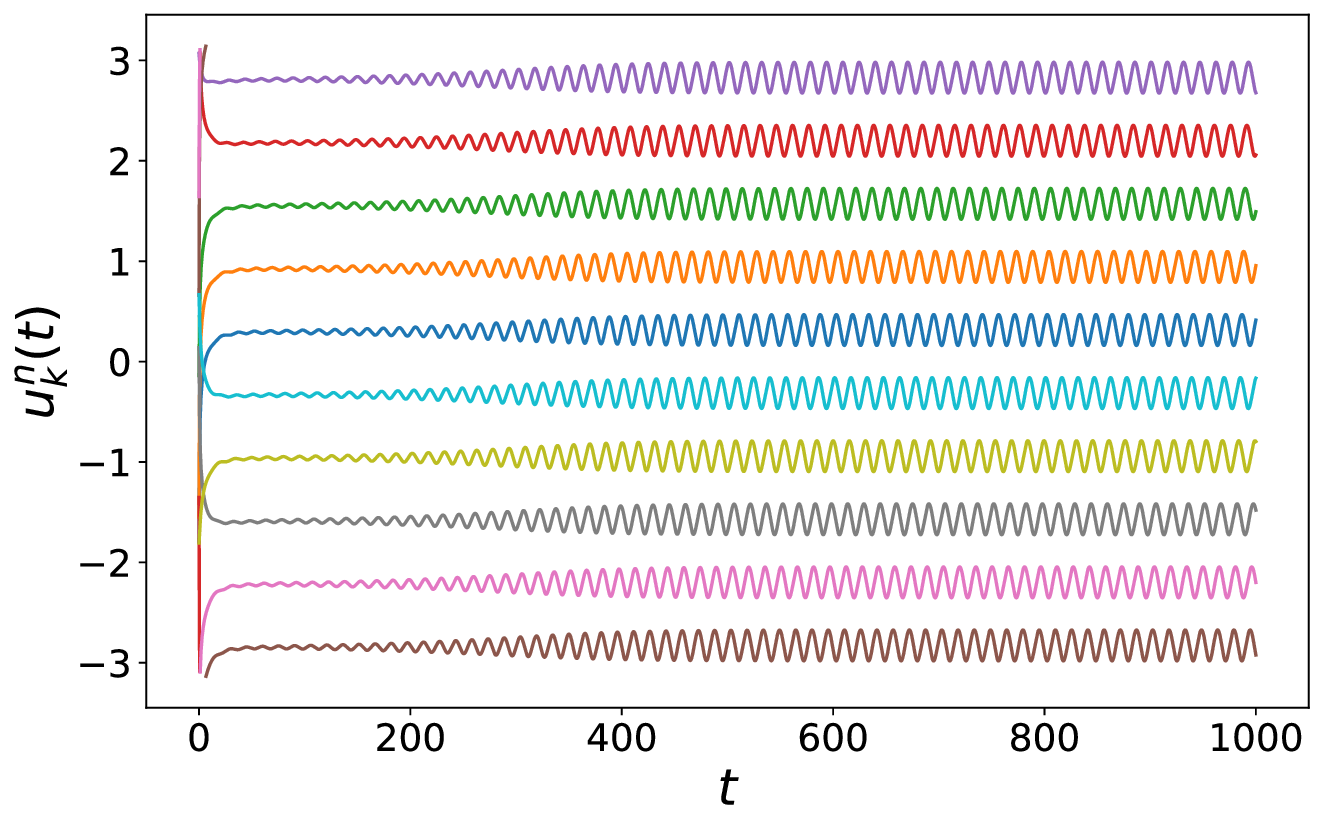}\\[-1ex]
{\footnotesize(b)}
\end{center}
\end{minipage}
\vspace*{0.5ex}

\begin{minipage}[t]{0.495\textwidth}
\begin{center}
\includegraphics[scale=0.265]{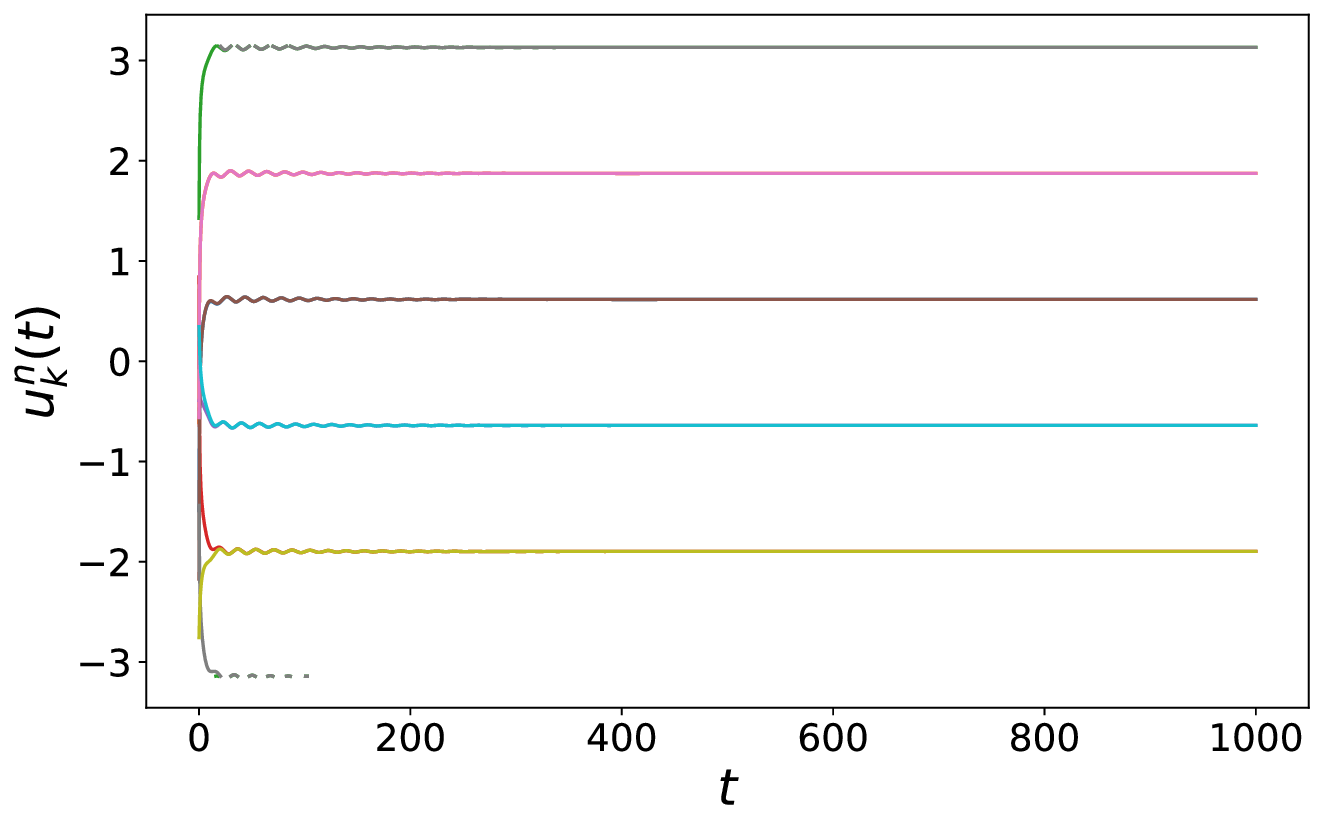}\\[-1ex]
{\footnotesize(c)}
\end{center}
\end{minipage}
\begin{minipage}[t]{0.495\textwidth}
\begin{center}
\includegraphics[scale=0.265]{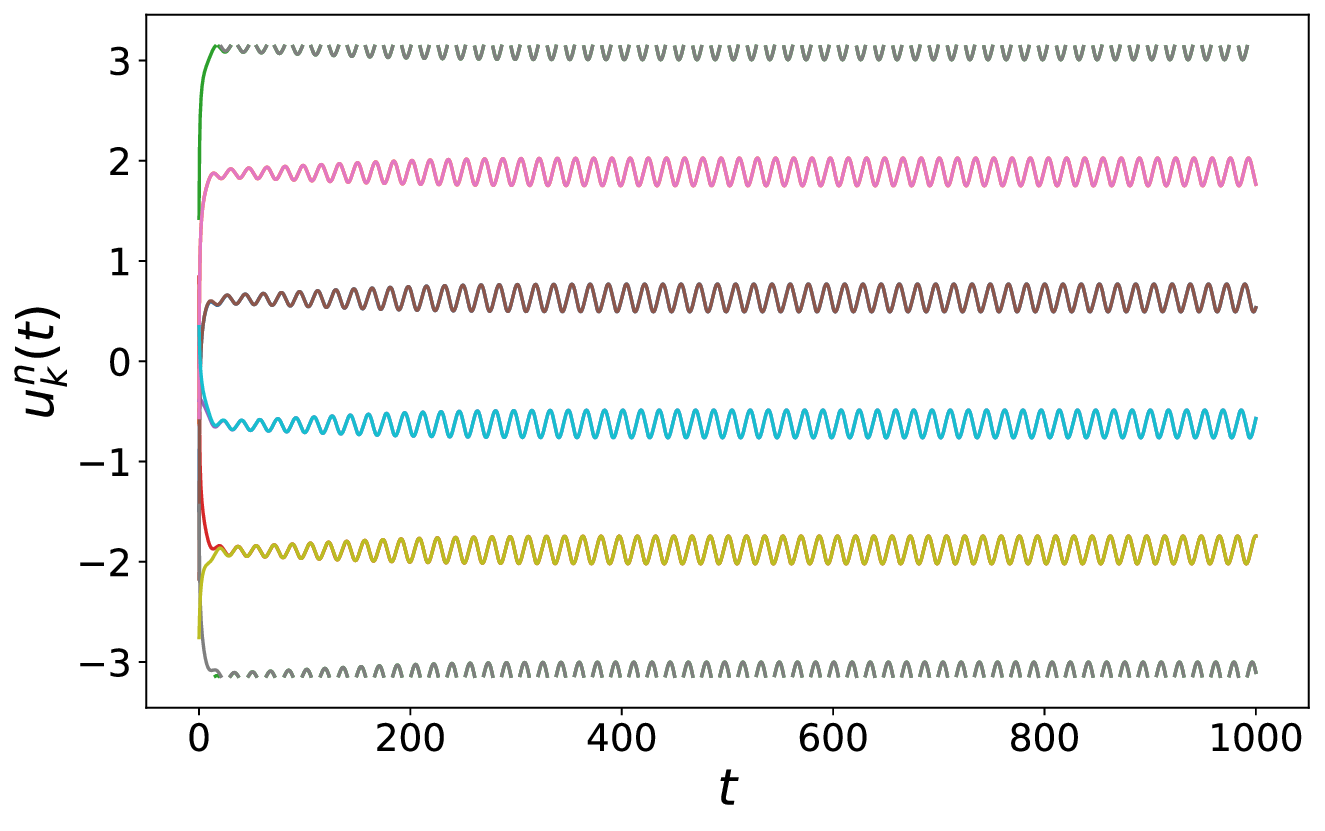}\\[-1ex]
{\footnotesize(d)}
\end{center}
\end{minipage}
\vspace*{0.5ex}

\begin{minipage}[t]{0.495\textwidth}
\begin{center}
\includegraphics[scale=0.265]{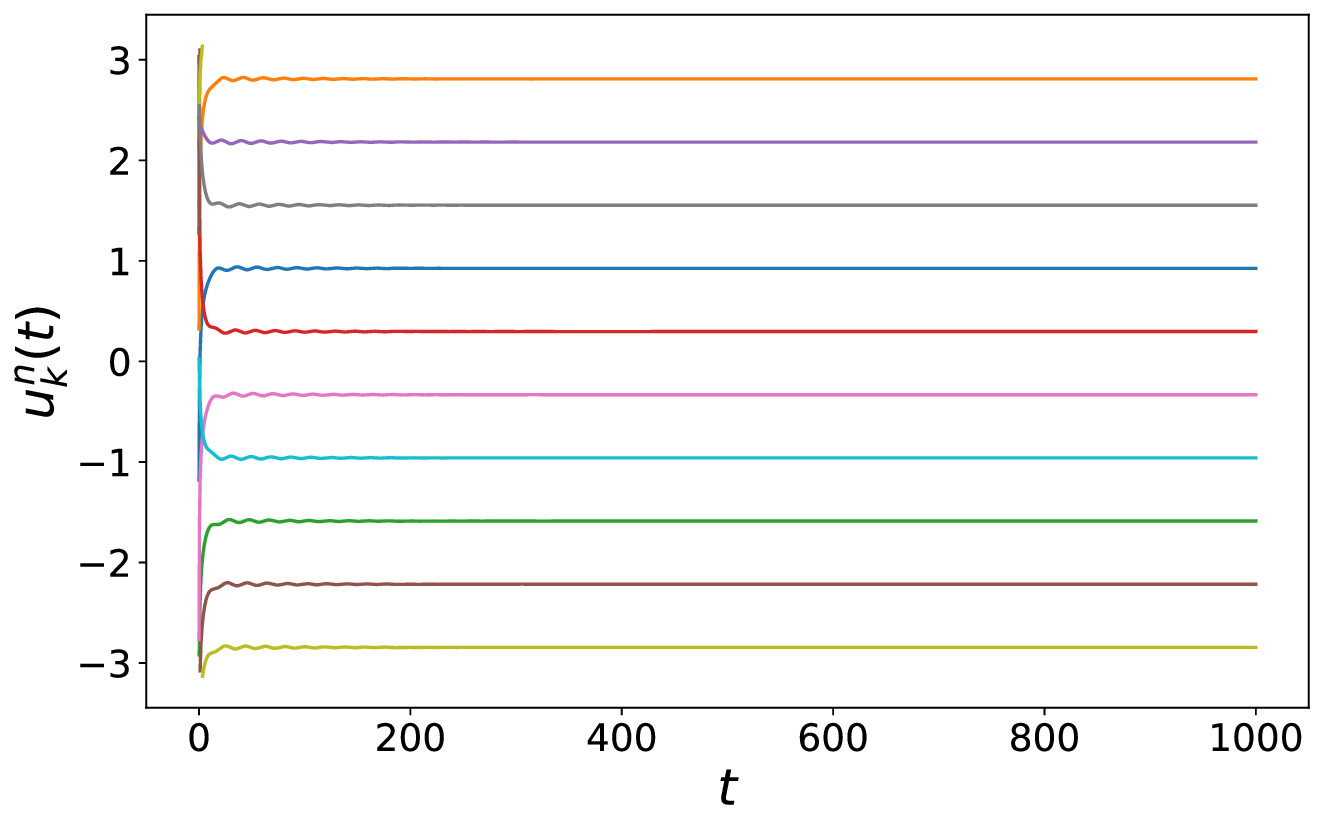}\\[-1ex]
{\footnotesize(e)}
\end{center}
\end{minipage}
\begin{minipage}[t]{0.495\textwidth}
\begin{center}
\includegraphics[scale=0.265]{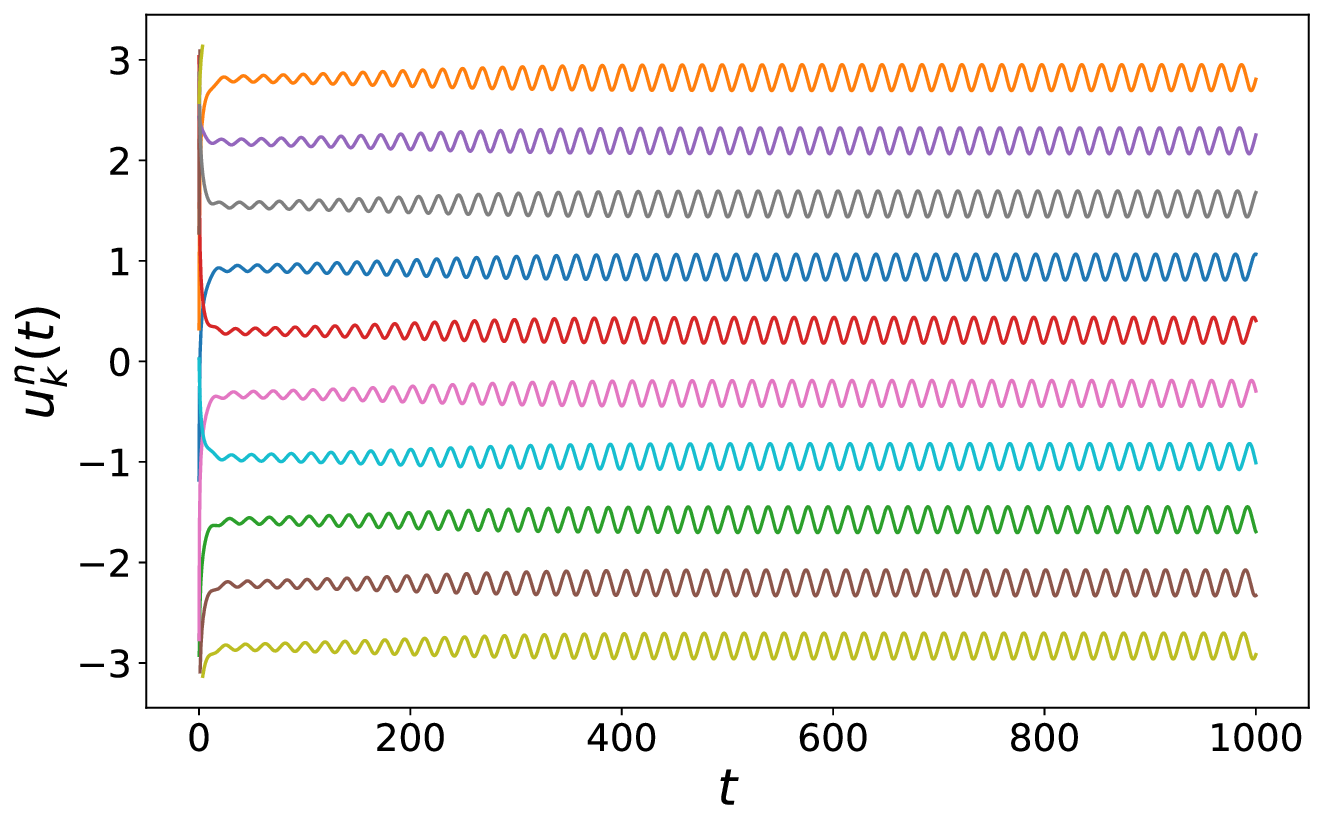}\\[-1ex]
{\footnotesize(f)}
\end{center}
\end{minipage}
\vspace*{0.5ex}

\begin{minipage}[t]{0.495\textwidth}
\begin{center}
\includegraphics[scale=0.265]{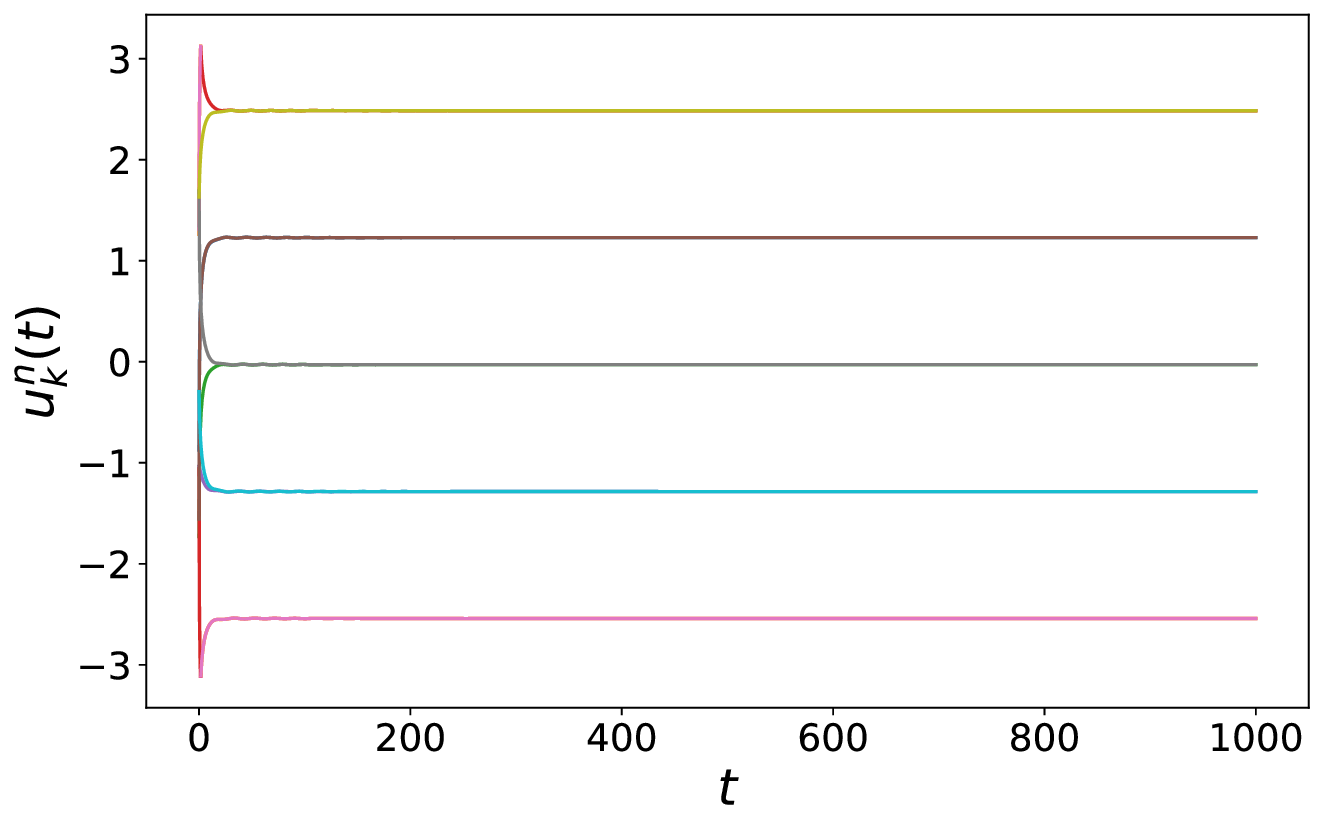}\\[-1ex]
{\footnotesize(g)}
\end{center}
\end{minipage}
\begin{minipage}[t]{0.495\textwidth}
\begin{center}
\includegraphics[scale=0.265]{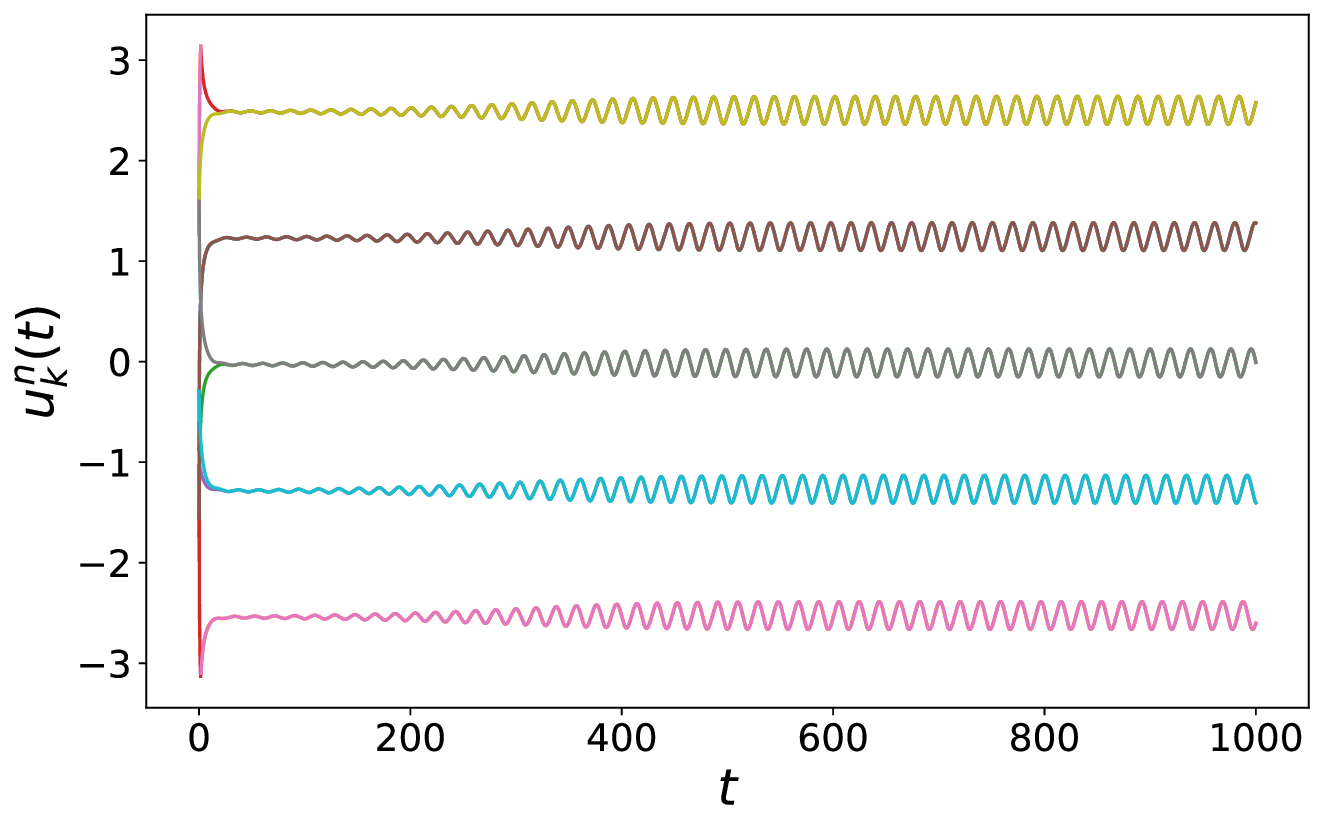}\\[-1ex]
{\footnotesize(h)}
\end{center}
\end{minipage}
\caption{Numerical simulation results for the KM \eqref{eqn:dsys}
with $n=1000$, $\kappa=0.4$, $\sigma=\pi/3$ and $b_3=0.5$:
(a) $(q,b_1)=(1,0.08)$; (b) $(1,0.06)$;
(c) $(2,0.275)$; (d) $(2,0.255)$;
(e) $(3,0.17)$; (f) $(3,0.15)$;
(g) $(4,0.245)$; (h) $(4,0.225)$.
See also the caption of Fig.~\ref{fig:5a1}.
\vspace*{5mm}
}
\label{fig:5b1}
\end{figure}

We carried out numerical simulations for the KM \eqref{eqn:dsys},
 using the DOP853 solver \cite{HNW93}, for $q\in[4]$.
We took $n=1000$
 and chose the initial values $u_k^n(0)$, $k\in[n]$, independently randomly
 according to the uniform distribution on the intervals $[-\pi+2\pi q k/n,\pi+2\pi q k/n]$
 centered at the $q$-twisted state \eqref{eqn:ts} with $\Omega_\D=0$.
So if there is an asymptotically stable that are different
 from the twisted and  modulated or oscillating twisted states,
 then the responses of \eqref{eqn:dsys} may converge to it as $t\to\infty$.
We also considered two cases $\sigma=0$ and $\pi/3$ for the phase lag,
 and $\kappa=0.4$ and $0.5$ for the neighbor size.
Recall that the $\kappa$-nearest neighbor graph
 reduces to a complete simple one when $\kappa=0.5$.

We begin with numerical results
 for $\kappa$-nearest neighbor graphs with $\kappa=0.4$,
 for which the $q$-twisted solution \eqref{eqn:tsol} is unstable
 in the uncontrolled CL \eqref{eqn:csys} with $b_1,b_3=0$,
 since condition~\eqref{eqn:ls} does not hold for $\ell=q\in[4]$ when $b_1=0$
 as seen from Table~\ref{tbl:4a}.
The results for $\kappa=0.5$ are provided in the next section.

Figures~\ref{fig:5a1} and \ref{fig:5b1}
 show the time-histories of every $100$th node (from 50th to 950th)
 for $\sigma=0$ and $\pi/3$, respectively.
The values of $b_1$ in the left and right columns of each figure
 were chosen such that they are larger and smaller, respectively,
 than the bifurcation points,
 which are approximated by $b_{1q}$, $q\in[4]$, (see Eq.~\eqref{eqn:b1q})
 for the $q$-twisted solutions \eqref{eqn:tsol}
 in the CL \eqref{eqn:csys}.
We see that the responses of the KM \eqref{eqn:dsys}
 converge to their steady states rapidly
 although their initial values are randomly distributed on the wide interval.
Moreover, oscillations occur for the smaller values of $b_1$
 when $\sigma=\pi/3$, in the right column of Fig.~\ref{fig:5b1},
 as detected by Theorem~\ref{thm:4b} for the CL \eqref{eqn:csys}.
We also notice that no rotation is observed even when $\sigma=\pi/3$,
 like the $q$-twisted solution \eqref{eqn:tsol} in the CL \eqref{eqn:csys},
 in contrast to observations for the uncontrolled KM \eqref{eqn:dsys}
 with $b_1,b_3=0$ in \cite{Y24c}.

\begin{figure}[t]
\begin{minipage}[t]{0.495\textwidth}
\begin{center}
\includegraphics[scale=0.265]{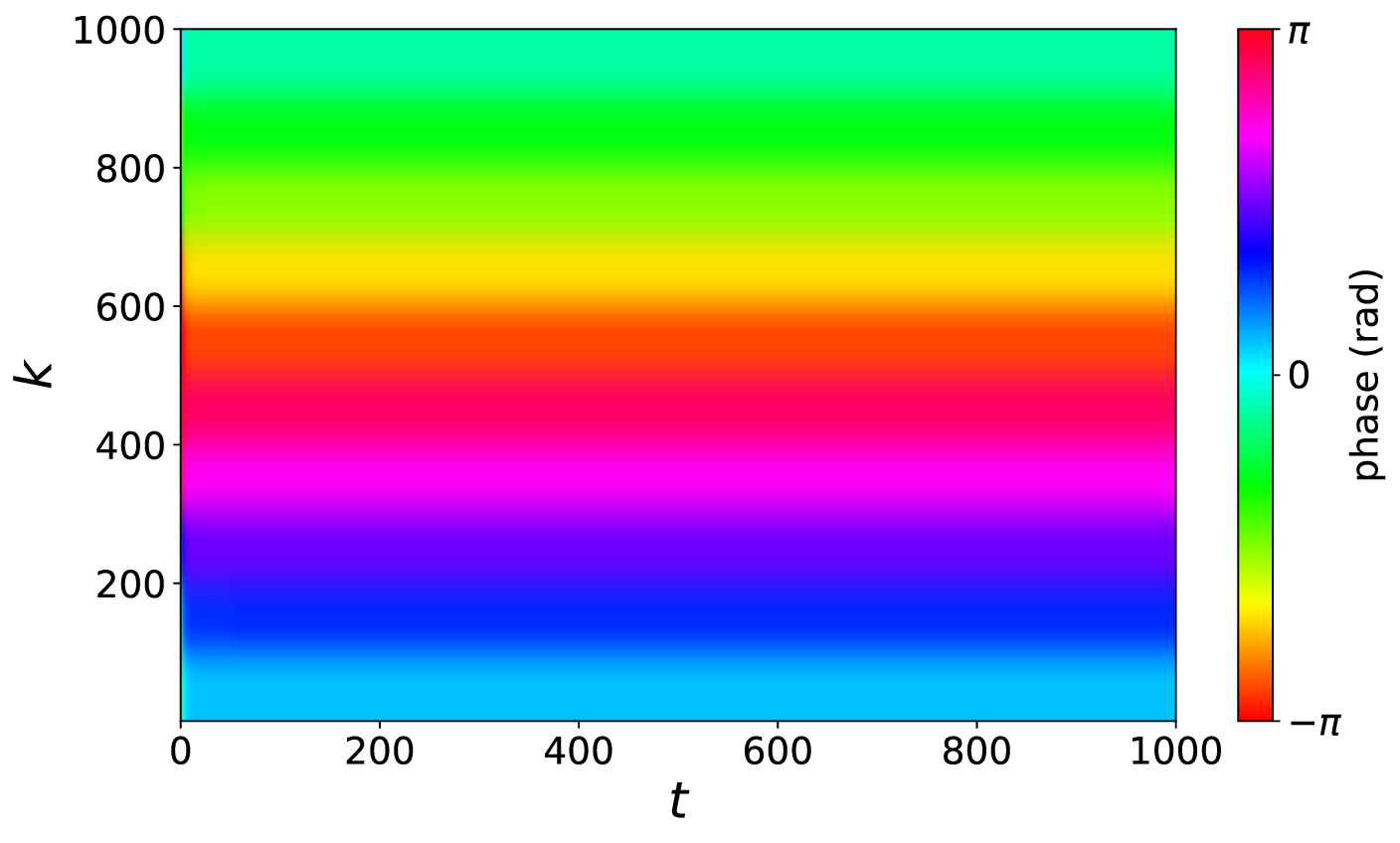}\\[-1ex]
{\footnotesize(a)}
\end{center}
\end{minipage}
\begin{minipage}[t]{0.495\textwidth}
\begin{center}
\includegraphics[scale=0.265]{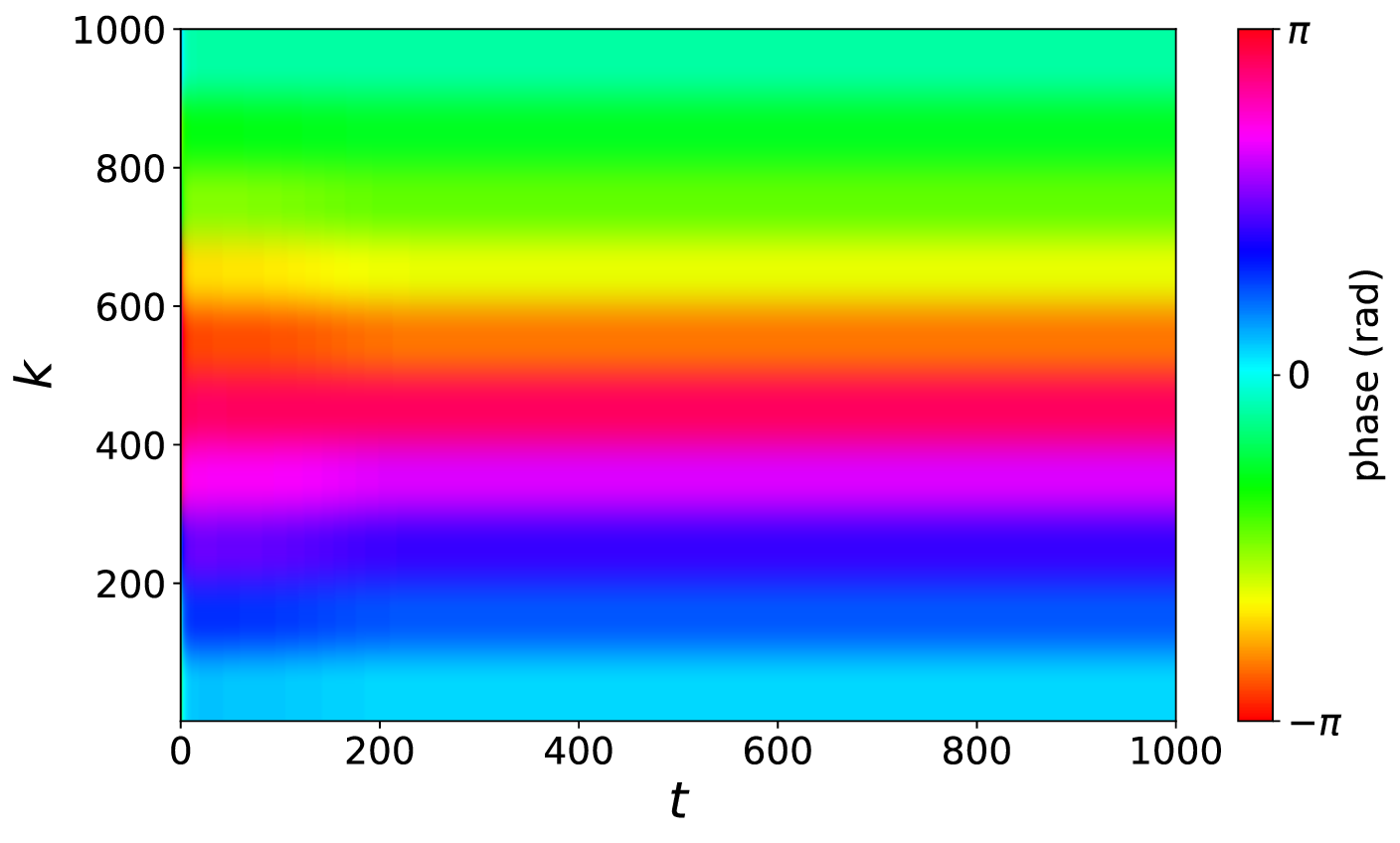}\\[-1ex]
{\footnotesize(b)}
\end{center}
\end{minipage}

\begin{minipage}[t]{0.495\textwidth}
\begin{center}
\includegraphics[scale=0.265]{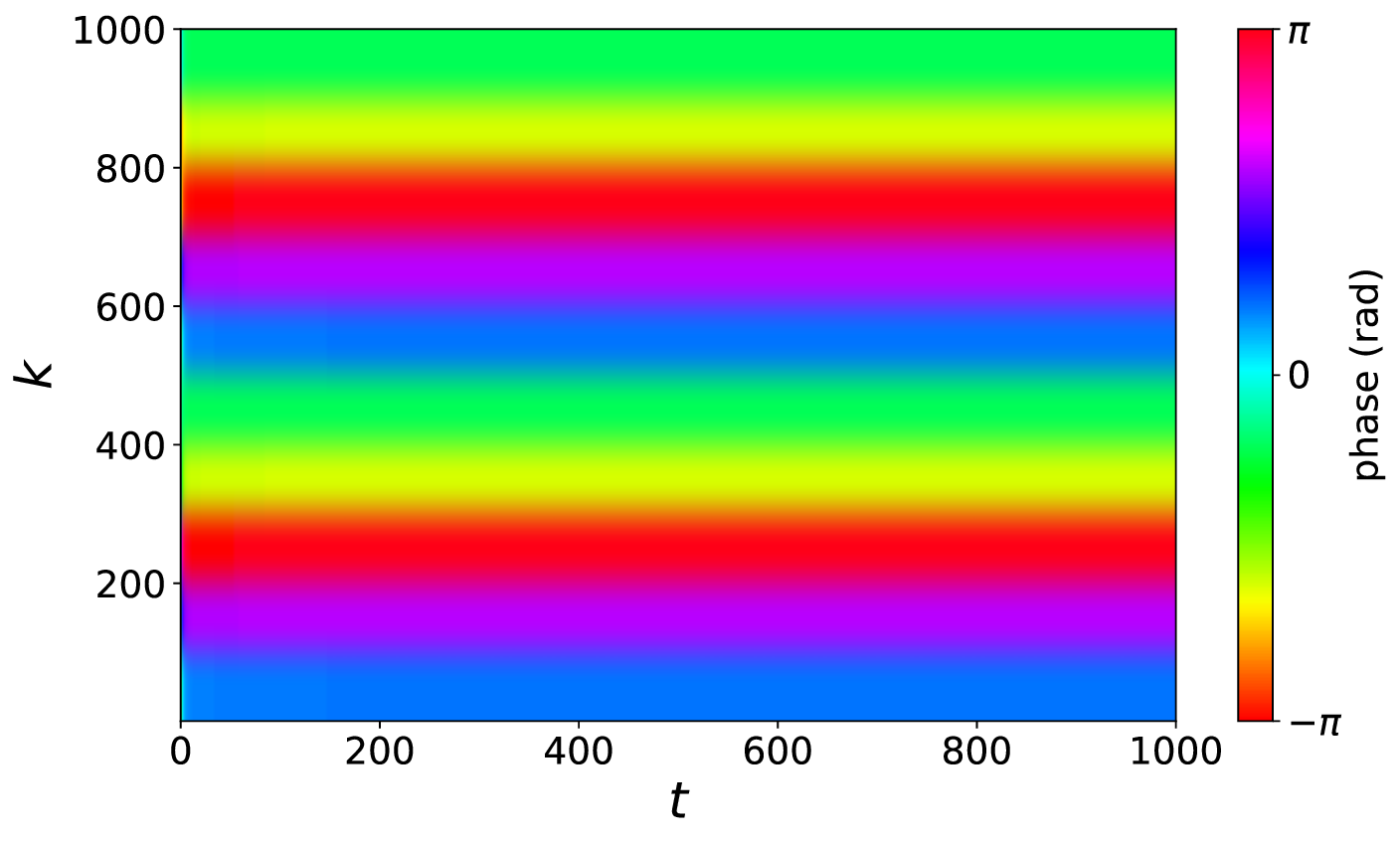}\\[-1ex]
{\footnotesize(c)}
\end{center}
\end{minipage}
\begin{minipage}[t]{0.495\textwidth}
\begin{center}
\includegraphics[scale=0.265]{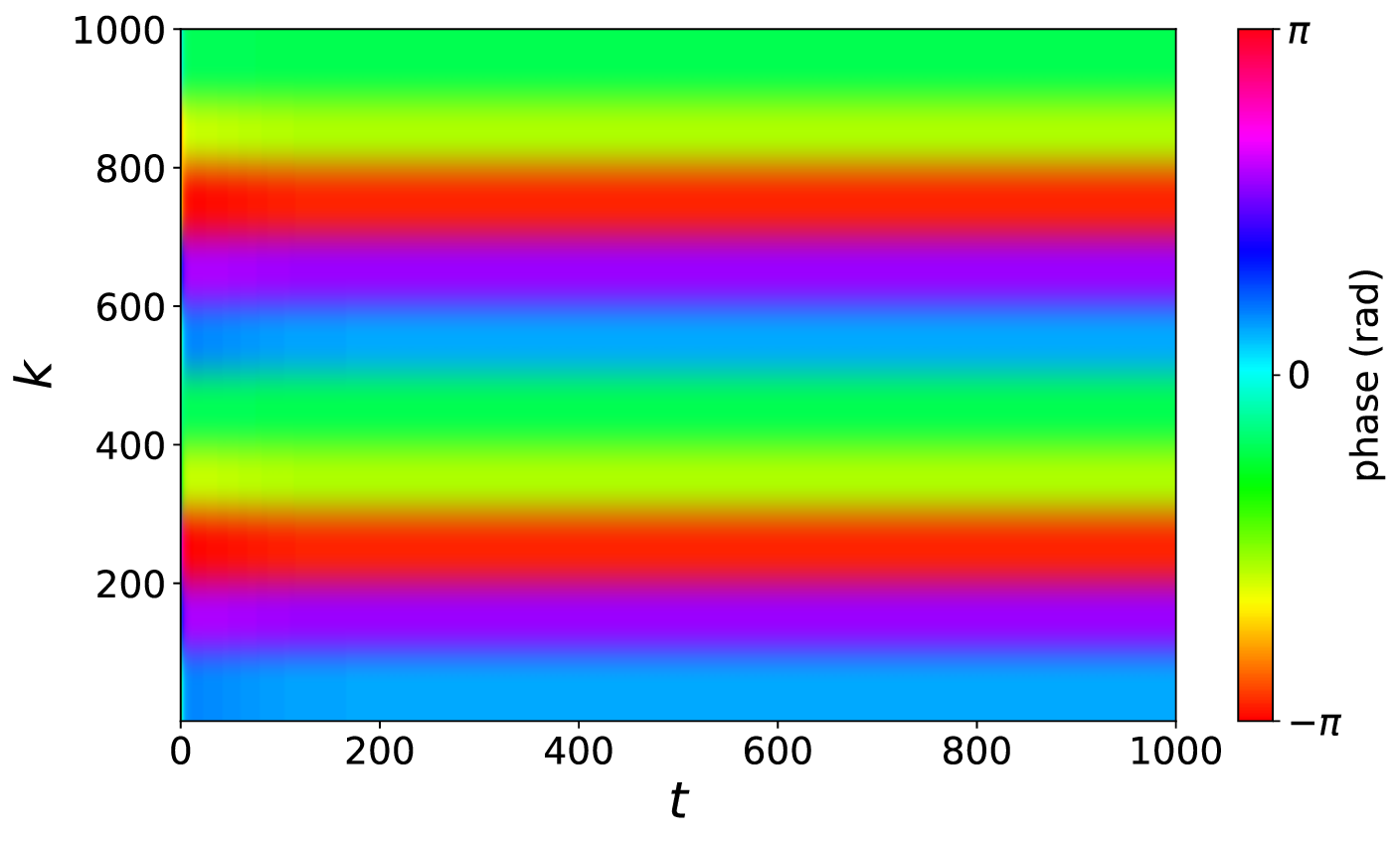}\\[-1ex]
{\footnotesize(d)}
\end{center}
\end{minipage}

\begin{minipage}[t]{0.495\textwidth}
\begin{center}
\includegraphics[scale=0.265]{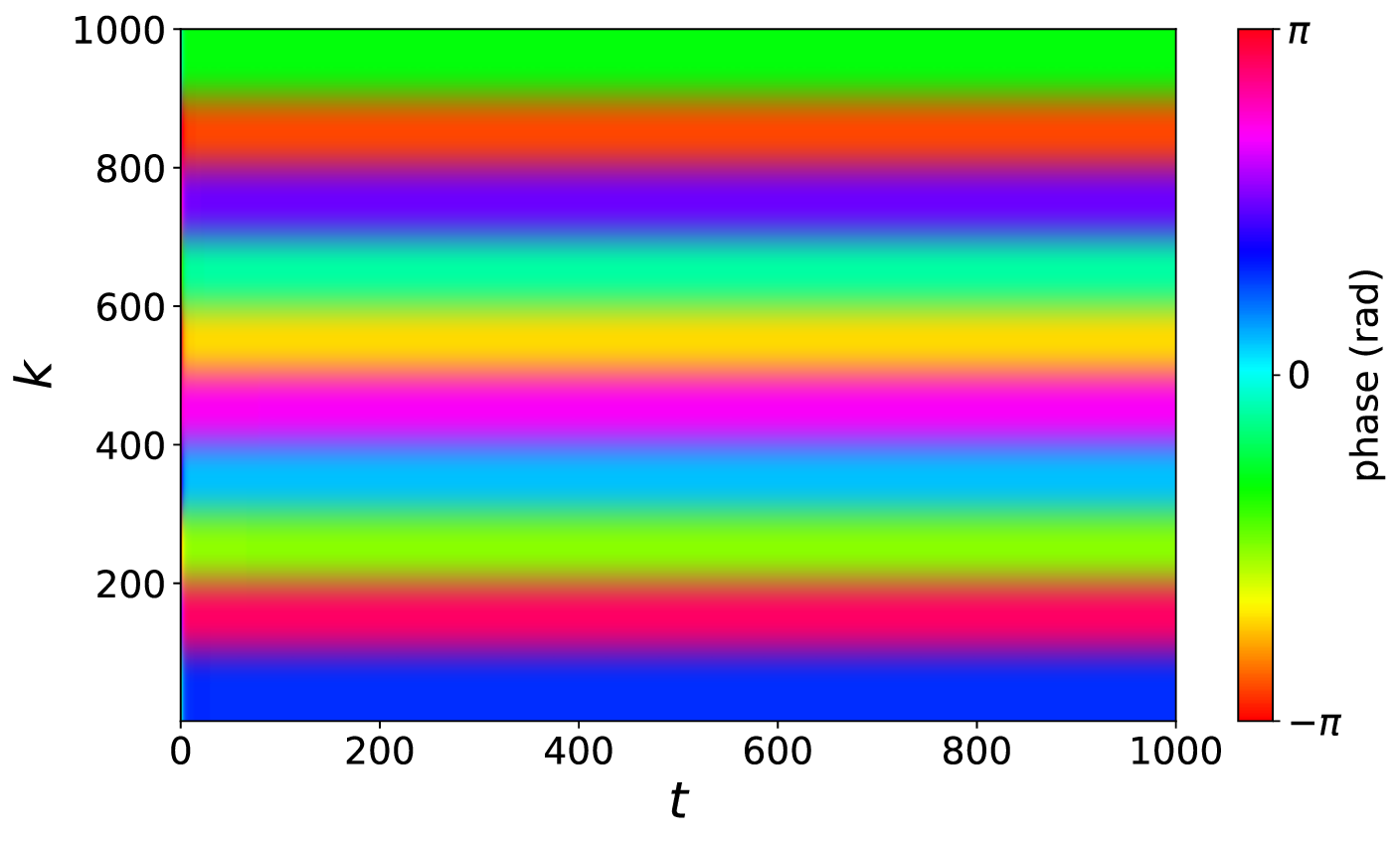}\\[-1ex]
{\footnotesize(e)}
\end{center}
\end{minipage}
\begin{minipage}[t]{0.495\textwidth}
\begin{center}
\includegraphics[scale=0.265]{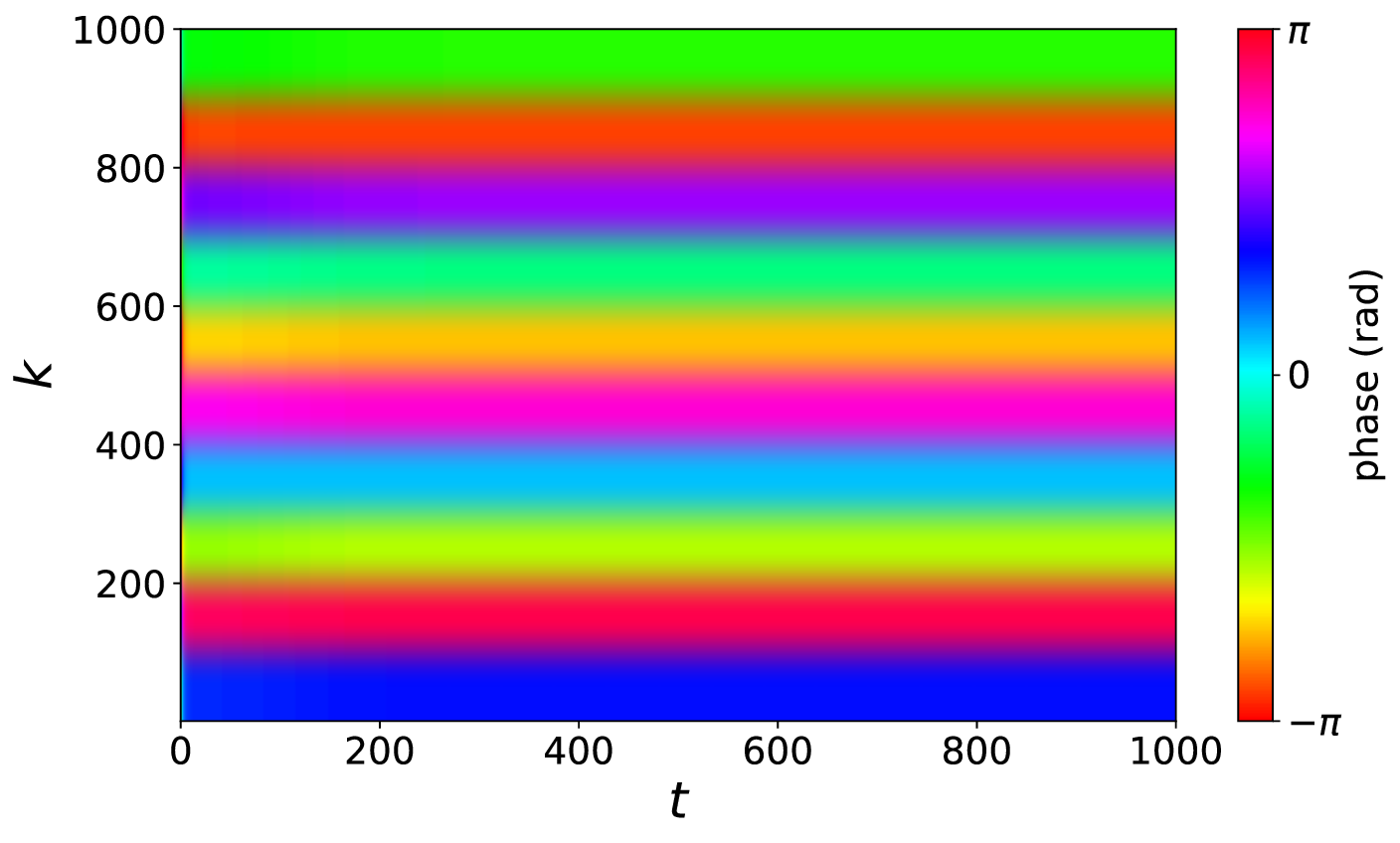}\\[-1ex]
{\footnotesize(f)}
\end{center}
\end{minipage}

\begin{minipage}[t]{0.495\textwidth}
\begin{center}
\includegraphics[scale=0.265]{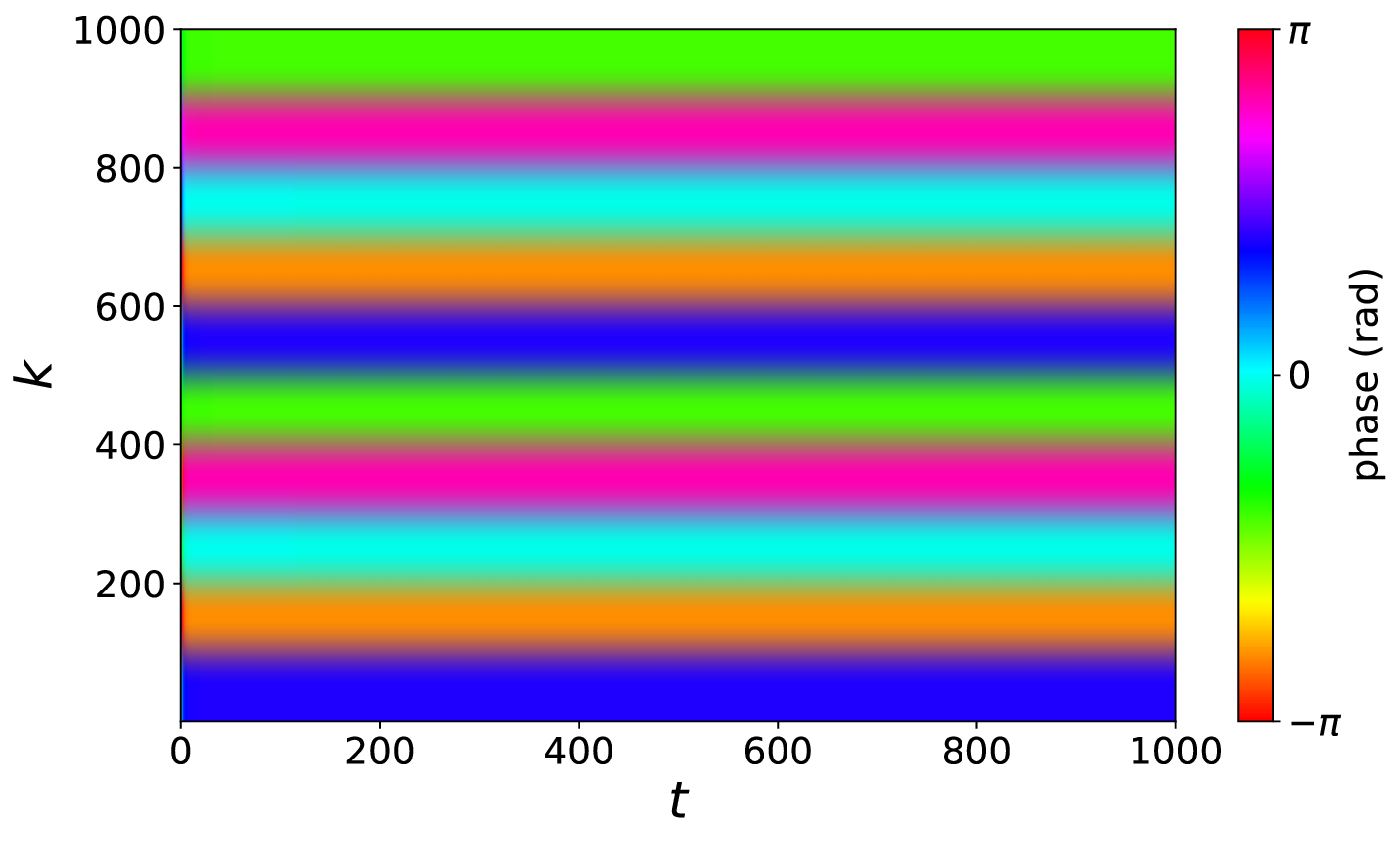}\\[-1ex]
{\footnotesize(g)}
\end{center}
\end{minipage}
\begin{minipage}[t]{0.495\textwidth}
\begin{center}
\includegraphics[scale=0.265]{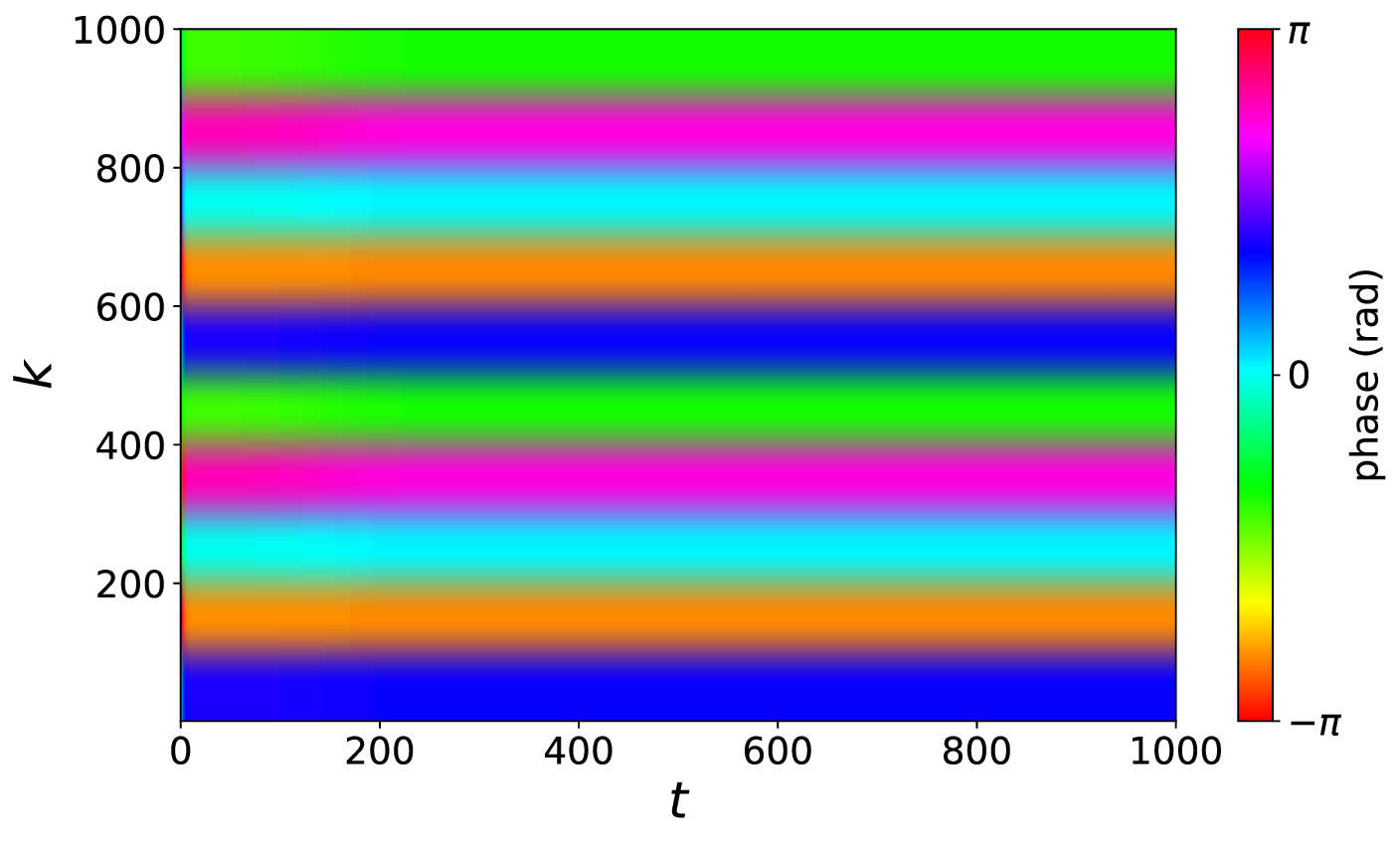}\\[-1ex]
{\footnotesize(h)}
\end{center}
\end{minipage}
\caption{Space-time plots of all oscillator phases $u_k^n(t)$, $k\in[n]$,
 in the KM \eqref{eqn:dsys} with $n=1000$, $\kappa=0.4$ and $\sigma=0$:
(a) $(q,b_1,b_3)=(1,0.16,1)$; (b) $(1,0.12,1)$;
(c) $(2,0.275,0.5)$; (d) $(2,0.255,0.5)$;
(e) $(3,0.17,0.5)$; (f) $(3,0.15,0.5)$;
(g) $(4,0.245,0.5)$; (h) $(4,0.225,0.5)$.
}
\label{fig:5a1'}
\end{figure}

\begin{figure}[t]
\begin{minipage}[t]{0.495\textwidth}
\begin{center}
\includegraphics[scale=0.265]{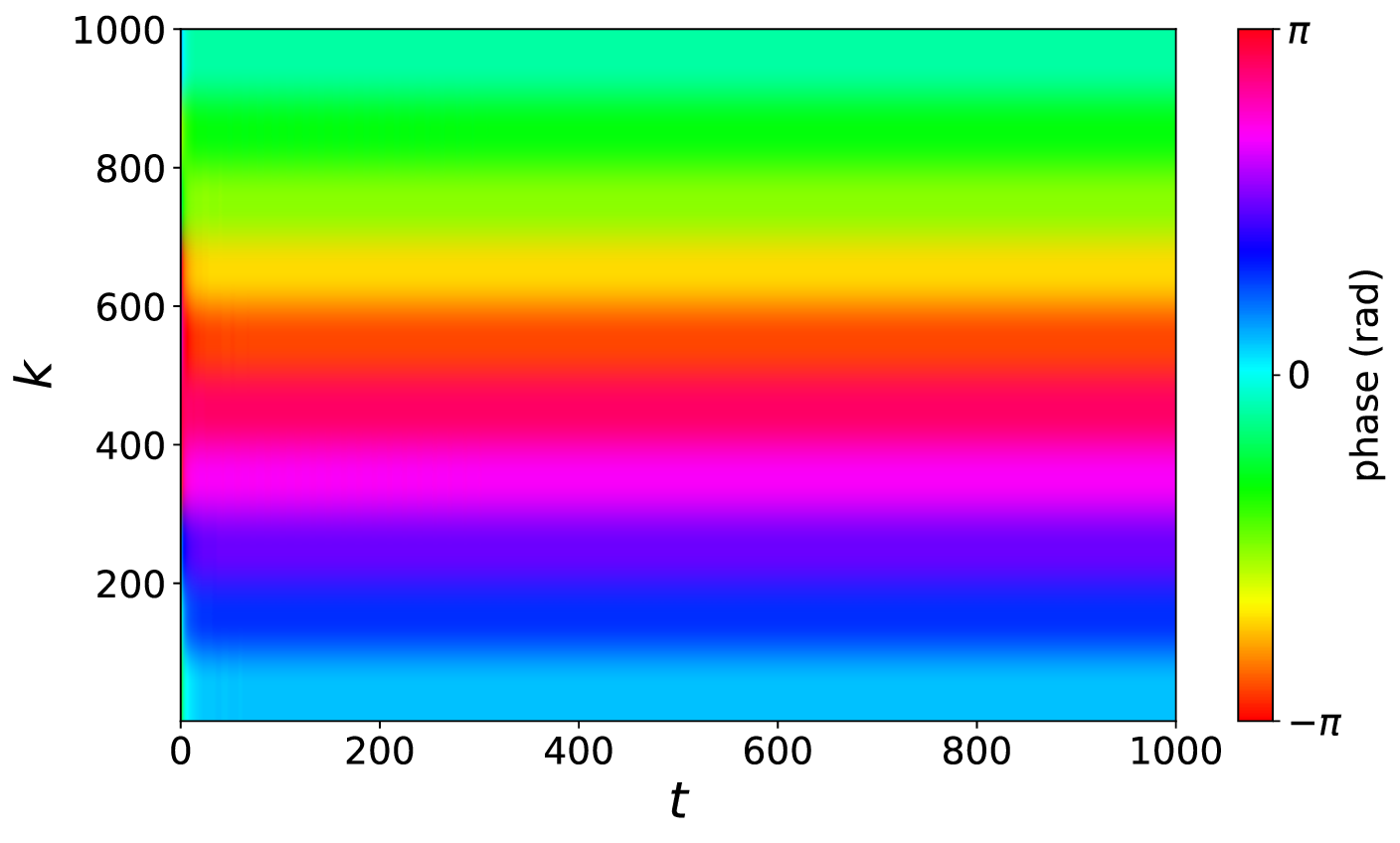}\\[-1ex]
{\footnotesize(a)}
\end{center}
\end{minipage}
\begin{minipage}[t]{0.495\textwidth}
\begin{center}
\includegraphics[scale=0.265]{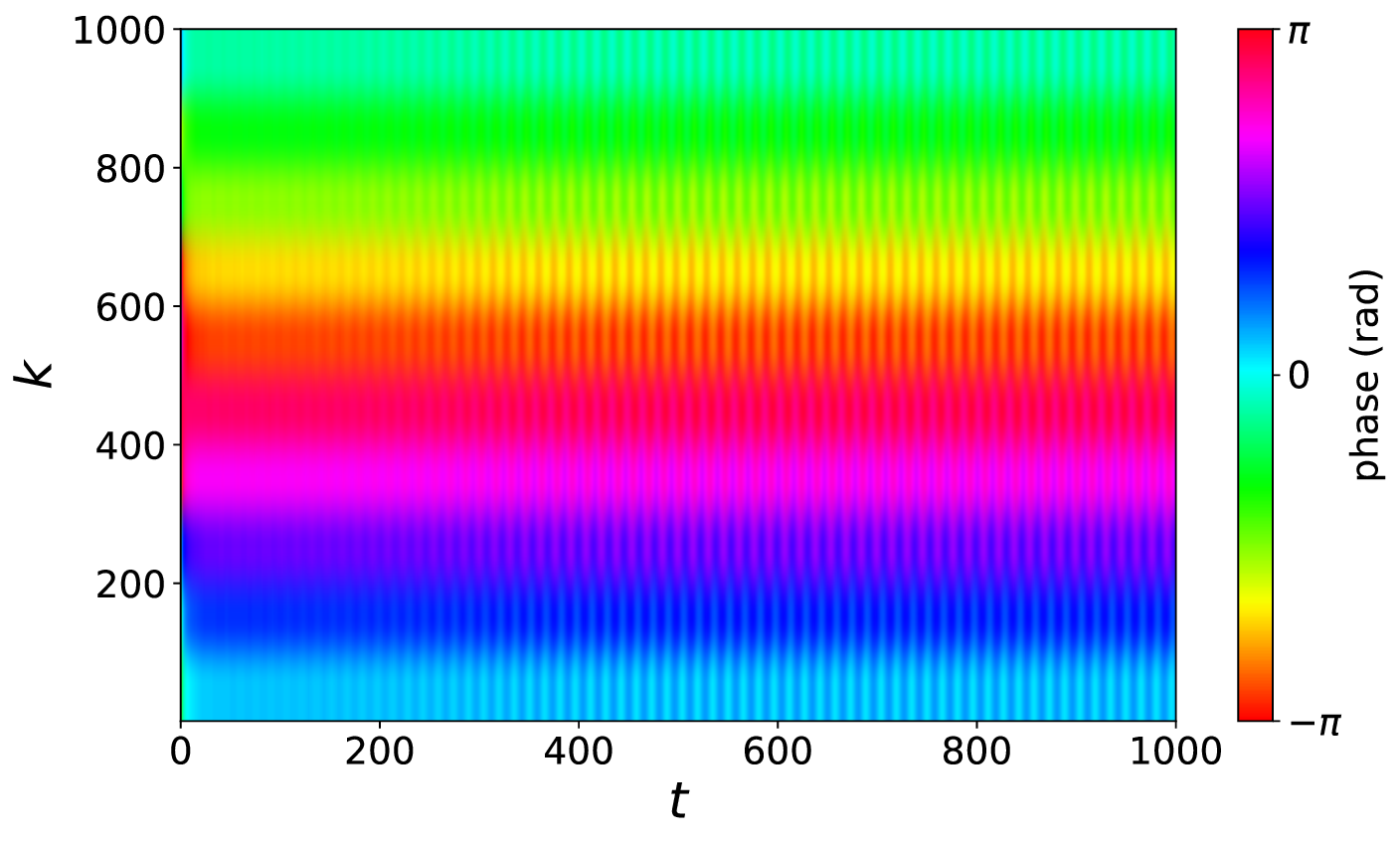}\\[-1ex]
{\footnotesize(b)}
\end{center}
\end{minipage}

\begin{minipage}[t]{0.495\textwidth}
\begin{center}
\includegraphics[scale=0.265]{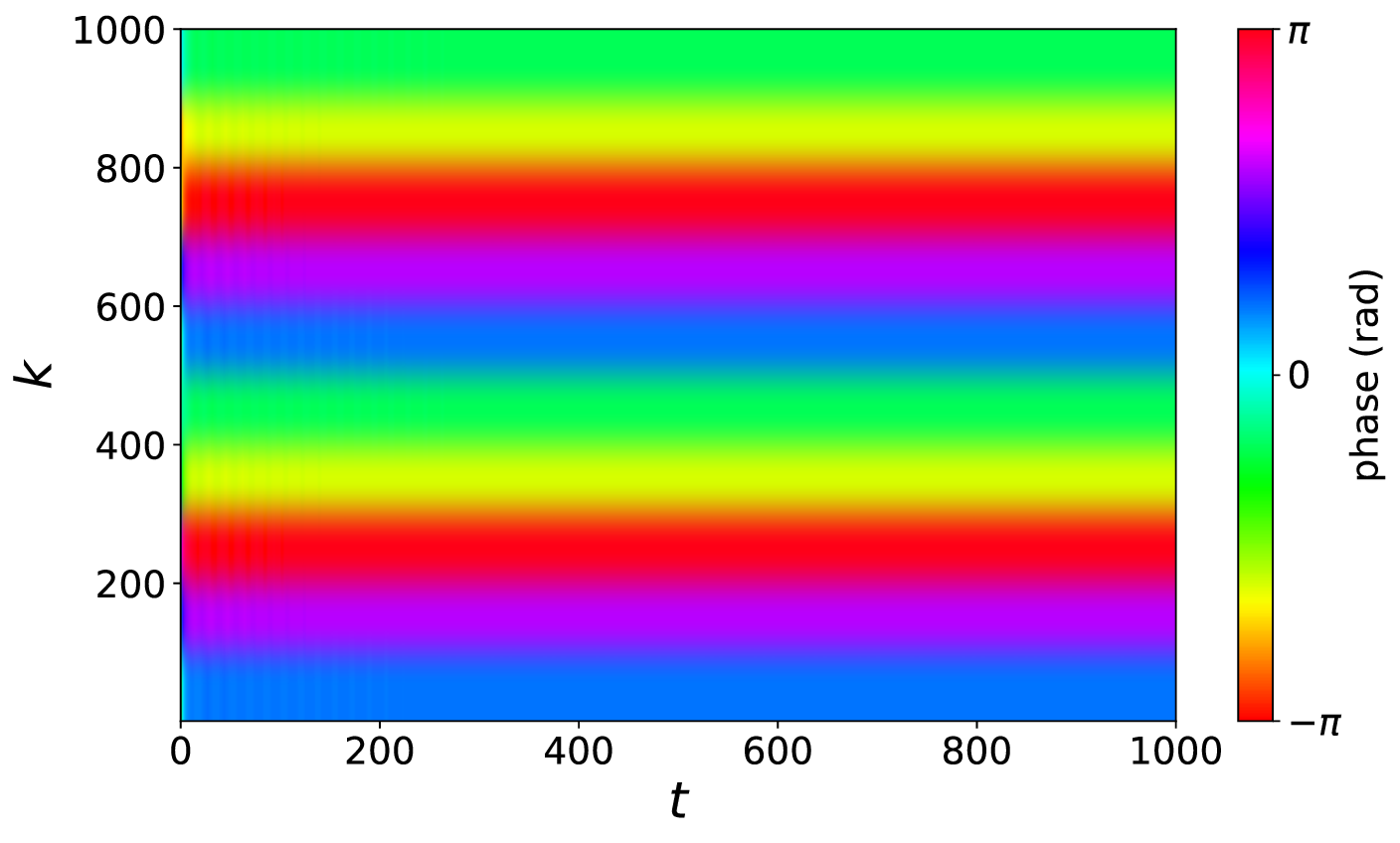}\\[-1ex]
{\footnotesize(c)}
\end{center}
\end{minipage}
\begin{minipage}[t]{0.495\textwidth}
\begin{center}
\includegraphics[scale=0.265]{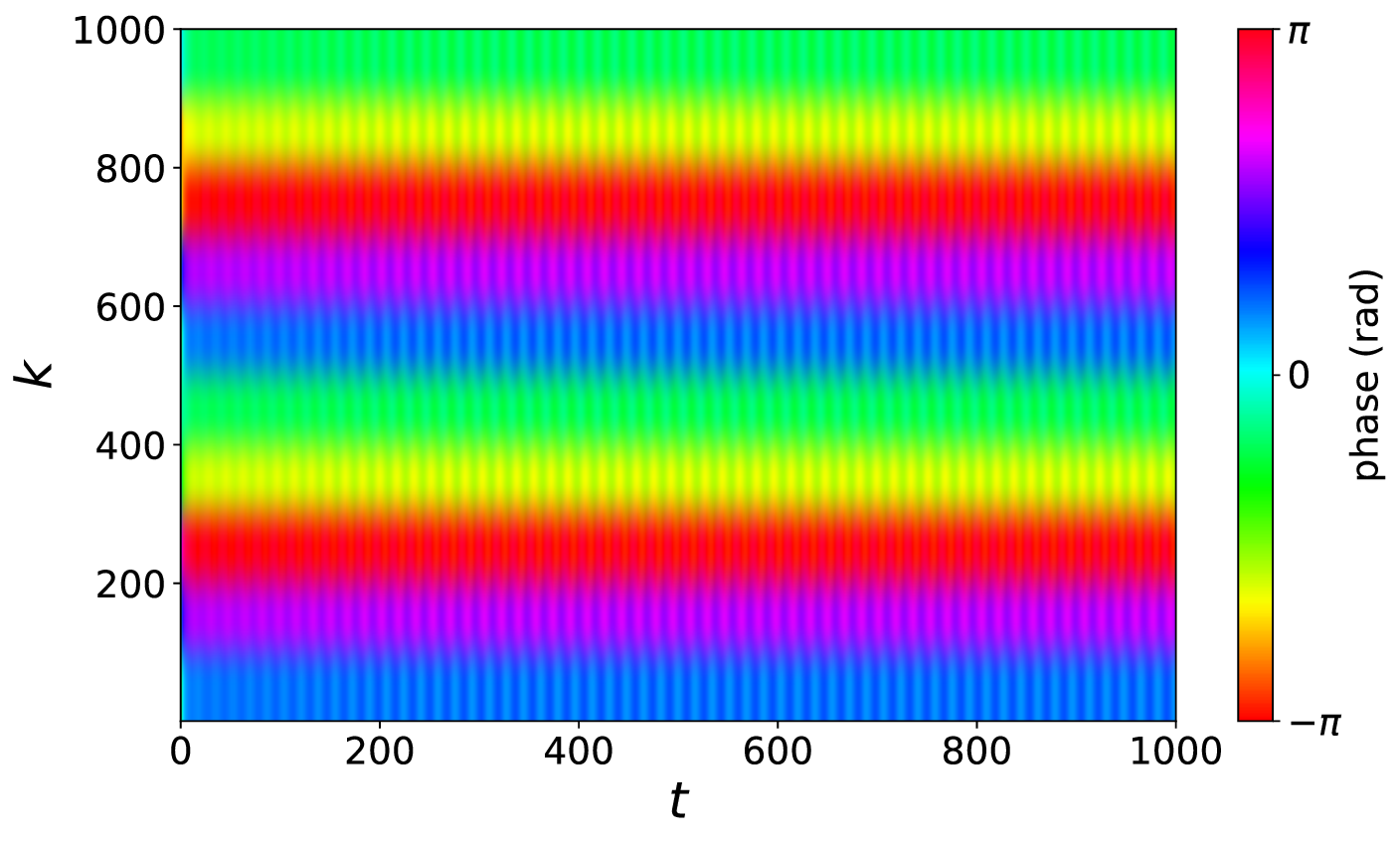}\\[-1ex]
{\footnotesize(d)}
\end{center}
\end{minipage}

\begin{minipage}[t]{0.495\textwidth}
\begin{center}
\includegraphics[scale=0.265]{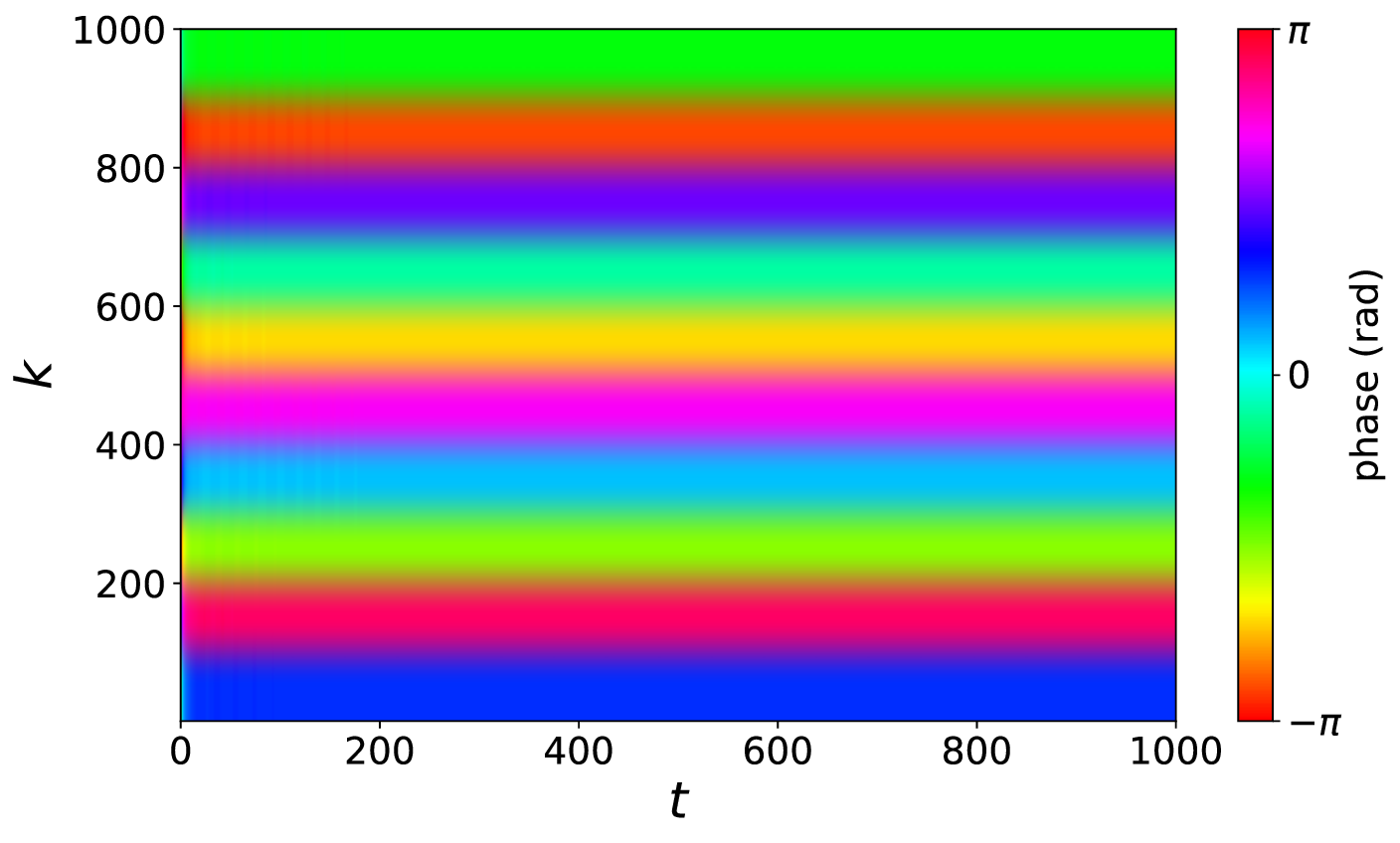}\\[-1ex]
{\footnotesize(e)}
\end{center}
\end{minipage}
\begin{minipage}[t]{0.495\textwidth}
\begin{center}
\includegraphics[scale=0.265]{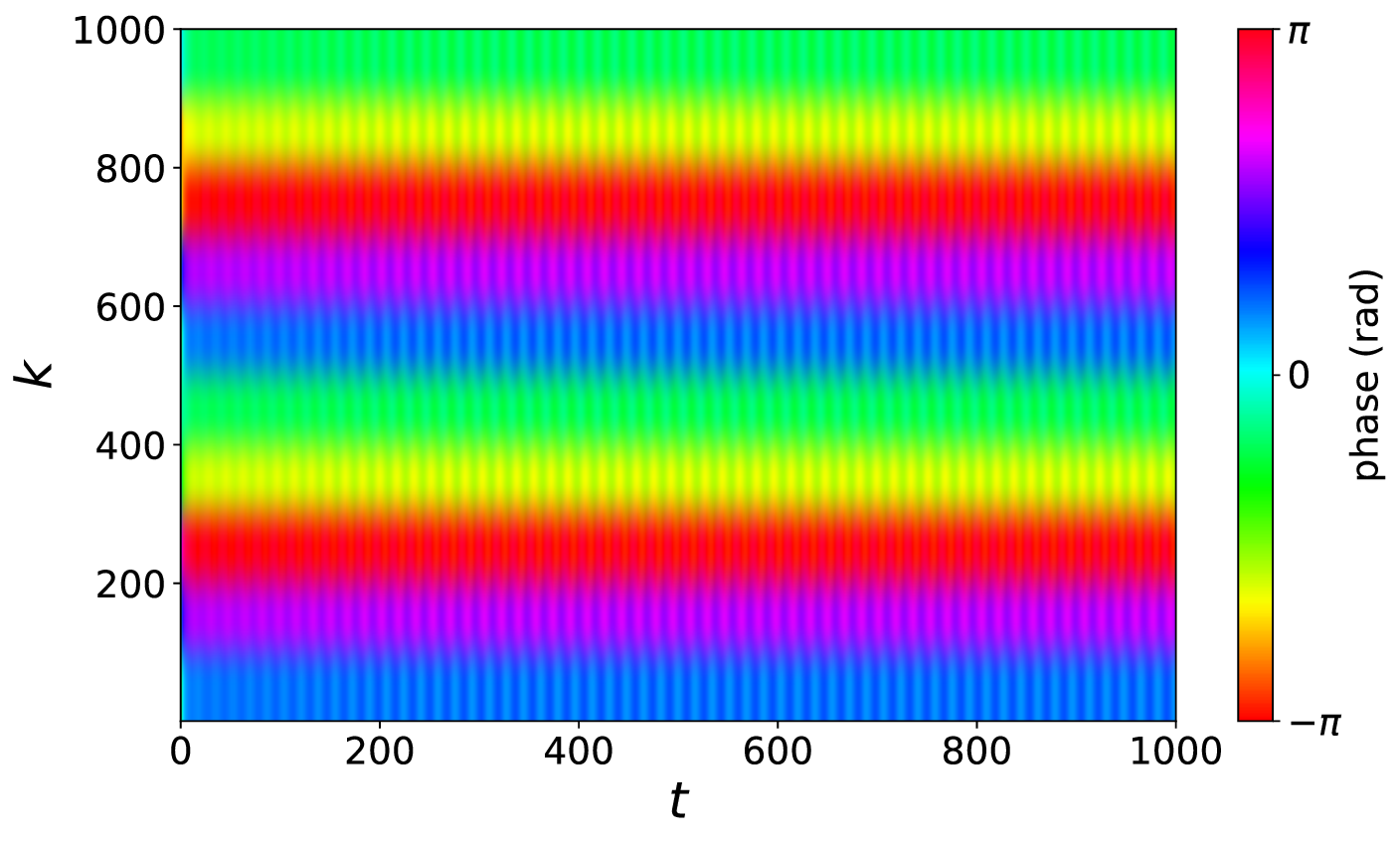}\\[-1ex]
{\footnotesize(f)}
\end{center}
\end{minipage}

\begin{minipage}[t]{0.495\textwidth}
\begin{center}
\includegraphics[scale=0.265]{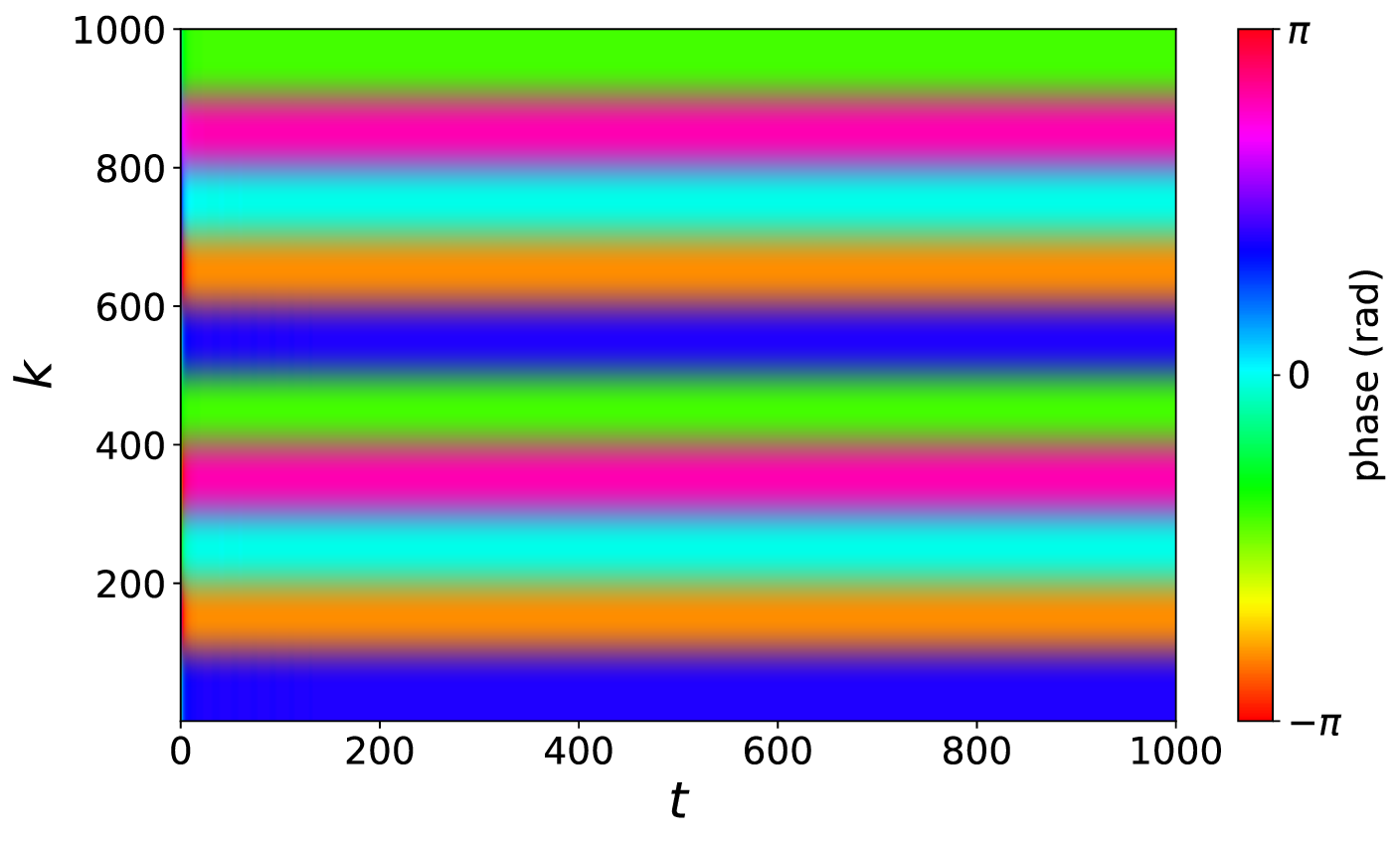}\\[-1ex]
{\footnotesize(g)}
\end{center}
\end{minipage}
\begin{minipage}[t]{0.495\textwidth}
\begin{center}
\includegraphics[scale=0.265]{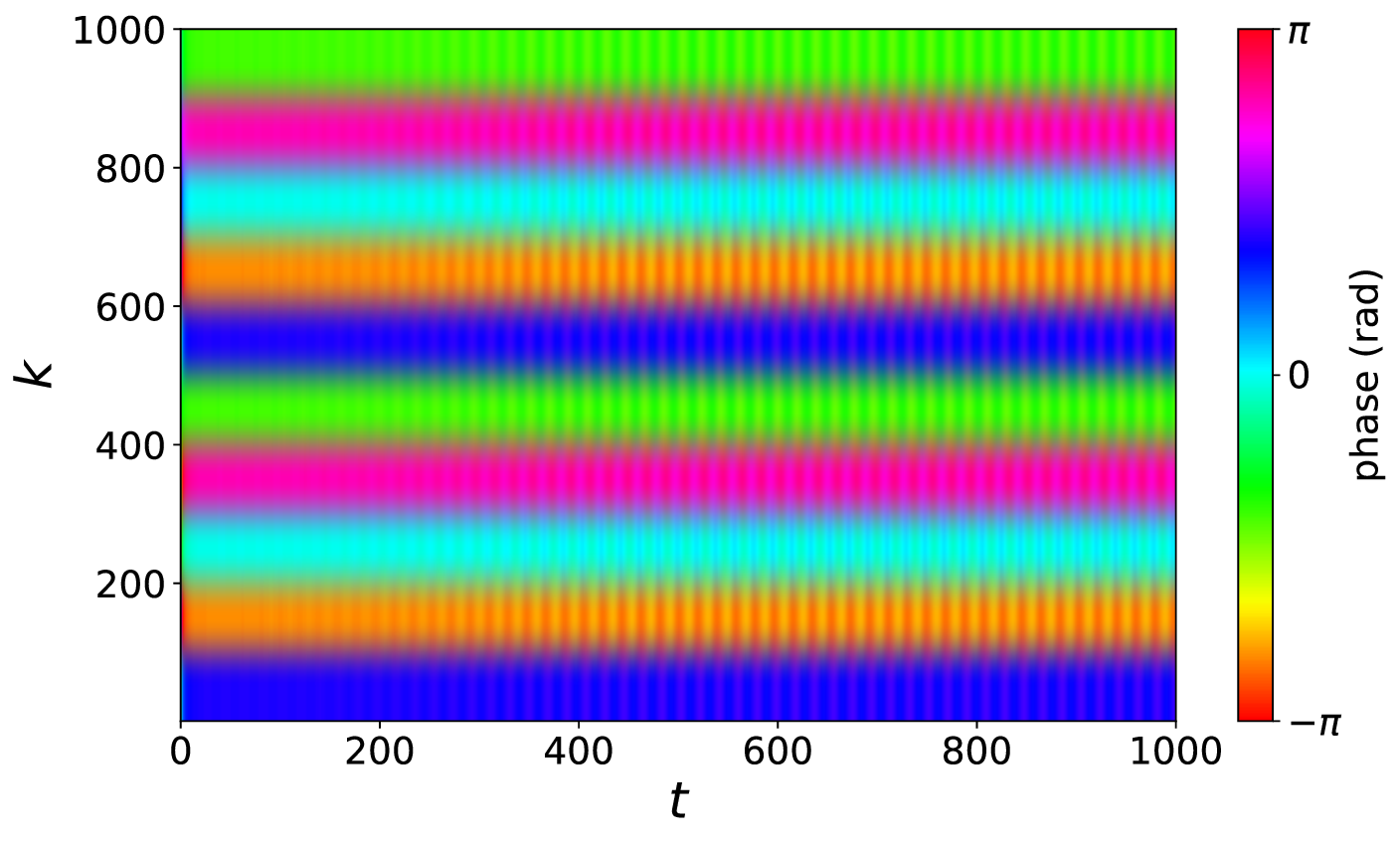}\\[-1ex]
{\footnotesize(h)}
\end{center}
\end{minipage}
\caption{Space-time plots of all oscillator phases $u_k^n(t)$, $k\in[n]$,
 in the KM \eqref{eqn:dsys} with $n=1000$, $\kappa=0.4$, $\sigma=\pi/3$ and $b_3=0.5$:
(a) $(q,b_1)=(1,0.08)$; (b) $(1,0.06)$;
(c) $(2,0.275)$; (d) $(2,0.255)$;
(e) $(3,0.17)$; (f) $(3,0.15)$;
(g) $(4,0.245)$; (h) $(4,0.225)$.
See also the caption of Fig.~\ref{fig:5a1}.
\vspace*{5mm}
}
\label{fig:5b1'}
\end{figure}

{\color{black}
Figures \ref{fig:5a1'} and \ref{fig:5b1'} show
 space-time plots of all oscillator phases $u_k^n(t)$, $k\in[n]$,
 for $\sigma=0$ and $\pi/3$, respectively.
Here the same values of $b_1$ and $u_k^n(0)$, $k\in[n]$,
 as in Figs.~\ref{fig:5a1} and \ref{fig:5b1} were used.
We see that all phases of the KM \eqref{eqn:dsys}
 converge to their steady states rapidly, as in Figs.~\ref{fig:5a1} and \ref{fig:5b1},
 and that they exhibit oscillations for the smaller values of $b_1$
 when $\sigma=\pi/3$, in the right column of Fig.~\ref{fig:5b1'},
 whereas only small differences between the larger and smaller values of $b_1$
 are observed in Fig.~\ref{fig:5a1'}.
Moreover, in both Figs.~\ref{fig:5a1'} and \ref{fig:5b1'},
 the variation in $u_k^n(t)$ from $k=1$ to $n$ increases as $q$ increases.}
 
\begin{figure}[t]
\begin{minipage}[t]{0.495\textwidth}
\begin{center}
\includegraphics[scale=0.245]{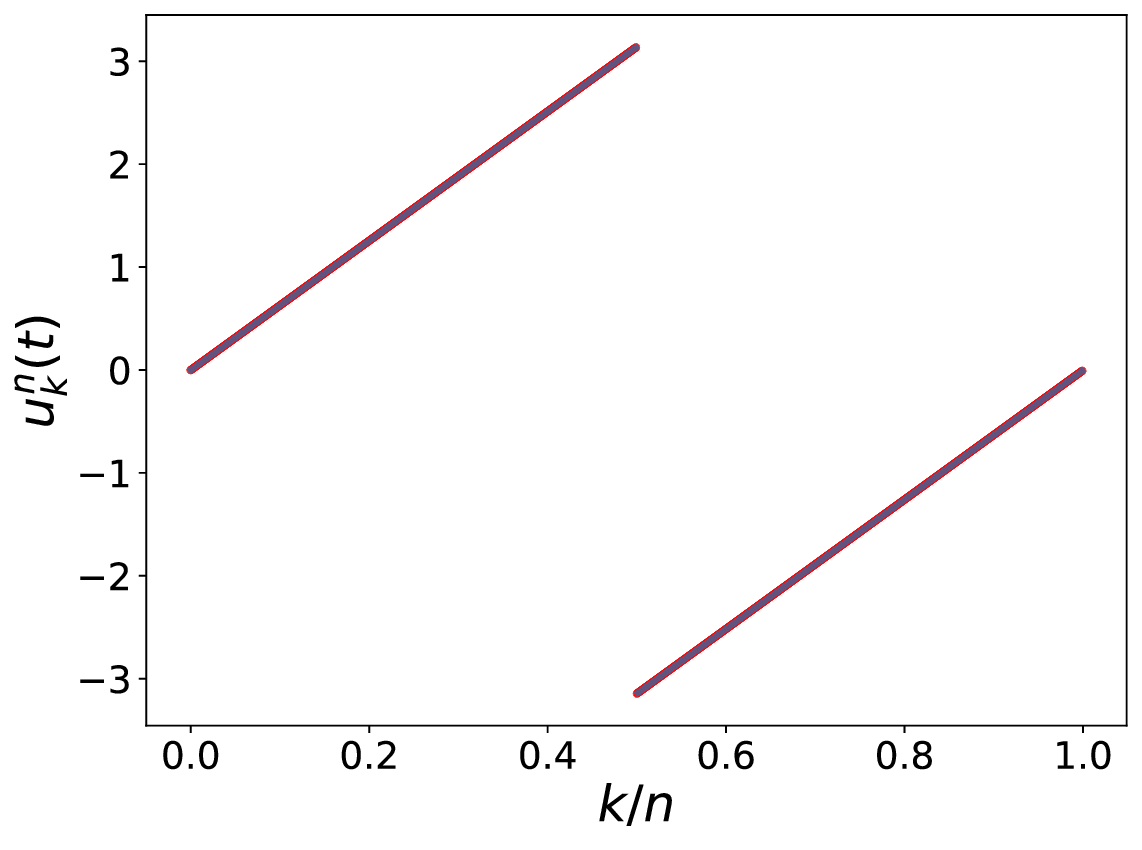}\\[-1ex]
{\footnotesize(a)}
\end{center}
\end{minipage}
\begin{minipage}[t]{0.495\textwidth}
\begin{center}
\includegraphics[scale=0.245]{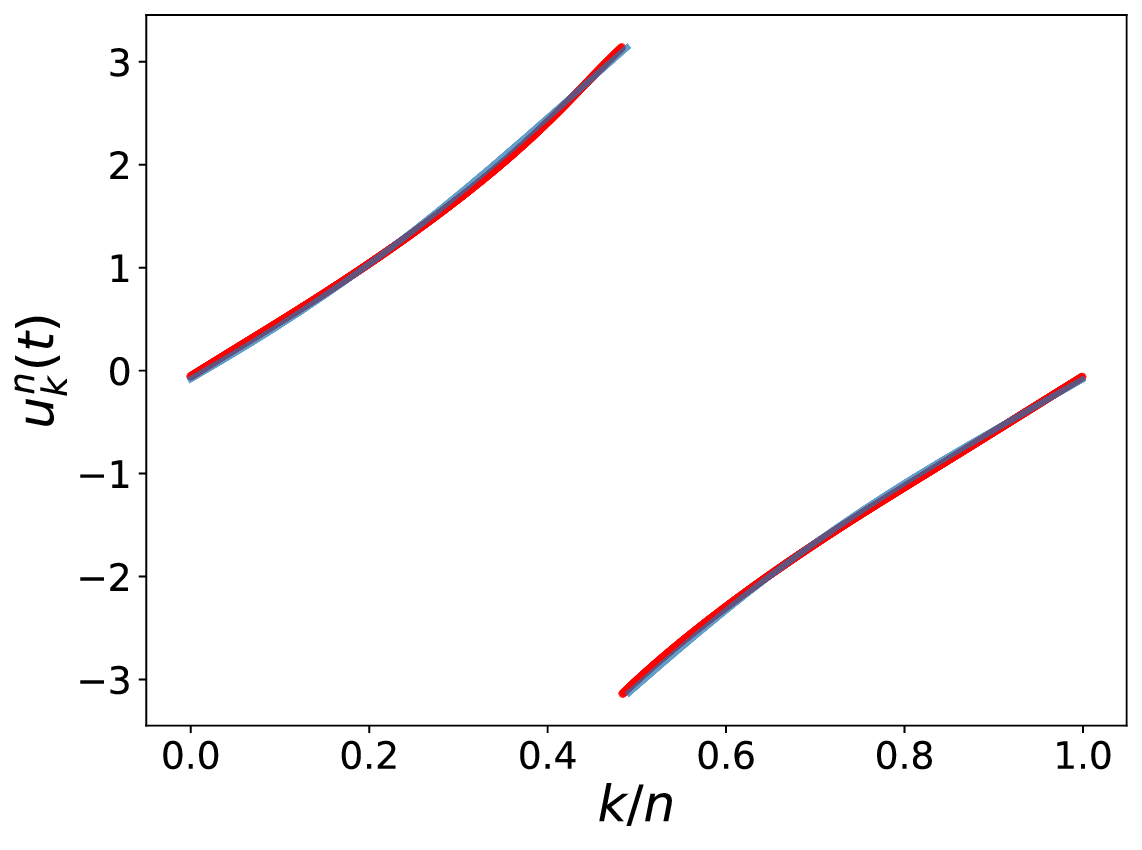}\\[-1ex]
{\footnotesize(b)}
\end{center}
\end{minipage}
\vspace*{0.5ex}

\begin{minipage}[t]{0.495\textwidth}
\begin{center}
\includegraphics[scale=0.245]{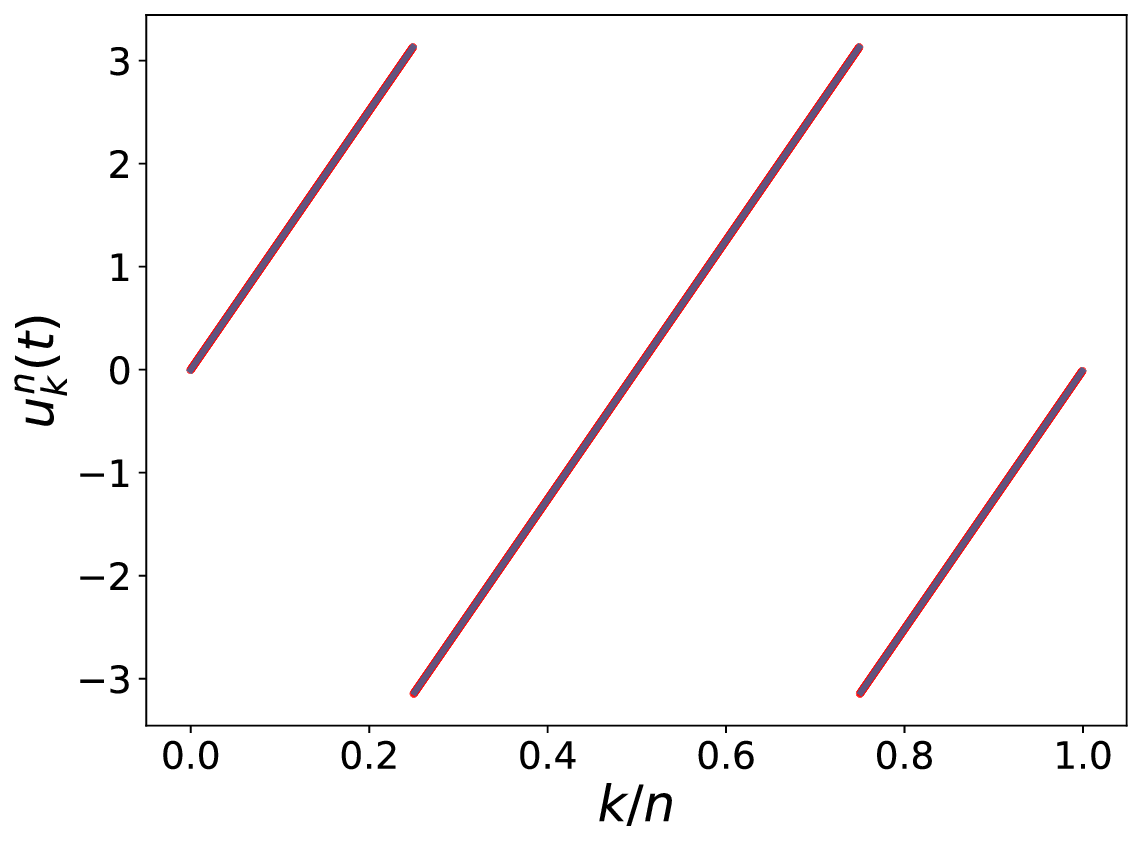}\\[-1ex]
{\footnotesize(c)}
\end{center}
\end{minipage}
\begin{minipage}[t]{0.495\textwidth}
\begin{center}
\includegraphics[scale=0.245]{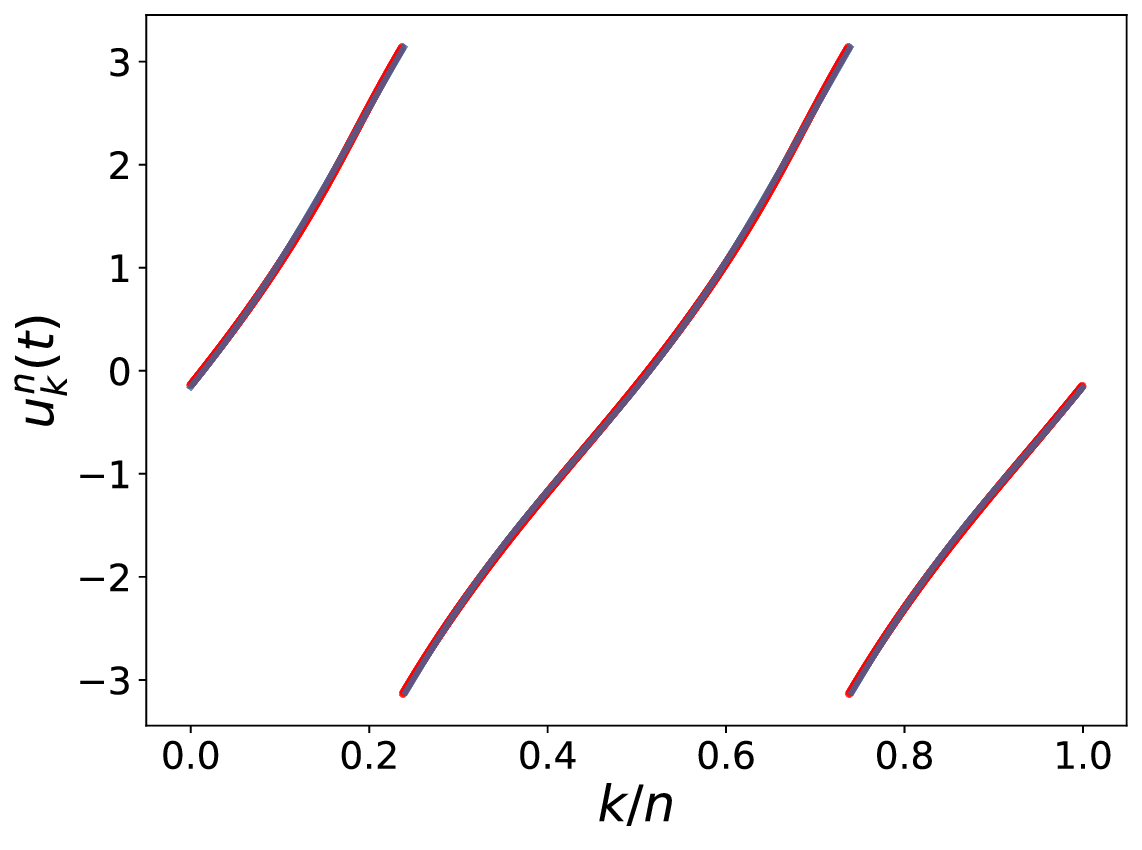}\\[-1ex]
{\footnotesize(d)}
\end{center}
\end{minipage}
\vspace*{0.5ex}

\begin{minipage}[t]{0.495\textwidth}
\begin{center}
\includegraphics[scale=0.245]{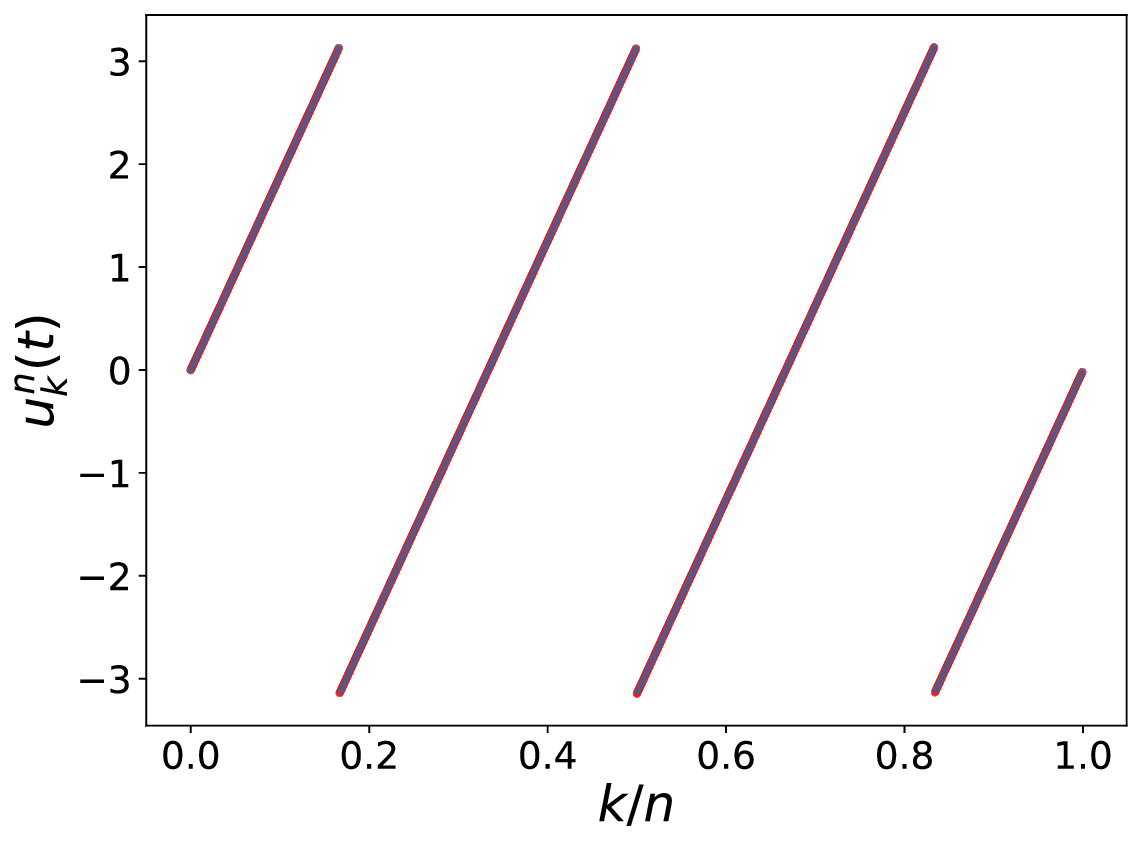}\\[-1ex]
{\footnotesize(e)}
\end{center}
\end{minipage}
\begin{minipage}[t]{0.495\textwidth}
\begin{center}
\includegraphics[scale=0.245]{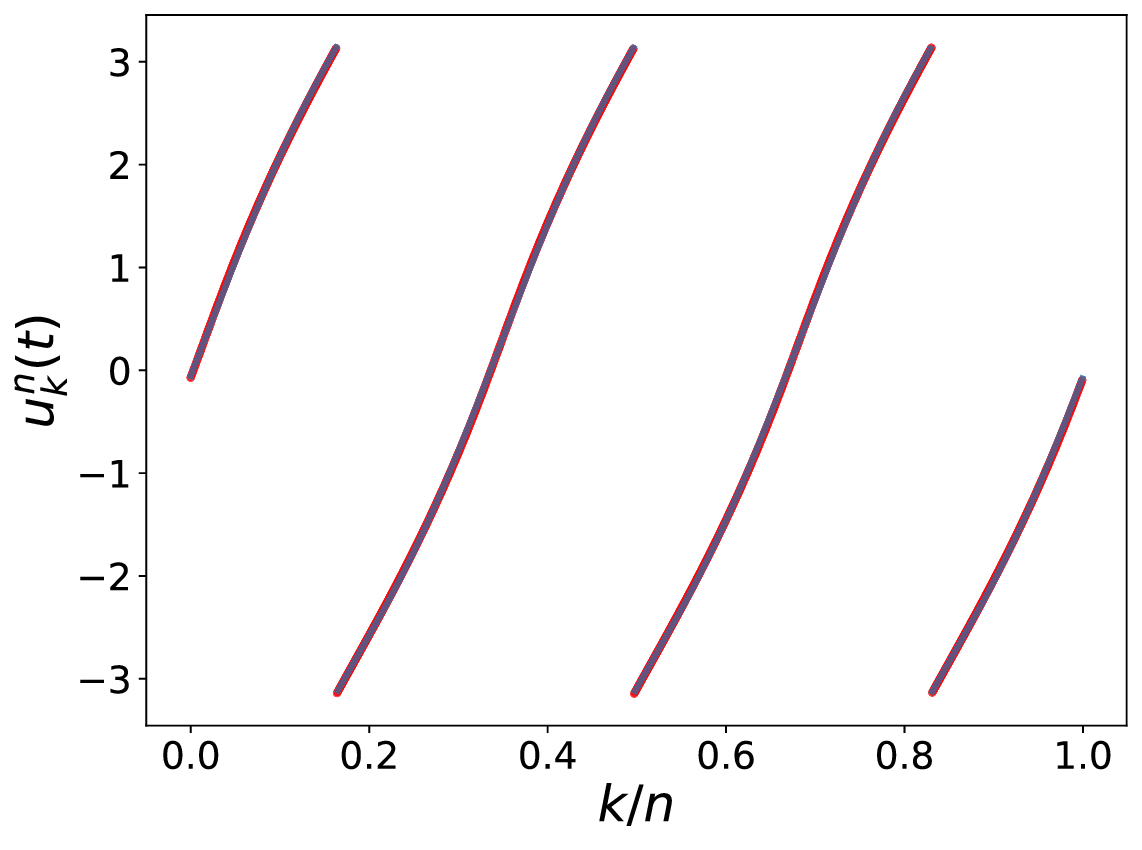}\\[-1ex]
{\footnotesize(f)}
\end{center}
\end{minipage}
\vspace*{0.5ex}

\begin{minipage}[t]{0.495\textwidth}
\begin{center}
\includegraphics[scale=0.245]{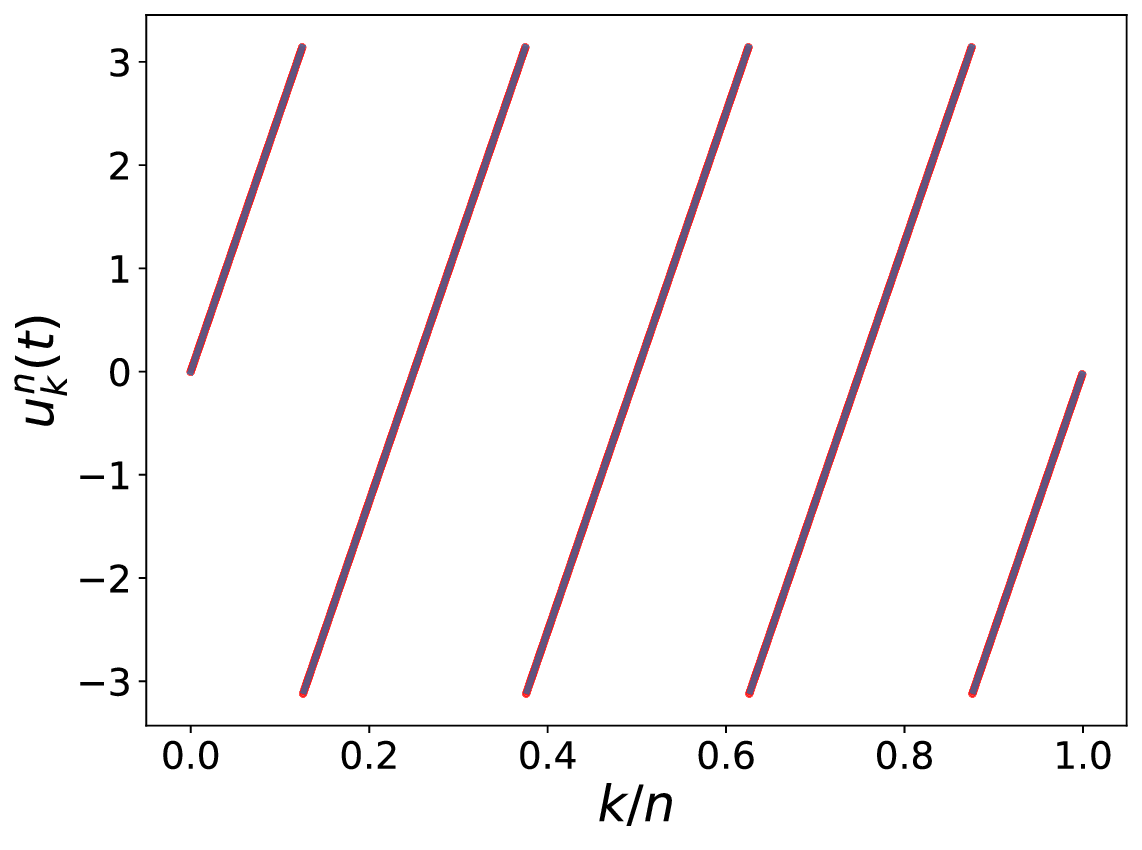}\\[-1ex]
{\footnotesize(g)}
\end{center}
\end{minipage}
\begin{minipage}[t]{0.495\textwidth}
\begin{center}
\includegraphics[scale=0.245]{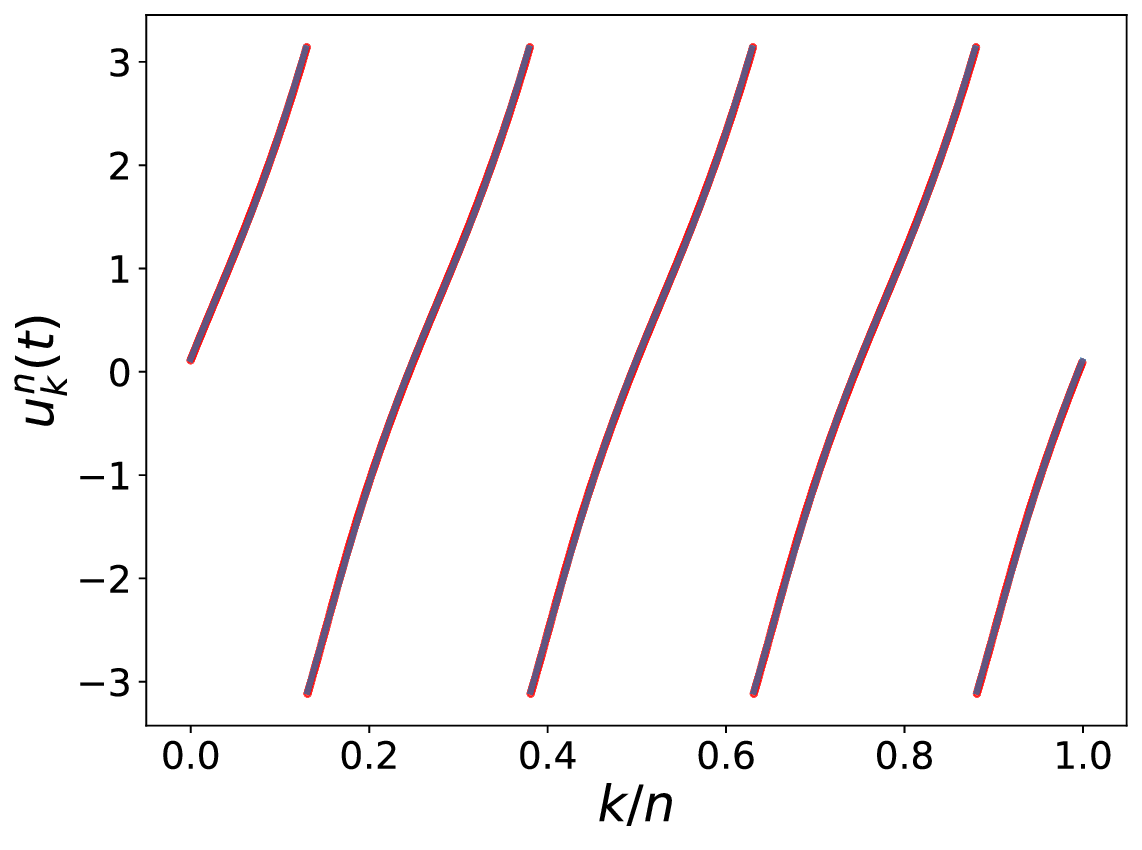}\\[-1ex]
{\footnotesize(h)}
\end{center}
\end{minipage}
\caption{Steady states of the KM \eqref{eqn:dsys}
 with $n=1000$ and $\sigma=0$ at $t=1000$:
(a) $(q,b_1,b_3)=(1,0.16,1)$; (b) $(1,0.12,1)$;
(c) $(2,0.55,0.5)$; (d) $(2,0.51,0.5)$;
(e) $(3,0.34,0.5)$; (f) $(3,0.3,0.5)$;
(g) $(4,0.49,0.5)$; (h) $(4,0.45,0.5)$.
The values of $u_k^n(t)\mod 2\pi$, $k\in[n]$, are plotted as the ordinates.
The simulation results are plotted as small red disks
 and the most probable leading terms in \eqref{eqn:thm4a} and \eqref{eqn:thm4b}
 estimated from them as blue lines
 although they coincide almost completely.}
\label{fig:5a2}
\end{figure}

\begin{figure}[t]
\begin{minipage}[t]{0.495\textwidth}
\begin{center}
\includegraphics[scale=0.245]{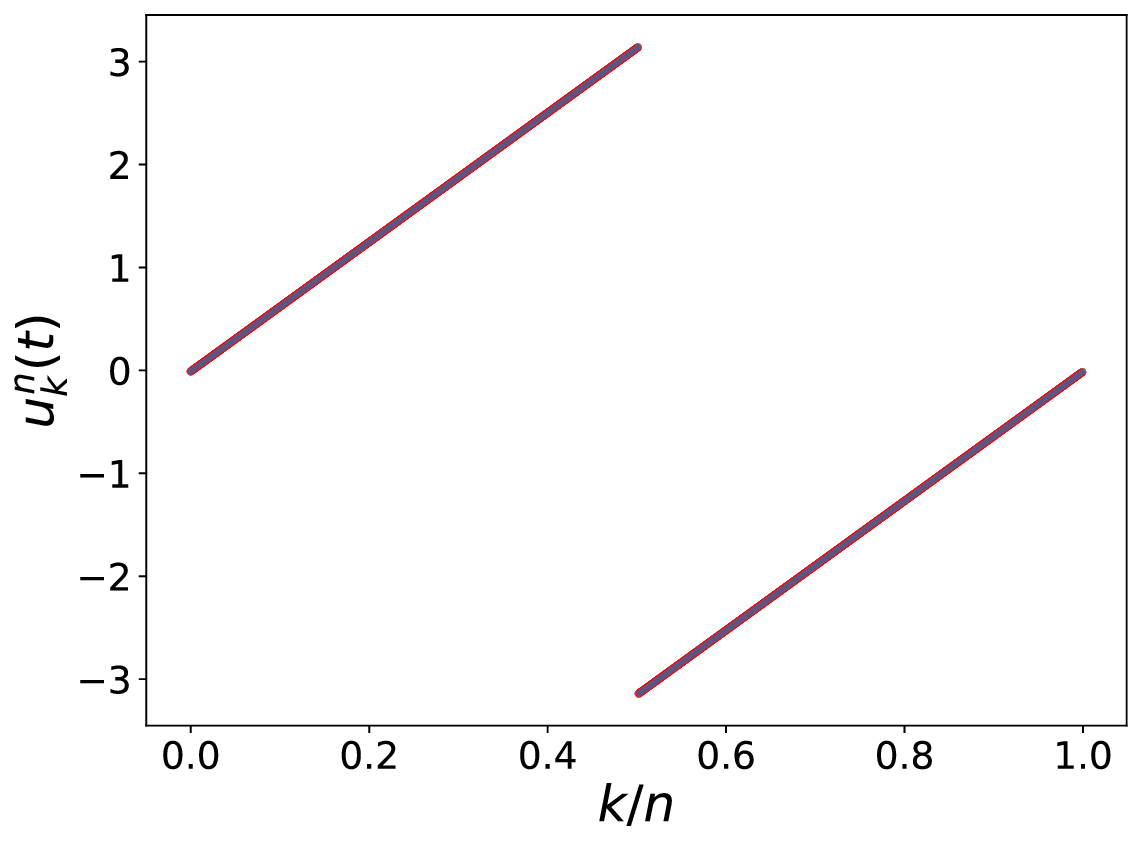}\\[-1ex]
{\footnotesize(a)}
\end{center}
\end{minipage}
\begin{minipage}[t]{0.495\textwidth}
\begin{center}
\includegraphics[scale=0.245]{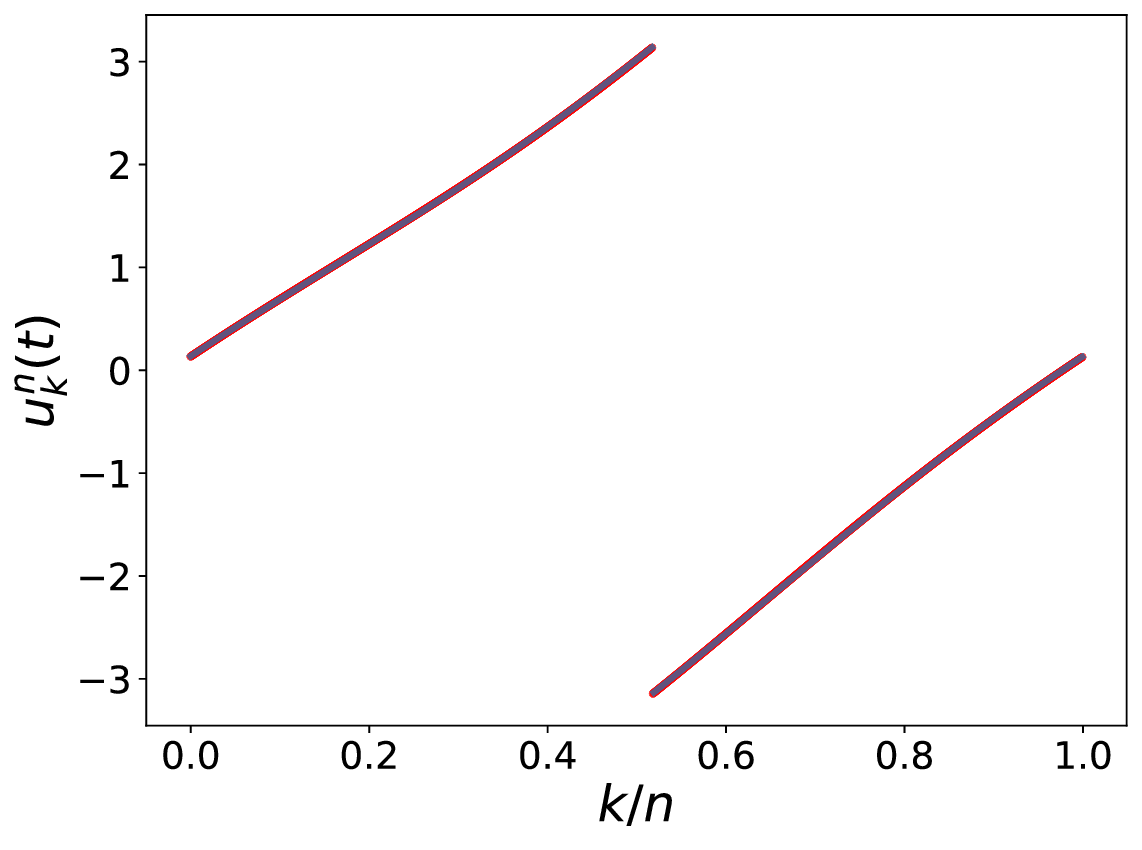}\\[-1ex]
{\footnotesize(b)}
\end{center}
\end{minipage}
\vspace*{0.5ex}

\begin{minipage}[t]{0.495\textwidth}
\begin{center}
\includegraphics[scale=0.245]{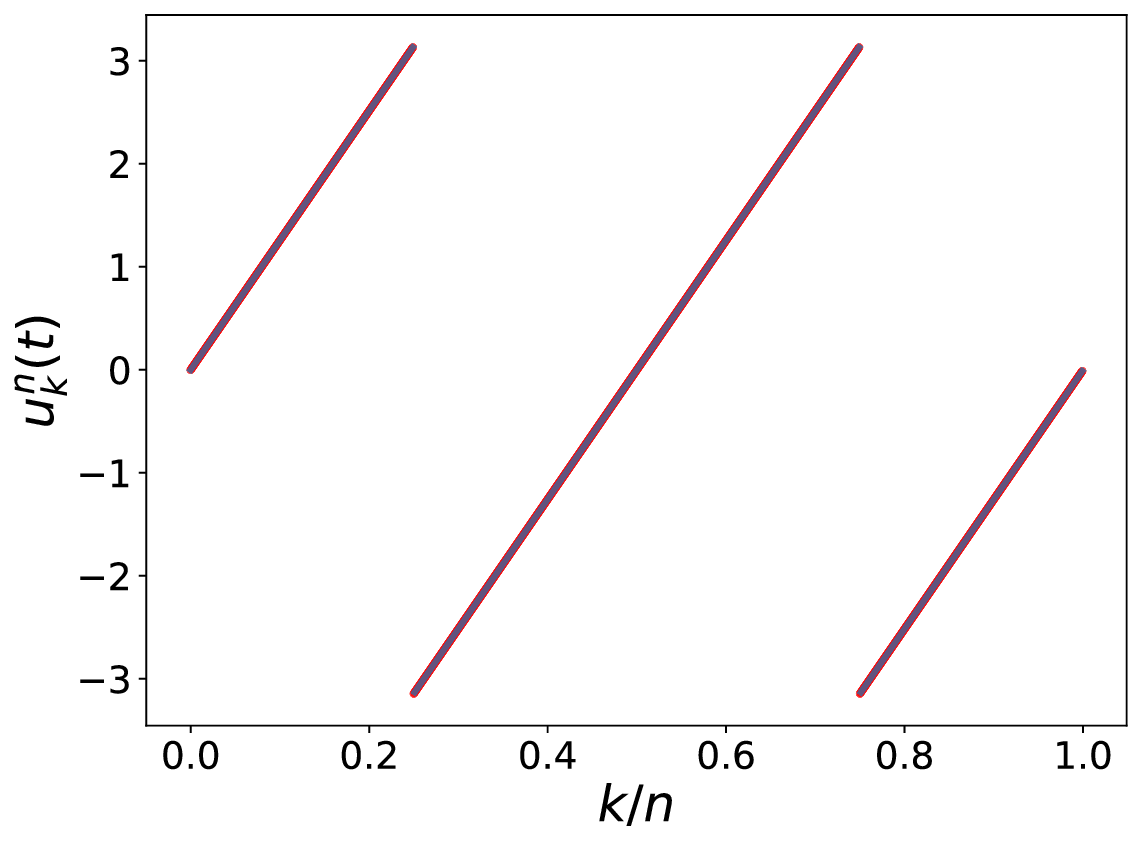}\\[-1ex]
{\footnotesize(c)}
\end{center}
\end{minipage}
\begin{minipage}[t]{0.495\textwidth}
\begin{center}
\includegraphics[scale=0.245]{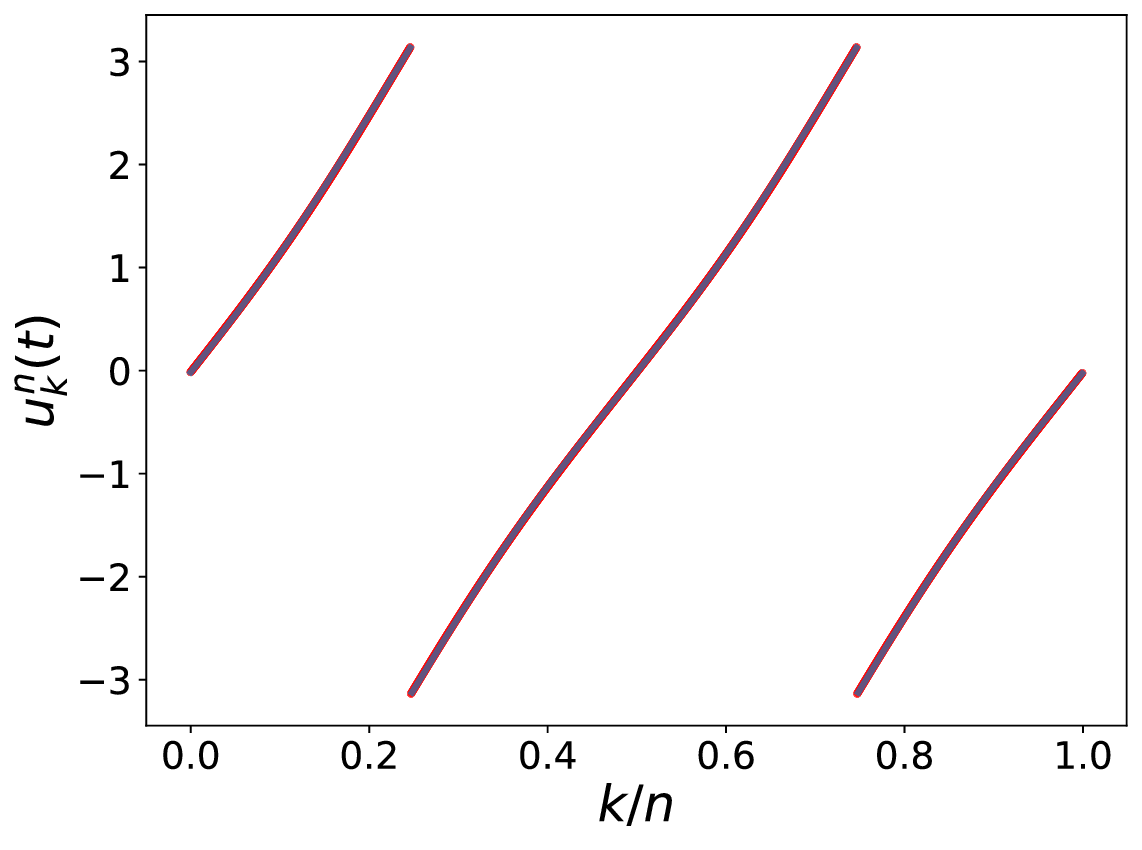}\\[-1ex]
{\footnotesize(d)}
\end{center}
\end{minipage}
\vspace*{0.5ex}

\begin{minipage}[t]{0.495\textwidth}
\begin{center}
\includegraphics[scale=0.245]{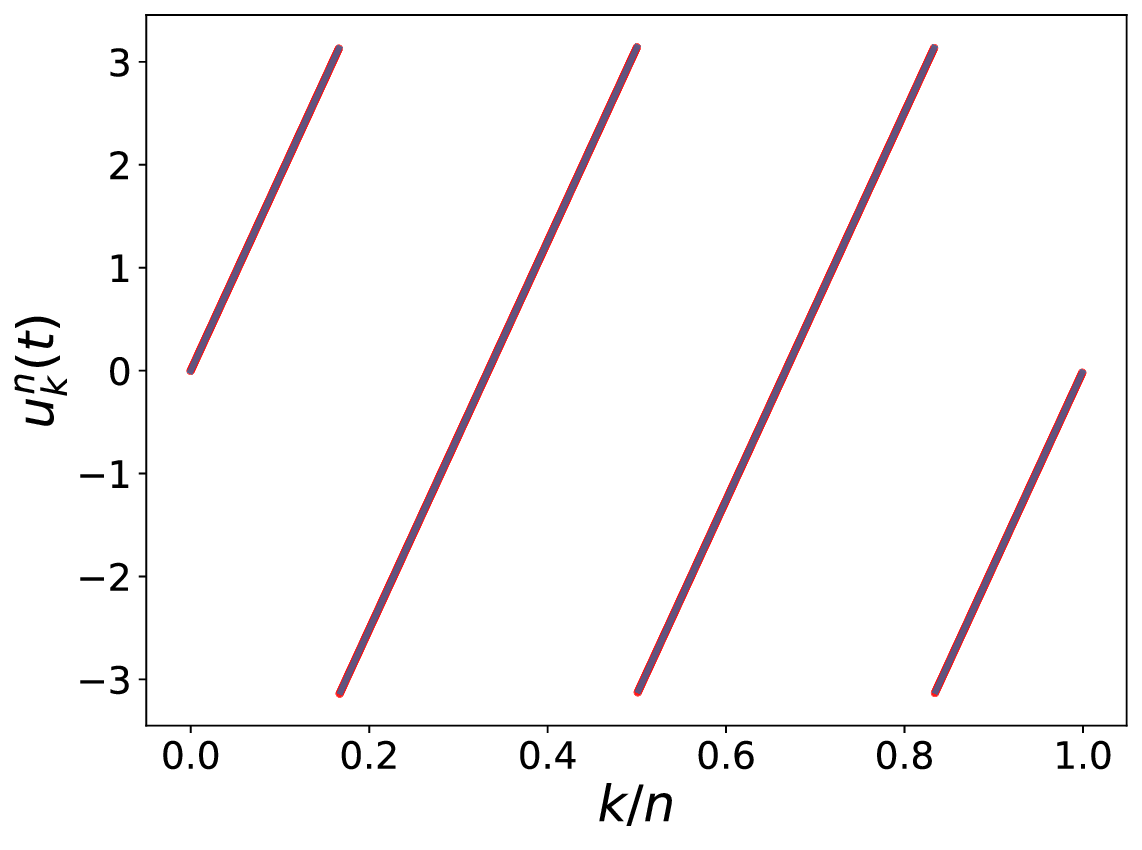}\\[-1ex]
{\footnotesize(e)}
\end{center}
\end{minipage}
\begin{minipage}[t]{0.495\textwidth}
\begin{center}
\includegraphics[scale=0.245]{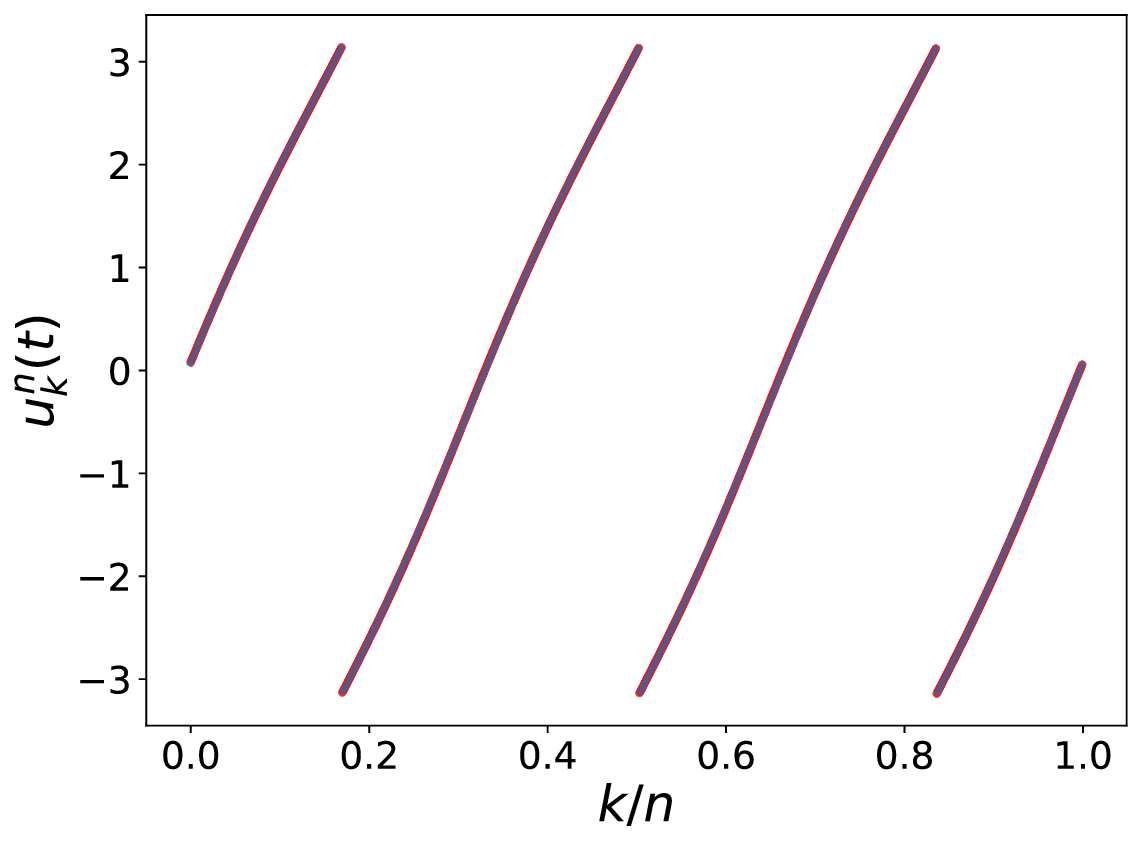}\\[-1ex]
{\footnotesize(f)}
\end{center}
\end{minipage}
\vspace*{0.5ex}

\begin{minipage}[t]{0.495\textwidth}
\begin{center}
\includegraphics[scale=0.245]{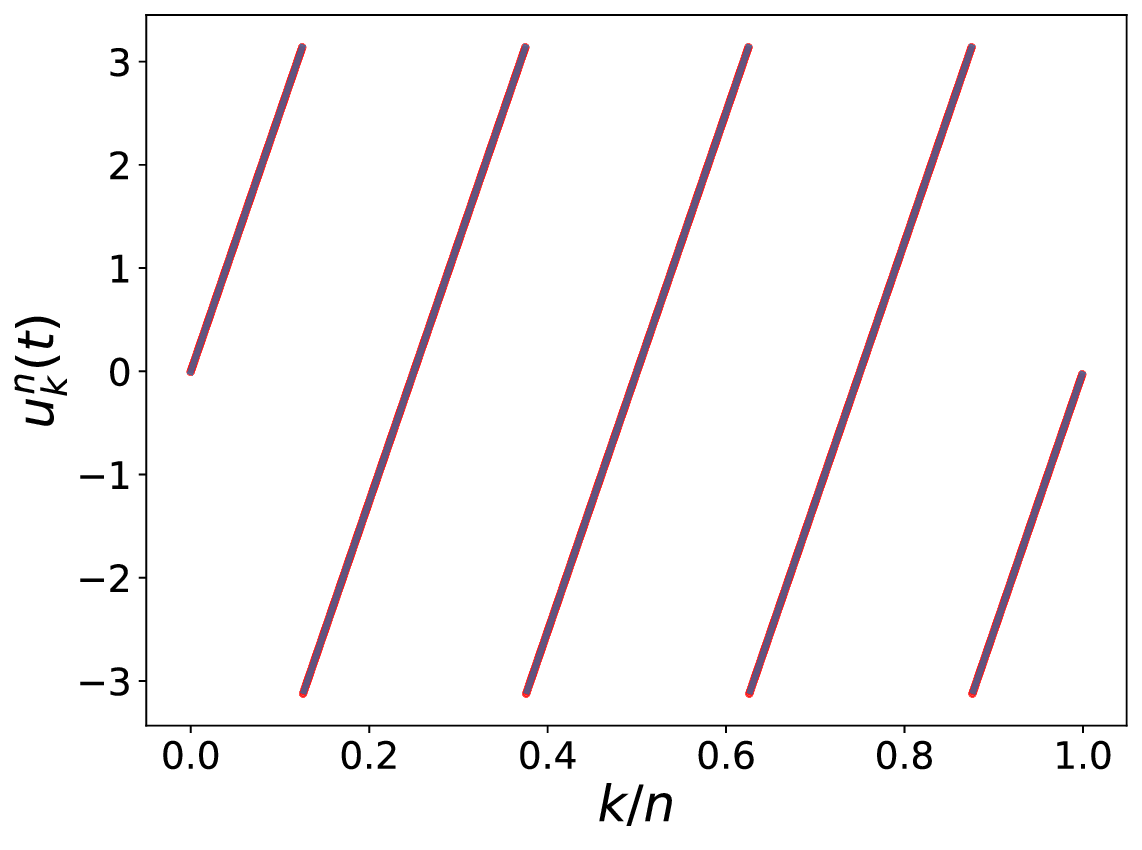}\\[-1ex]
{\footnotesize(g)}
\end{center}
\end{minipage}
\begin{minipage}[t]{0.495\textwidth}
\begin{center}
\includegraphics[scale=0.245]{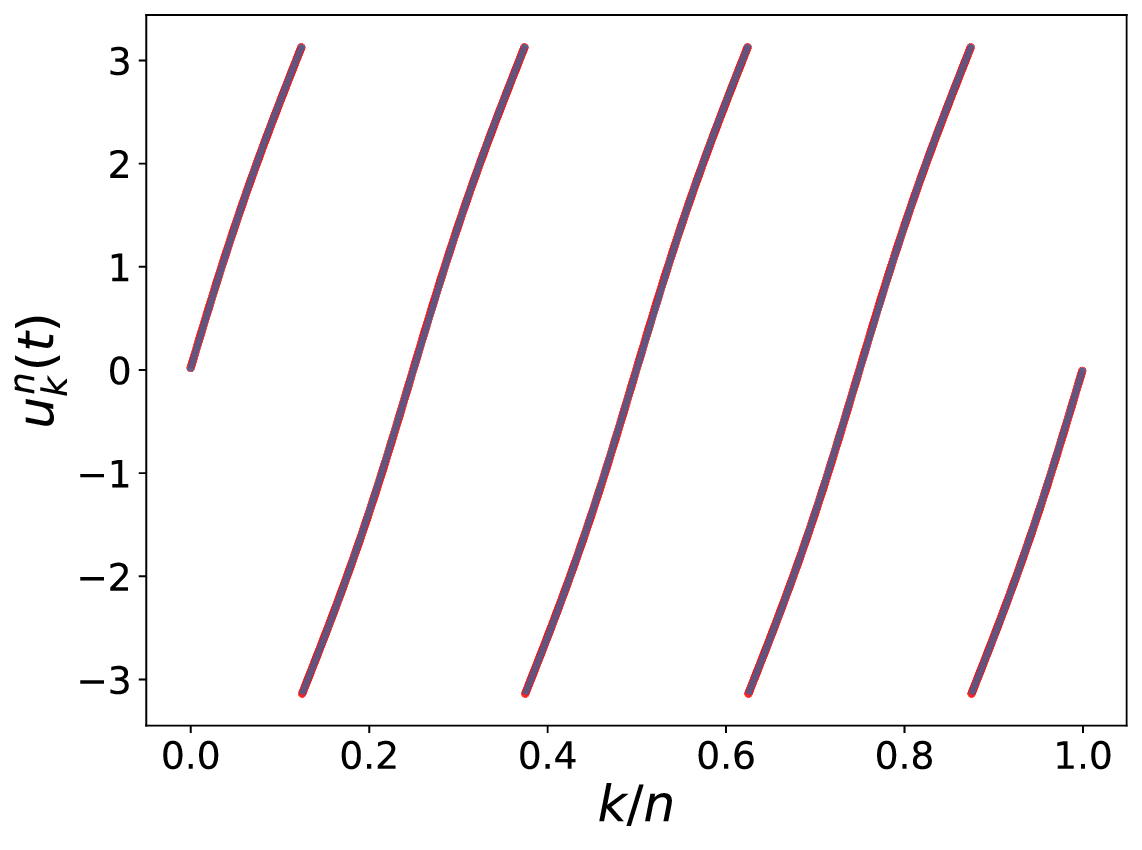}\\[-1ex]
{\footnotesize(h)}
\end{center}
\end{minipage}
\caption{Steady states of the KM \eqref{eqn:dsys}
 with $n=1000$, $\sigma=\pi/3$ and $b_3=0.5$ at $t=1000$:
(a) $(q,b_1)=(1,0.08)$; (b) $(1,0.06)$;
(c) $(2,0.275)$; (d) $(2,0.255)$;
(e) $(3,0.17)$; (f) $(3,0.15)$;
(g) $(4,0.245)$; (h) $(4,0.225)$.
See also the caption of Fig.~\ref{fig:5a2}.}
\label{fig:5b2}
\end{figure}

In Figs.~\ref{fig:5a2} and \ref{fig:5b2},
 $u_k^n(t)$, $k\in[n]$, at $t=1000$,
 which may be regarded as the steady states
 from Figs.~\ref{fig:5a1} and \ref{fig:5b1},
 are plotted as small red disks  for $\sigma=0$ and $\pi/3$, respectively.
Here the same values of $b_1$ and $u_k^n(0)$, $k\in[n]$,
 as in Figs.~\ref{fig:5a1} and \ref{fig:5b1} were used.
We observe that the responses of the KM \eqref{eqn:dsys}
 converge to the twisted and modulated or oscillating twisted states,
 respectively, for the larger and smaller values of $b_1$,
 as predicted by Theorems~\ref{thm:4a} and \ref{thm:4b}
 with the assistance of Corollary~\ref{cor:2a} and Theorem~\ref{thm:2e}.
Indeed, we confirmed that
 the deviation from the twisted state is about $10^{-7}$ at most
 in the left column of each figure for the larger values of $b_1$.
In particular, the target state \eqref{eqn:ts}
 is accomplished there.

The most probably leading term,
\begin{equation}
u(x)=2\pi qx+r(t)\sin(2\pi qx+\psi(t))+\Omega t,
\label{eqn:ss}
\end{equation}
in the modulated and oscillating twisted solutions \eqref{eqn:thm4a} and \eqref{eqn:thm4b}
 was estimated from the numerical simulation results for each cases
 by using the least mean square method as
\begin{equation}
\Omega t=\frac{1}{n}\sum_{k=1}^nv_k^n(t),\quad
r(t)=2\sqrt{c(t)^2+s(t)^2}
\label{eqn:ss1}
\end{equation}
and
\begin{align*}
\psi(t)=&\arctan\frac{s(t)}{c(t)}\quad
\left(\mbox{resp. }\arctan\frac{s(t)}{c(t)}+pi\mbox{ or } 
\arctan\frac{s(t)}{c(t)}-\pi\right)
\end{align*}
for $c(t)>0$ (resp. $c(t)<0$ and  $s(t)>0$ or  $s(t)<0$) with
\[
v_k^n(t)=u_k^n(t)-\frac{2\pi q k}{n}
\]
and
\[
c(t)=\frac{1}{n}\sum_{k=1}^n v_k^n(t)\cos\frac{2\pi k}{n},\quad
s(t)=\frac{1}{n}\sum_{k=1}^n v_k^n(t)\sin\frac{2\pi k}{n},
\]
and it is plotted as a blue line in each figure.
The computed results coincide with the simulation results
 for the KM \eqref{eqn:dsys} almost completely,
 as detected by Theorems~\ref{thm:4a} and \ref{thm:4b}
 for the CL \eqref{eqn:csys}.
Here the estimated values of $\Omega t$ in \eqref{eqn:ss1} were very small,
 and more precisely about $10^{-12}$ and $10^{-2}$ at most
 for $\sigma=0$ and $\pi/3$, respectively, even when $t=1000$.
 
\begin{figure}[t]
\begin{minipage}[t]{0.495\textwidth}
\begin{center}
\includegraphics[scale=0.245]{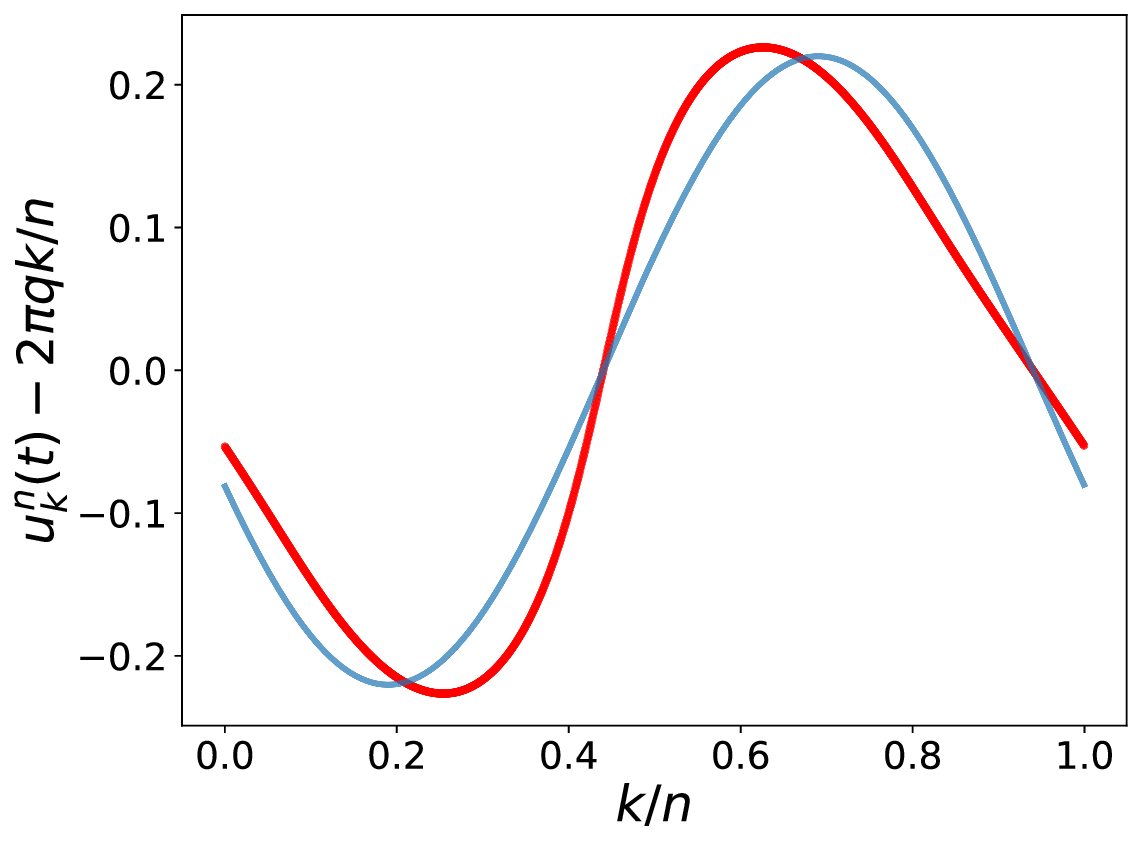}\\[-1ex]
{\footnotesize(a)}
\end{center}
\end{minipage}
\begin{minipage}[t]{0.495\textwidth}
\begin{center}
\includegraphics[scale=0.245]{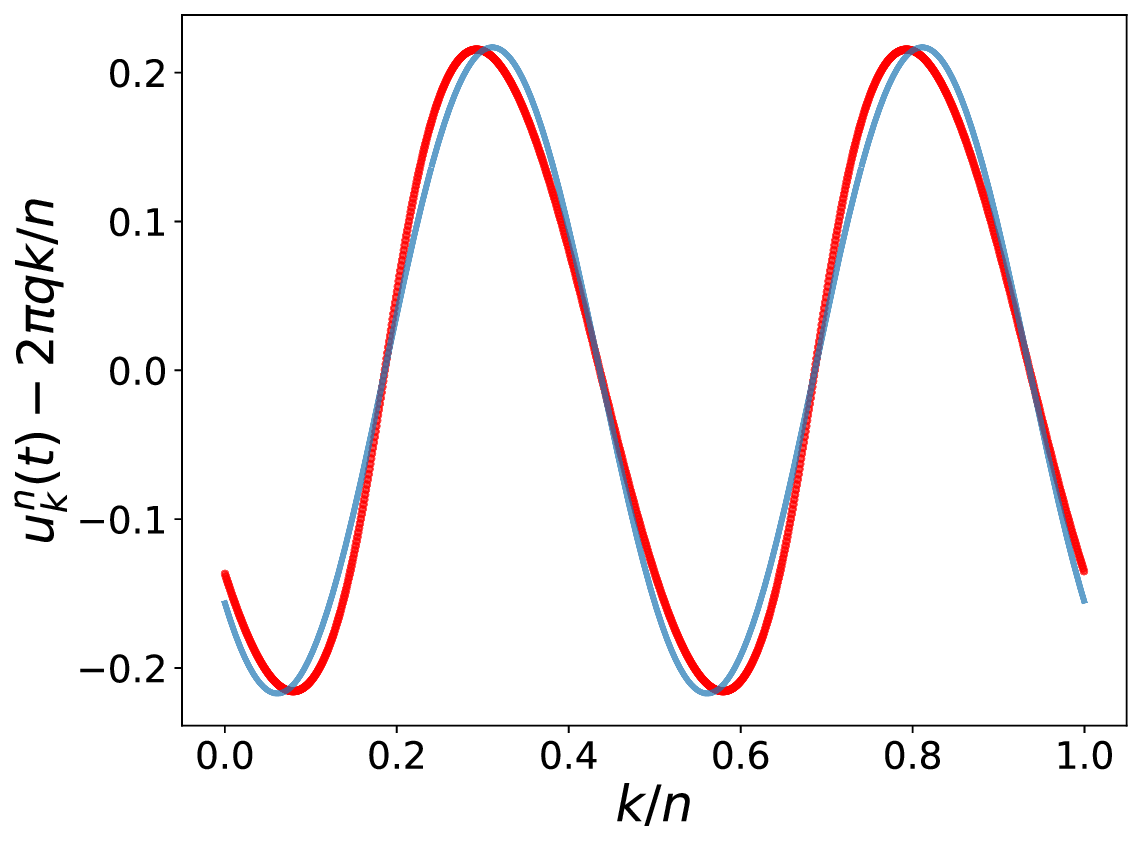}\\[-1ex]
{\footnotesize(b)}
\end{center}
\end{minipage}
\vspace*{0.5ex}

\begin{minipage}[t]{0.495\textwidth}
\begin{center}
\includegraphics[scale=0.245]{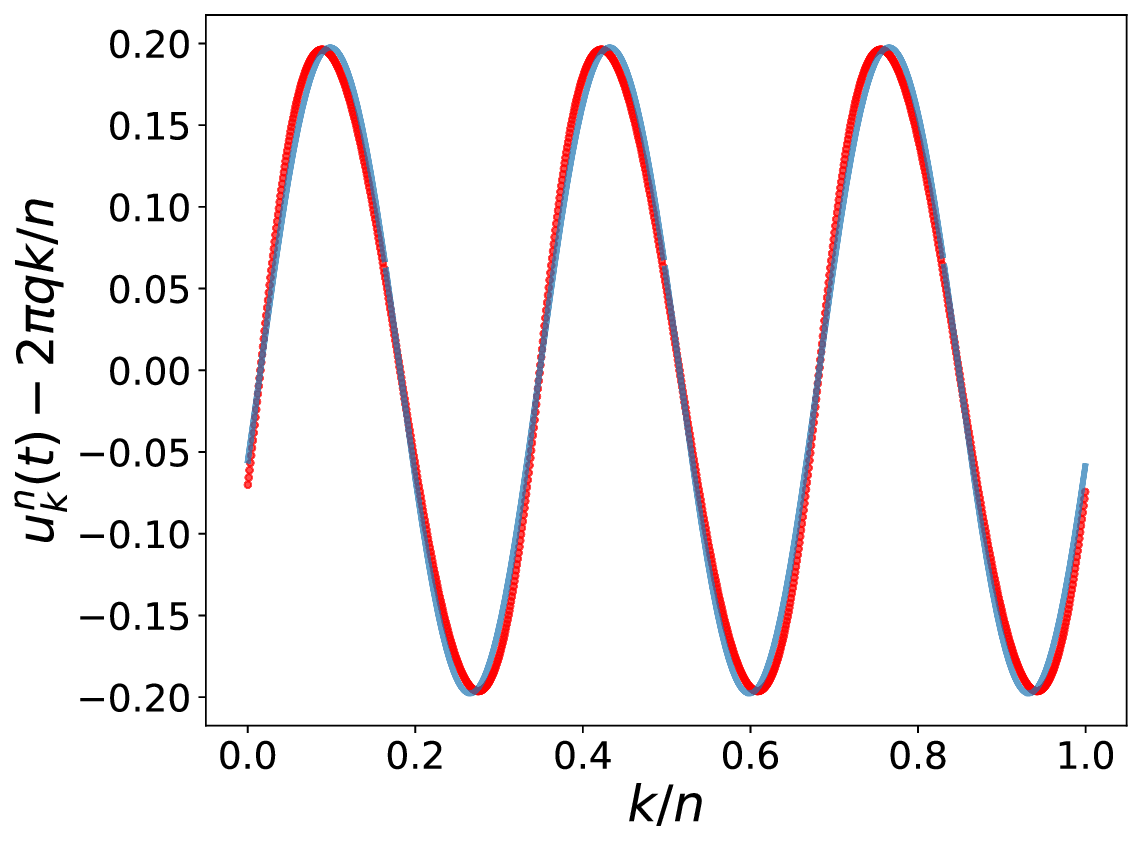}\\[-1ex]
{\footnotesize(c)}
\end{center}
\end{minipage}
\begin{minipage}[t]{0.495\textwidth}
\begin{center}
\includegraphics[scale=0.245]{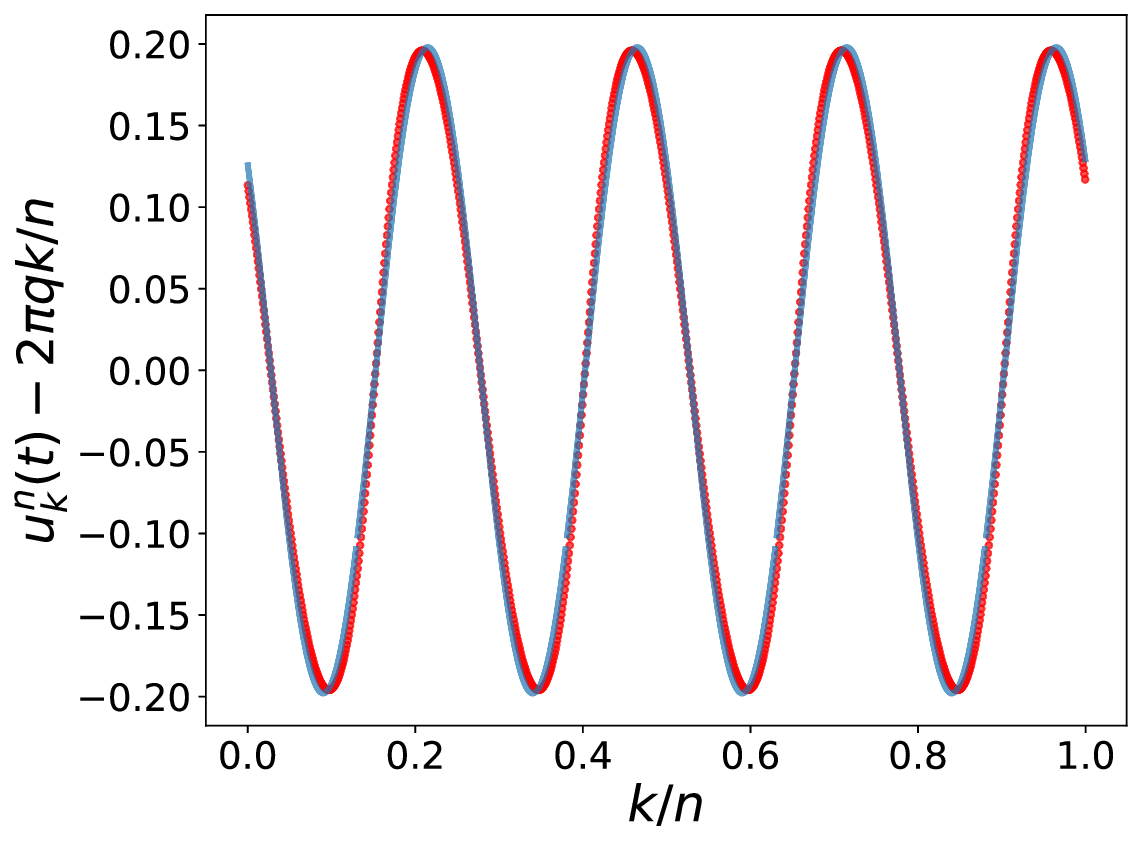}\\[-1ex]
{\footnotesize(d)}
\end{center}
\end{minipage}
\caption{Deviation from the $q$-twisted states
 in the steady states of the KM \eqref{eqn:dsys}
 with $n=1000$ and $\sigma=0$ at $t=1000$:
(a) $(q,b_1,b_3)=(1,0.12,1)$;
(b) $(2,0.51,0.5)$;
(c) $(3,0.3,0.5)$;
(d) $(4,0.45,0.5)$.
The simulation results are plotted as small red disks,
 and the estimates from their most probably leading terms
 given by \eqref{eqn:ss} are plotted as blue lines.}
\label{fig:5a3}
\end{figure}

\begin{figure}[t]
\begin{minipage}[t]{0.495\textwidth}
\begin{center}
\includegraphics[scale=0.245]{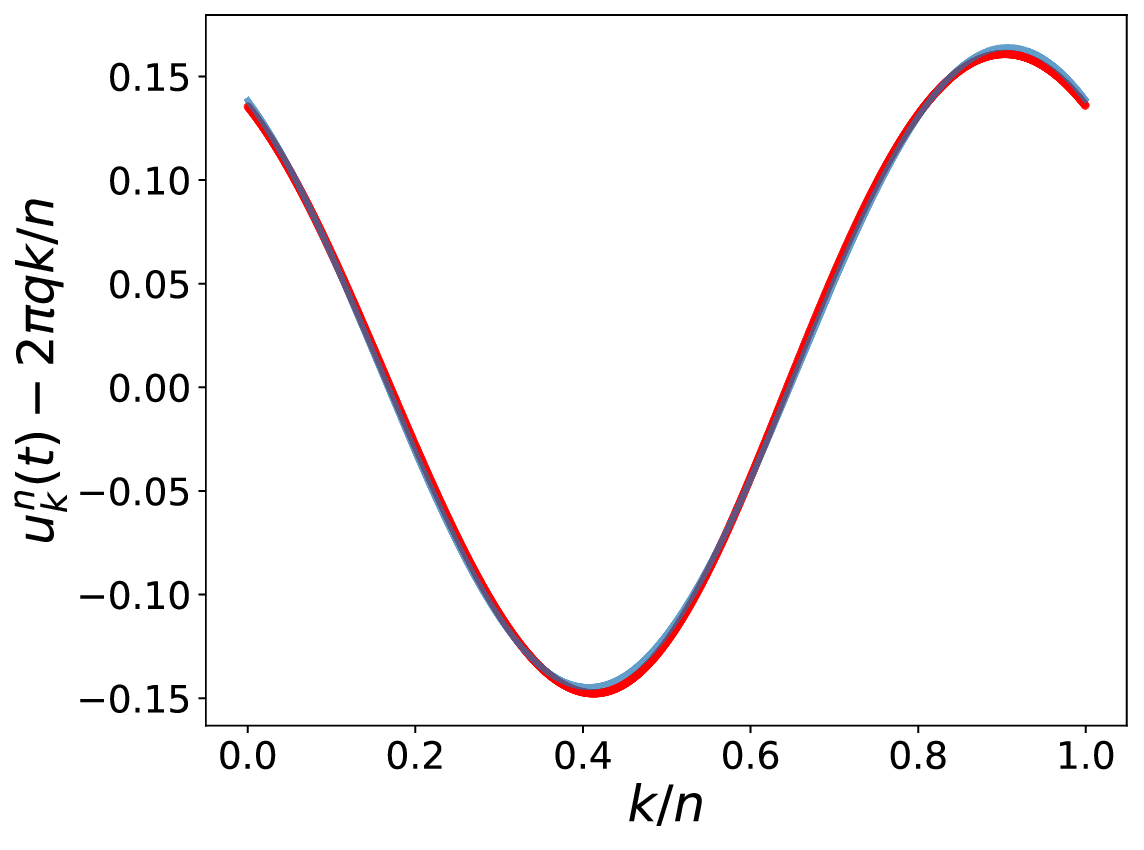}\\[-1ex]
{\footnotesize(a)}
\end{center}
\end{minipage}
\begin{minipage}[t]{0.495\textwidth}
\begin{center}
\includegraphics[scale=0.245]{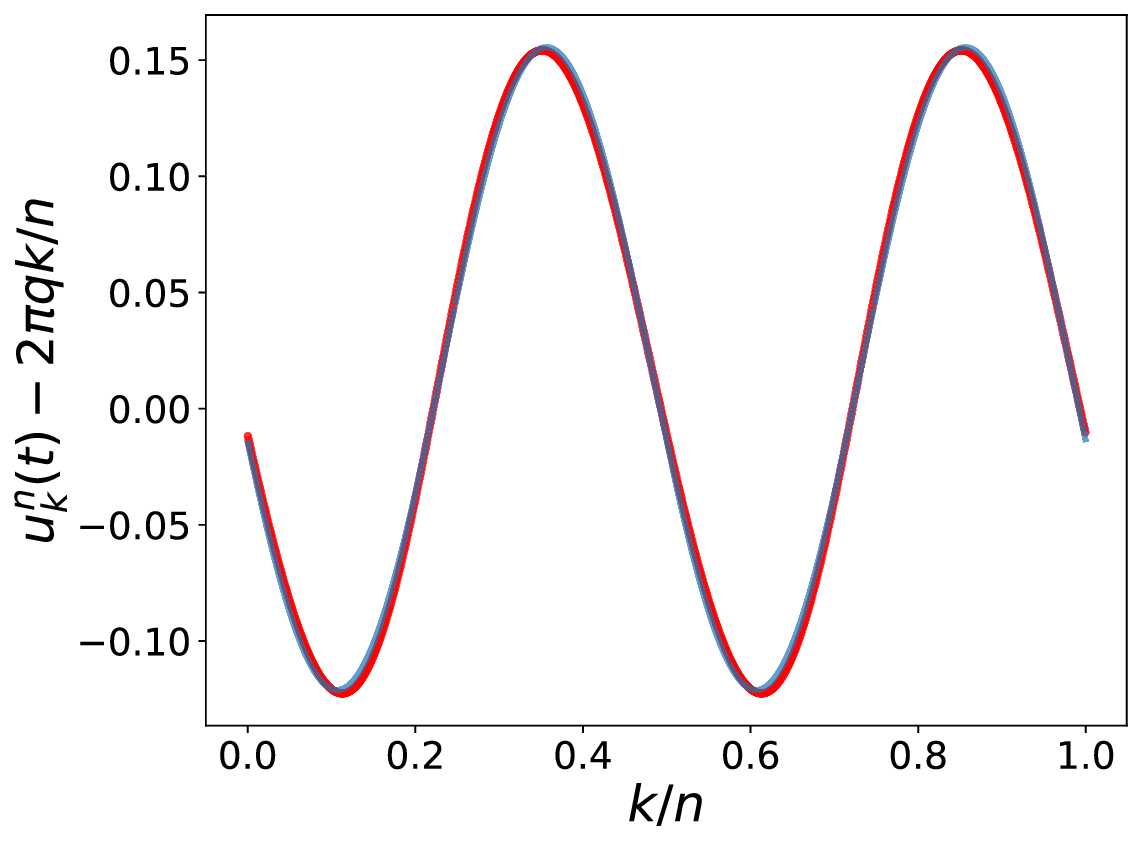}\\[-1ex]
{\footnotesize(b)}
\end{center}
\end{minipage}
\vspace*{0.5ex}

\begin{minipage}[t]{0.495\textwidth}
\begin{center}
\includegraphics[scale=0.245]{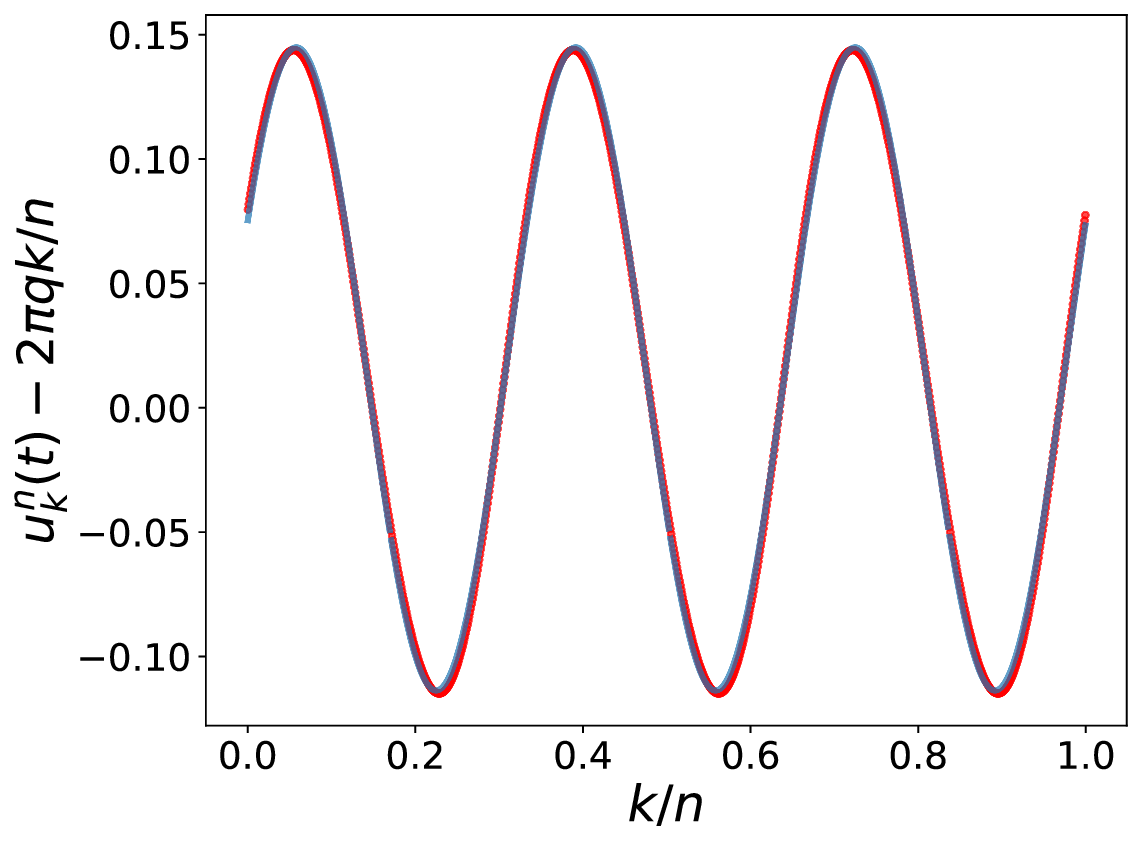}\\[-1ex]
{\footnotesize(c)}
\end{center}
\end{minipage}
\begin{minipage}[t]{0.495\textwidth}
\begin{center}
\includegraphics[scale=0.245]{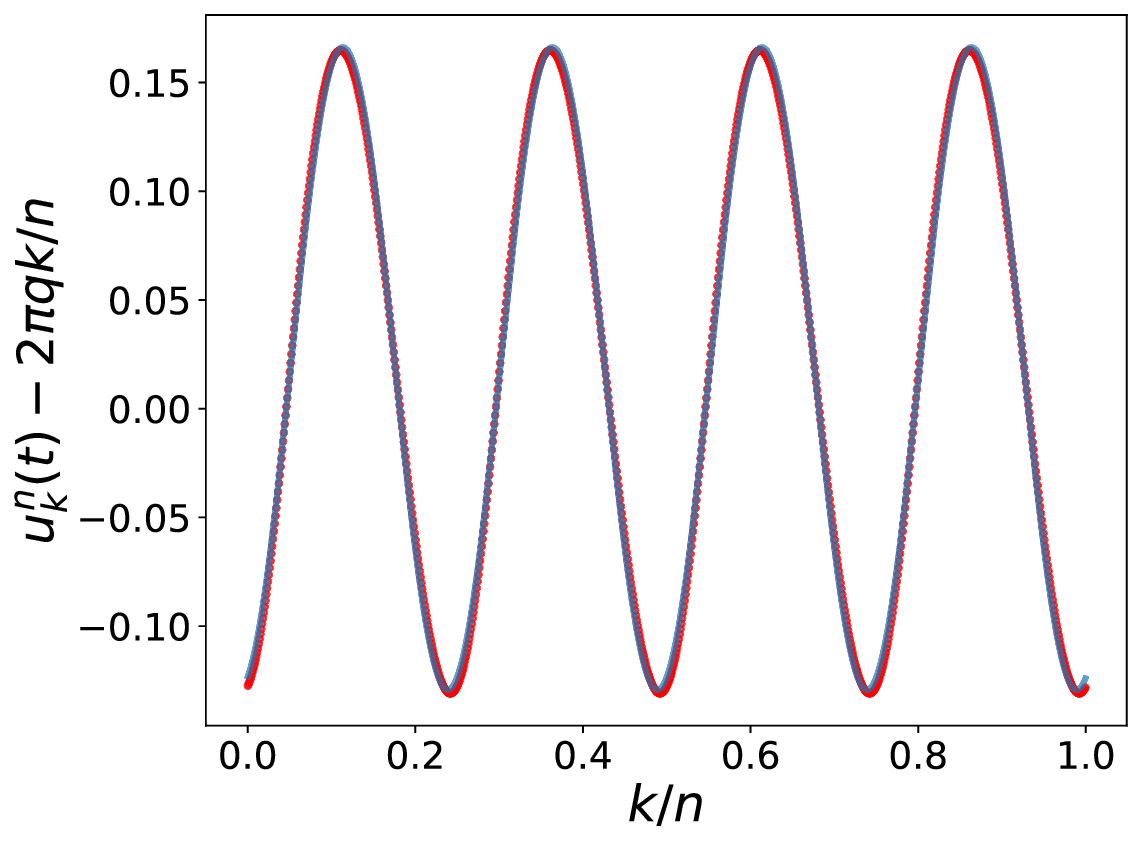}\\[-1ex]
{\footnotesize(d)}
\end{center}
\end{minipage}
\caption{Deviation from the $q$-twisted states
 in the steady states of the KM \eqref{eqn:dsys}
 with $n=1000$, $\sigma=\pi/3$ and $b_3=0.5$ at $t=1000$:
(a) $(q,b_1)=(1,0.06)$;
(b) $(2,0.255)$;
(c) $(3,0.15)$;
(d) $(4,0.225)$.
See also the caption of Fig.~\ref{fig:5a3}.}
\label{fig:5b3}
\end{figure}

In Figs.~\ref{fig:5a3} and \ref{fig:5b3},
 the deviation, $u_k^n(t)-2\pi q k/n$, $k\in[n]$,
 of the steady state in the right columns
 of Figs.~\ref{fig:5a2} and \ref{fig:5b2}
 from the $q$-twisted state \eqref{eqn:ts} in the KM \eqref{eqn:dsys}
 for $\sigma=0$ and $\sigma=\pi/3$, respectively,
 when $b_1$ is considered to be smaller than the bifurcation point,
 is plotted as small red disks.
It was also estimated from the most probably leading term
 displayed in Figs.~\ref{fig:5a2} and \ref{fig:5b2}
 and is plotted as a blue line.
The agreement between both results is fine
  except in Fig.~\ref{fig:5a3}(a) for $q=1$ and $\kappa=0.4$.
The reason for their disagreement in Fig.~\ref{fig:5a3}(a)
 is considered to be that the absolute value of 
 $\mu_{2q}=(\chi_1(2q,q)-\chi_1(q,q))\cos\sigma$ is small
 and the $2q$-oscillation mode is easily to be excited
 when $q=1$ and $\kappa=0.4$  (see Fig.~\ref{fig:3a}),
 compared with the other cases.

\begin{figure}[t]
\begin{minipage}[t]{0.495\textwidth}
\begin{center}
\includegraphics[scale=0.4]{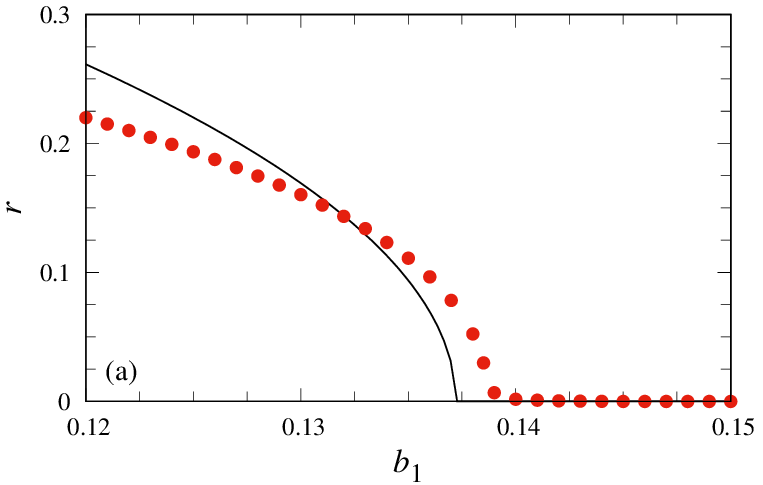}
\end{center}
\end{minipage}
\begin{minipage}[t]{0.495\textwidth}
\begin{center}
\includegraphics[scale=0.4]{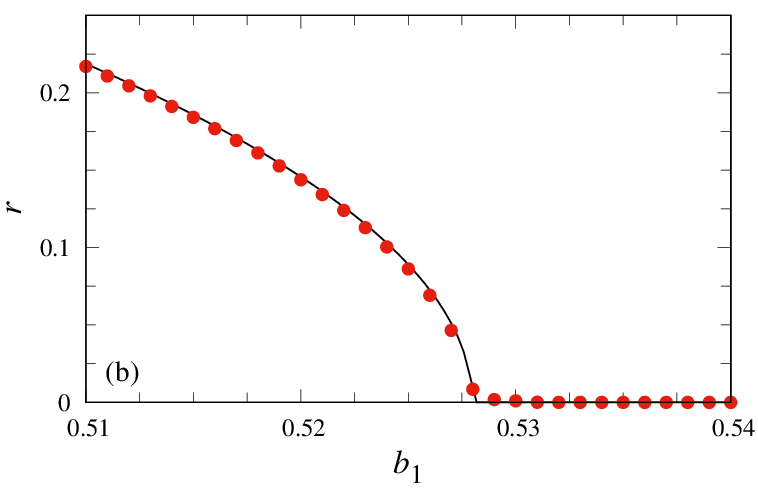}
\end{center}
\end{minipage}
\vspace*{0.5ex}

\begin{minipage}[t]{0.495\textwidth}
\begin{center}
\includegraphics[scale=0.4]{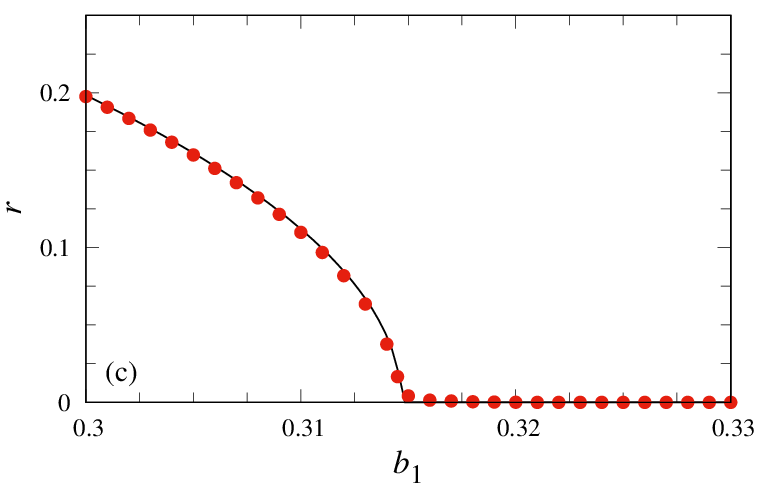}
\end{center}
\end{minipage}
\begin{minipage}[t]{0.495\textwidth}
\begin{center}
\includegraphics[scale=0.4]{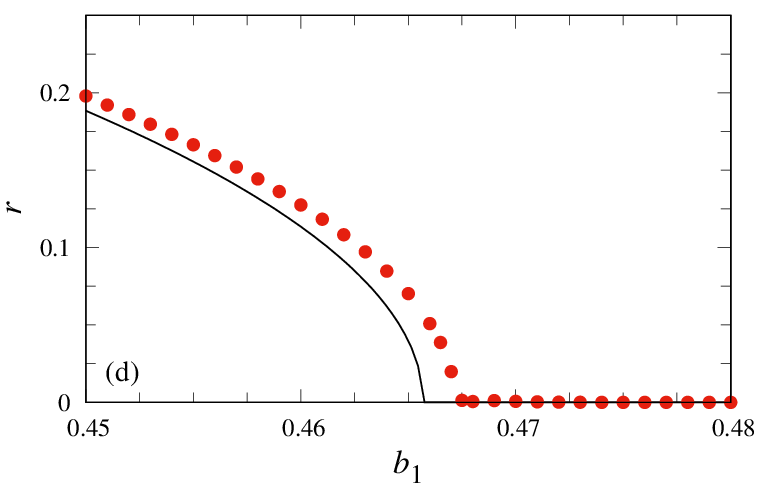}
\end{center}
\end{minipage}
\caption{Bifurcation diagram of the steady states in  the KM \eqref{eqn:dsys}
 with  $n=1000$, $\kappa=0.4$ and $\sigma=0$:
(a) $(q,b_3)=(1,1)$;
(b) $(2,0.5)$;
(c) $(3,0.5)$;
(d) $(4,0.5)$.
The amplitude $r$ in 
 \eqref{eqn:ss} estimated from the simulation results and theoretical predictions
 (see the text for more details)
 are plotted as small red disks and black solid lines, respectively.}
\label{fig:5a4}
\end{figure}

\begin{figure}[t]
\begin{minipage}[t]{0.495\textwidth}
\begin{center}
\includegraphics[scale=0.4]{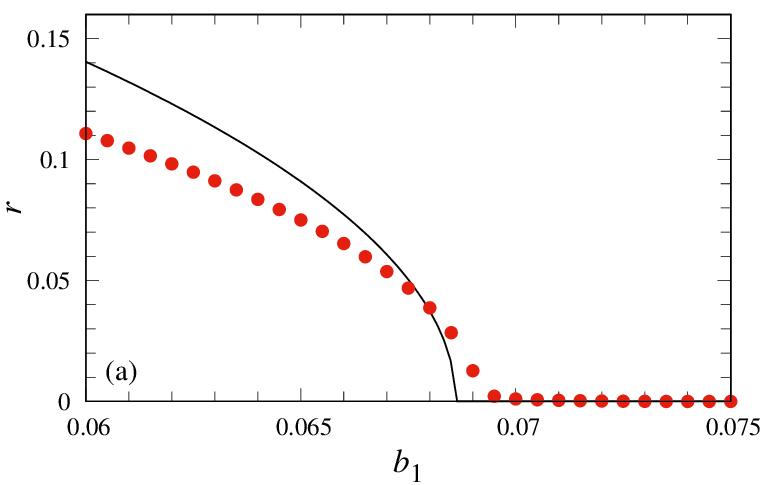}
\end{center}
\end{minipage}
\begin{minipage}[t]{0.495\textwidth}
\begin{center}
\includegraphics[scale=0.4]{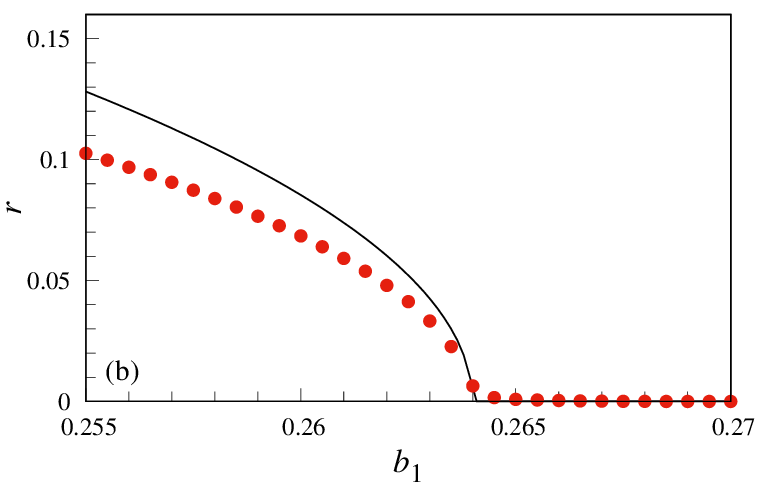}
\end{center}
\end{minipage}
\vspace*{0.5ex}

\begin{minipage}[t]{0.495\textwidth}
\begin{center}
\includegraphics[scale=0.4]{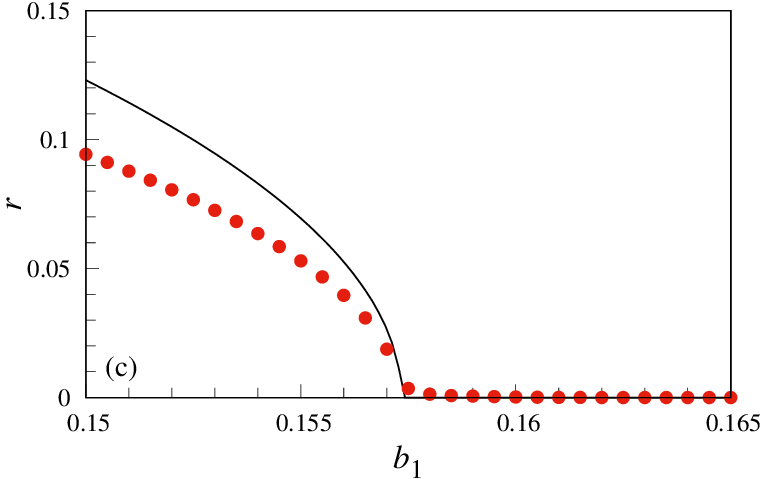}
\end{center}
\end{minipage}
\begin{minipage}[t]{0.495\textwidth}
\begin{center}
\includegraphics[scale=0.4]{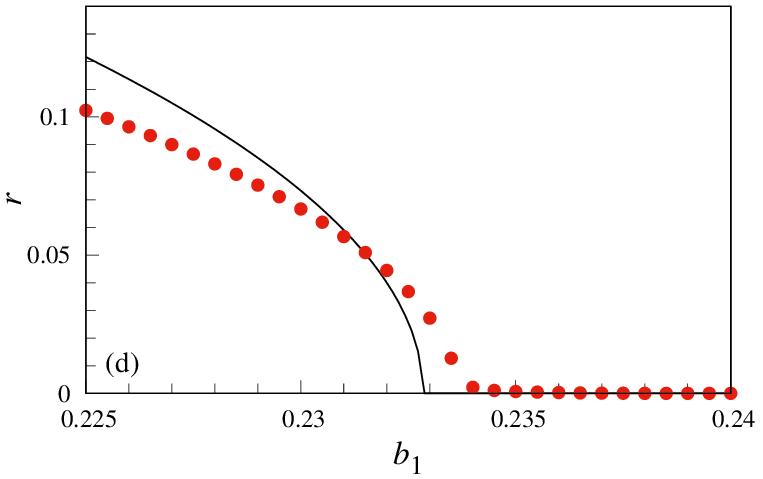}
\end{center}
\end{minipage}
\caption{Bifurcation diagram of the steady states in  the KM \eqref{eqn:dsys}
 with  $n=1000$, $\kappa=0.4$. $\sigma=\pi/3$  and $b_3=0.5$:
(a) $q=1$;
(b) $2$;
(c) $3$;
(d) $4$.
See also the caption of Fig.~\ref{fig:5a4}.}
\label{fig:5b4}
\end{figure}

Finally, we present numerically computed bifurcation diagrams
 for $\sigma=0$ and $\pi/3$ in Figs.~\ref{fig:5a4} and \ref{fig:5b4}, respectively.
The amplitude $r$ of the expression \eqref{eqn:ss}
 estimated from the numerical simulation results for the steady states
 as in Figs.~\ref{fig:5a2}-\ref{fig:5b3}
 are plotted as small red disks, and the theoretical predictions,
\begin{equation}
\sqrt{-\frac{b_1-b_{1q}}{\beta_0}}\quad\mbox{and}\quad
\sqrt{-\frac{b_1-b_{1q}}{\beta_\sigma}},
\label{eqn:rth}
\end{equation}
obtained from Theorems~\ref{thm:4a} and \ref{thm:4b}
 are plotted as black solid lines for $\sigma=0$ and $\sigma=\pi/3$, respectively,
 where $\beta_0$ and $\beta_\sigma$
 are given by \eqref{eqn:beta1} and  \eqref{eqn:beta1s}.
Good agreement between both results is found,
 especially in Figs.~\ref{fig:5a4}(b) and (c),
 although slight differences are seen in the other figures.


\section{Numerical Simulations: Complete Simple Graphs}

We next give numerical results for complete simple graphs,
 i.e., $\kappa=\tfrac{1}{2}$,
 for which the $q$-twisted solutions given by \eqref{eqn:tsol0}
 still exists but is unstable in the uncontrolled CL \eqref{eqn:csys} with $b_1,b_3=0$.

\begin{figure}[t]
\begin{minipage}[t]{0.495\textwidth}
\begin{center}
\includegraphics[scale=0.265]{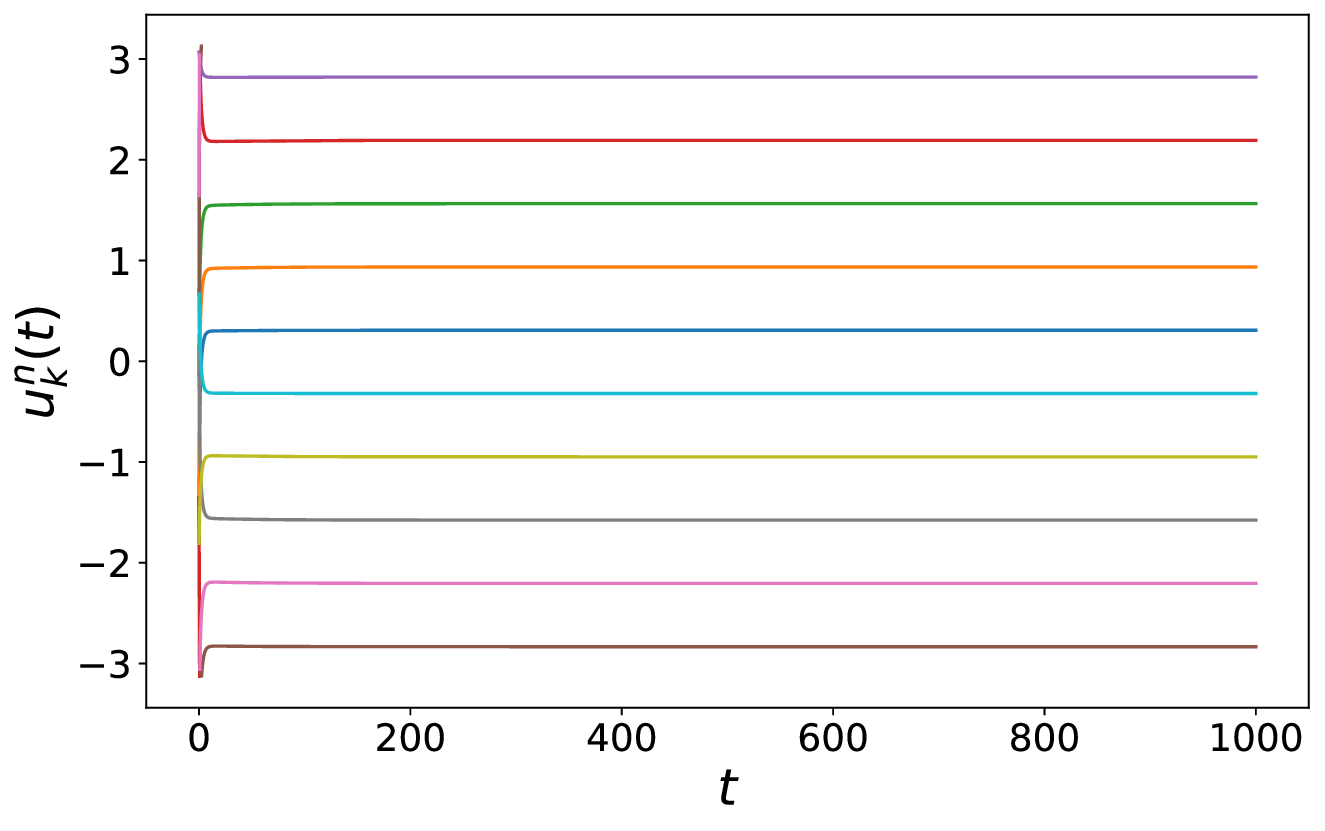}\\[-1ex]
{\footnotesize(a)}
\end{center}
\end{minipage}
\begin{minipage}[t]{0.495\textwidth}
\begin{center}
\includegraphics[scale=0.265]{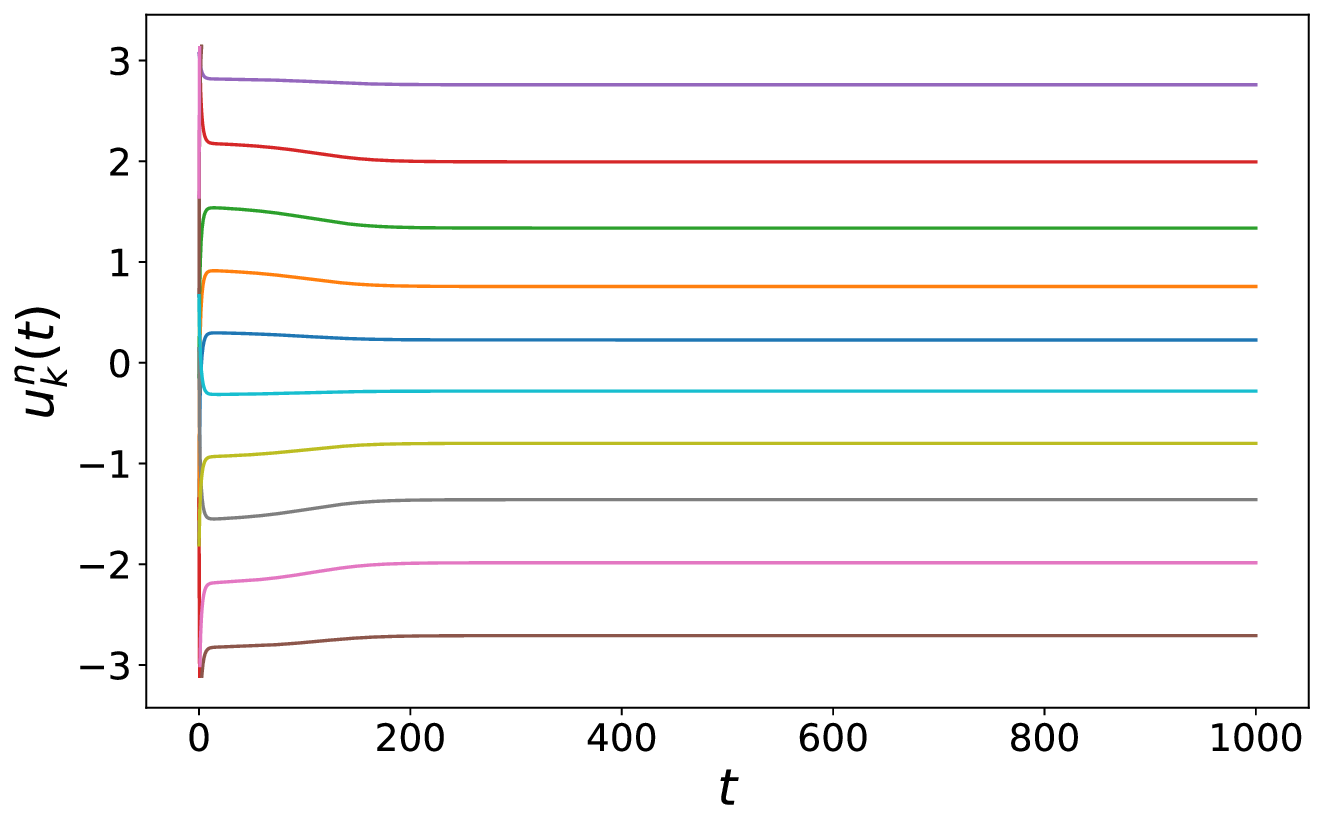}\\[-1ex]
{\footnotesize(b)}
\end{center}
\end{minipage}
\vspace*{0.5ex}

\begin{minipage}[t]{0.495\textwidth}
\begin{center}
\includegraphics[scale=0.265]{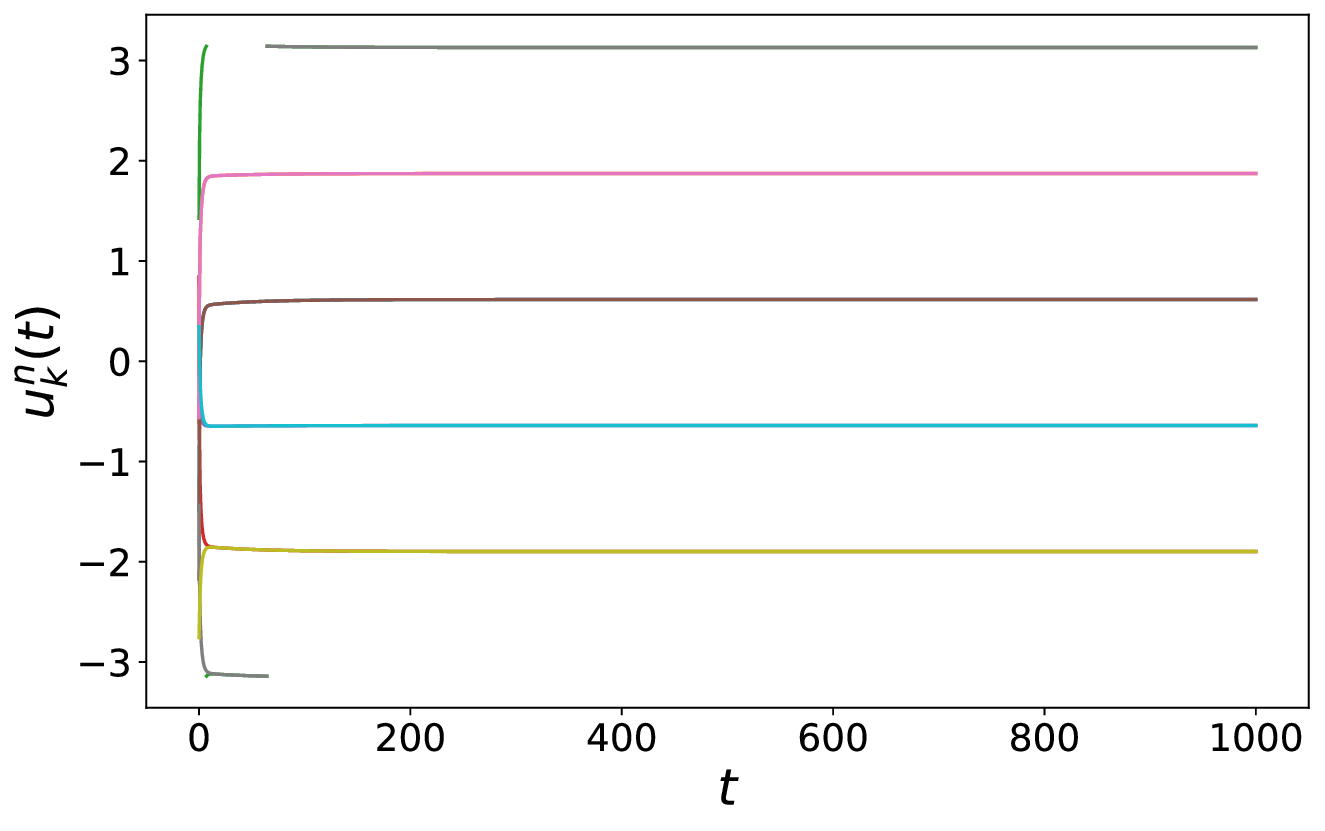}\\[-1ex]
{\footnotesize(c)}
\end{center}
\end{minipage}
\begin{minipage}[t]{0.495\textwidth}
\begin{center}
\includegraphics[scale=0.265]{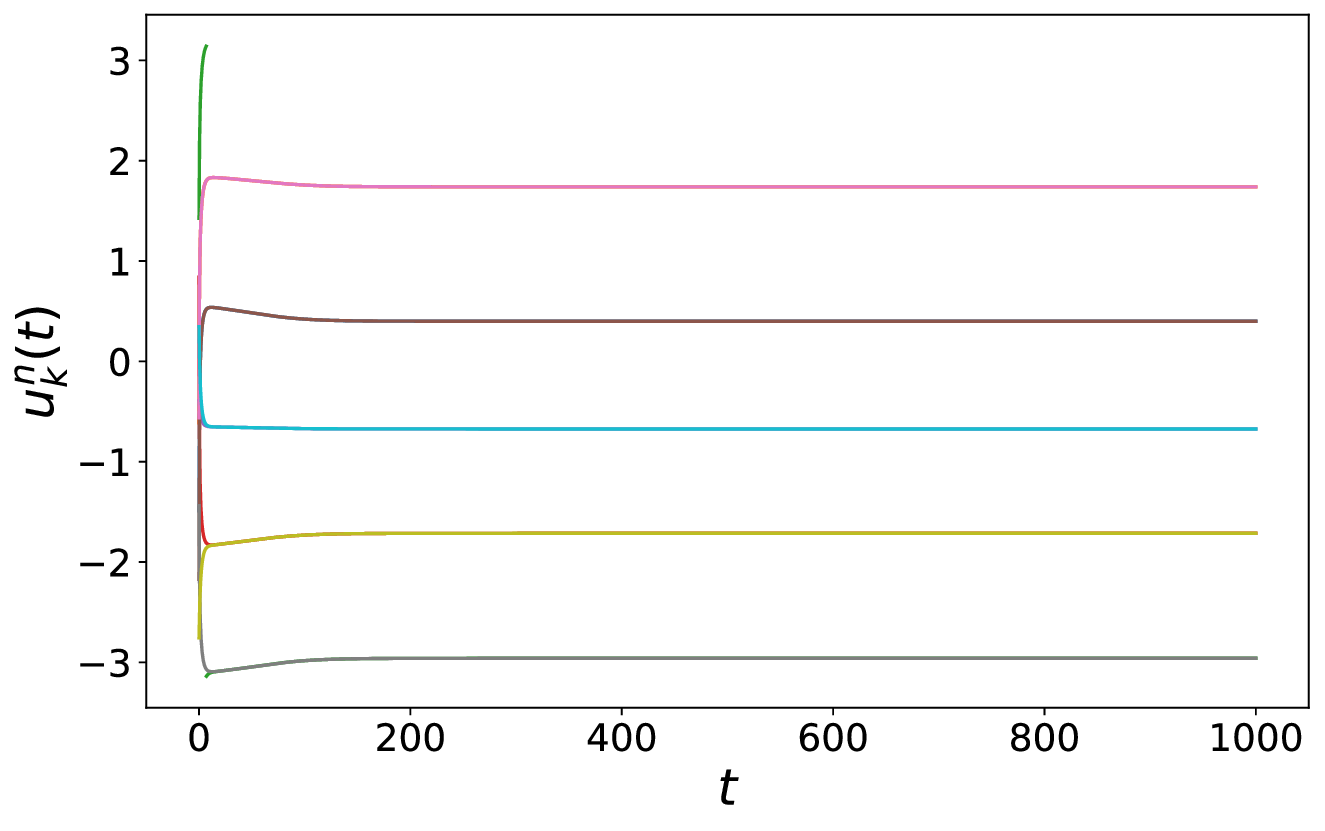}\\[-1ex]
{\footnotesize(d)}
\end{center}
\end{minipage}
\vspace*{0.5ex}

\begin{minipage}[t]{0.495\textwidth}
\begin{center}
\includegraphics[scale=0.265]{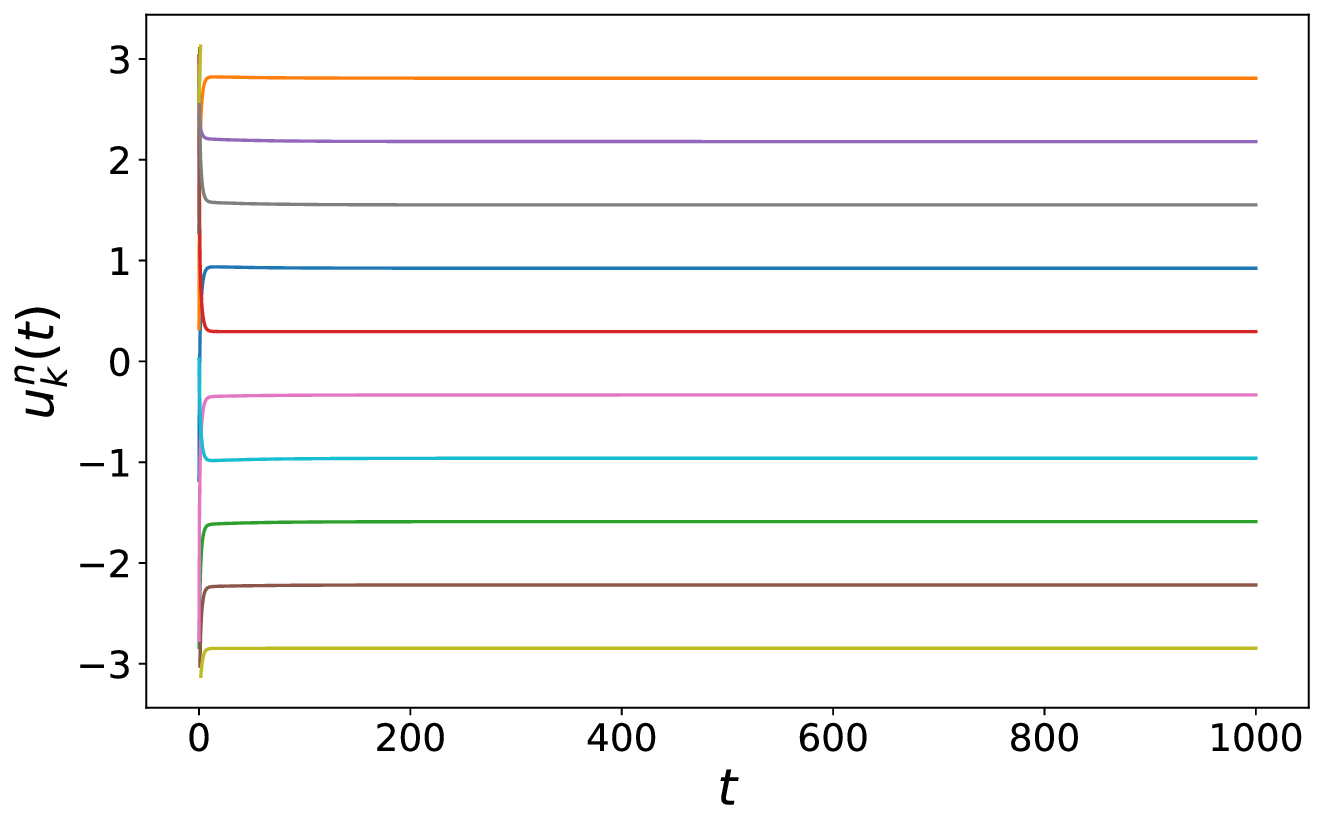}\\[-1ex]
{\footnotesize(e)}
\end{center}
\end{minipage}
\begin{minipage}[t]{0.495\textwidth}
\begin{center}
\includegraphics[scale=0.265]{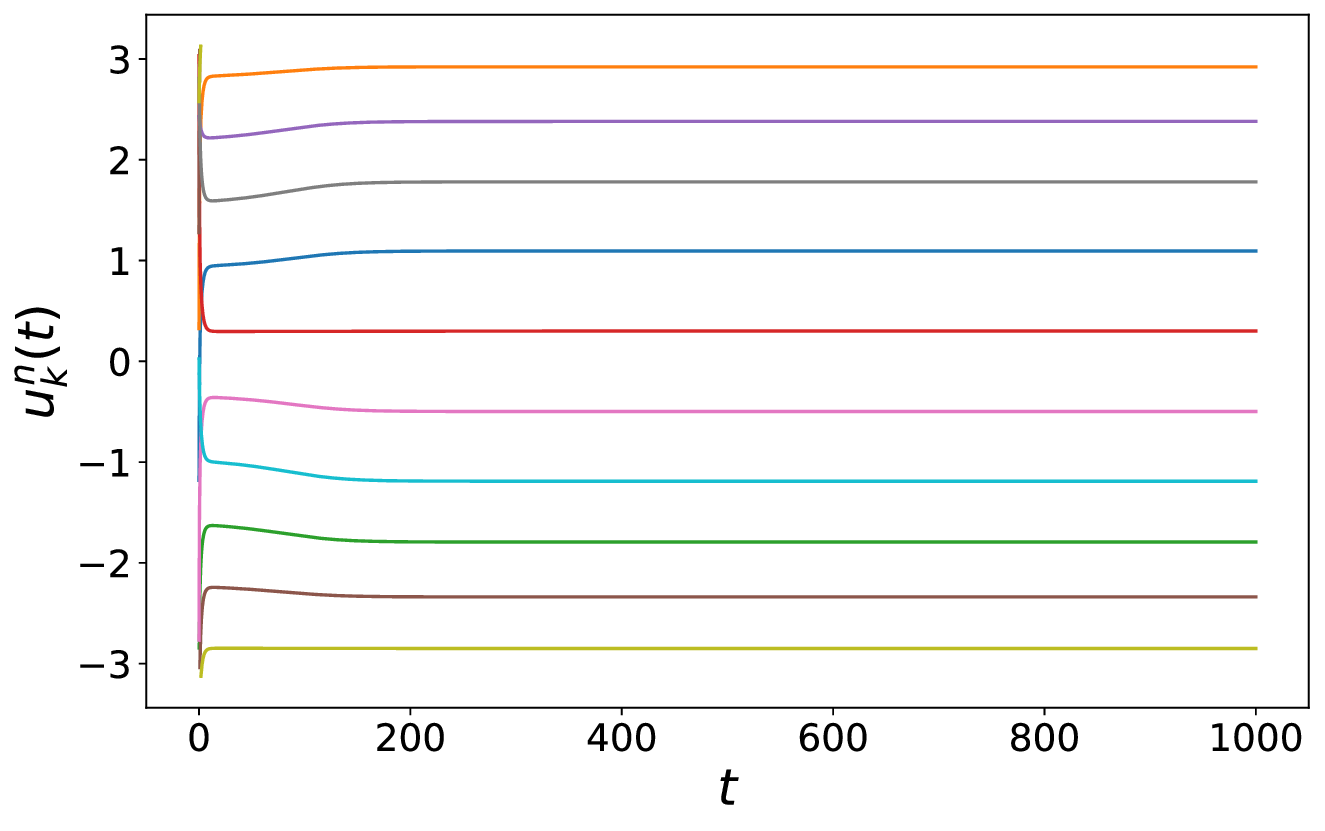}\\[-1ex]
{\footnotesize(f)}
\end{center}
\end{minipage}
\vspace*{0.5ex}

\begin{minipage}[t]{0.495\textwidth}
\begin{center}
\includegraphics[scale=0.265]{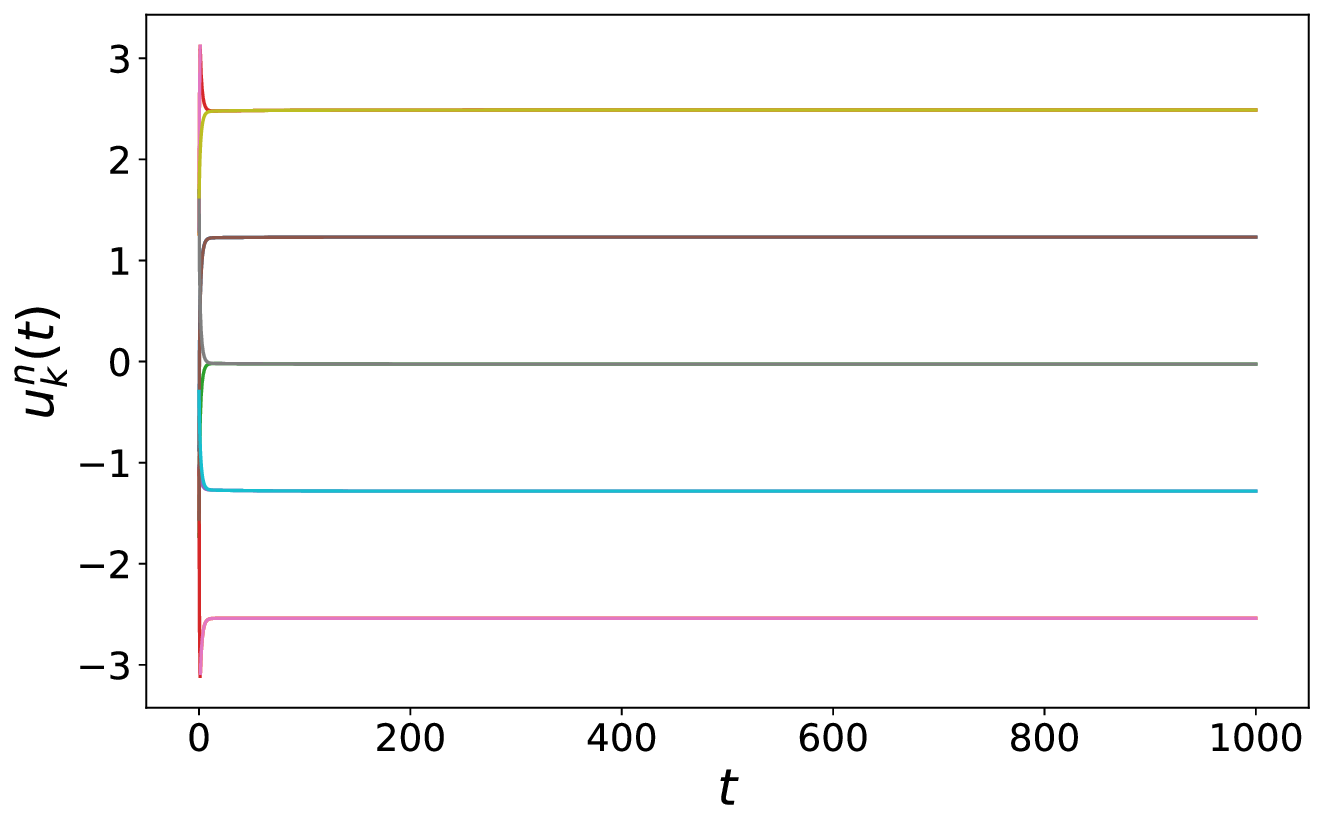}\\[-1ex]
{\footnotesize(g)}
\end{center}
\end{minipage}
\begin{minipage}[t]{0.495\textwidth}
\begin{center}
\includegraphics[scale=0.265]{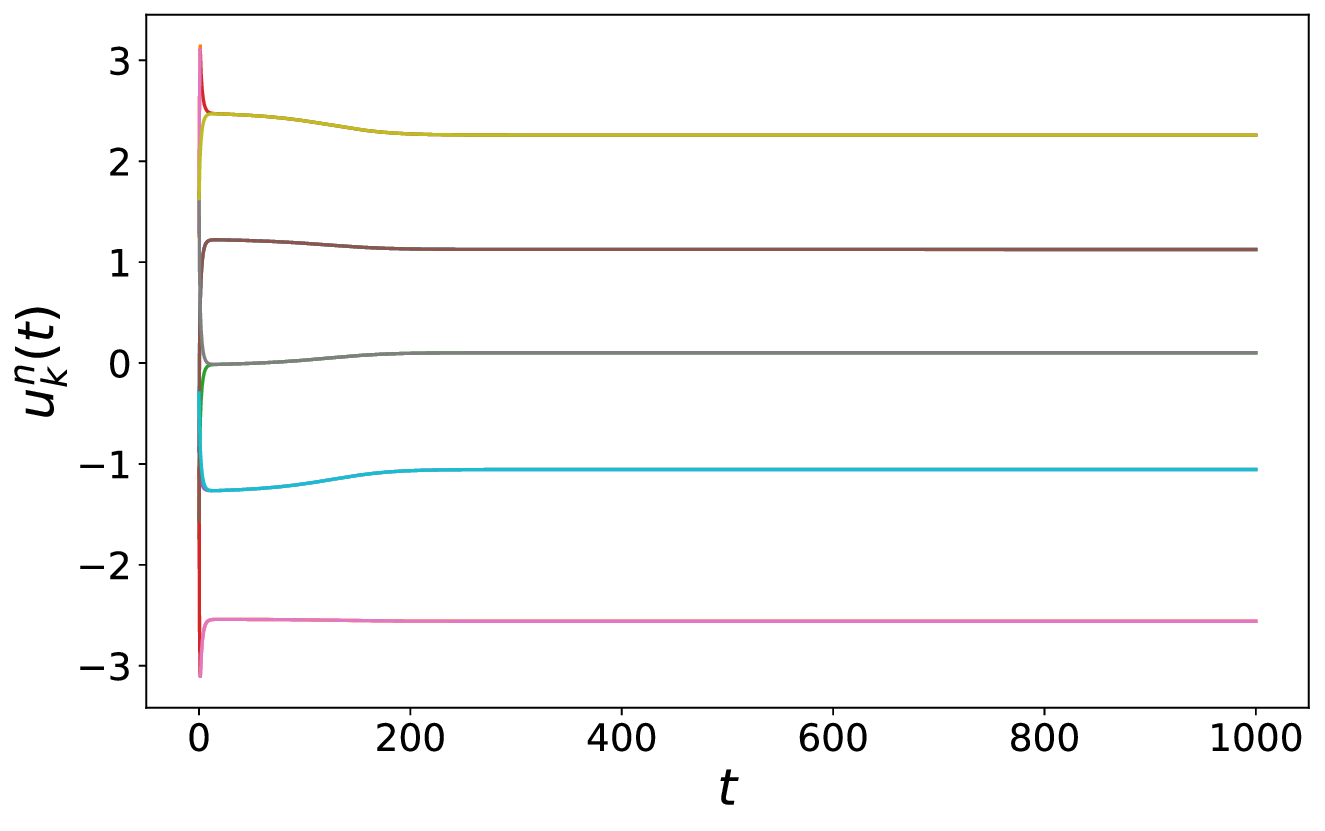}\\[-1ex]
{\footnotesize(h)}
\end{center}
\end{minipage}
\caption{Numerical simulation results for the KM \eqref{eqn:dsys}
 with $n=1000$, $\kappa=0.5$, $\sigma=0$ and $b_3=0.5$:
(a) $(q,b_1)=(1,0.52)$;
(b) $(1,0.48)$;
(c) $(2,0.52)$;
(d) $(2,0.48)$;
(e) $(3,0.52)$;
(f) $(3,0.48)$;
(g) $(4,0.52)$;
(h) $(4,0.48)$.
See also the caption of Fig.~\ref{fig:5a1}.
\vspace*{5mm}}
\label{fig:6a1}
\end{figure}

\begin{figure}[t]
\begin{minipage}[t]{0.495\textwidth}
\begin{center}
\includegraphics[scale=0.265]{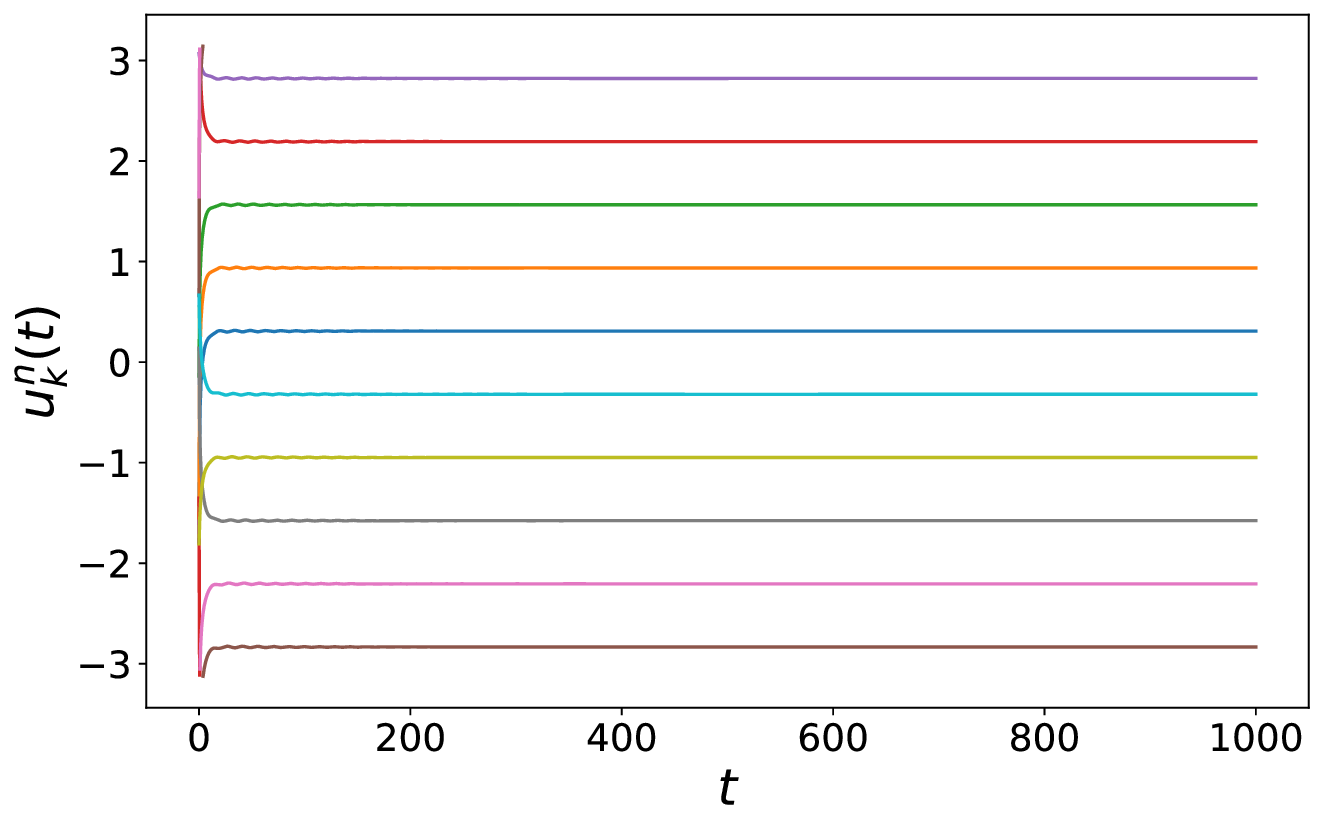}\\[-1ex]
{\footnotesize(a)}
\end{center}
\end{minipage}
\begin{minipage}[t]{0.495\textwidth}
\begin{center}
\includegraphics[scale=0.265]{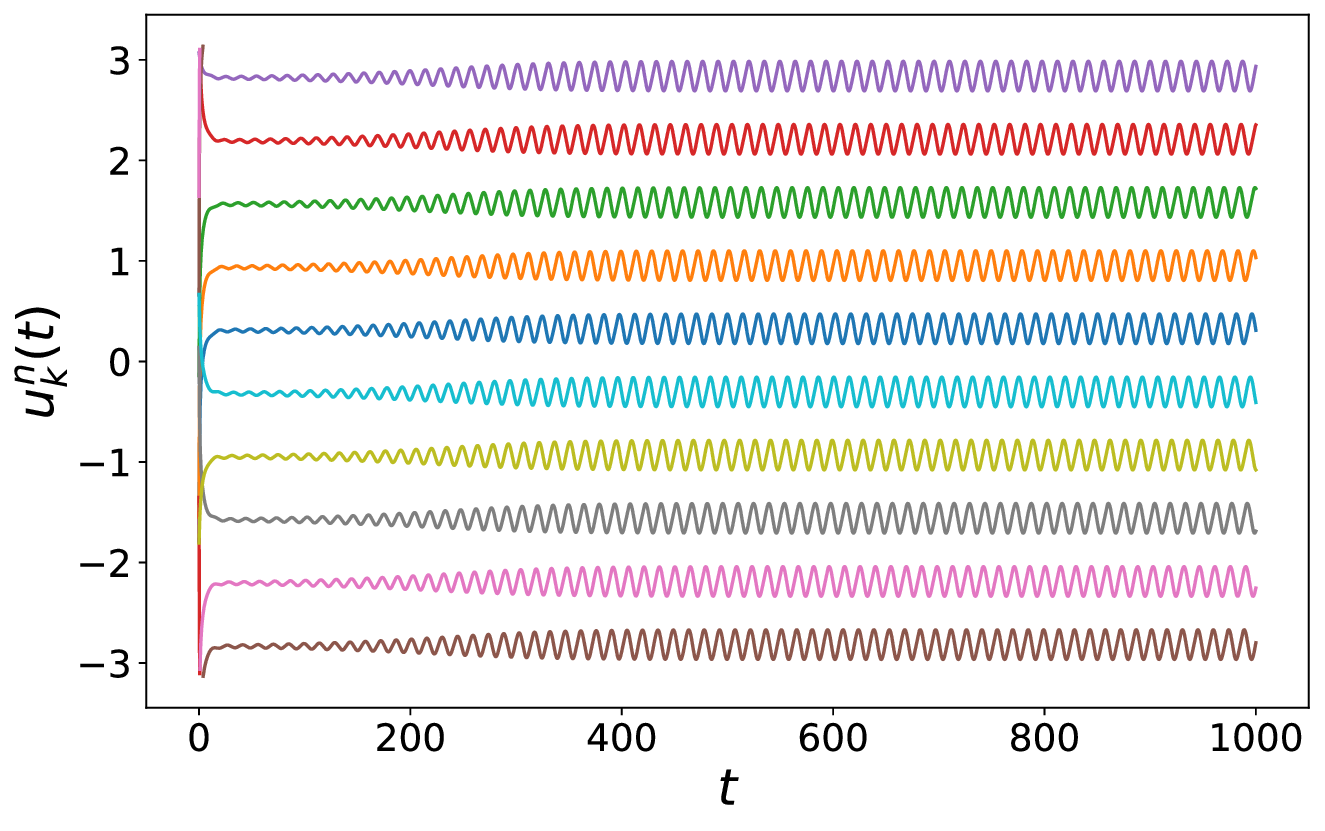}\\[-1ex]
{\footnotesize(b)}
\end{center}
\end{minipage}
\vspace*{0.5ex}

\begin{minipage}[t]{0.495\textwidth}
\begin{center}
\includegraphics[scale=0.265]{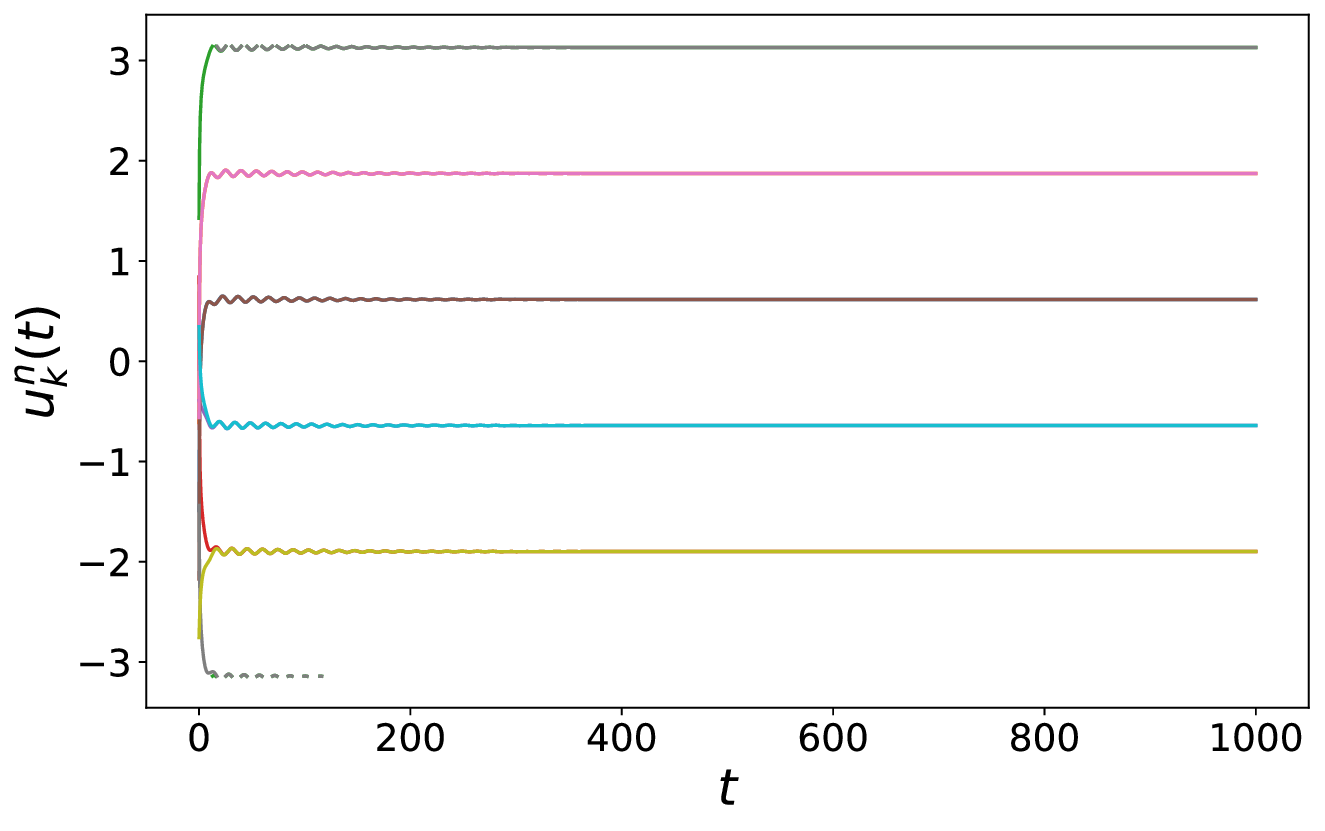}\\[-1ex]
{\footnotesize(c)}
\end{center}
\end{minipage}
\begin{minipage}[t]{0.495\textwidth}
\begin{center}
\includegraphics[scale=0.265]{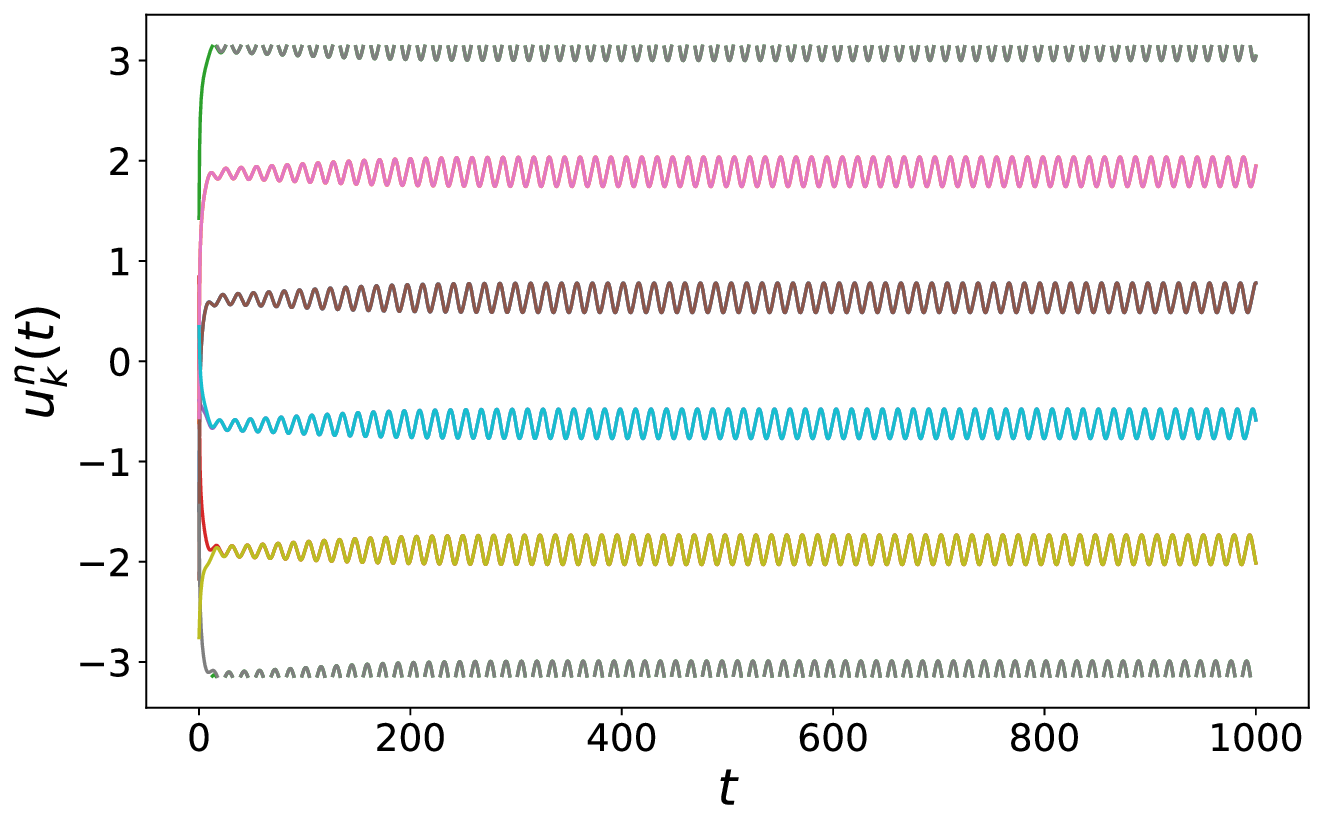}\\[-1ex]
{\footnotesize(d)}
\end{center}
\end{minipage}
\vspace*{0.5ex}

\begin{minipage}[t]{0.495\textwidth}
\begin{center}
\includegraphics[scale=0.265]{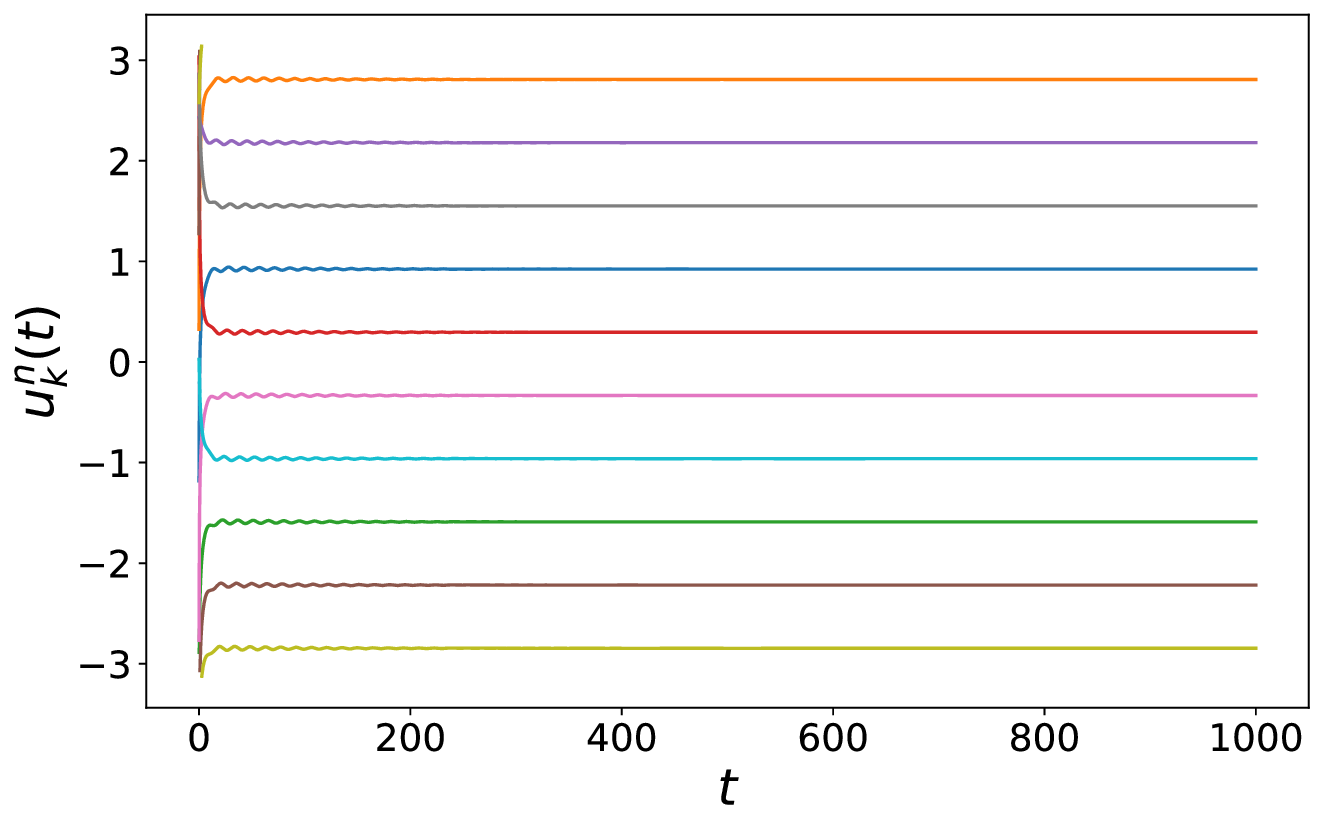}\\[-1ex]
{\footnotesize(e)}
\end{center}
\end{minipage}
\begin{minipage}[t]{0.495\textwidth}
\begin{center}
\includegraphics[scale=0.265]{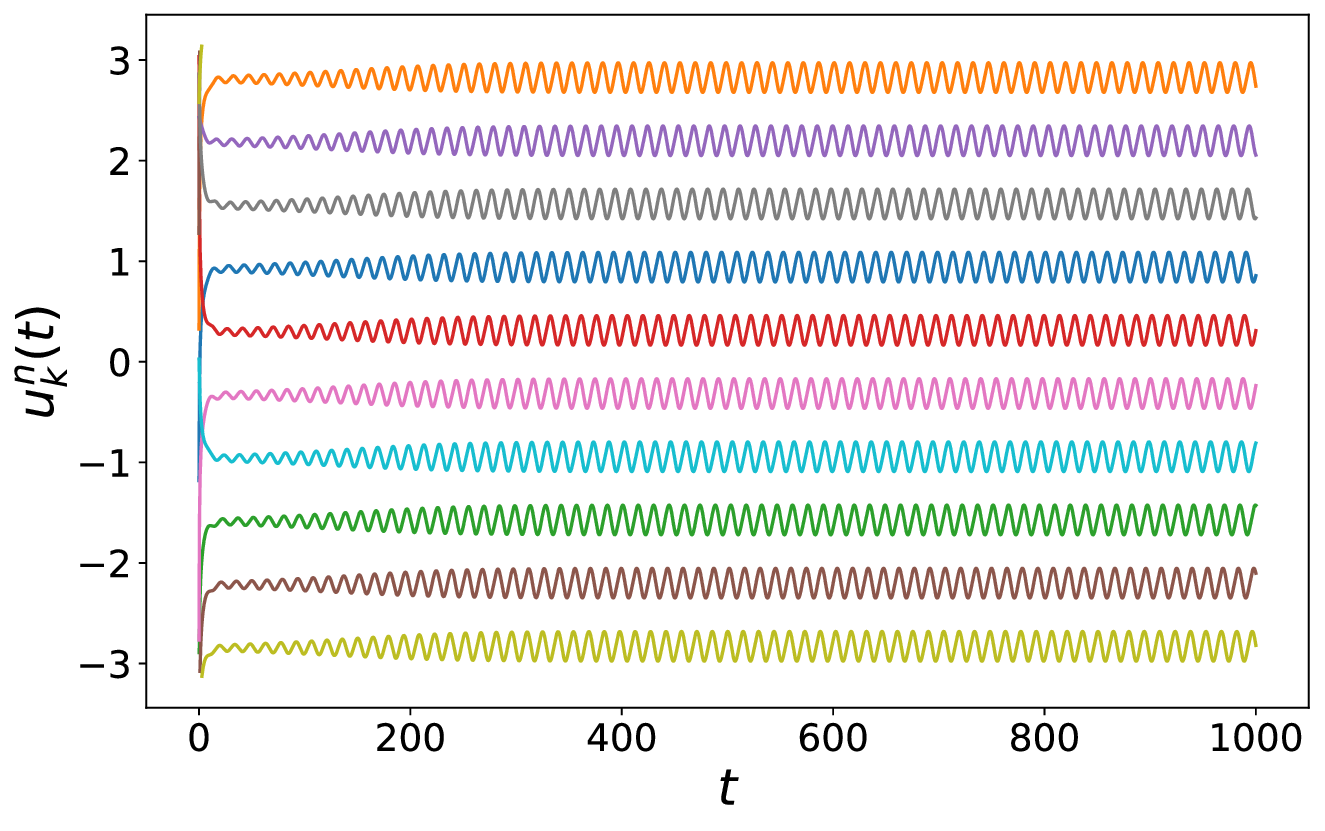}\\[-1ex]
{\footnotesize(f)}
\end{center}
\end{minipage}
\vspace*{0.5ex}

\begin{minipage}[t]{0.495\textwidth}
\begin{center}
\includegraphics[scale=0.265]{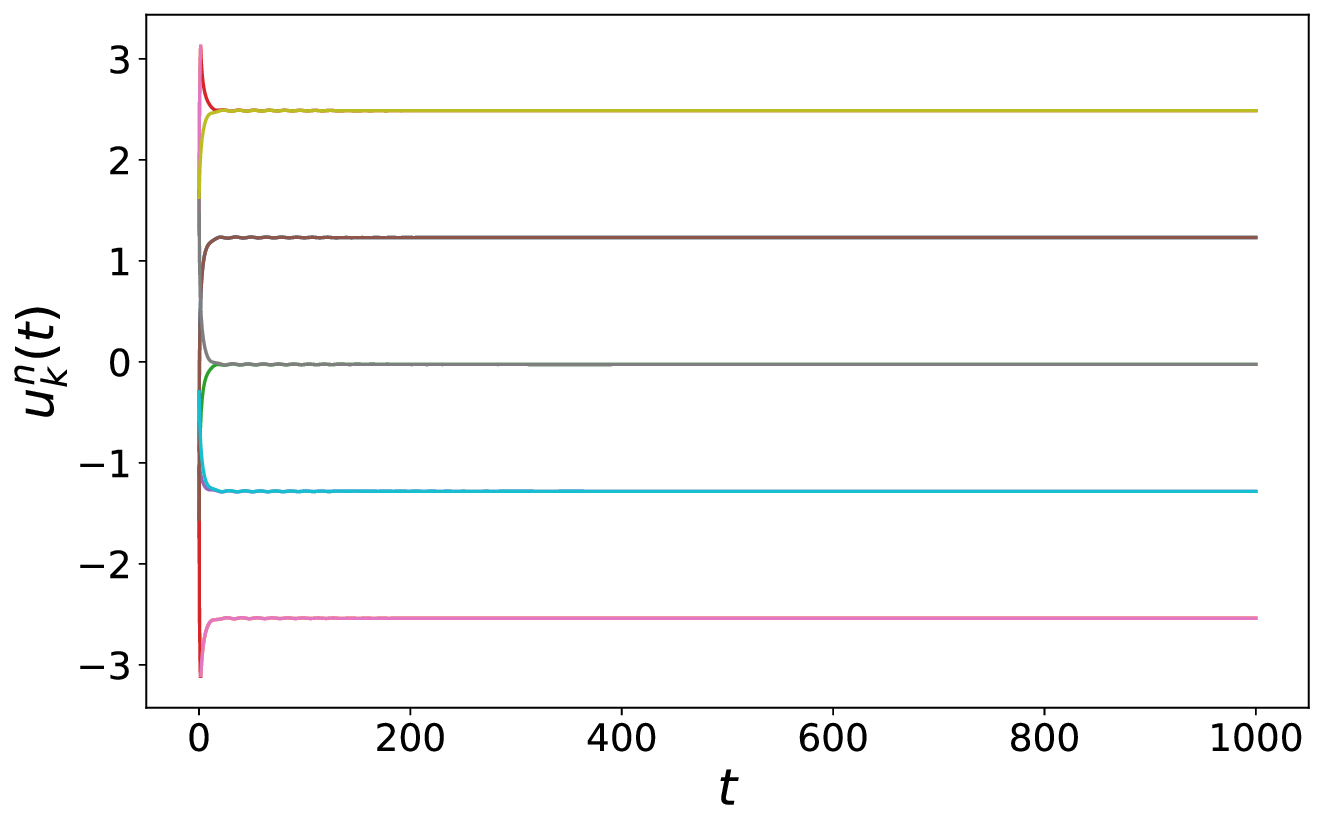}\\[-1ex]
{\footnotesize(g)}
\end{center}
\end{minipage}
\begin{minipage}[t]{0.495\textwidth}
\begin{center}
\includegraphics[scale=0.265]{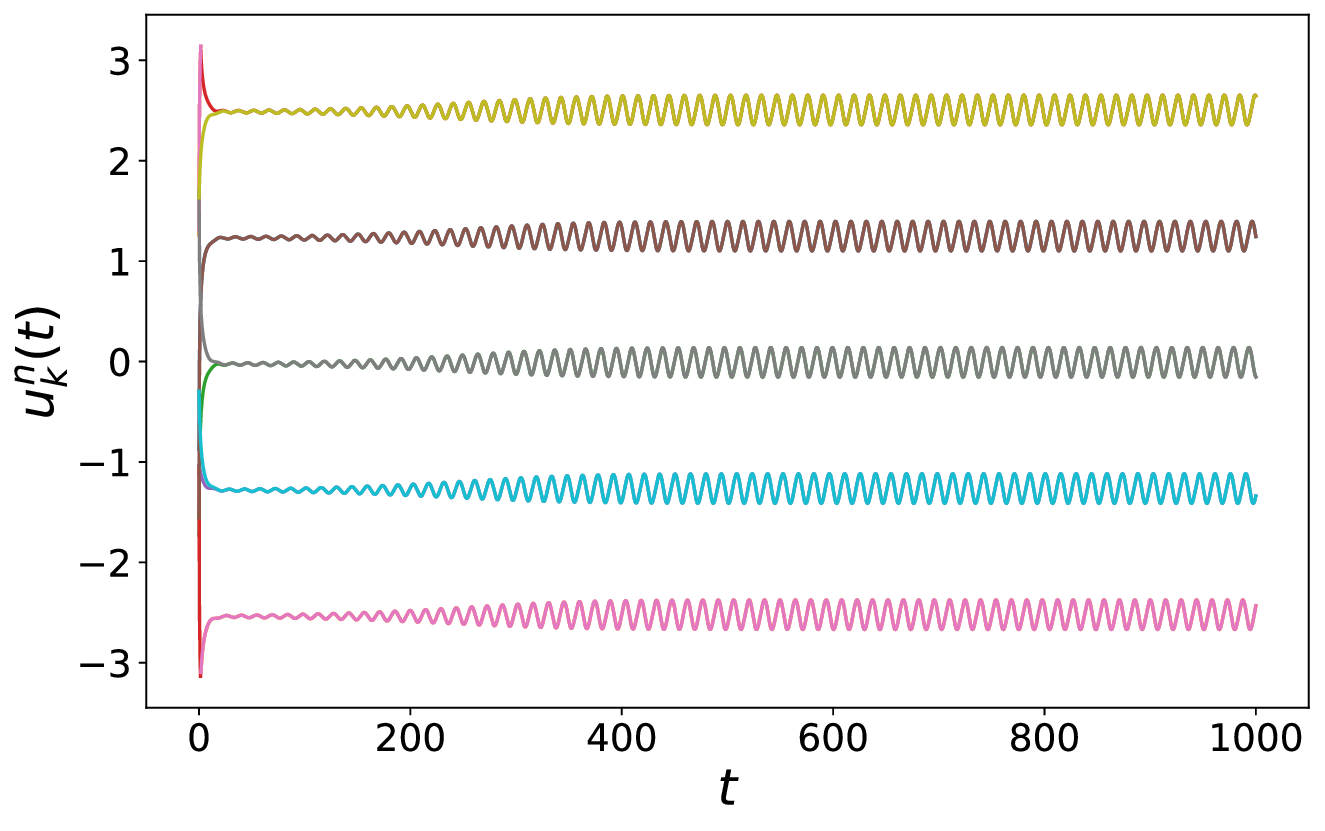}\\[-1ex]
{\footnotesize(h)}
\end{center}
\end{minipage}
\caption{Numerical simulation results for the KM \eqref{eqn:dsys}
 with $n=1000$, $\kappa=0.5$, $\sigma=\pi/3$ and $b_3=0.5$:
(a) $(q,b_1)=(1,0.26)$;
(b) $(1,0.24)$;
(c) $(2,0.26)$;
(d) $(2,0.24)$;
(e) $(3,0.26)$;
(f) $(3,0.24)$;
(g) $(4,0.26)$;
(h) $(4,0.24)$.
See also the caption of Fig.~\ref{fig:5a1}.
\vspace*{5mm}}
\label{fig:6b1}
\end{figure}

\begin{figure}[t]
\begin{minipage}[t]{0.495\textwidth}
\begin{center}
\includegraphics[scale=0.245]{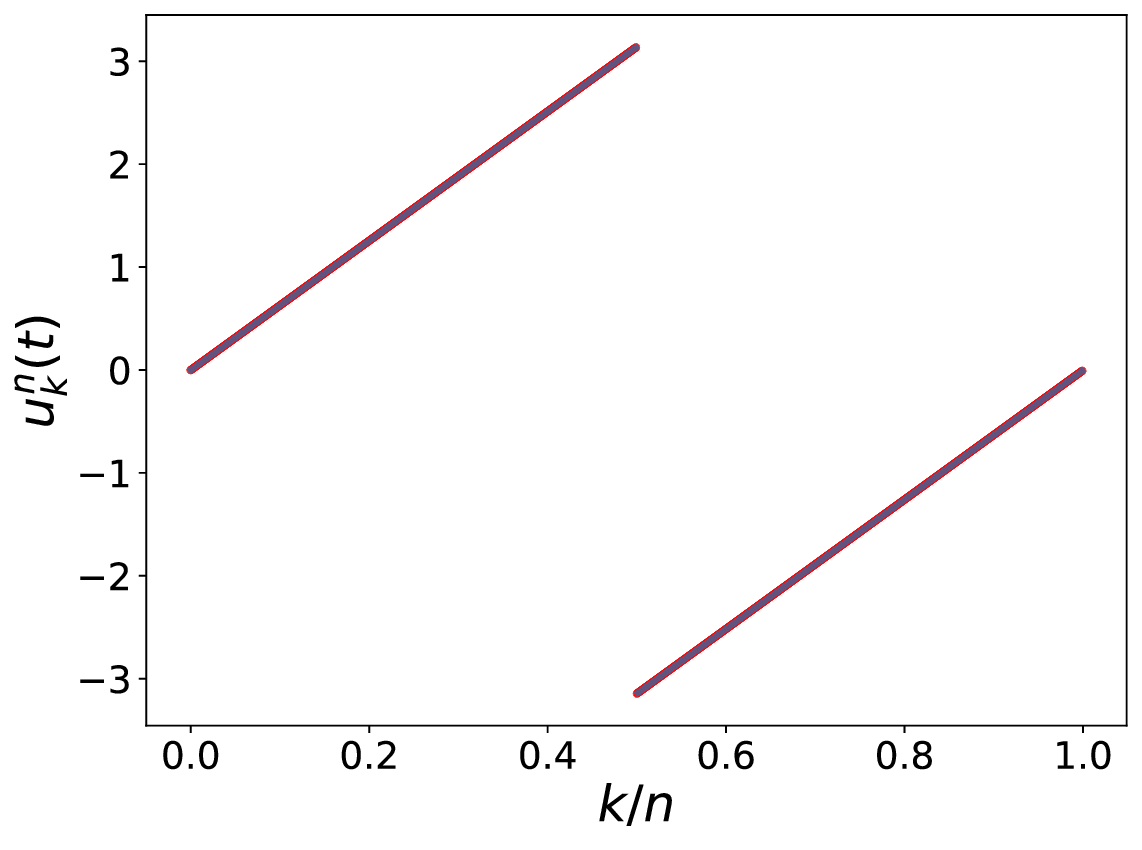}\\[-1ex]
{\footnotesize(a)}
\end{center}
\end{minipage}
\begin{minipage}[t]{0.495\textwidth}
\begin{center}
\includegraphics[scale=0.245]{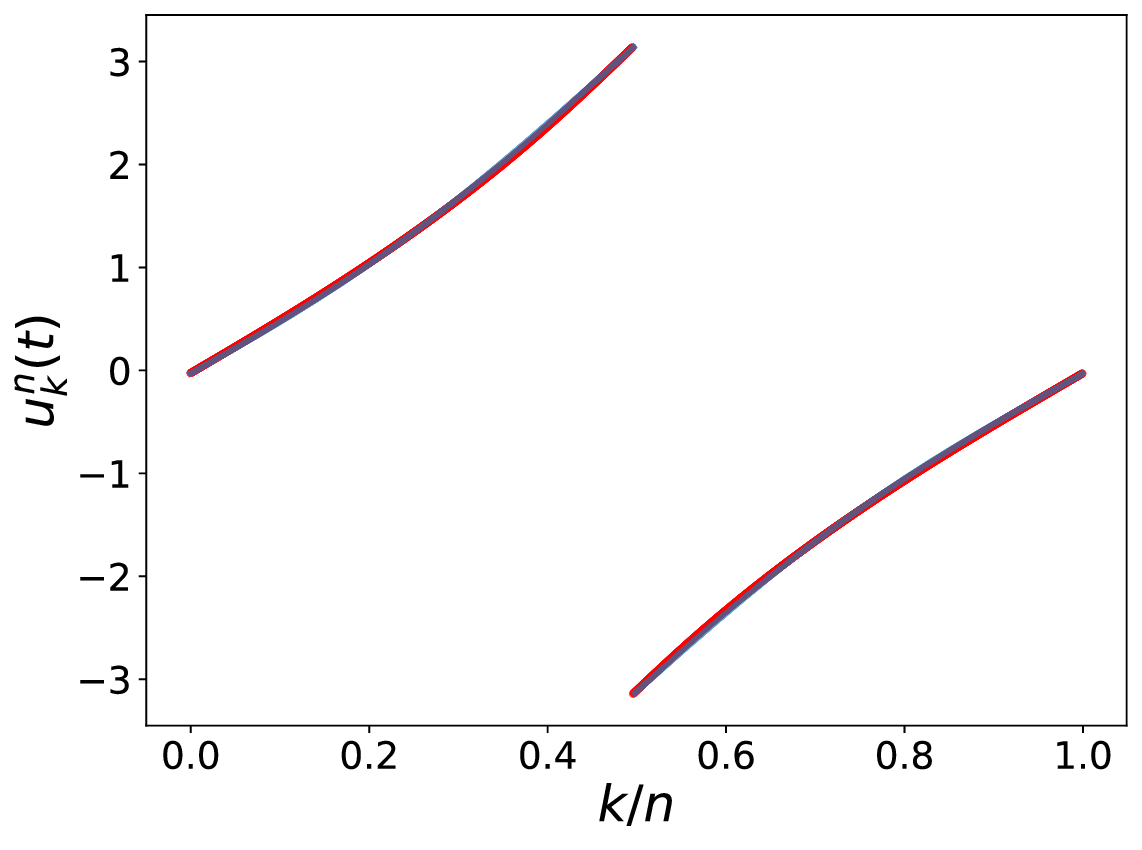}\\[-1ex]
{\footnotesize(b)}
\end{center}
\end{minipage}
\vspace*{0.5ex}

\begin{minipage}[t]{0.495\textwidth}
\begin{center}
\includegraphics[scale=0.245]{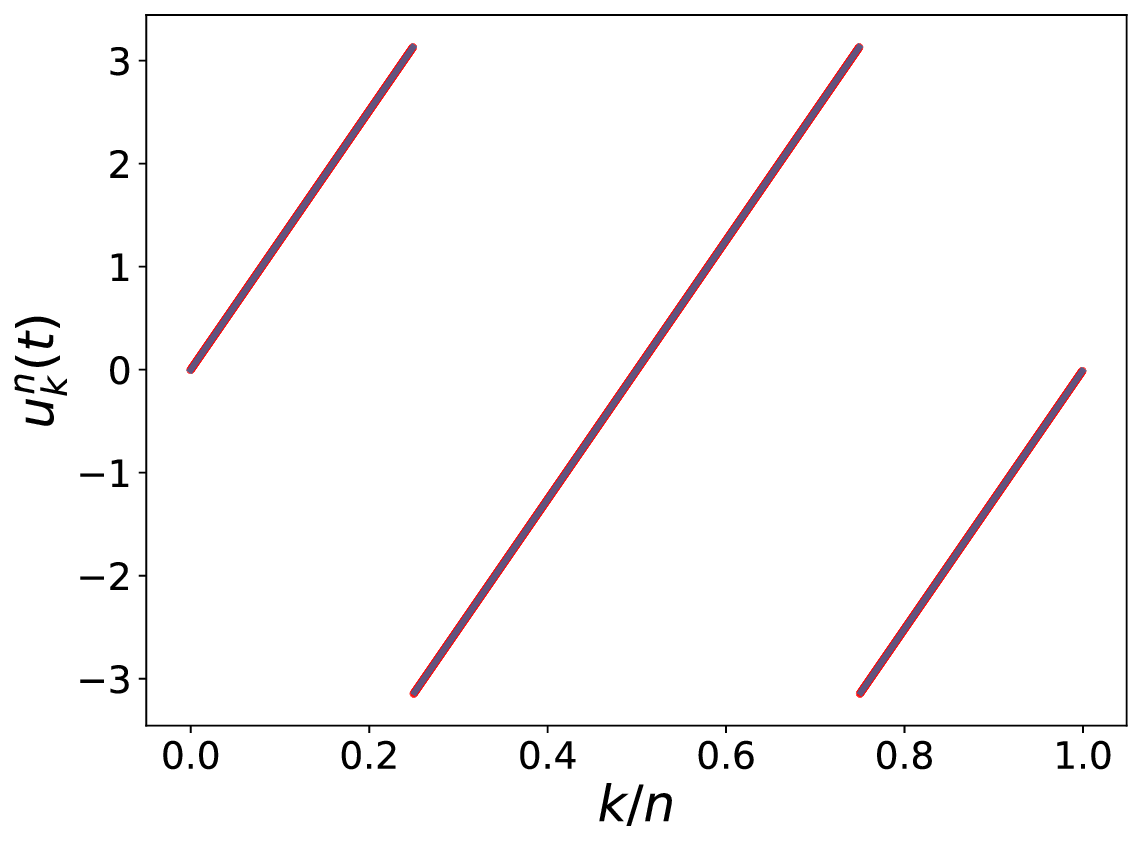}\\[-1ex]
{\footnotesize(c)}
\end{center}
\end{minipage}
\begin{minipage}[t]{0.495\textwidth}
\begin{center}
\includegraphics[scale=0.245]{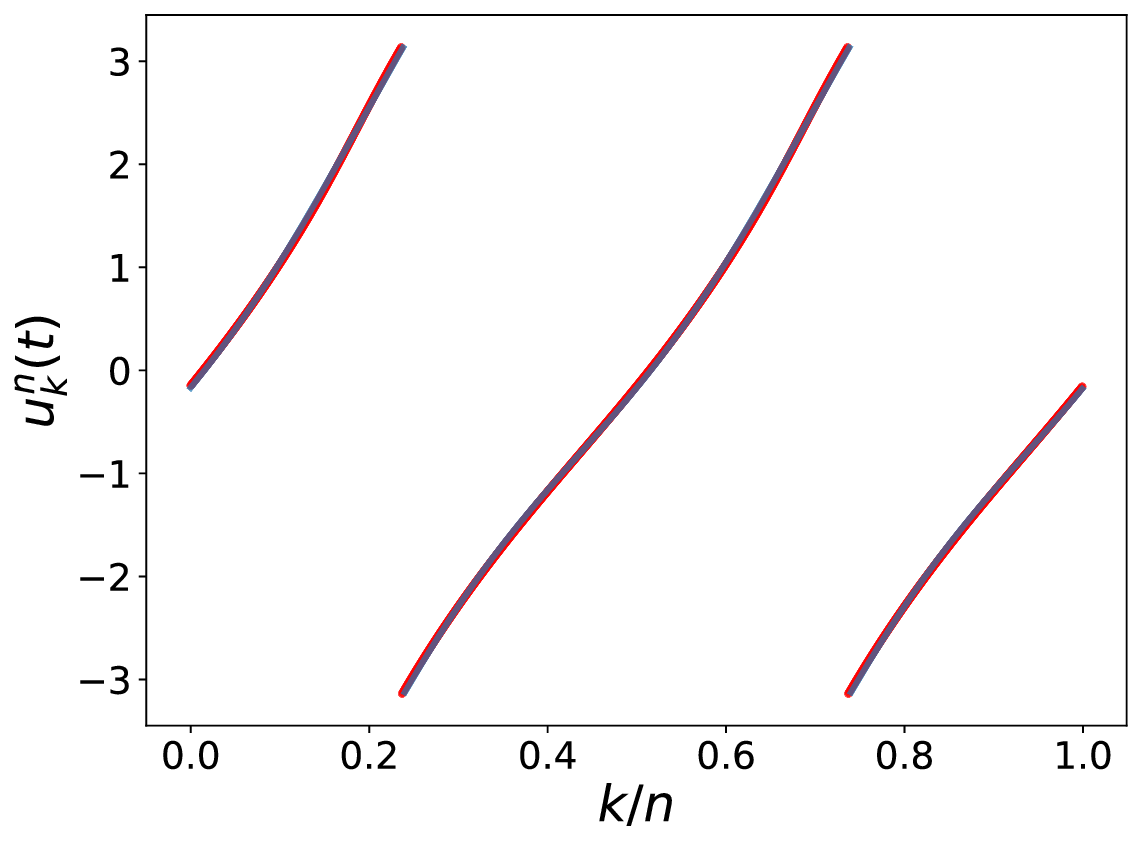}\\[-1ex]
{\footnotesize(d)}
\end{center}
\end{minipage}
\vspace*{0.5ex}

\begin{minipage}[t]{0.495\textwidth}
\begin{center}
\includegraphics[scale=0.245]{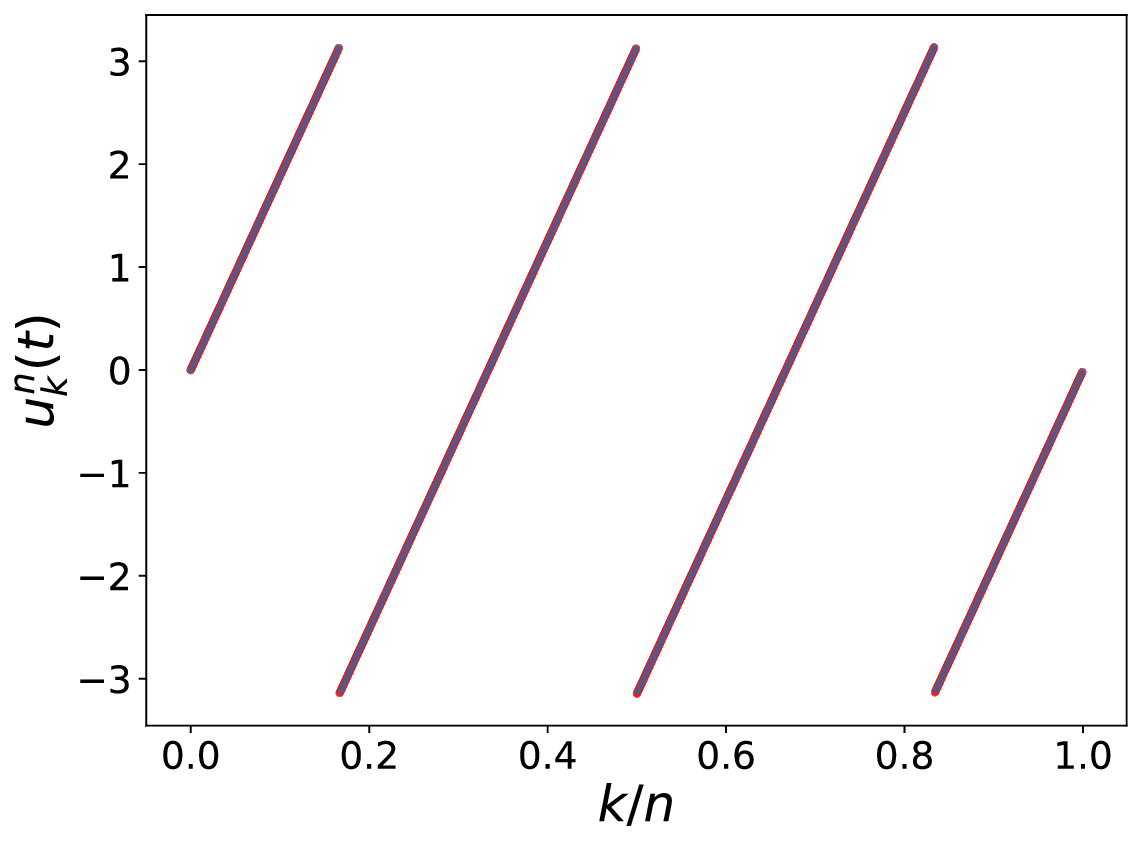}\\[-1ex]
{\footnotesize(e)}
\end{center}
\end{minipage}
\begin{minipage}[t]{0.495\textwidth}
\begin{center}
\includegraphics[scale=0.245]{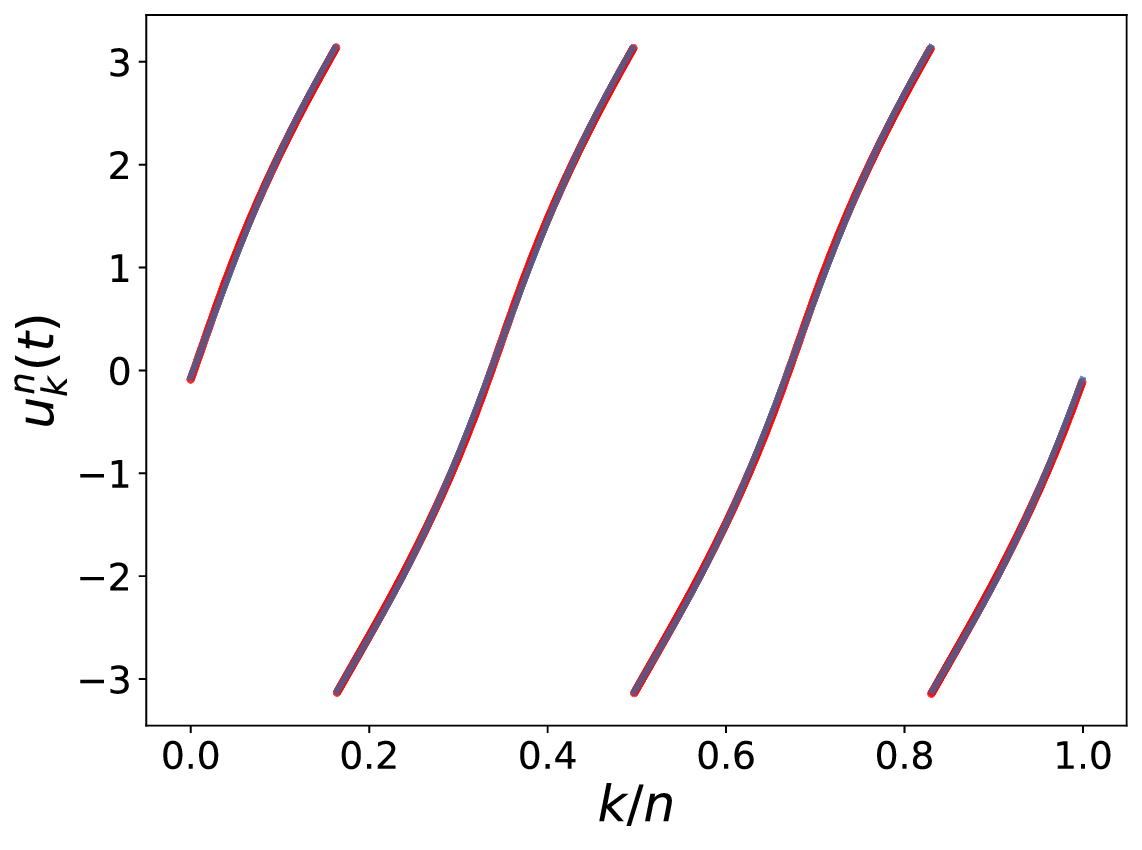}\\[-1ex]
{\footnotesize(f)}
\end{center}
\end{minipage}
\vspace*{0.5ex}

\begin{minipage}[t]{0.495\textwidth}
\begin{center}
\includegraphics[scale=0.245]{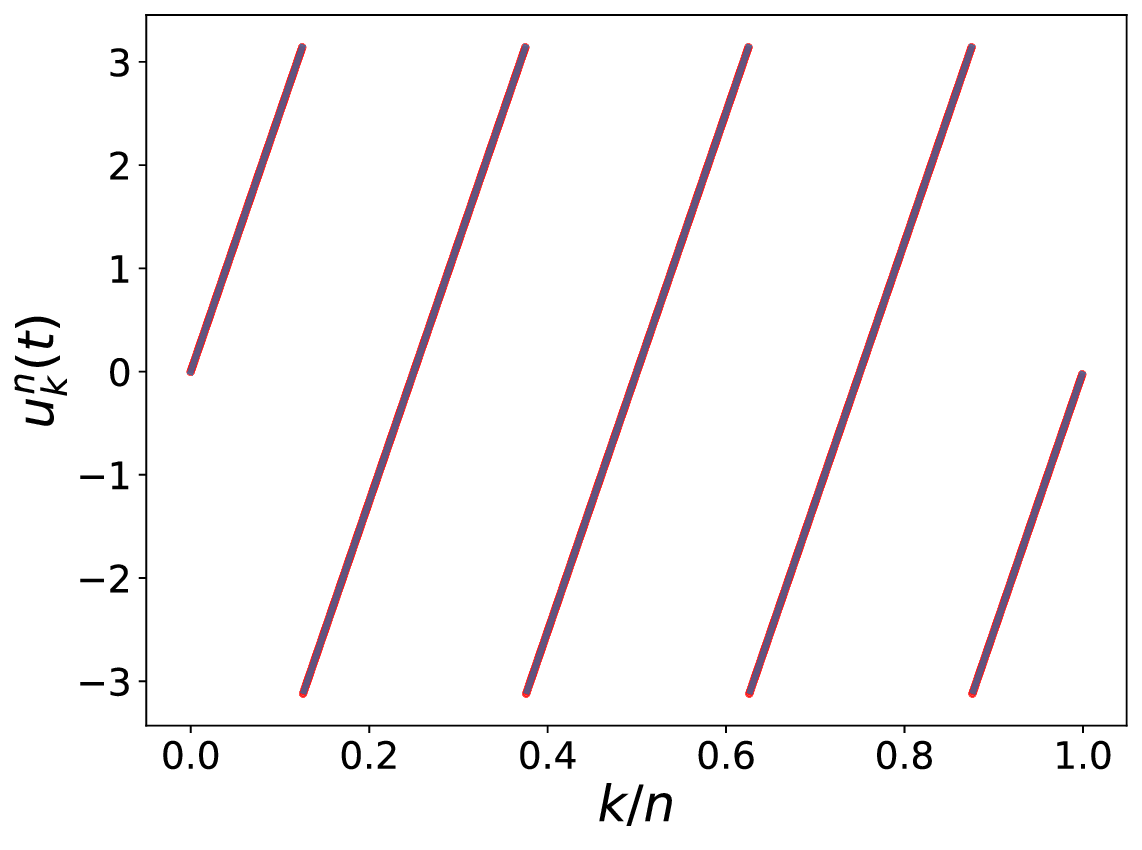}\\[-1ex]
{\footnotesize(g)}
\end{center}
\end{minipage}
\begin{minipage}[t]{0.495\textwidth}
\begin{center}
\includegraphics[scale=0.245]{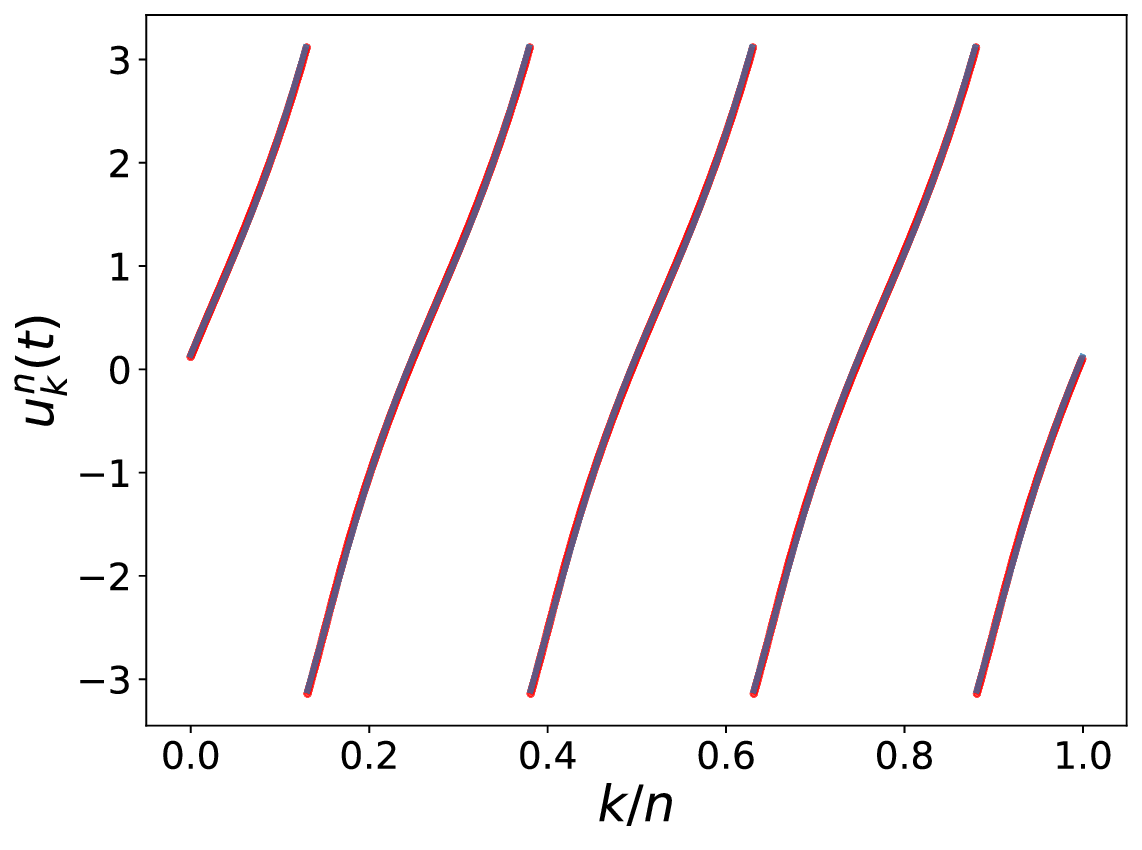}\\[-1ex]
{\footnotesize(h)}
\end{center}
\end{minipage}
\caption{Steady states of the KM \eqref{eqn:dsys}
 with $n=1000$, $\kappa=0.5$, $\sigma=0$ and $b_3=0.5$ at $t=1000$:
(a) $(q,b_1)=(1,0.52)$;
(b) $(1,0.48)$;
(c) $(2,0.52)$;
(d) $(2,0.48)$;
(e) $(3,0.52)$;
(f) $(3,0.48)$;
(g) $(4,0.52)$;
(h) $(4,0.48)$.
See also the caption of Fig.~\ref{fig:5a2}.}
\label{fig:6a2}
\end{figure}

\begin{figure}[t]
\begin{minipage}[t]{0.495\textwidth}
\begin{center}
\includegraphics[scale=0.245]{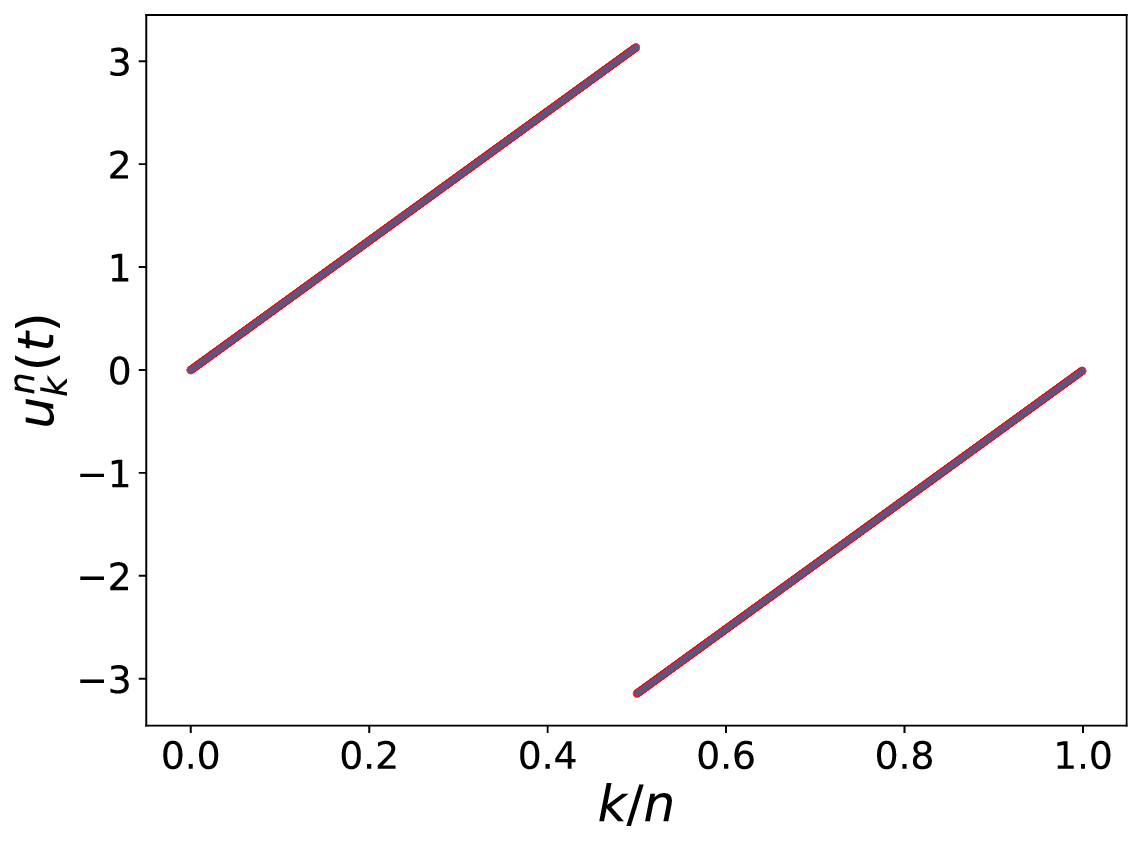}\\[-1ex]
{\footnotesize(a)}
\end{center}
\end{minipage}
\begin{minipage}[t]{0.495\textwidth}
\begin{center}
\includegraphics[scale=0.245]{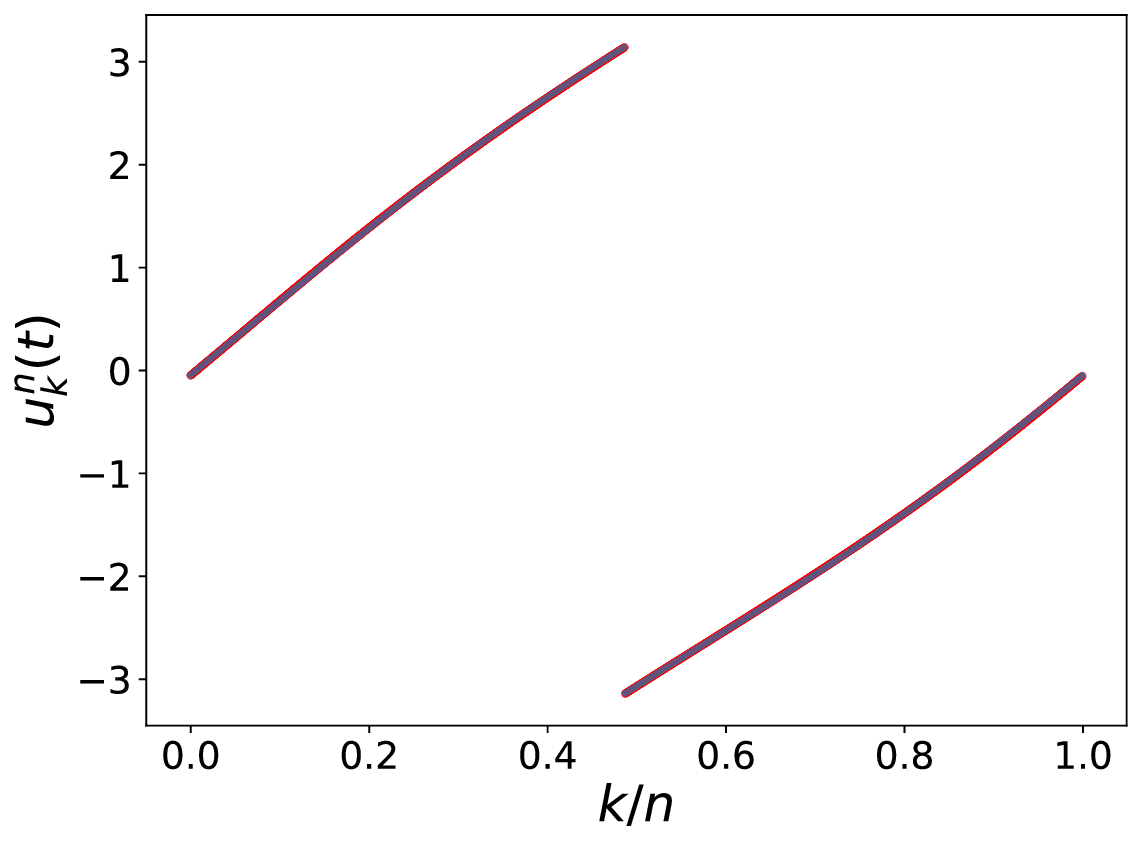}\\[-1ex]
{\footnotesize(b)}
\end{center}
\end{minipage}
\vspace*{0.5ex}

\begin{minipage}[t]{0.495\textwidth}
\begin{center}
\includegraphics[scale=0.245]{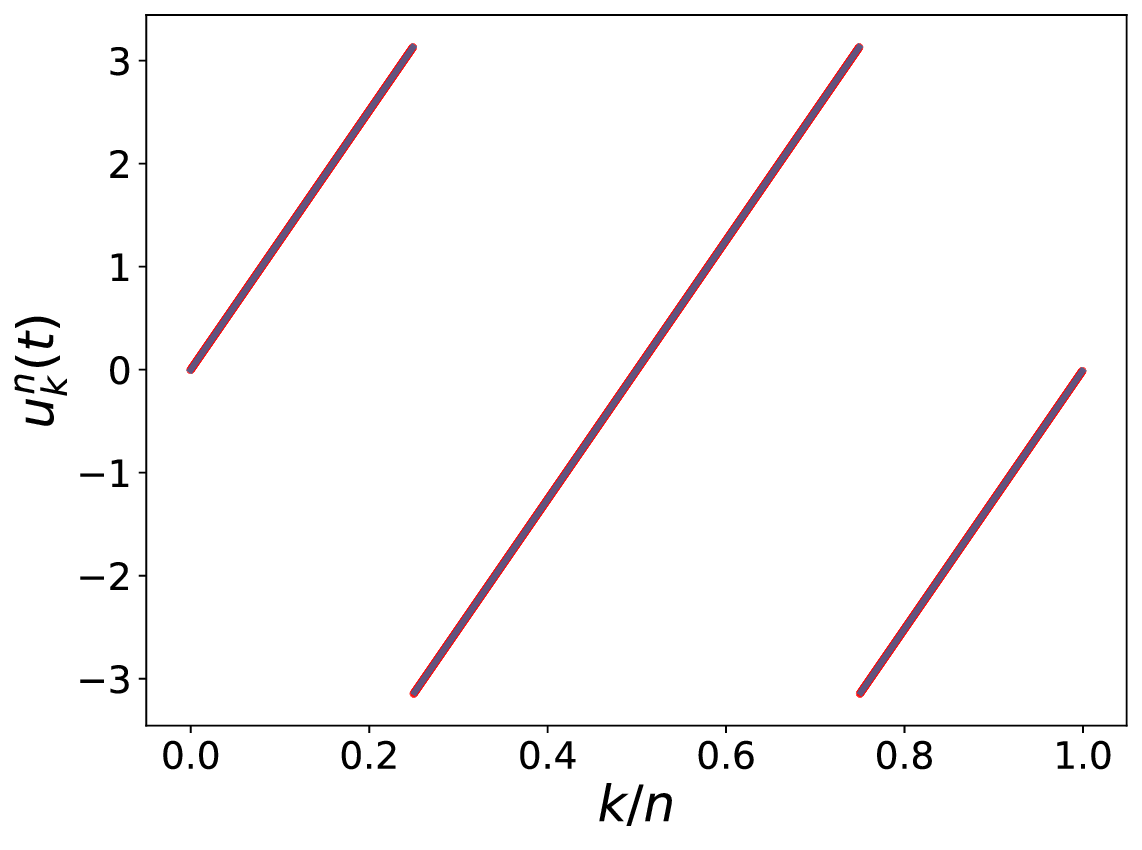}\\[-1ex]
{\footnotesize(c)}
\end{center}
\end{minipage}
\begin{minipage}[t]{0.495\textwidth}
\begin{center}
\includegraphics[scale=0.245]{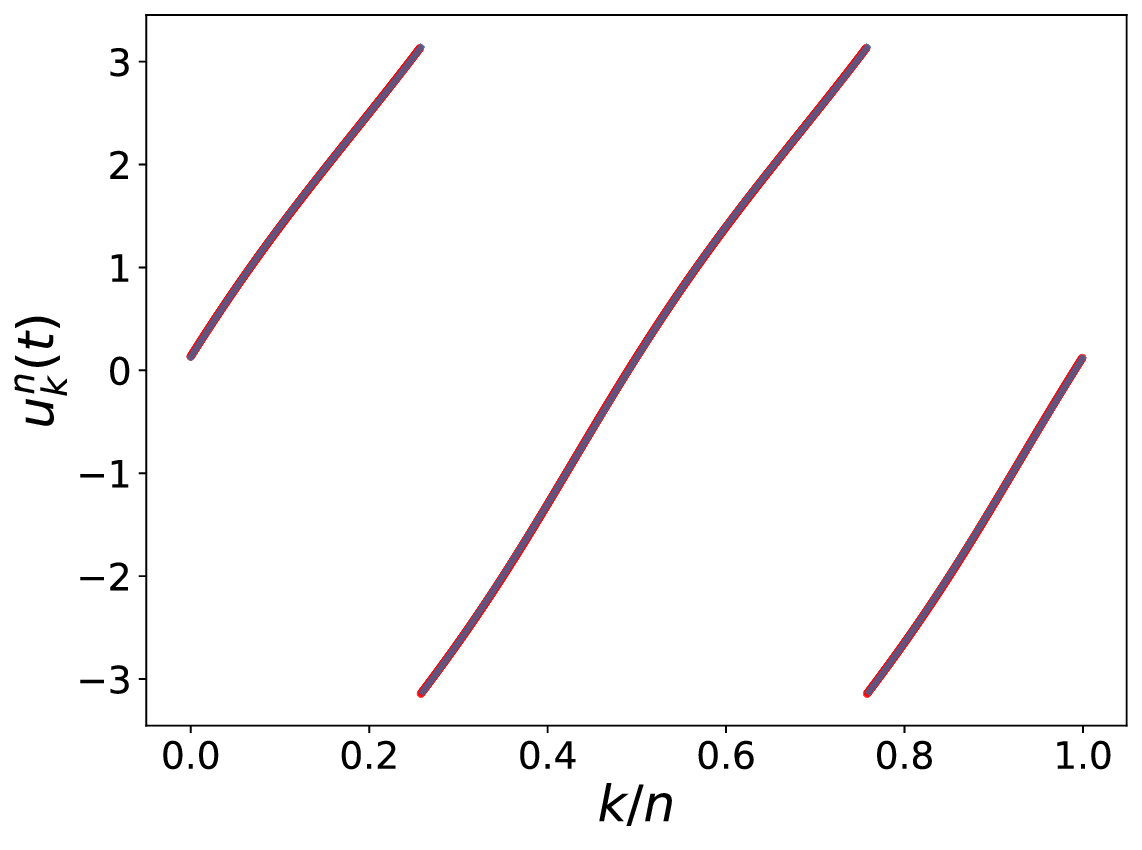}\\[-1ex]
{\footnotesize(d)}
\end{center}
\end{minipage}
\vspace*{0.5ex}

\begin{minipage}[t]{0.495\textwidth}
\begin{center}
\includegraphics[scale=0.245]{s3c1_fig3.eps}\\[-1ex]
{\footnotesize(e)}
\end{center}
\end{minipage}
\begin{minipage}[t]{0.495\textwidth}
\begin{center}
\includegraphics[scale=0.245]{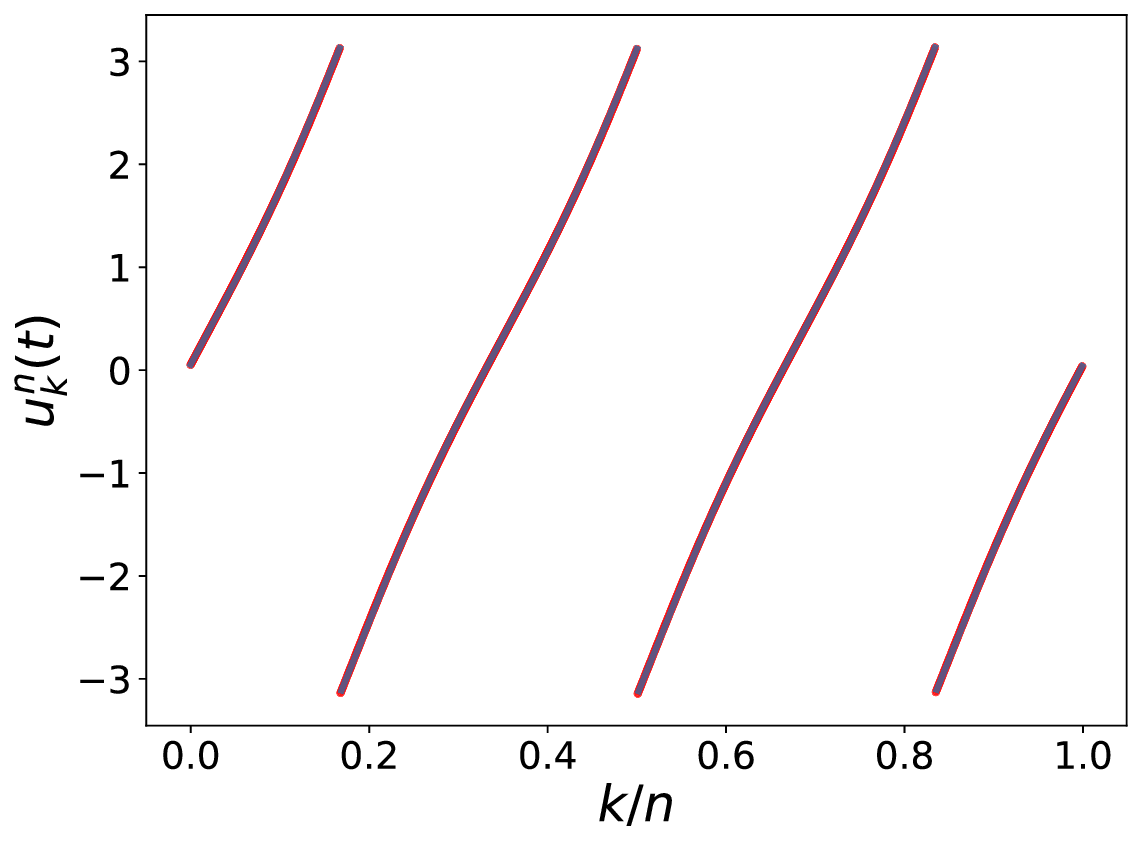}\\[-1ex]
{\footnotesize(f)}
\end{center}
\end{minipage}
\vspace*{0.5ex}

\begin{minipage}[t]{0.495\textwidth}
\begin{center}
\includegraphics[scale=0.245]{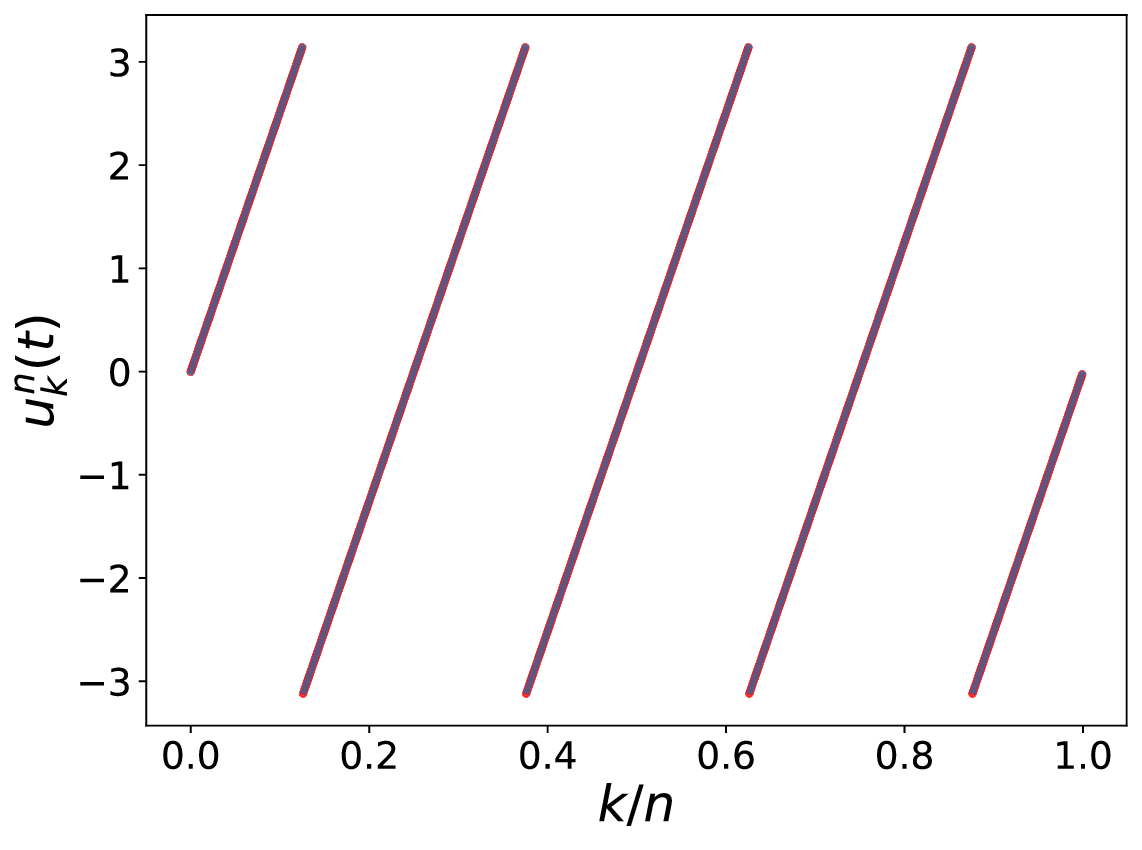}\\[-1ex]
{\footnotesize(g)}
\end{center}
\end{minipage}
\begin{minipage}[t]{0.495\textwidth}
\begin{center}
\includegraphics[scale=0.245]{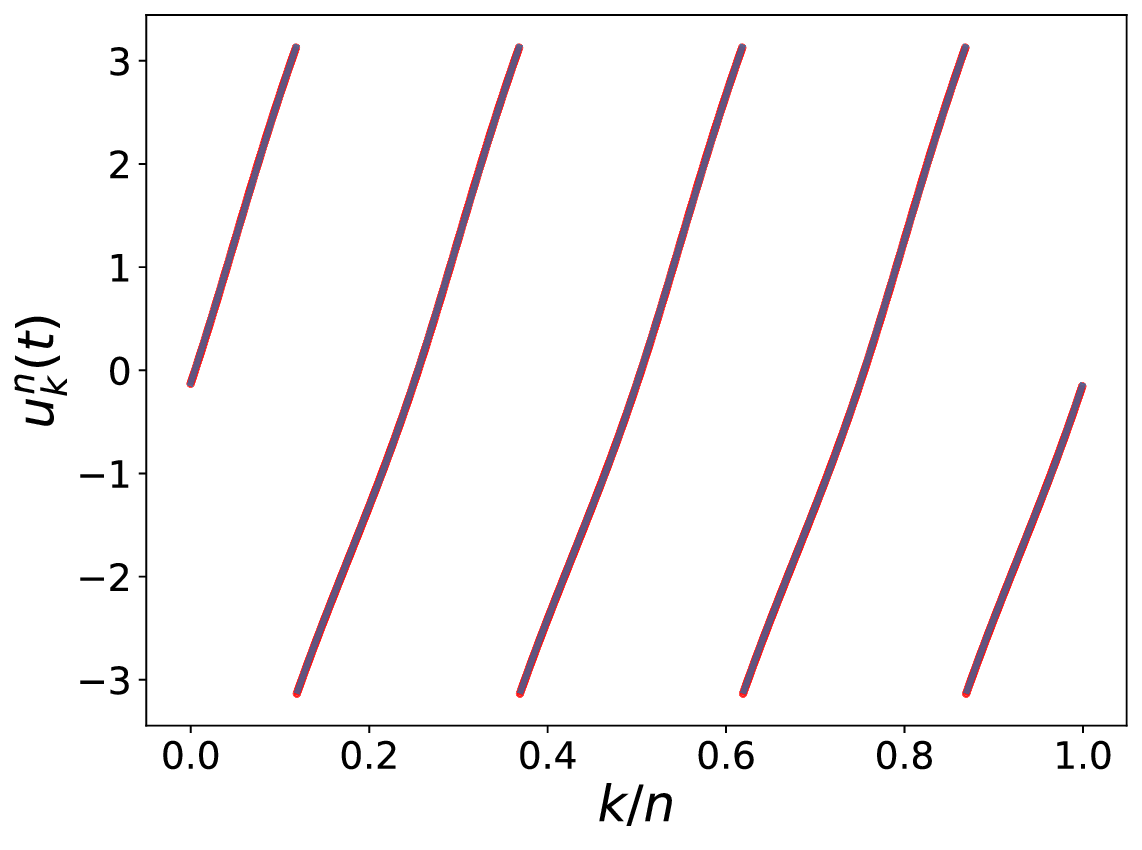}\\[-1ex]
{\footnotesize(h)}
\end{center}
\end{minipage}
\caption{Steady states of  the KM \eqref{eqn:dsys}
 with $n=1000$, $\kappa=0.5$, $\sigma=\pi/3$ and $b_3=0.5$ at $t=1000$:
(a) $(q,b_1)=(1,0.26)$;
(b) $(1,0.24)$;
(c) $(2,0.26)$;
(d) $(2,0.24)$;
(e) $(3,0.26)$;
(f) $(3,0.24)$;
(g) $(4,0.26)$;
(h) $(4,0.24)$.
See also the caption of Fig.~\ref{fig:5a2}.}
\label{fig:6b2}
\end{figure}

\begin{figure}[t]
\begin{minipage}[t]{0.495\textwidth}
\begin{center}
\includegraphics[scale=0.245]{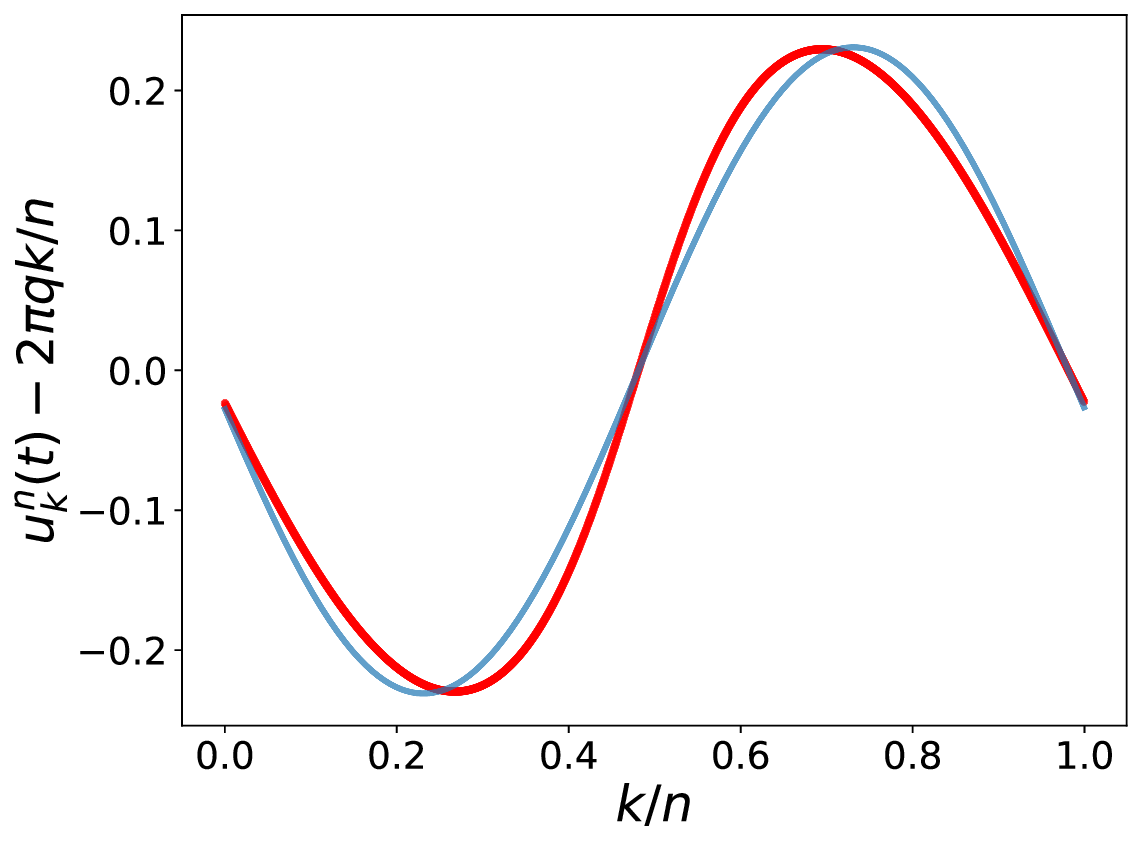}\\[-1ex]
{\footnotesize(a)}
\end{center}
\end{minipage}
\begin{minipage}[t]{0.495\textwidth}
\begin{center}
\includegraphics[scale=0.245]{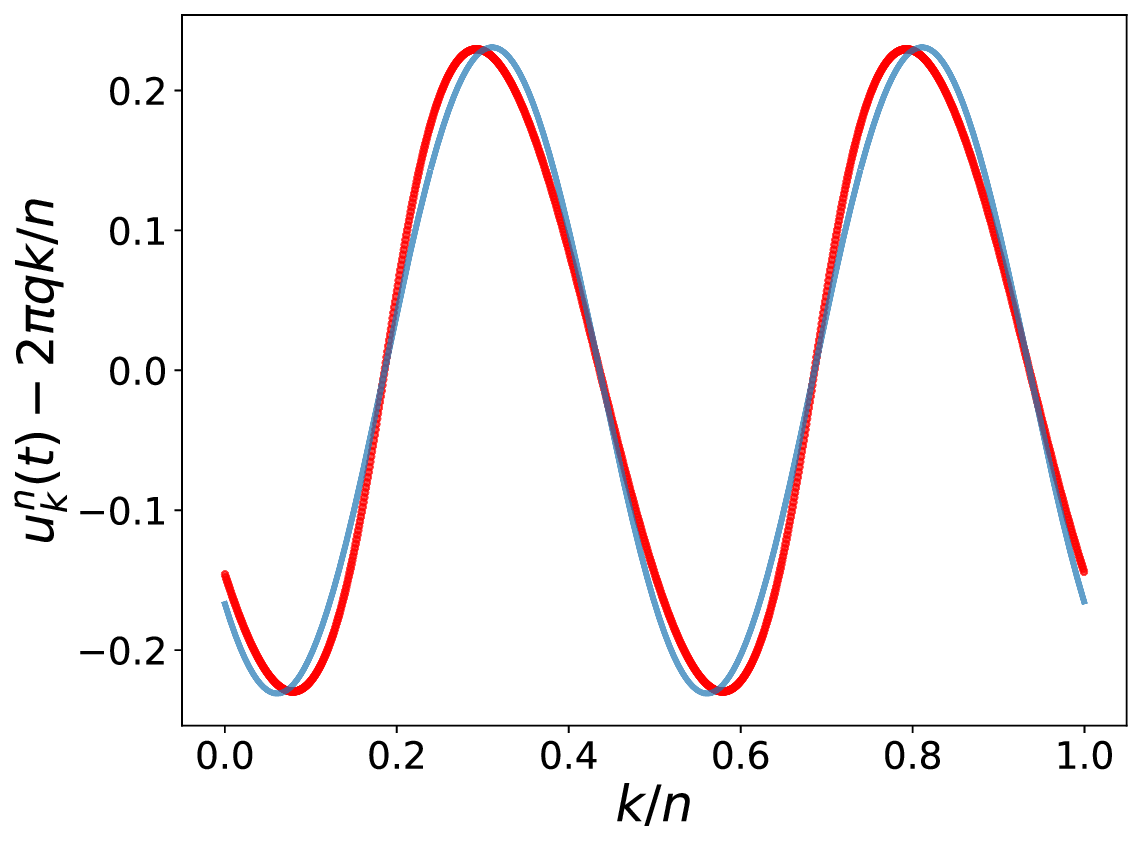}\\[-1ex]
{\footnotesize(b)}
\end{center}
\end{minipage}
\vspace*{0.5ex}

\begin{minipage}[t]{0.495\textwidth}
\begin{center}
\includegraphics[scale=0.245]{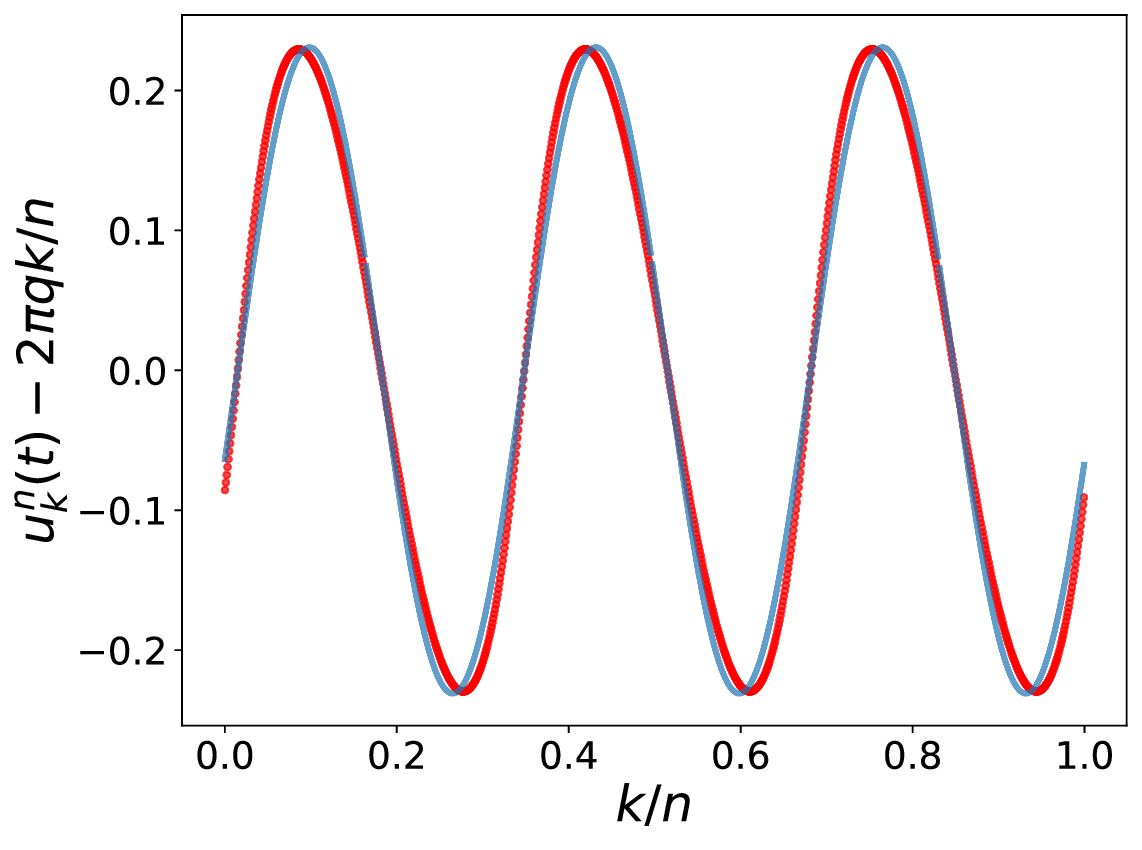}\\[-1ex]
{\footnotesize(c)}
\end{center}
\end{minipage}
\begin{minipage}[t]{0.495\textwidth}
\begin{center}
\includegraphics[scale=0.245]{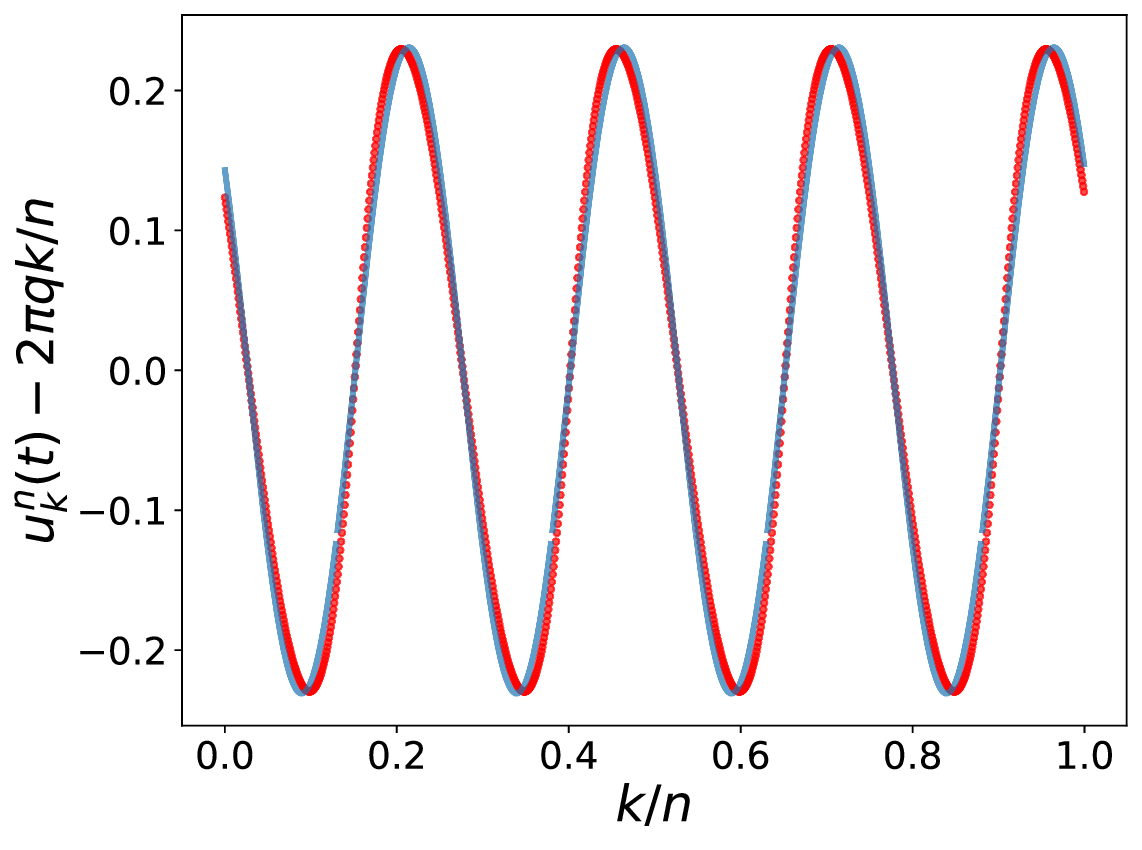}\\[-1ex]
{\footnotesize(d)}
\end{center}
\end{minipage}
\caption{Deviation from the $q$-twisted states in the steady states
 in the KM \eqref{eqn:dsys}
 with $n=1000$, $\kappa=0.5$, $\sigma=0$ and $b_1=0.48$ and $b_3=0.5$:
(a) $q=1$;
(b) $2$;
(c) $3$;
(d) $4$.
See also the caption of Fig.~\ref{fig:5a3}.}
\label{fig:6a3}
\end{figure}

\begin{figure}[t]
\begin{minipage}[t]{0.495\textwidth}
\begin{center}
\includegraphics[scale=0.245]{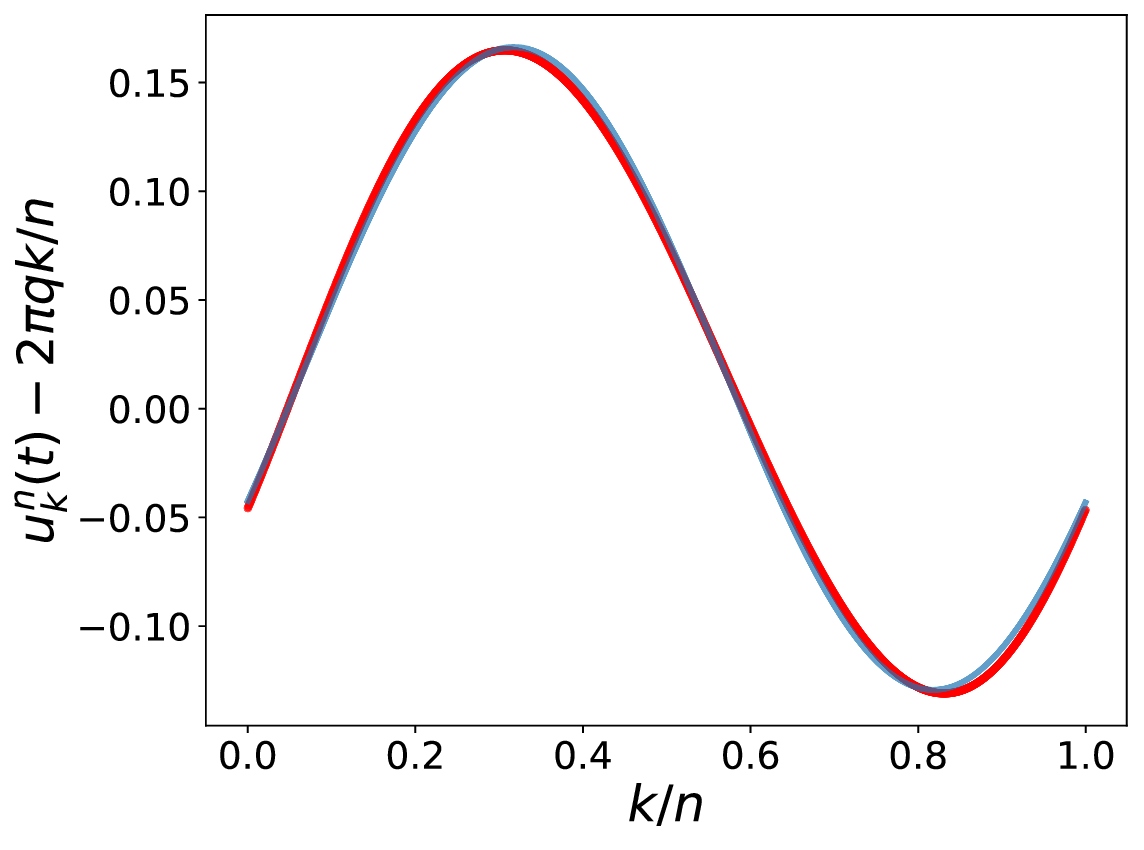}\\[-1ex]
{\footnotesize(a)}
\end{center}
\end{minipage}
\begin{minipage}[t]{0.495\textwidth}
\begin{center}
\includegraphics[scale=0.245]{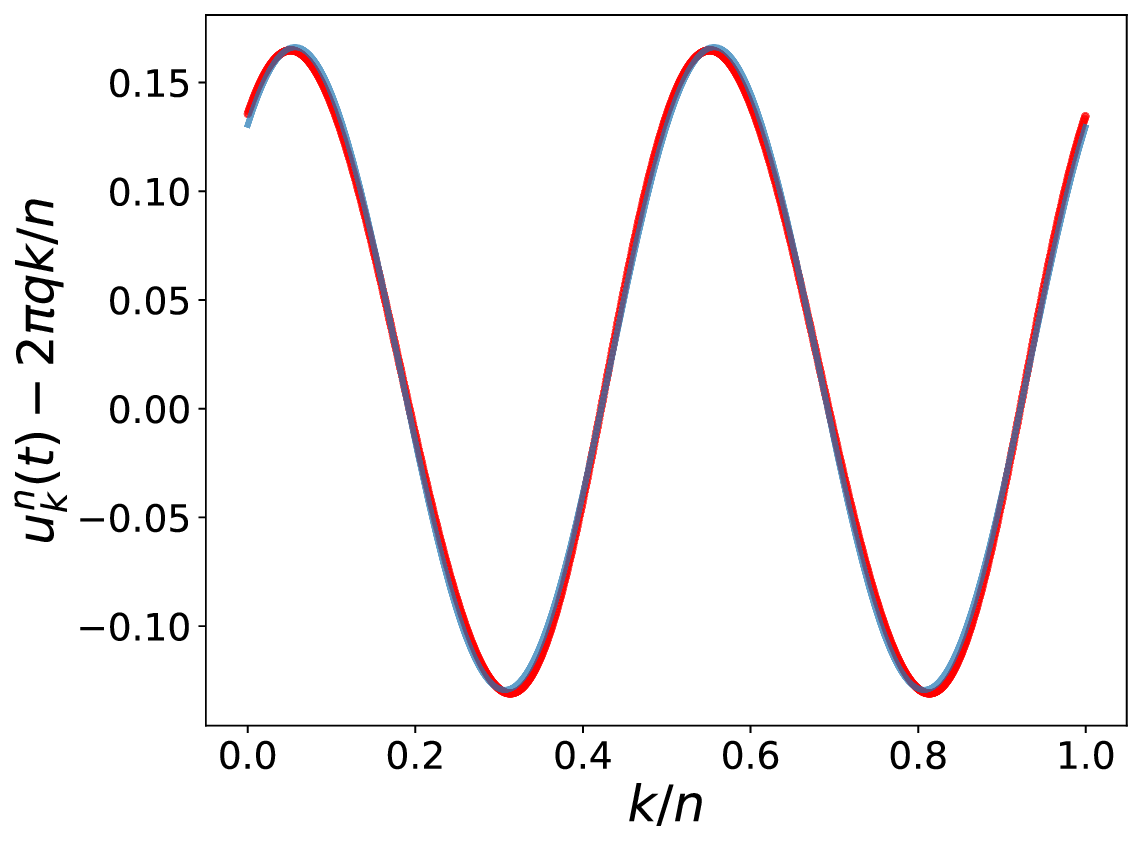}\\[-1ex]
{\footnotesize(b)}
\end{center}
\end{minipage}
\vspace*{0.5ex}

\begin{minipage}[t]{0.495\textwidth}
\begin{center}
\includegraphics[scale=0.245]{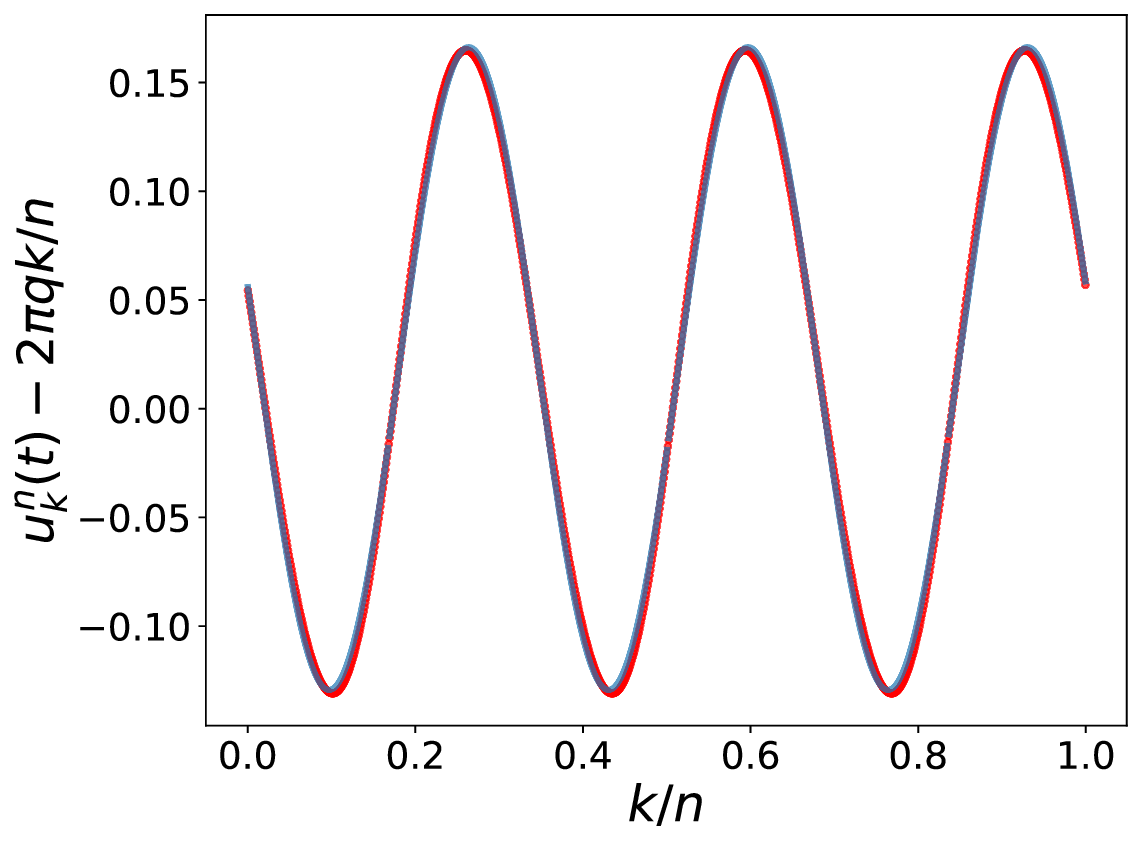}\\[-1ex]
{\footnotesize(c)}
\end{center}
\end{minipage}
\begin{minipage}[t]{0.495\textwidth}
\begin{center}
\includegraphics[scale=0.245]{s4c2_fig4.eps}\\[-1ex]
{\footnotesize(d)}
\end{center}
\end{minipage}
\caption{Deviation from the $q$-twisted states in the steady states
 in  the KM \eqref{eqn:dsys}
 with  $n=1000$, $\kappa=0.5$, $\sigma=\pi/3$,  $b_1=0.24$ and $b_3=0.5$:
(a) $q=1$;
(b) $2$;
(c) $3$;
(d) $4$.
See also the caption of Fig.~\ref{fig:5a3}.}
\label{fig:6b3}
\end{figure}

\begin{figure}[t]
\begin{minipage}[t]{0.495\textwidth}
\begin{center}
\includegraphics[scale=0.4]{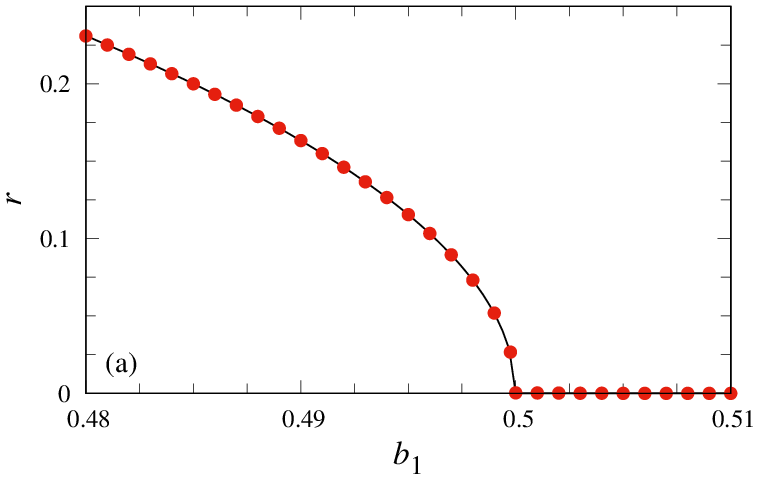}
\end{center}
\end{minipage}
\begin{minipage}[t]{0.495\textwidth}
\begin{center}
\includegraphics[scale=0.4]{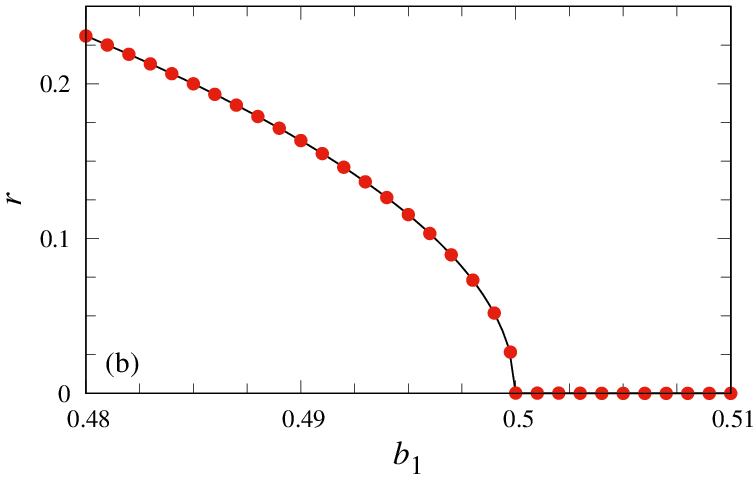}
\end{center}
\end{minipage}

\begin{minipage}[t]{0.495\textwidth}
\begin{center}
\includegraphics[scale=0.4]{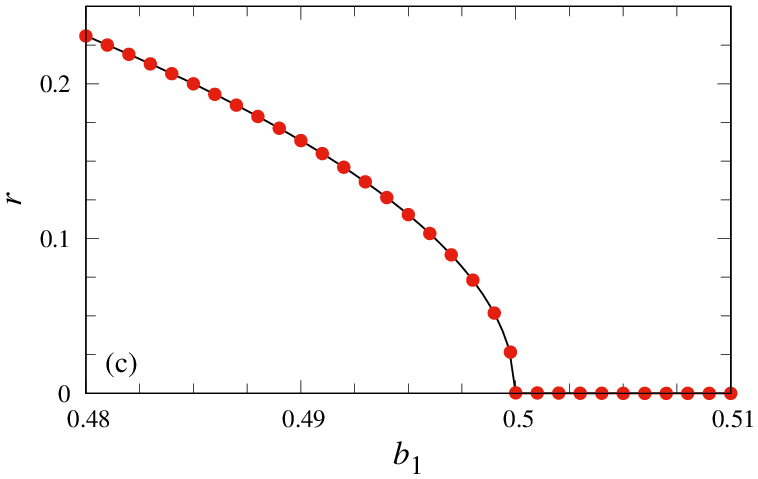}
\end{center}
\end{minipage}
\begin{minipage}[t]{0.495\textwidth}
\begin{center}
\includegraphics[scale=0.4]{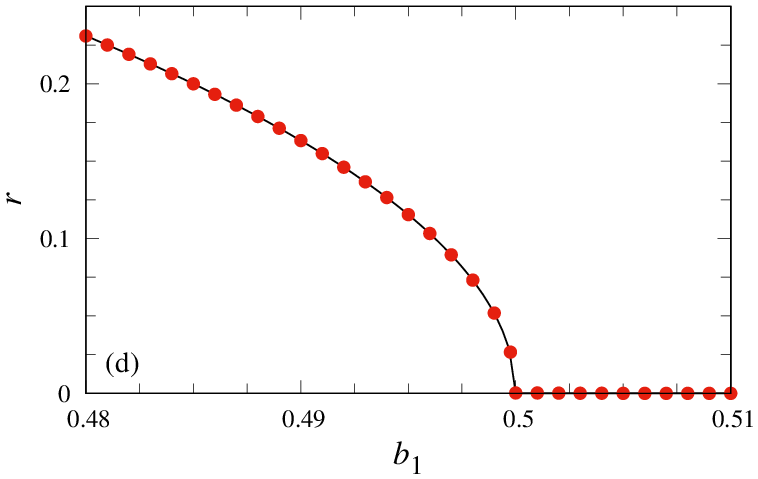}
\end{center}
\end{minipage}
\caption{Bifurcation diagram of the steady states in  the KM \eqref{eqn:dsys}
 with  $n=1000$, $\kappa=0.5$, $\sigma=0$ and $b_3=0.5$:
(a) $q=1$;
(b) $2$;
(c) $3$;
(d) $4$.
See also the caption of Fig.~\ref{fig:5a4}.}
\label{fig:6a4}
\end{figure}

Figures~\ref{fig:6a1} and \ref{fig:6b1}
 show the time-histories of every $100$th node (from 50th to 950th)
 for $\sigma=0$ and $\pi/3$, respectively, like Figs.~\ref{fig:5a1} and \ref{fig:5b1}.
Here the values of $u_k^n(t)\mod 2\pi$, $k\in[n]$, are plotted as the ordinates.
The values of $b_1=0.52$ and $0.48$ (resp. $b_1=0.26$ and $0.24$)
 were chosen in the left and right columns of Fig.~\ref{fig:6a1}
 (resp. Fig.~\ref{fig:6b1}), respectively,
 and they are larger and smaller than the bifurcation point
 approximated by $b_{1q}=0.5$ (resp. $b_{1q}=0.25$) (see Table.~\ref{tbl:4a}).
We see that the responses converge to the steady states rapidly,
 and oscillations occur in the right column of Fig.~\ref{fig:6b1}
 for $\sigma=\pi/3$ and $b_1=0.24$,
 as detected by Theorem~\ref{thm:4b} for the CL \eqref{eqn:csys}.

In Figs.~\ref{fig:6a2} and \ref{fig:6b2},
 $u_k^n(t)$, $k\in[n]$, at $t=1000$,
 which may be regarded as the steady states
 from the results of Figs.~\ref{fig:6a1} and \ref{fig:6b1},
 are plotted as small red disks  for $\sigma=0$ and $\pi/3$, respectively.
Here the same values of $b_1$ and $u_k^n(0)$, $k\in[n]$,
 as in Figs.~\ref{fig:6a1} and \ref{fig:6b1} were used.
We observe that the responses of the KM \eqref{eqn:dsys}
 converge to the twisted and modulated or oscillating twisted states
 for $b_1=0.52$ or $0.26$ and $0.48$ or $0.24$
 in the left and right columns of each figure, respectively,
 as predicted by Theorems~\ref{thm:4a} and \ref{thm:4b}
 with the assistance of Corollary~\ref{cor:2a} and Theorem~\ref{thm:2e}.
Indeed, we confirmed that
 the deviation from the twisted state is about $10^{-11}$ and $10^{-6}$ at most
 in the left columns of Figs.~\ref{fig:6a2} and \ref{fig:6b2}
 for $b_1=0.52$ and $0.26$, respectively.
In particular, the target state \eqref{eqn:ts} is accomplished there.
The most probably leading term \eqref{eqn:ss}
 in the modulated and oscillating twisted solutions \eqref{eqn:thm4a} and \eqref{eqn:thm4b}
 was also estimated from the numerical simulation results
 by using the least mean square method
 and is plotted as a blue line in each figure,
 as in Figs.~\ref{fig:5a2} and \ref{fig:5b2}.
Both results coincide almost completely,
 as detected by Theorems~\ref{thm:4a} and \ref{thm:4b}
 for the CL \eqref{eqn:csys}.

In Figs.~\ref{fig:6a3} and \ref{fig:6b3},
 the deviation, $u_k^n(t)-2\pi q k/n$, $k\in[n]$,
 of the steady state in the right columns
 of Figs.~\ref{fig:6a2} and \ref{fig:6b2}
 from the desired $q$-twisted one in the KM \eqref{eqn:dsys}
 with $\sigma=0$ and $\sigma=\pi/3$, respectively,
 is plotted as small red disks.
Estimates obtained from the most probably leading terms
 displayed in Figs.~\ref{fig:6a2} and \ref{fig:6b2}
 are also plotted as blue lines, as in Figs.~\ref{fig:5a3} and \ref{fig:5b3}.
Both results coincide almost completely.

\begin{figure}[t]
\begin{minipage}[t]{0.495\textwidth}
\begin{center}
\includegraphics[scale=0.4]{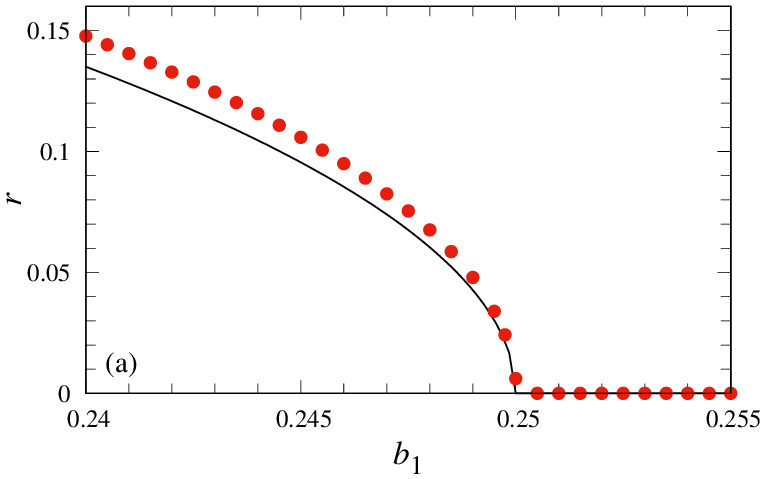}
\end{center}
\end{minipage}
\begin{minipage}[t]{0.495\textwidth}
\begin{center}
\includegraphics[scale=0.4]{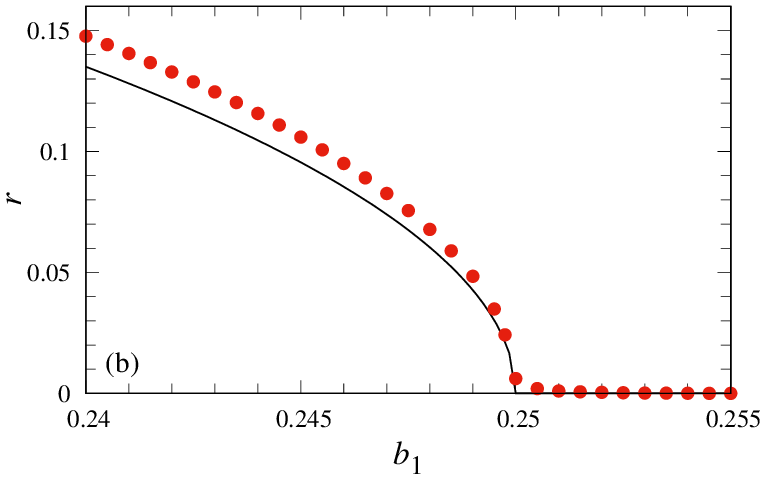}
\end{center}
\end{minipage}

\begin{minipage}[t]{0.495\textwidth}
\begin{center}
\includegraphics[scale=0.4]{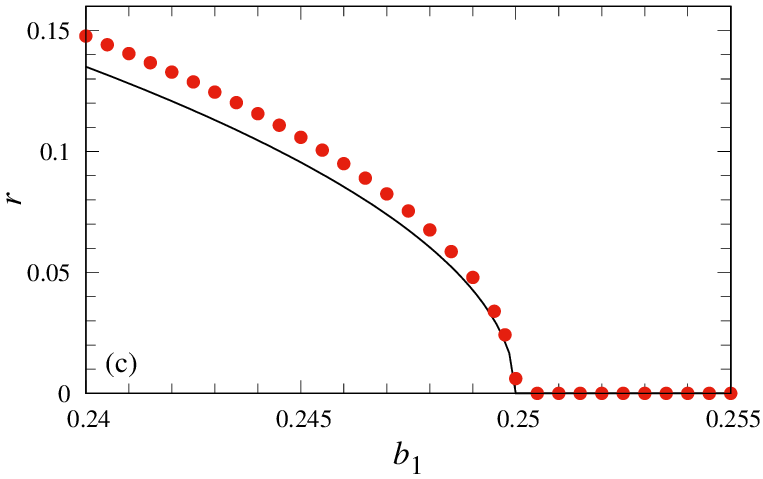}
\end{center}
\end{minipage}
\begin{minipage}[t]{0.495\textwidth}
\begin{center}
\includegraphics[scale=0.4]{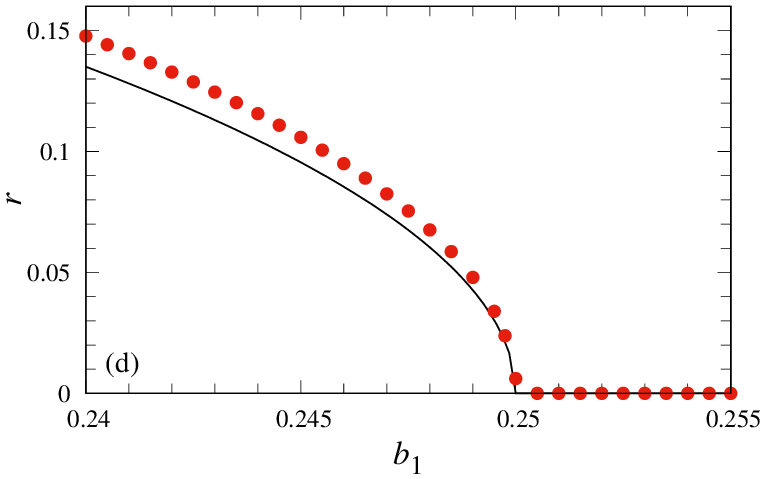}
\end{center}
\end{minipage}
\caption{Bifurcation diagram of the steady states in  the KM \eqref{eqn:dsys}
 with  $n=1000$, $\kappa=0.5$, $\sigma=\pi/3$ and $b_3=0.5$:
(a) $q=1$;
(b) $2$;
(c) $3$;
(d) $4$.
See also the caption of Fig.~\ref{fig:5a4}.}
\label{fig:6b4}
\end{figure}

Finally, we present numerically computed bifurcation diagrams
 for $\sigma=0$ and $\pi/3$ in Figs.~\ref{fig:6a4} and \ref{fig:6b4}, respectively,
 as in Figs.~\ref{fig:5a4} and \ref{fig:5b4}.
The amplitude $r$ of \eqref{eqn:ss}
 estimated from the numerical simulation results
 are plotted as small red disks,
 and the theoretical predictions given by \eqref{eqn:rth}
 are plotted as black solid lines.
Good agreement between both results is found,
 especially in Fig.~\ref{fig:6a4},
 although slight differences are seen in Fig.~\ref{fig:6b4}.

\section{Concluding Remarks}
We summarize this paper as follows:
We studied feedback control of twisted states in the KM \eqref{eqn:dsys}
 of identical oscillators defined on deterministic nearest neighbor graphs
 containing complete simple ones when it may have phase-lag.
Using the center manifold reduction technique \cite{HI11},
 we analyzed the stability and bifurcations 
 of twisted solutions in the CL \eqref{eqn:csys}
 for the KM \eqref{eqn:dsys} subjected to feedback control.
In particular, it was shown that the twisted solutions exist and can be stabilized
 not only for nearest neighbor graphs but also for complete simple graphs.
Moreover, the CL \eqref{eqn:csys} was shown to suffer bifurcations
 at which the twisted solutions becomes unstable
 and a stable one-parameter family of modulated or oscillating twisted solutions
 is born, depending on whether the phase-lag is zero or not.
We demonstrated the theoretical results by numerical simulations
 for the feedback controlled KM \eqref{eqn:dsys}
 on deterministic nearest neighbor and complete simple graphs.

{\color{black}
From an applied viewpoint, the present work suggests that
 such a simple feedback law as the linear-cubic terms in \eqref{eqn:dsys}
 can be used to realize and maintain prescribed twisted phase patterns
 in coupled oscillator networks by tuning only a few scalar gains.
Kuramoto-type phase descriptions are widely used as reduced models across disciplines,
 including power-system and microgrid settings
 where the stability of phase-locked operating states
 and their control under disturbances are central
 \cite{FNP08,DB12,GZLWY21,AA22,SA15},
 and biological and neural oscillator networks
 where spatiotemporal rhythms and network level oscillatory patterns are of interest
 \cite{BHD10,BGLM20,SNKJZBPPB20,YLYWQZH25,ZLQWWL22}.
For the KM \eqref{eqn:dsys},
 the CL analysis provides explicit stability regions and bifurcations
 of twisted solutions,
 thereby offering analytically tractable guidance for gain tuning and mode selection.
More broadly, the results fit into the general theme of controlling collective behavior
 in complex networks under practical constraints \cite{CKM13,M15,DNM22,SA16}.
Although our main analysis is carried out in the CL framework,
 it is reasonable to expect that it remains at least qualitatively valid
 for moderate-size graphs (e.g., $n\sim 50$).
The CL predictions can then be used 
 to indicate parameter regions where qualitative transitions occur,
 which can be tested numerically on such moderate-size graphs.
}

Finally, we give some comments for future work.
The KM and its generalization with phase-lag
 was studied for different solutions from twisted ones
 in \cite{OW16,C17,O18,BGLM20,MM22,MRSKG24,LXYZLLJ25}.
In particular, chimera states were discussed in \cite{O18,BGLM20,MM22}.
The theory reviewed in Section~2 is also applicable to these cases
 and may be useful to uncover their dynamics.
The KM with time delay,
 which possesses important applications in neuroscience \cite{CEVB97,BNBDMSM23,EK09},
 has often been a subject of research
 \cite{ABPRS05,ADKMZ08,PR15,RPJK16,MRSKG24}.
It will be another next target to extend the theory of Section~2
 to the KM with time delay and its CL.

\section*{Acknowledgements}
This work was partially supported by the JSPS KAKENHI Grant Number JP23K22409.
{\color{black}
The author thanks the anonymous referees
 for their helpful comments and constructive suggestions.
In particular, by one of them,
 he could realize an error in Theorem~\ref{thm:2c} of the original manuscript
 and correct it.}
 

\appendix

\renewcommand{\theequation}{\Alph{section}.\arabic{equation}}

\section{Derivation of \eqref{eqn:ifex}}

We first rewrite the CL \eqref{eqn:csys}
 in the rotational frame with the rotational speed $\Omega$ as
\begin{align}
\frac{\partial}{\partial t}u(t,x)
=&\omega-\Omega+p\left(\cos u(t,x)\int_{x-\kappa}^{x+\kappa}\sin u(t,y)\d y\right.\notag\\
&\qquad\qquad
\left. -\sin u(t,x)\int_{x-\kappa}^{x+\kappa}\cos u(t,y)\d y\right)\cos\sigma\notag\\
&
+p\left(\sin u(t,x)\int_{x-\kappa}^{x+\kappa}\sin u(t,y)\d y\right.\notag\\
&\qquad\qquad
\left. +\cos u(t,x)\int_{x-\kappa}^{x+\kappa}\cos u(t,y)\d y\right)\sin\sigma\notag\\
& -b_1(\bar{u}(t,x)-u(t,x))-b_3(\bar{u}(t,x)-u(t,x))^3.
\label{eqn:csys2}
\end{align}
Letting \eqref{eqn:solex} with $\Omega=0$, we have
\begin{align*}
\cos u(t,x)
=&\cos 2\pi qx-\sin 2\pi qx
\biggl(\xi_0+\sum_{j=1}^\infty(\xi_j\cos 2\pi jx+\eta_j\sin2\pi jx)\biggr)\\
&
-\cos 2\pi qx\bigl(\tfrac{1}{4}((\xi_q^2+\eta_q^2)
 +(\xi_q^2-\eta_q^2)\cos 4\pi qx+2\xi_q\eta_q\sin 4\pi qx)\\
&\quad
+\xi_0(\xi_q\cos2\pi q x+\eta_q\sin2\pi q x)\\
&\quad
+\tfrac{1}{2}\sum_{j\neq q}((\xi_q\xi_j+\eta_q\eta_j)\cos2\pi(q-j)x
 -(\xi_q\eta_j-\xi_j\eta_q)\sin2\pi(q-j)x\\
&\qquad
+(\xi_q\xi_j-\eta_q\eta_j)\cos2\pi(q+j)x
 +(\xi_q\eta_j+\xi_j\eta_q)\sin2\pi(q+j)x)\bigr)\\
&
+\sin 2\pi qx\bigl(
 \tfrac{1}{8}(\xi_q^2+\eta_q^2)(\xi_q\cos2\pi qx+\eta_q\sin 2\pi qx)\\
&\quad
+\tfrac{1}{24}((\xi_q^2-3\eta_q^2)\xi_q\cos 6\pi qx
 +(3\xi_q^2-\eta_q^2)\eta_q\sin 6\pi qx)\bigr)+\cdots
\end{align*}
and
\begin{align*}
\sin u(t,x)
=& \sin 2\pi qx+\cos 2\pi qx
\biggl(\xi_0+\sum_{j=1}^\infty(\xi_j\cos 2\pi jx+\eta_j\sin2\pi jx)\biggr)\\
&
-\sin 2\pi qx\bigl(\tfrac{1}{4}((\xi_q^2+\eta_q^2)
 +(\xi_q^2-\eta_q^2)\cos 4\pi  qx+2\xi_q\eta_q\sin 4\pi qx)\\
&\quad
+\xi_0(\xi_q\cos2\pi q x+\eta_q\sin2\pi q x)\\
&\quad
+\tfrac{1}{2}\sum_{j\neq q}((\xi_q\xi_j+\eta_q\eta_j)\cos2\pi(q-j)x
 -(\xi_q\eta_j-\xi_j\eta_q)\sin2\pi(q-j)x\\
&\qquad
+(\xi_q\xi_j-\eta_q\eta_j)\cos2\pi(q+j)x
 +(\xi_q\eta_j+\xi_j\eta_q)\sin2\pi(q+j)x)\bigr)\\
&
-\cos 2\pi qx\bigl(
 \tfrac{1}{8}(\xi_q^2+\eta_q^2)(\xi_q\cos2\pi qx+\eta_q\sin 2\pi qx)\\
&\quad
+\tfrac{1}{24}((\xi_q^2-3\eta_q^2)\xi_q\cos 6\pi qx
 +(3\xi_q^2-\eta_q^2)\eta_q\sin 6\pi qx)\bigr)+\cdots,
\end{align*}
where `$\cdots$' represents higher-order terms of
\[
O\left(\xi_q^4+\eta_q^4+\xi_0^2
 +\sum_{j=1,j\neq q}^\infty(\xi_j^2+\eta_j^2)\right).
\]

We compute the integrals in \eqref{eqn:csys2} as
\begin{align*}
&
\int_{x-\kappa}^{x+\kappa}\cos u(t,y)\d y\hspace*{6em}\\
&
=-a_2(q,0)\cos2\pi qx+a_2(q,0)\xi_0\sin2\pi qx\\
&\quad
 -\sum_{j=1}^\infty(a_1(q,j)(\xi_j\sin2\pi jx-\eta_j\cos2\pi jx)\cos2\pi qx\\
&\qquad\quad
-a_2(q,j)(\xi_j\cos2\pi j+\eta_j\sin2\pi jx)\sin2\pi qx)\\
&\quad
 +\tfrac{1}{4}a_2(q,0)(\xi_q^2+\eta_q^2)\cos 2\pi qx\\
&\quad
+\tfrac{1}{4}a_1(q,2q)((\xi_q^2-\eta_q^2)\sin4\pi qx-2\xi_q\eta_q\cos4\pi qx)\sin2\pi q x\\
&\quad
 +\tfrac{1}{4}a_2(q,2q)((\xi_q^2-\eta_q^2)\cos4\pi qx+2\xi_q\eta_q\sin4\pi qx)\cos2\pi qx
\hspace*{5em}
\end{align*}
\begin{align*}
&\quad
+\tfrac{1}{8}a_1(q,q)(\xi_q^2+\eta_q^2)(\xi_q\sin2\pi qx-\eta_q\cos2\pi qx)\cos2\pi qx\\
&\quad
-\tfrac{1}{8}a_2(q,q)(\xi_q^2+\eta_q^2)(\xi_q\cos2\pi qx+\eta_q\sin2\pi qx)\sin2\pi qx\\
&\quad
 +\tfrac{1}{24}a_1(q,3q)
 ((\xi_q^2-3\eta_q^2)\xi_q\sin6\pi qx-(3\xi_q^2-\eta_q^2)\eta_q\cos6\pi qx)\cos2\pi qx\\
&\quad
 -\tfrac{1}{24}a_2(q,3q)
 ((\xi_q^2-3\eta_q^2)\xi_q\cos6\pi qx+(3\xi_q^2-\eta_q^2)\eta_q\sin6\pi qx)\sin2\pi qx\\
&\quad
+a_1(q,q)(\xi_0\xi_q\sin2\pi qx-\xi_0\eta_q\cos2\pi qx)\sin2\pi qx\\
&\quad
+a_2(q,q)(\xi_0\xi_q\cos2\pi qx+\xi_0\eta_q\sin2\pi qx)\cos2\pi qx\\
&\quad
 +\tfrac{1}{2}\sum_{j\neq q}\bigl(a_1(q,q-j)((\xi_q\xi_j+\eta_q\eta_j)\sin2\pi(q-j)x\\
& \qquad\quad
 +(\xi_q\eta_j-\xi_j\eta_q)\cos2\pi(q-j)x)\sin2\pi qx\\
&\qquad
 +a_2(q,q-j)((\xi_q\xi_j+\eta_q\eta_j)\cos2\pi(q-j)x\\
 & \qquad\quad
  -(\xi_q\eta_j-\xi_j\eta_q)\sin2\pi(q-j)x)\cos2\pi q x\\
&\qquad
 +a_1(q,q+j)((\xi_q\xi_j-\eta_q\eta_j)\sin2\pi(q+j)x\\
& \qquad\quad
 -(\xi_q\eta_j+\xi_j\eta_q)\cos2\pi(q+j)x)\sin2\pi q x\\
&\qquad
 +a_2(q,q+j)((\xi_q\xi_j-\eta_q\eta_j)\cos2\pi(q+j)x\\
 & \qquad\quad
 +(\xi_q\eta_j+\xi_j\eta_q)\sin2\pi(q+j)x)\cos2\pi q x\bigr)
 +\cdots
\end{align*}
and
\begin{align*}
&
\int_{x-\kappa}^{x+\kappa}\sin u(t,y)\d y\\
&
=-a_2(q,0)\sin2\pi qx-a_2(q,0)\xi_0\cos2\pi qx\\
&\quad
-\sum_{j=1}^\infty(a_1(q,j)(\xi_j\sin2\pi jx-\eta_j\cos2\pi jx)\sin2\pi qx\\
&\qquad\quad
+a_2(q,j)(\xi_j\cos2\pi j+\eta_j\sin2\pi jx)\cos2\pi qx)\\
&\quad
 +\tfrac{1}{4}a_2(q,0)(\xi_q^2+\eta_q^2)\sin 2\pi qx\\
&\quad
-\tfrac{1}{4}a_1(q,2q)((\xi_q^2-\eta_q^2)\sin4\pi qx-2\xi_q\eta_q\cos4\pi qx)\cos2\pi q x\\
&\quad
+\tfrac{1}{4}a_2(q,2q)((\xi_q^2-\eta_q^2)\cos4\pi qx+2\xi_q\eta_q\sin4\pi qx)\sin2\pi qx\\
&\quad
+\tfrac{1}{8}a_1(q,q)(\xi_q^2+\eta_q^2)(\xi_q\sin2\pi qx-\eta_q\cos2\pi qx)\sin2\pi qx\\
&\quad
+\tfrac{1}{8}a_2(q,q)(\xi_q^2+\eta_q^2)(\xi_q\cos2\pi qx+\eta_q\sin2\pi qx)\cos2\pi qx\\
&\quad
+\tfrac{1}{24}a_1(q,3q)
 ((\xi_q^2-3\eta_q^2)\xi_q\sin6\pi qx-(3\xi_q^2-\eta_q^2)\eta_q\cos6\pi qx)\sin2\pi qx\\
&\quad
+\tfrac{1}{24}a_2(q,3q)
 ((\xi_q^2-3\eta_q^2)\xi_q\cos6\pi qx+(3\xi_q^2-\eta_q^2)\eta_q\sin6\pi qx)\cos2\pi qx\\
&\quad
-a_1(q,q)(\xi_0\xi_q\sin2\pi qx-\xi_0\eta_q\cos2\pi qx)\cos2\pi qx\\
&\quad
+a_2(q,q)(\xi_0\xi_q\cos2\pi q+\xi_0\eta_q\sin2\pi qx)\sin2\pi qx\\
&\quad
 -\tfrac{1}{2}\sum_{j\neq q}\bigl(a_1(q,q-j)((\xi_q\xi_j+\eta_q\eta_j)\sin2\pi(q-j)x\\
&\qquad\quad
 +(\xi_q\eta_j-\xi_j\eta_q)\cos2\pi(q-j)x)\cos2\pi q x\\
&\qquad
 -a_2(q,q-j)((\xi_q\xi_j+\eta_q\eta_j)\cos2\pi(q-j)x\\
&\qquad\quad
 -(\xi_q\eta_j-\xi_j\eta_q)\sin2\pi(q-j)x)\sin2\pi q x
\end{align*}
\begin{align*}
&\qquad
 +a_1(q,q+j)((\xi_q\xi_j-\eta_q\eta_j)\sin2\pi(q+j)x\\
&\qquad\quad
 -(\xi_q\eta_j+\xi_j\eta_q)\cos2\pi(q+j)x)\cos2\pi q x\\
&\qquad
 -a_2(q,q+j)((\xi_q\xi_j-\eta_q\eta_j)\cos2\pi(q+j)x\\
&\qquad\quad
 +(\xi_q\eta_j+\xi_j\eta_1)\sin2\pi(q+j)x)\sin2\pi q x\bigr)+\cdots.
\hspace*{10em}
\end{align*}
We substitute \eqref{eqn:solex} into \eqref{eqn:csys2},
 integrate the resulting equation with respect to $x$ from $0$ to $1$
 after multiplying it with $\cos 2\pi j$ or $\sin 2\pi j$, $j\in\Nset$.
Thus, we obtain \eqref{eqn:ifex} for $q\in[4]$ after lengthy calculations.


\end{document}